\documentclass[10pt,a4paper]{amsart}
\usepackage{amssymb,amsmath,amsthm,amscd}
\usepackage{leftidx,caption,enumitem,mathtools,booktabs,scrextend,array,hhline,tablefootnote}
\usepackage{graphicx,mathrsfs,bbm}
\usepackage{color}
\usepackage{microtype}
\usepackage{hyperref}
\usepackage[usenames,dvipsnames,svgnames,table]{xcolor}
\usepackage{tikz}
\usepackage{comment}
\usepackage{xfrac,faktor}
\usepackage{epstopdf} 
\usepackage{float}
\usepackage{adjustbox}
\usepackage{tikz-cd}
\usepackage[all]{xy}
\newtheorem{lemma}{Lemma}[section]
\newtheorem{theorem}[lemma]{Theorem}
\newtheorem{corollary}[lemma]{Corollary}
\newtheorem{proposition}[lemma]{Proposition}

\newtheorem{problem}[lemma]{Problem}
\theoremstyle{definition}
\newtheorem{definition}[lemma]{Definition}
\theoremstyle{example}
\newtheorem{example}[lemma]{Example}
\theoremstyle{remark}
\newtheorem{remark}[lemma]{Remark}
\newcommand{\C}{\mathbb{C}}
\newcommand{\D}{\mathbb{D}}
\newcommand{\N}{\mathbb{N}}
\newcommand{\Q}{\mathbb{Q}}
\newcommand{\R}{\mathbb{R}}
\newcommand{\Z}{\mathbb{Z}}
\newcommand{\HP}{\mathbb{H}}
\newcommand{\mate}{\bot \!\! \! {\bot}}

\newcommand{\cC}{\mathcal{C}}
\newcommand{\cK}{\mathcal{K}}
\newcommand{\cS}{\mathcal{S}}

\newcommand{\cN}{\mathcal{N}}
\newcommand{\cF}{\mathcal{F}}
\newcommand{\cT}{\mathcal{T}}

\newcommand*\circcled[1]{\tikz[baseline=(char.base)]{
\node[shape=circle,draw,inner sep=1.3pt] (char) {#1};}}
\newcommand{\ciq}{\circcled{?}}

\makeatletter
\DeclareFontFamily{U}{tipa}{}
\DeclareFontShape{U}{tipa}{m}{n}{<->tipa10}{}
\newcommand{\arc@char}{{\usefont{U}{tipa}{m}{n}\symbol{62}}}

\newcommand{\arc}[1]{\mathpalette\arc@arc{#1}}

\newcommand{\arc@arc}[2]{
  \sbox0{$\m@th#1#2$}
  \vbox{
    \hbox{\resizebox{\wd0}{\height}{\arc@char}}
    \nointerlineskip
    \box0
  }
}
\makeatother

\DeclareMathOperator{\re}{Re}

\DeclareMathOperator{\Int}{int}

\renewcommand{\epsilon}{\varepsilon}
\renewcommand{\phi}{\varphi}

\DeclareMathOperator{\Deg}{deg}

\DeclareMathOperator{\ind}{ind}

\date{\today}

\begin{document}

\title[Mirrors of conformal dynamics]{Mirrors of conformal dynamics:\\ {\footnotesize Interplay between anti-rational maps, reflection groups,}\\ {\footnotesize Schwarz reflections, and correspondences}}

\begin{author}[M.~Lyubich]{Mikhail Lyubich}
\address{Institute for Mathematical Sciences, Stony Brook University, 100 Nicolls Rd, Stony Brook, NY 11794-3660, USA}
\email{mlyubich@math.stonybrook.edu}
\end{author}
\thanks{M.L. was partially supported by NSF grants DMS-1901357 and 2247613, a fellowship from the Hagler Institute for Advanced Study, and the Clay fellowship.}

\begin{author}[S.~Mukherjee]{Sabyasachi Mukherjee}
\address{School of Mathematics, Tata Institute of Fundamental Research, 1 Homi Bhabha Road, Mumbai 400005, India}
\email{sabya@math.tifr.res.in}
\thanks{S.M. was supported by the Department of Atomic Energy, Government of India, under project no.12-R\&D-TFR-5.01-0500, an endowment of the Infosys Foundation, and SERB research project grants SRG/2020/000018 and MTR/2022/000248.}
\end{author}

\begin{abstract}
The goal of this survey is to present intimate interactions between four branches of conformal dynamics: iterations of anti-rational maps, actions of Kleinian reflection groups, dynamics generated by Schwarz reflections in quadrature domains, and algebraic correspondences. We start with several examples of Schwarz reflections as well as algebraic correspondences obtained by matings between anti-rational maps and reflection groups, and examples of Julia set realizations for limit sets of reflection groups (including classical Apollonian-like gaskets). 
We follow up these examples with dynamical relations between explicit Schwarz reflection parameter spaces and parameter spaces of anti-rational maps and of reflection groups. 
These are complemented by a number of general results and illustrations of important technical tools, such as David surgery and straightening techniques. We also collect several analytic applications of the above theory. 
\end{abstract}

\maketitle

\setcounter{tocdepth}{1}
\tableofcontents

\section{Overview}\label{intro_sec}

\subsection*{The dictionary}
In his pioneering work on Fuchsian groups, Poincar{\'e} studied discrete groups of isometries of the hyperbolic plane generated by reflections \cite{Poi82}. Decades later, Coxeter \cite{Cox34} and Vinberg \cite{Vin67} studied discrete groups generated by reflections in much bigger generality.

In a seemingly unrelated world, Fatou and Julia laid the foundation of the theory of dynamics of holomorphic maps, particularly the dynamics of rational maps on the Riemann sphere, in the first quarter of the twentieth century \cite{fatou-1919,fatou-1920a,fatou-1920b,fatou-1926,julia-1918,julia-1922}. These developments drove Fatou to observe similarities between the dynamics of rational maps and that of Kleinian groups \cite[p. 22]{Fatou29}. In the 1980s, this philosophical analogy was set on a firm footing by Sullivan who introduced quasiconformal techniques in the study of rational dynamics and paved the way for the discovery of various deep connections between these two branches of conformal dynamics \cite{sullivan-dict}. Further contributions to this dictionary between Kleinian groups (and the associated theory of $3$-manifolds) and rational dynamics were subsequently made by McMullen-Sullivan \cite{MS98}\footnote{Refining a 1983 IHES preprint by Sullivan.}, McMullen \cite{McM95,McM98b}, Lyubich-Minsky \cite{LM97}, Pilgrim \cite{Pil03}, and others.

Inspired by the Fatou-Sullivan dictionary between Kleinian groups and complex dynamics, it is natural to think of iterations of antiholomorphic rational maps (anti-rational for short) on the Riemann sphere $\widehat{\C}$ as the complex dynamics counterpart of actions of Kleinian reflection groups.

In this survey, we will expound several recent results that advance the above theme in the Fatou-Sullivan dictionary. These results reveal certain explicit and somewhat surprising connections between the dynamics of anti-rational maps and Kleinian reflection groups. Moreover, these novel links between the two branches of conformal dynamics have given rise to a fresh class of conformal dynamical systems on the Riemann sphere generated by Schwarz reflection maps associated with quadrature domains. 

\subsection*{Schwarz reflection dynamics}
A domain in the complex plane with piecewise analytic boundary is called a \emph{quadrature domain} if the Schwarz reflection map with respect to its boundary extends anti-meromorphically to its interior. Such domains were first investigated by Davis \cite{Dav74}, and independently by Aharonov and Shapiro \cite{AS73, AS76, AS78}. Since then, quadrature domains have played an important role in various areas of complex analysis and fluid dynamics (see \cite{QD} and the references therein).

\begin{figure}[h!]
\captionsetup{width=0.96\linewidth}
\begin{center}
\includegraphics[width=0.36\linewidth]{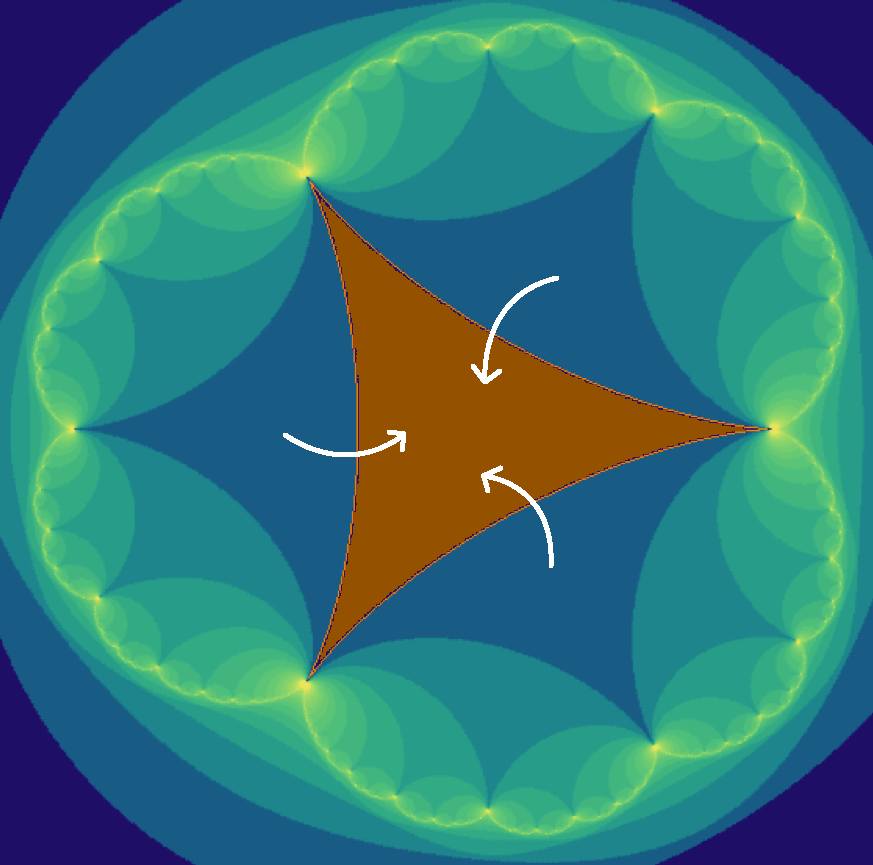}\qquad \includegraphics[width=0.36\linewidth]{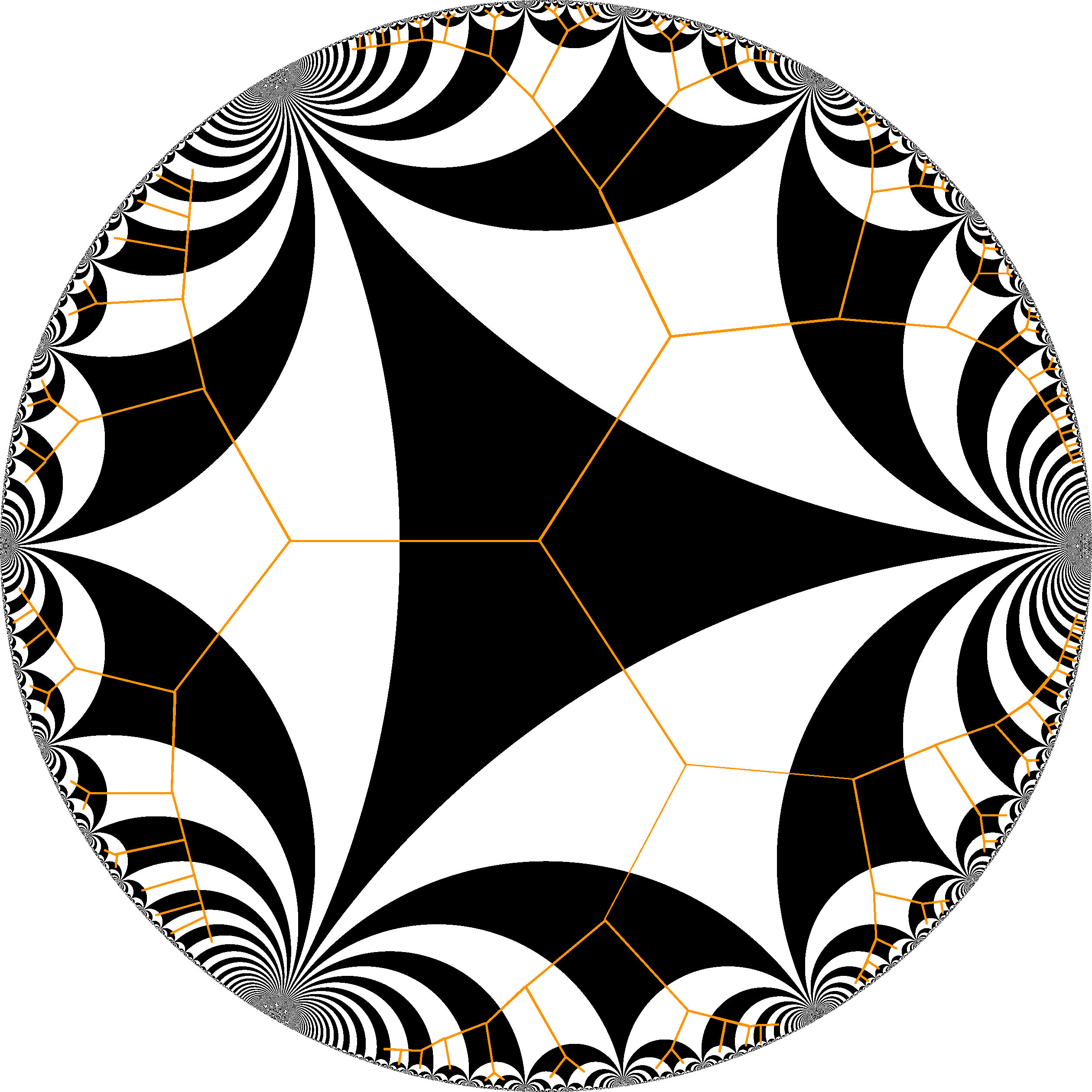}
\end{center}
\caption{Left: The dynamical plane of the Schwarz reflection map with respect to a deltoid curve. The interior (respectively, exterior) of the green Jordan curve is the escaping/tiling set (respectively, the non-escaping/filled Julia set), where the map behaves like the ideal triangle reflection group (respectively, like the quadratic map $\overline{z}^2$). Right: The tessellation of the unit disk under the ideal triangle reflection group and the corresponding Cayley tree are depicted.}
\label{deltoid_intro_fig}
\end{figure}

Iterations of Schwarz reflections was first studied by Seung-Yeop Lee and Nikolai Makarov in \cite{LM} to address some questions of interest in statistical physics concerning topology and singular points of quadrature domains. Subsequently, a systematic exploration of Schwarz reflection dynamics was launched in \cite{LLMM1,LLMM2,LLMM3,LMM1,LMM2}...,\footnote{To indicate interrelations among the papers that this survey is based on, we include the references to the original arXiv versions in the bibliography.} which demonstrated that Schwarz dynamics can combine features of dynamics of anti-rational maps and Kleinian reflection groups in a common dynamical plane. More precisely, the dynamical plane of a Schwarz reflection map admits an invariant partition into the \emph{escaping/tiling} and \emph{non-escaping/filled Julia} sets, which often parallel the action of a reflection group and of an anti-rational map, respectively. The simplest instance of this combination phenomenon (also called \emph{mating} in the dynamical setting) is displayed in Figure~\ref{deltoid_intro_fig}.
It also transpired through these studies that the parameter spaces of Schwarz reflections are intimately related to parameter loci of anti-rational maps and reflection~groups.

\subsection*{Antiholomorphic correspondences}
Another important role in this survey is played by antiholomorphic correspondences on the Riemann sphere; i.e., multi-valued maps on $\widehat{\C}$ with antiholomorphic local branches. The phenomenon of mating or combining quadratic rational maps with the modular group was discovered by Bullett and Penrose in the context of iterated holomorphic correspondences \cite{BP}, and was studied comprehensively by Bullett and Lomonaco in recent years \cite{BuLo1, BuLo2,BuLo3}.
It turns out that the study of Schwarz reflection dynamics can be used profitably to 
construct in a regular way antiholomorphic analogs of the Bullett-Penrose algebraic correspondences and to generalize them to arbitrary degree (where the modular group is replaced with anti-conformal analogs of Hecke groups, called \emph{anti-Hecke groups}). In \cite{LLMM3,LMM3,LLM23}, certain Schwarz reflection maps were constructed as hybrid dynamical systems and they were lifted to produce antiholomorphic correspondences whose dynamics combine anti-rational maps with the entire structure of anti-Hecke groups or ideal polygon reflection~groups.

\subsection*{Combination theorems}
The phenomena described above fits into the larger story of combination theorems, which has a long and rich history in groups, geometry and dynamics. Roughly speaking, the aim of a combination procedure is to take two compatible objects, and combine them to produce a richer and more general object that retains some of the essential features of the initial objects. Important examples of such constructions include the Klein Combination Theorem for two Kleinian groups \cite{Klein}, the Bers Simultaneous Uniformization Theorem that combines two surfaces (or equivalently, two Fuchsian groups) \cite{Bers60}, the Thurston Double Limit Theorem that allows one to combine two projective measured laminations (or equivalently, two groups on the boundary of the corresponding Teichm{\"u}ller space) \cite{Thu86,Otal98}, the Bestvina-Feighn Combination Theorem for Gromov-hyperbolic groups \cite{BF92}, etc. Douady and Hubbard transferred the notion of a combination theorem from the world of groups to that of holomorphic dynamics by designing the theory of \emph{polynomial mating} \cite{Dou83,Hub12}. Some of these classical combination theorems can be regarded as the underlying motivation and driving principles for the mating results that will be discussed in this survey.

While various explicit Schwarz reflection maps (and the associated correspondences) can be recognized as matings in the sense described above, the converse task of `interbreeding' reflection groups with anti-rational maps without a priori knowledge of the answer (i.e., constructing Schwarz reflection maps or correspondences as matings of given anti-rational maps and reflection groups) presents substantial technical challenges. The first obstruction comes from the inherent mismatch between invertible dynamical systems and non-invertible ones. This is circumvented by replacing a reflection group with its so-called \emph{Nielsen map}; i.e., a piecewise anti-M{\"o}bius non-invertible map that is orbit equivalent to the group (similar maps in the holomorphic setting are often called \emph{Bowen-Series} maps, cf. \cite{BS79}). It turns out that these Nielsen maps can often be topologically mated with antiholomorphic polynomials (anti-polynomials for short) along the lines of Douady--Hubbard mating of polynomials. The next difficulty lies in upgrading such topological hybrid dynamical systems to conformal ones. The lack of availability of Thurston-type realization theorems (which are used to construct rational maps as matings of two polynomials) for partially defined dynamical systems and the existence of parabolic elements in reflection groups cause serious impediments to the desired uniformization of topological matings. A general and flexible technique of David surgery (generalizations of quasiconformal surgery) was devised in \cite{LMMN} to surmount the above hurdles and to construct Schwarz reflections as combinations of large classes of anti-polynomials and Nielsen maps associated with reflection groups. It should be mentioned that such a surgery procedure first appeared in the work of Ha{\"i}ssinsky in the context of complex polynomials \cite{Hai98,Hai00}, and in the work of Bullett and Ha{\"i}ssinsky in the problem of mating quadratic polynomials with the modular group \cite{BH07}.

\begin{figure}[h!]
\captionsetup{width=0.96\linewidth}
\includegraphics[width=0.41\linewidth]{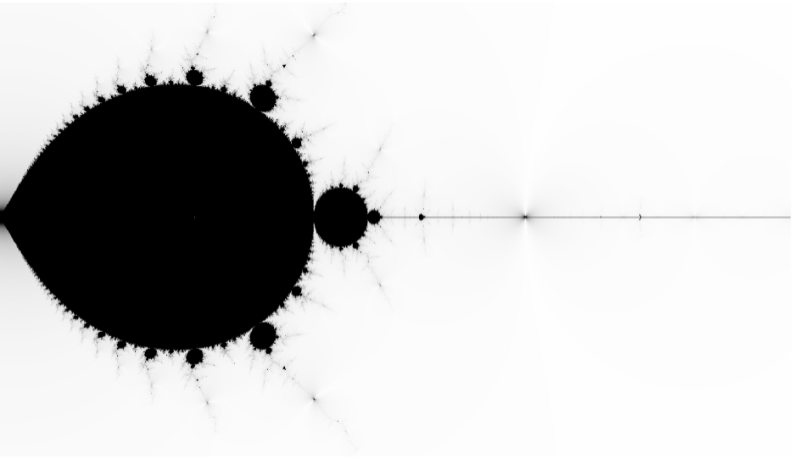}\ \includegraphics[width=0.5\linewidth]{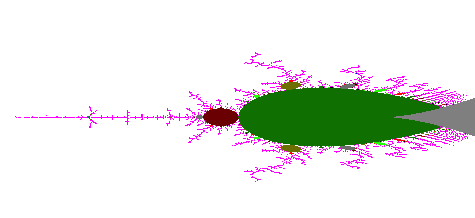}
\caption{Left: The connectedness locus of the Circle-and-Cardioid family of Schwarz reflections. Right: A part of the Tricorn.}
\label{c_and_c_para_fig}
\end{figure}

\subsection*{Parameter spaces of Schwarz reflections, anti-rational maps, and reflection groups}
The presence of common traits between the dynamics of anti-rational maps and Schwarz reflections manifests itself in the parameter spaces of Schwarz reflection maps as well. This was first observed numerically by Lee and Makarov. Their computer experiments showed that the connectedness locus of a family of quadratic Schwarz reflection maps (the so-called \emph{Circle-and-Cardioid} or \emph{C\&C} family) looks identical to (a part of) the \emph{Tricorn}, the connectedness locus of quadratic anti-polynomials (see Figure~\ref{c_and_c_para_fig}). While the appearance of `copies' of polynomial connectedness loci in parameter spaces of various holomorphic maps can be justified using the theory of polynomial-like maps and the associated straightening maps (cf. \cite{DH2,IK}), the situation here is more subtle as the Schwarz reflections under consideration only exhibit \emph{pinched/degenerate} anti-polynomial-like restrictions that cannot be straightened to anti-polynomials using quasiconformal surgery tools. To bypass this issue, a combinatorial straightening route was adopted in \cite{LLMM2} to relate the connectedness locus of the Circle-and-Cardioid family to (a part of) the Tricorn. 
Due to certain quasiconformal flexibility properties of antiholomorphic maps (associated with parabolic dynamics), this combinatorial straightening map only yields a homeomorphism between combinatorial models of the two connectedness loci.

The results of \cite{LLMM2} have been sharpened and generalized to arbitrary degree in a recent work where the space of `polygonal Schwarz reflections' of degree $d$ was studied \cite{LLM23} (for $d=2$, this space reduces to the C\&C family and deltoid-like Schwarz reflections)\footnote{Polygonal Schwarz reflections are defined in tree-like quadrature multi-domains, and their external maps come from ideal polygon reflection groups.}. A combination of combinatorial straightening techniques and puzzle machinery was used to construct a \emph{dynamically natural} homeomorphism between purely repelling combinatorial classes of degree $d$ polygonal Schwarz reflections and degree $d$ anti-polynomials with connected limit/Julia set. A dynamically natural bijection between geometrically finite maps in these two families was also established, and this bijection was shown to be generically continuous but discontinuous at some places. Here, ``dynamically natural'' means that for any anti-polynomial $f$, the corresponding polygonal Schwarz reflection is a conformal mating of $f$ with the regular ideal polygon reflection group.

There are several other parameter spaces of Schwarz reflection maps that are closely related to parameter spaces of anti-rational maps. A prototypical example of such families is the \emph{cubic Chebyshev} family of Schwarz reflections, which arises from univalent restrictions of the cubic Chebyshev polynomial to appropriate round disks. While these Schwarz reflections also fall outside the scope of usual polynomial-like straightening theory, their pinched anti-polynomial-like restrictions are somewhat more tame than the ones furnished by polygonal Schwarz reflections. In fact, a classical theorem of Warschawski (on the boundary behavior of conformal maps of topological strips) was used in \cite{LLMM3} to quasiconformally straighten these pinched anti-polynomial-like restrictions (of Schwarz reflections in the cubic Chebyshev family) to quadratic parabolic anti-rational maps. This result was instrumental in producing quadratic antiholomorphic correspondences as matings of quadratic parabolic anti-rational maps and an anti-conformal analogue of the modular group. It enabled us to define a straightening map from the cubic Chebyshev family of Schwarz reflections (or equivalently, from the resulting space of correspondences) to a family of quadratic parabolic anti-rational maps. 

The results mentioned in the previous paragraph have higher degree analogs too. A generalization of the cubic Chebyshev family of Schwarz reflections was introduced and studied in \cite{LMM3}. This family of Schwarz reflections arise from the space of degree $d+1$ polynomials ($d\geq 2$) that are injective on the closed disk and have a unique critical point on the unit circle (equivalently, from cardioid-like quadrature domains with a unique cusp on their boundaries). As in the $d=2$ case, a quasiconformal straightening surgery was designed for these Schwarz reflections, and the corresponding parameter space was shown to be a close cousin of a family of parabolic anti-rational maps of degree $d$. Once again, this straightening surgery played a fundamental role in the proof of existence of bi-degree $d$:$d$ antiholomorphic correspondences that are matings of degree $d$ parabolic anti-rational maps and anti-Hecke groups.

Special real two-dimensional slices in the above space of Schwarz reflections are obtained when the underlying polynomials (that are injective on the disk) are \emph{Shabat polynomials}; i.e., they have two critical values in the plane. These can be seen as higher degree  one-parameter generalizations of the cubic Chebyshev family. The connectedness loci of such families of Schwarz reflections are \emph{combinatorially equivalent} to certain parameter spaces of \emph{Belyi} parabolic anti-rational maps (a Belyi anti-rational map is an anti-rational map with at most three critical values) \cite{LMM4}.

On the group side, let us recall that the index two Fuchsian subgroup of an ideal polygon reflection group uniformizes a punctured sphere. In general, such surfaces have moduli and pinching suitable closed geodesics on these surfaces allows one to study the deformations/degenerations of reflection groups. It turns out that certain classes of Schwarz reflection maps are also amenable to similar deformation techniques. Such deformations were used in \cite{LMM2} to produce a dynamically natural homeomorphism between a space of Schwarz reflections and the closure of the Bers slice of the ideal polygon reflection group. 

\subsection*{Interplay between the holomorphic and antiholomorphic cases}

In the 1990s, Shaun Bullett and Christopher Penrose \cite{BP} discovered that some quadratic polynomials and the modular group can co-exist in the same dynamical plane for a bi-degree $2$:$2$ algebraic correspondence, so this correspondence can be viewed as the {\em mating} of a quadratic map and the modular group. They conjectured that actually any quadratic polynomial can be mated with the modular group in this way and that the parameter space of the relevant correspondences is naturally homeomorphic to the Mandelbrot set. Important progress towards the resolution of this  conjecture was made by Bullett and Ha{\"i}ssinsky, who showed in \cite{BH07} that this space of correspondences contains matings of the modular group and a large class of quadratic polynomials with connected Julia set (including Collet-Eckmann maps). But the full conjecture remained out of reach for another decade.

In the early 2010s, Luna Lomonaco introduced in her thesis the class of {\em parabolic-like maps}, a parabolic version of polynomial-like maps.
Such a map is defined on a domain containing a parabolic point $\alpha$, but it lacks a polynomial-like quality near $\alpha$.
She proved that any parabolic-like map can be quasiconformally straightened to a parabolic rational map \cite{Lom15}. It was then suggested by Adam Epstein ``that parabolic-like mappings might be applied to the Bullett-Penrose family of correspondences'' \cite[p. 209]{BuLo1}. This idea was pursued by Bullett and Lomonaco:
they first showed that a certain pinched polynomial-like restriction of a Bullett-Penrose correspondence can be extended to a quasi-regular parabolic-like map, then straightened it to a quadratic parabolic rational map, and then moved on to conclude that the algebraic correspondence in question is the mating of that rational map with the modular group  \cite{BuLo1}.

The idea of a parabolic straightening was adapted in \cite{LLMM3} for the antiholomorphic setting to show that the Chebyshev family of Schwarz reflections is naturally bijectively equivalent to the \emph{parabolic Tricorn} (which is the connectedness locus of quadratic parabolic anti-rational maps). First we perform a surgery that replaces the Blaschke external map of a parabolic rational map with a ``Farey external map''\footnote{where all the maps in question are anti-holomorphic}
yielding a Schwarz reflection dynamics  (see Subsection~\ref{chebyshev_gen_subsec}).  This is a ``pinched version'' of  the classical straightening theorem.
Then we lift the Schwarz dynamics to a correspondence dynamics by means of the uniformizing Chebyshev polynomial.

In the course of this development (at the first step), a new version of the surgery machinery was designed for straightening pinched polynomial-like maps that may not necessarily admit a holomorphic extension around the parabolic point (see \cite[\S 5]{LLMM3} and Subsections~\ref{chebyshev_subsec} and~\ref{chebyshev_gen_subsec}). (This machinery makes use of classical  Warschawski's Theorem on the uniformization of topological strips, though on various special occasions it can be replaced with a hand-crafted construction.) It avoids an intermediate quasi-regular extension of the map yielding a good control of the dependence of the straightening on the parameters.
Altogether, it then directly implies a natural bijection between the parameter space of the correspondences in question and the parabolic Tricorn.

In a parallel development, Bullett and Lomonaco completed their proof of the fact that the connectedness locus of the space of Bullett-Penrose algebraic correspondences is homeomorphic to the parabolic Mandelbrot set \cite{BuLo3}. Finally, due to the work of Petersen and Roesch that appeared meanwhile \cite{PR21}, the parabolic Mandelbrot set turns out to be homeomorphic to the genuine Mandelbrot set. Altogether, it confirmed the Bullett-Penrose Conjecture.

Going beyond degree two, in \cite{BuFr} Bullett and Freiberger put forward a general Mating Conjecture that all polynomials with connected Julia set can be mated with Hecke groups as correspondences. By means of the above surgery machinery (based on Warschawski's Theorem), a parabolic version of this conjecture was established in the antiholomorphic setting (where anti-polynomials are replaced with parabolic anti-rational maps) in \cite{LMM3} (see Section~\ref{general_mating_corr_sec}).

To conclude this subsection, let us mention some recent combination theorems in the holomorphic world which were inspired by the study of Schwarz reflection dynamics and its connection with correspondences. In \cite{MM1}, Bowen-Series maps of Fuchsian punctured sphere groups were used as natural holomorphic analogues of Nielsen maps (these Bowen-Series maps are  expansive covering maps of the circle), and combination theorems for such Bowen-Series maps and hyperbolic complex polynomials were proved. 

To allow more general genus zero orbifold groups to be incorporated in the mating framework, Bowen-Series maps of punctured spheres were generalized to the class of \emph{factor Bowen-Series maps} in \cite{MM2}. These maps, which are associated with a large collection of genus zero orbifolds containing punctured spheres as well as Hecke surfaces as special cases, were also shown to be conformally mateable with hyperbolic complex polynomials. When the complex polynomials are \emph{real-symmetric}, the corresponding conformal matings turn out to be complex conjugates of Schwarz reflection maps. This fact was used in \cite{MM2} to produce, on the one hand, various new examples of holomorphic correspondences that are matings of genus zero orbifolds and complex polynomials, and on the other hand, complex-analytic embeddings of Teichm{\"u}ller spaces of genus zero orbifolds into the space of algebraic correspondences. 

Subsequently in \cite{LLM23}, a complete characterization of the conformal matings of factor Bowen-Series maps and complex polynomials were given by introducing a holomorphic analog of Schwarz reflection maps. These partially defined maps, which act as involutions on the boundary of their domain of definition (as opposed to the trivial action of Schwarz reflections on the boundaries of quadrature domains), were termed as \emph{B-involutions}. A detailed study of parameter spaces of B-involutions was carried out in \cite{LLM23}, and this was leveraged to construct algebraic correspondences as matings of generic complex polynomials in the connectedness loci (geometrically finite maps, and periodically repelling finitely renormalizable maps) with copious genus zero orbifolds (including punctured spheres and Hecke surfaces).

Finally, the parabolic version of the Bullett-Freiberger conjecture has been recently resolved in the holomorphic setting in a joint work of the authors with Shaun Bullett and Luna Lomonaco \cite{BLLM}. 
Some of the key ideas employed in the proof of this result have roots in the classical Douady-Hubbard theory of polynomial-like maps where such maps are realized as matings of hybrid and external classes (or external maps) \cite{DH2}, \cite[\S 3]{Lyu99}. 
To handle the existence of parabolic points for Hecke groups and rational maps, we extend the theory of polynomial-like maps to \emph{pinched polynomial-like maps} with parabolic external maps (such that the domain and the range are allowed to touch).
We apply this theory to construct pinched polynomial-like maps as matings of hybrid classes arising from parabolic rational maps and certain external maps associated with the Hecke group. 
%These external maps, which are called the \emph{Farey map}\footnote{The Farey map is the factor Bowen-Series map associated with the Hecke group.} and the \emph{Hecke map}, are conformally equivalent; and hence the resulting matings are conformally conjugate. 
Finally, these pinched polynomial-like maps are lifted to the desired algebraic correspondences on appropriate branched covering spaces. 
%This is done in a `geometric' way that uses basic covering theory for Riemann surfaces and an `analytic' way based upon quasiconformal surgery. 

\subsection*{When Julia sets look like Kleinian limit sets}
\begin{figure}[h!]
\captionsetup{width=0.96\linewidth}
\includegraphics[width=0.32\linewidth]{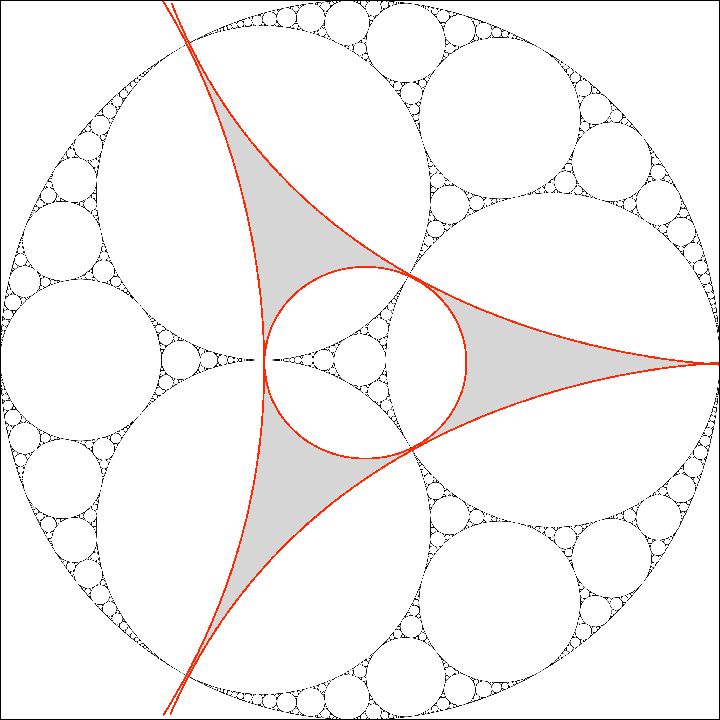}\ \hspace{4mm}\  \includegraphics[width=0.324\linewidth]{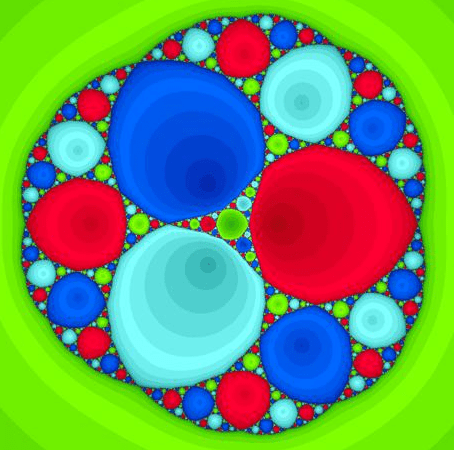}  
\caption{Two homeomorphic fractals: the classical Apollonian gasket and the Julia set of a cubic anti-rational map.}
\label{limit_julia_intro_fig}
\end{figure}

Another major theme of this survey concerns explicit dynamical relations between limit sets of Kleinian reflection groups and Julia sets of anti-rational maps. While many topological, analytic and measure-theoretic similarities between Kleinian limit sets and rational Julia sets have been known for long, no example of dynamically natural homeomorphism between such fractals was known until recently, to the best of our knowledge. The first non-trivial example of such a homeomorphism was produced in \cite{LLMM4} by turning the Nielsen map of the classical Apollonian gasket reflection group into a critically fixed anti-rational map (see Figure~\ref{limit_julia_intro_fig}). The main idea of this construction was to use a certain compatibility property between Nielsen maps and power maps to cook up a topological branched cover from the Nielsen map of the Apollonian group and then invoke the Thurston Realization Theorem to obtain the desired anti-rational map. This recipe was generalized in \cite{LLM1}, where dynamically natural homeomorphisms between Julia sets of critically fixed anti-rational maps and limit sets of Kleinian reflection groups (arising from finite circle packings) were manufactured. This dictionary has various dynamical consequences; for instance, the geodesic lamination models of the corresponding limit and Julia sets can be explicitly related, and the resulting bijection between critically fixed anti-rational maps and Kleinian reflection groups commutes with the operation of mating in the respective categories. From an analytic point of view, it is worth mentioning that these fractals are usually\footnote{It was recently shown that these sets can actually be quasiconformally equivalent; cf. \cite{LMM26}} not quasiconformally equivalent \cite{LZ23b} (although they may have isomorphic quasisymmetry groups, see \cite{LLMM4}), but these Julia sets can be mapped onto the corresponding limit sets by global David homeomorphisms \cite{LMMN}.

The above correspondence between critically fixed anti-rational maps and Kleinian reflection groups (arising from finite circle packings) also has fundamental parameter space implications, which were investigated in \cite{LLM2} using rescaling limit techniques. This revealed striking similarities between the parameter spaces of anti-rational maps and reflection groups. In particular, it was shown that the quasiconformal deformation space of a kissing reflection group is bounded if and only if a suitable deformation space of the corresponding critically fixed anti-rational map is bounded (which can be seen as an analogue of Thurston's Boundedness Theorem in the context of anti-rational maps), and that the bifurcation structures of these deformation spaces have the same combinatorial patterns. Further, it was demonstrated that the union of suitable deformation spaces of critically fixed anti-rational maps admits a monodromy representation onto the mapping class group of a punctured sphere, which is in harmony with a result of Hatcher and Thurston on the global topological complexity of parameter spaces of reflection groups~\cite{HT}.

\subsection*{Applications to analytic problems}
The development of iteration theory of Schwarz reflection maps has interesting consequences to certain questions of purely analytic origin. In fact, the characterization of (simply connected) quadrature domains as univalent images of the disk under rational maps gives abundant examples of Schwarz reflection maps, and connects the study of Schwarz reflection dynamics to the classical theory of univalent functions in geometric function theory. The intimate links between Schwarz reflections, quadrature domains, and univalent rational maps have been utilized to study the topology of quadrature domains and answer related questions with statistical physics motivation, to solve extremal problems for suitable classes of univalent maps, and to study domains of univalence of complex polynomials. 
The crux of the matter is to translate the above analytic problems to questions regarding the dynamics of naturally associated Schwarz reflection maps, and then apply the dynamical theory of Schwarz reflections to obtain the desired solutions \cite{LM,LM1,LMM1,LMM2,LMM4}. 

In the same vein, David surgery tools developed in the mating theory described above also have several analytic applications. Indeed, David extension theorems and associated surgery techniques were used in \cite{LMMN} to prove, among other things, that many naturally arising non-quasisymmetric circle homeomorphisms are still welding maps. It is used to show that some interesting Julia sets (called \emph{pine trees}, see Figure~\ref{pine_tree_fig}) and limit sets of necklace reflection groups are conformally removable (the latter was a well-known open problem among experts working in the area).

\subsection*{Structure of the survey}

We begin the survey (Section~\ref{interplay_sec}) with a quick and somewhat informal introduction to the main mathematical objects, their interconnections, and how this interplay leads to various new results in the antiholomorphic chapter of the Fatou-Sullivan dictionary. Section~\ref{antiholo_background_sec} covers the necessary preliminaries: here we discuss various elementary properties of Kleinian reflection groups, anti-rational dynamics, and Schwarz reflection maps. Sections~\ref{quadratic_examples_sec} and~\ref{cubic_examples_sec} illustrate various features of Schwarz reflection dynamics, new mating phenomena, straightening techniques, topological and analytic connections between limit and Julia sets, etc. through a number of concrete examples. Section~\ref{schwarz_para_space_sec} expounds the parameter space structure of some special families of Schwarz reflections and their relations with appropriate spaces of anti-rational maps and reflection groups. More general straightening and mating results require recently developed David surgery tools, which are discussed in Section~\ref{david_surgery_sec}. Section~\ref{new_line_dict_sec} explicates dynamically natural homeomorphisms between Julia sets of anti-rational maps and limit sets of reflection groups, and parameter space consequences of this connection. The next two Sections,~\ref{mating_anti_poly_nielsen_sec} and~\ref{mating_para_space_sec}, describe a general mating theory for Nielsen maps of reflection groups and anti-polynomials, and relate the parameter spaces of the associated Schwarz reflections to connectedness loci of anti-polynomials. Section~\ref{general_mating_corr_sec} is devoted to the construction of antiholomorphic generalizations of Bullett-Penrose correspondences and the underlying straightening surgery. Finally, some of the analytic applications of the theory are recorded in Section~\ref{anal_app_sec}.

\medskip

\section{The main characters, their interplay, and some applications at a glance}\label{interplay_sec}

\subsection{Four models for external dynamics}

To unify some of the principal players of this survey, we will refer to the dynamics of an anti-rational map on a (marked) completely invariant Fatou component and to the dynamics of a Schwarz reflection map on its escaping set as their \emph{external dynamics}. More generally, this term applies to the dynamics of anti-polynomial-like maps (or their degenerate analogs, often called pinched anti-polynomial-like maps) on their escaping sets. In accordance with the classical theory of polynomial-like maps (cf. \cite{DH2}), such external dynamics can be modeled by appropriate piecewise anti-analytic (i.e., real-analytic and orientation-reversing) covering maps of the circle, called \emph{external maps}.

Many of the families of anti-rational maps and Schwarz reflections that we will be concerned with have fixed external dynamics. Relations between the corresponding external maps (i.e., piecewise anti-analytic circle coverings) lie at the core of all fundamental connections between reflection groups and anti-rational maps.

\vspace{1mm}

\noindent \textbf{i) Power map.} The first and the most well-known of them is the power map $\overline{z}^d$. Monic, centered anti-polynomials of degree $d$ with connected Julia set have $\overline{z}^d$ as the conformal model of their external dynamics. We denote the connectedness locus of monic, centered degree $d$ anti-polynomials by $\mathscr{C}_d$.
\smallskip

\noindent \textbf{ii) Parabolic anti-Blaschke product.} The unicritical antiholomorphic Blaschke (anti-Blaschke for brevity) product 
$$
B_d~=~\frac{(d+1)\overline z^d + (d-1)}{(d-1)\overline z^d + (d+1)}
$$ 
is topologically conjugate to $\overline{z}^d$ on the unit circle $\mathbb{S}^1$; however, unlike the expanding endomorphism $\overline{z}^d$, the map $B_d$ has a parabolic fixed point on the circle. The space of anti-rational maps admitting $B_d$ as their external dynamics can be thought of as the parabolic counterpart of the connectedness locus of degree $d$ anti-polynomials. This space is denoted by $\pmb{\mathcal{B}}_d$ (see Subsection~\ref{para_anti_rat_gen_subsubsec}).
\smallskip

\begin{figure}[h!]
\captionsetup{width=0.96\linewidth}
\begin{center}
\includegraphics[width=0.36\linewidth]{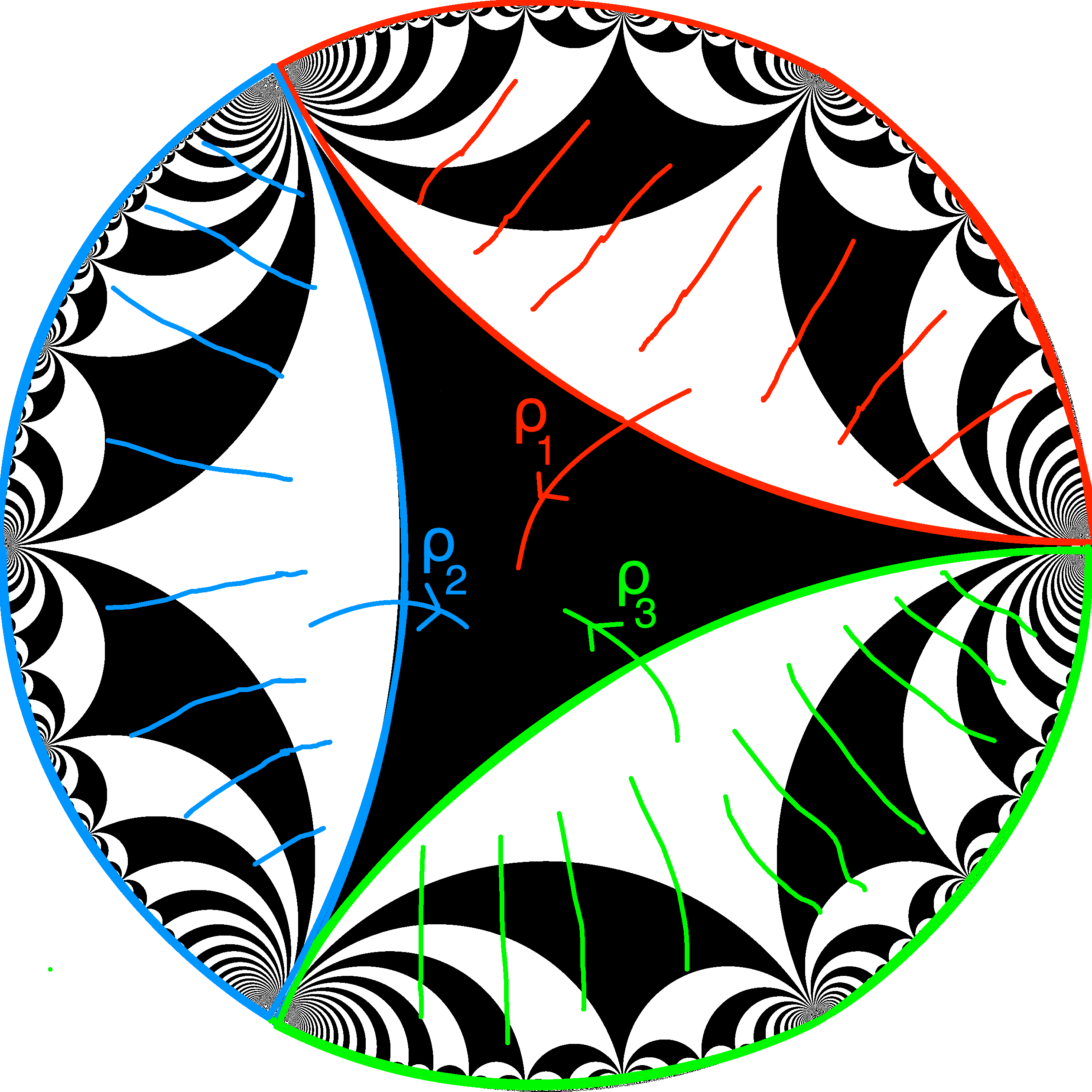}
\end{center}
\caption{The ideal triangle reflection group is generated by the reflections $\rho_1,\rho_2,\rho_3$ in the sides of an ideal hyperbolic triangle. These generators define the corresponding Nielsen map $\pmb{\cN}_2$ as a piecewise anti-M{\"o}bius map on the shaded regions.}
\label{itg_nielsen_fig}
\end{figure}

\noindent \textbf{iii) Nielsen map.} An analog of the map $\overline{z}^d$ in the reflection group world is given by the \emph{Nielsen map} $\pmb{\cN}_d$ associated with the group $\pmb{G}_d$ generated by reflections in the sides of a regular ideal $(d+1)-$gon in the hyperbolic plane. It is a piecewise anti-M{\"o}bius map topologically conjugate to $\overline{z}^d$ on $\mathbb{S}^1$ (see Figure~\ref{itg_nielsen_fig} and Subsection~\ref{nielsen_map_subsubsec}). The map $\pmb{\cN}_d$ has $d+1$ parabolic fixed points on $\mathbb{S}^1$ at the ideal vertices of the regular ideal $(d+1)-$gon, so the above conjugacy is not quasisymmetric. The collection $\cS_{\pmb{\cN}_d}$ of antiholomorphic maps with $\pmb{\cN}_d$ as their external map coincides with the connectedness locus of \emph{regular polygonal} Schwarz reflections; i.e., a certain class of degree~$d$ piecewise Schwarz reflection maps associated with tree-like quadrature multi-domains (see Subsection~\ref{c_and_c_general_subsec} and Section~\ref{mating_para_space_sec}).

\noindent \textbf{iv) Anti-Farey map.} Yet another external class arising from reflection groups is obtained as a factor of $\pmb{\cN}_d$. In fact, since $\pmb{\cN}_d$ commutes with rotation by $2\pi/(d+1)$, it descends to a piecewise anti-analytic, degree $d$, orientation-reversing covering map $\pmb{\cF}_d:\mathbb{S}^1\to\mathbb{S}^1$ with a unique parabolic fixed point on the circle (see Subsections~\ref{nielsen_first_return_external_map_subsubsec},~\ref{anti_farey_subsubsec}). The map $\pmb{\cF}_d$ is also topologically conjugate to the above three maps; but more importantly, it is \emph{quasisymmetrically} conjugate to the parabolic anti-Blaschke product $B_d$. Although the external map $\pmb{\cF}_d$ is obtained as a factor of the piecewise M{\"o}bius Nielsen map $\pmb{\cN}_d$, it has a fully ramified critical point. Antiholomorphic maps having $\pmb{\cF}_d$ as the external class can be described as the connectedness locus of Schwarz reflections associated with cardioid-like quadrature domains with a unique critical point (of local degree $d+1$) that escapes in one iterate. We denote this space of Schwarz reflections by $\cS_{\pmb{\cF}_d}$ (see Subsection~\ref{upgrade_step_subsubsec}).

\subsection{Relation between $\overline{z}^d$ and $\pmb{\cN}_d$, and its implications to the dictionary}

As mentioned before, the Nielsen map $\pmb{\cN}_{d}$ is topologically conjugate to the power map $\overline{z}^d$ on $\mathbb{S}^1$ via a circle homeomorphism $\pmb{\mathcal{E}}_d$. Due to its close relation with the classical \emph{Minkowski question mark function}, we call the conjugating map $\pmb{\mathcal{E}}_d$ the \emph{$d$-th Minkowski circle homeomorphism} (see Subsection~\ref{question_mark_subsubsec}). This topological compatibility between $\overline{z}^d$ and $\pmb{\cN}_d$ serves as a bridge between Kleinian reflection groups and anti-rational maps. Specifically, the existence of this map is instrumental in the addition of the following entries in the Fatou-Sullivan dictionary.
\smallskip

\noindent \textbf{i) Mating anti-polynomials with necklace reflection groups, and parameter spaces of Schwarz reflections.} There exists a large class of Schwarz reflection maps which are matings of degree $d$ anti-polynomials with connected Julia set and necklace reflection groups (i.e., Kleinian reflection groups in the closure of the Bers slice of $\pmb{G}_d$, see Subsection~\ref{necklace_subsubsec}). Various instances of this mating phenomena are explicated in Subsections~\ref{deltoid_subsec},~\ref{c_and_c_center_subsec},~\ref{talbot_subsec},~\ref{c_and_c_general_subsec},~\ref{sigma_d_subsec}, and general existence results are collected in Section~\ref{mating_anti_poly_nielsen_sec}. In its simplest avatar, the above combination results can be thought of as a fusion of the Bers Simultaneous Uniformization Theorem for Fuchsian groups and simultaneous uniformization of a pair of Blaschke products. More sophisticated versions of matings of anti-polynomials and reflection groups run along the lines of the Douady-Hubbard mating theory for polynomials and the Thurston Double Limit Theorem.

The parameter spaces of these families of Schwarz reflections bear strong resemblance with the connectedness locus $\mathscr{C}_d$ and the Bers slice closure of the reflection group $\pmb{G}_d$ (see Subsections~\ref{c_and_c_general_subsec},~\ref{sigma_d_subsec} and Section~\ref{mating_para_space_sec}).
\smallskip

\noindent \textbf{ii) Equivariant homeomorphisms between Julia and limit sets, and deformation space analogies.} Limit sets of Kleinian reflection groups arising from circle packings (including the classical Apollonian gasket) are homeomorphic to Julia sets of critically fixed anti-rational maps in a dynamically natural fashion (see Subsections~\ref{apollo_group_map_schwarz_subsec},~\ref{sigma_d_subsec} for examples of this connection and Subsection~\ref{new_line_dict_dyn_subsec} for a general result). 

Moreover, the bifurcation structure, boundedness properties, and global topologies of the deformation spaces of kissing reflection groups and critically fixed anti-rational maps have stark similarities (see Subsection~\ref{new_line_dict_para_subsec}).

\subsection{David surgery as a key technical tool}

A key technical ingredient in the proof of simultaneous uniformization of Blaschke products (and in many other important surgery techniques in holomorphic dynamics) is the Ahlfors-Beurling Extension Theorem, which states that a quasisymmetric homeomorphism of the circle extends continuously to a quasiconformal homeomorphism of the disk. Since the Minkowski circle homeomorphism $\ciq$ conjugates parabolic dynamics to the  hyperbolic one, it is not quasisymmetric. 
However, it was shown in \cite{LLMM4} that its inverse admits a David extension to the disk.  Indeed,  by a direct number theoretic analysis,
It was verified that $\ciq$ satisfies Chen-Chen-He's \cite{CCH96} and Zakeri's \cite{Zak04} distortion property that is sufficient for the David extendability. This result was then extended, by dynamical means,   to the general class of circle homeomorphisms (and their local counterparts) conjugating hyperbolic dynamics to the parabolic one. It laid down a foundation for a {\em general David surgery machinery} that facilitates construction of parabolic conformal dynamical systems from hyperbolic  ones.
This novel machinery is elaborated in Section~\ref{david_surgery_sec} 
(see Section~\ref{qc_david_appendix} for background on David homeomorphisms). Not only is it indispensable for matings of anti-polynomials with necklace groups, but it also yields a direct passage from subhyperbolic (anti-)rational maps to geometrically finite rational maps (generalizing the work of Ha{\"i}ssinsky from the 1990s) and to kissing reflection groups, thereby shedding new light on the analytic geometry of such Julia and limit sets (see Subsection~\ref{conf_removable_subsec}).

\subsection{Quasisymmetric compatibility of $B_d$ and $\pmb{\cF}_d$, and antiholomorphic paradigm for Bullett-Penrose correspondences}

The existence of a quasisymmetric conjugacy between the parabolic anti-Blaschke product $B_d$ and the anti-Farey map $\pmb{\cF}_d$ (on the circle) enables one to mate the filled Julia dynamics of maps in $\pmb{\mathcal{B}}_d$ with the external map $\pmb{\cF}_d$ and realize the matings as Schwarz reflection maps in $\cS_{\pmb{\cF}_d}$. The fact that the anti-Farey map $\pmb{\cF}_d$ has a fully ramified critical point implies that the corresponding quadrature domains are uniformized by degree $d+1$ polynomials. These uniformizing polynomials can be used to lift Schwarz reflections in $\cS_{\pmb{\cF}_d}$ to construct antiholomorphic correspondences of bi-degree $d$:$d$ that are matings of parabolic anti-rational maps in $\pmb{\mathcal{B}}_d$ and an anti-conformal version of the classical Hecke group. This gives a regular framework for producing correspondences which are antiholomorphic counterparts (for arbitrary degree) of the Bullett-Penrose degree two holomorphic correspondences (see Subsections~\ref{chebyshev_subsec},~\ref{chebyshev_gen_subsec} for $d=2$ examples and Section~\ref{general_mating_corr_sec} for the general realization theorem).

\subsection{The emergence of pinched polynomial-like maps and novel straightening techniques}

While polynomial-like maps enjoy a central role in holomorphic dynamics, there are certain situations where one encounters degenerate or pinched versions of polynomial-like maps. Roughly speaking, a degenerate polynomial-like map is a proper holomorphic map from a pinched topological disk onto a larger pinched topological disk with finitely many touching points between the domain and co-domain. Such objects naturally appear in the study of maps with parabolic external dynamics (cf. \cite{Lom15,PR21}).

The study of Schwarz reflection maps (especially those that combine the actions of anti-polynomials and necklace groups) brings degenerate polynomial-like maps to the fore. In fact, all maps in $\cS_{\pmb{\cN}_d}$ and $\cS_{\pmb{\cF}_d}$ (i.e., with Nielsen and anti-Farey external maps) admit degenerate polynomial-like structures that cannot be upgraded to actual polynomial-like maps (see Subsections~\ref{deltoid_degeneration_subsubsec},~\ref{c_and_c_basilica_pinched_anti_quad_subsubsec},~\ref{chebyshev_center_hybrid_conj_subsubsec} for examples). Straightening degenerate polynomial-like maps to rational maps encounters various subtleties. In fact, for degenerate anti-polynomials-like restrictions arising from maps in $\cS_{\pmb{\cN}_d}$, there are analytic obstructions to quasiconformal straightening (as $\pmb{\cN}_d$ has more than one parabolic fixed points) which makes the study of this space more difficult. As mentioned in the previous section, this compels one to apply combinatorial techniques and puzzle machinery to relate spaces of polygonal Schwarz reflections to connectedness loci of anti-polynomials (see Subsection~\ref{c_and_c_general_subsec} and Section~\ref{mating_para_space_sec}). On the other hand, the pinched anti-polynomial-like restrictions of maps in $\cS_{\pmb{\cF}_d}$ are amenable to quasiconformal straightening (since their external dynamics have a unique parabolic fixed point with controlled geometry), which allows one to relate the parameter space of $\cS_{\pmb{\cF}_d}$ to the parameter space of the parabolic anti-rational family $\pmb{\mathcal{B}}_d$ (see Subsections~\ref{cubic_cheby_qc_straightening_subsubsec},~\ref{mating_regularity_subsec}).

\subsection{General conjectures and questions}

Let us briefly mention some open questions regarding the general structure of parameter spaces of Schwarz reflection maps and antiholomorphic correspondences.

\subsubsection{Product structure in spaces of Schwarz reflections and combinatorial rigidity}

As mentioned before, the dynamical plane of a piecewise Schwarz reflection map can be decomposed into two invariant subsets: the non-escaping set and the escaping/tiling set. In the mating locus, Schwarz reflections behave like anti-rational maps on their non-escaping sets and exhibit features of necklace reflection groups on their tiling sets.
Thus, freezing the dynamics of a Schwarz reflection map on its non-escaping set, and deforming its
dynamics on the tiling set should give rise to a `copy' of the Teichm{\"u}ller space of a necklace group in the Schwarz parameter space. On the other hand, fixing the conformal class of the dynamics on the tiling set, and changing it on the non-escaping set should produce a `copy' of an anti-rational parameter space in the Schwarz parameter space.
The above heuristics suggest that the parameter space of appropriate families of Schwarz reflection maps should have a local product structure; i.e. locally they should be products of anti-rational parameter spaces and Teichm{\"u}ller spaces of reflection groups. While the existence of such \emph{Bers slices} has been justified in several special cases (see Subsections~\ref{c_and_c_general_subsec},~\ref{sigma_d_subsec} and Section~\ref{mating_para_space_sec}), the general picture demands further investigation. This would require a good understanding of puzzle structures and combinatorial rigidity properties of Schwarz reflection maps.

\subsubsection{Degenerations, and Double Limit Theorems}
In the spirit of degenerations of rational maps and Kleinian groups, it is natural to study degenerations of Schwarz reflections that arise as conformal matings of anti-polynomials and necklace groups. This is particularly interesting when the anti-polynomials tend to the boundaries hyperbolic components and the groups go to the boundary of their quasiconformal deformation spaces. 
In certain cases, this produces a `phase transition'; i.e., such a degenerating sequence of matings converges to a limiting piecewise Schwarz reflection dynamical system, but
at least one of the quadrature domains gets pinched into a disjoint collection of quadrature domains (or equivalently, the Carath{\'e}odory limits of some of the uniformizing rational maps of the associated quadrature domains undergo degree drop).  This is an entirely new degeneration phenomenon in conformal dynamics that deserves to be better understood.

In general, one conjectures that there should be analogues of the Thurston Double Limit Theorem for this setup which would describe the dynamics of the limiting map as a quotient of the dynamics of the degenerating sequence of conformal matings.

\subsubsection{Discreteness locus and mating locus in the space of correspondences}

The antiholomorphic correspondences that arise as matings of anti-rational maps and anti-Hecke groups sit inside a larger space of correspondences generated by deck transformations of polynomials and a circular reflection (see Section~\ref{general_mating_corr_sec}). It would be quite interesting to understand the dynamics of these more general correspondences. For instance, in the spirit of J{\"o}rgensen's inequality for Kleinian groups, one can ask when such correspondences exhibit suitable discreteness properties (e.g., when do they act discretely on some part of the sphere?). It is also natural to ask for an intrinsic characterization of the \emph{mating locus} in this bigger space of correspondences. Satisfactory answers to these questions would involve exploring uncharted territories.

\section{Background on antiholomorphic dynamics}\label{antiholo_background_sec}

\subsection{Kleinian reflection groups}\label{ref_group_subsec}

\subsubsection{Circle packings}\label{circle_pack_subsubsec}
An oriented circle in $\widehat{\C}$ encloses an open round disk, which we will refer to as the interior of the circle.
A {\em circle packing} $\mathcal{P}$ is a connected finite collection of at least three oriented circles in $\widehat{\C}$ with disjoint interiors.
The combinatorial configuration of a circle packing can be encoded by its {\em contact graph} $\Gamma$, which has a vertex associated with each circle, and an edge connecting two vertices if and only if the two associated circles touch.
The embedding of $\mathcal{P}$ in $\widehat\C$ endows its contact graph with a plane structure (equivalently, a cyclic order of the edges meeting at a vertex). Clearly, the contact graph of a circle packing is simple. In fact, this is the only constraint on the graph (See \cite[Chapter~13]{Thurston78}).

\begin{theorem}[Circle Packing Theorem]\label{circle_packing_thm}
Every connected, simple, plane graph  is isomorphic to the contact graph of some circle packing.
\end{theorem}

\begin{definition}\label{graph_term_def}
Let $\Gamma$ be a finite connected graph.
\noindent\begin{enumerate}\upshape
\item $\Gamma$ is said to be {\em $k$-connected} if $\Gamma$ contains more than $k$ vertices and remains connected if any $k-1$ vertices and their corresponding incident edges are removed.

\item $\Gamma$ is called {\em polyhedral} if $\Gamma$ is isomorphic to the $1$-skeleton of a convex polyhedron. Equivalently, $\Gamma$ is {\em polyhedral} if it is planar and $3$-connected.

\item $\Gamma$ is said to be {\em outerplanar} if it has a planar drawing for which all vertices lie on the boundary of some face.

\item $\Gamma$ is called {\em Hamiltonian} if there exists a Hamiltonian cycle, i.e., a closed path in $\Gamma$ visiting each of its vertices exactly once. 
\end{enumerate}
\end{definition}
\noindent According to a theorem of Steinitz, a graph is polyhedral if and only if it is $3$-connected and planar.
Given a polyhedral graph, we have a stronger version of the Circle Packing Theorem \cite{Sch92}.

\begin{theorem}[Circle Packing Theorem for polyhedral graphs]\label{thm:gcpt}
Suppose $\Gamma$ is a polyhedral graph, then there is a pair of circle packings whose contact graphs are isomorphic to $\Gamma$ and its planar dual.
Moreover, the two circle packings intersect orthogonally at their points of tangency.
\smallskip

This pair of circle packings is unique up to M\"obius transformations.
\end{theorem}

\subsubsection{Kissing reflection groups}\label{kissing_group_subsubsec}

Let $\Gamma$ be a connected simple plane graph.
By the Circle Packing Theorem, $\Gamma$ is (isomorphic to) the contact graph of some circle packing 
$$
\mathcal{P}=\{C_1,..., C_{d+1}\}.
$$
We define the {\em kissing reflection group} associated with this circle packing $\mathcal{P}$ as
$$
G_\mathcal{P} := \langle \rho_1,..., \rho_{d+1}\rangle,
$$
where $\rho_i$ is the reflection along the circle $C_i$. As an abstract group, $G_\mathcal{P}$ is the free product of $d+1$ copies of $\Z/2\Z$.

Since a kissing reflection group is a discrete subgroup of $\textrm{Aut}^\pm(\widehat{\C})$ (the group of all M{\"o}bius and anti-M{\"o}bius automorphisms of $\widehat{\C}$) \cite[Part~II, Chapter~5, Theorem~1.2]{VS93}, definitions of limit set and domain of discontinuity easily carry over to kissing reflection groups (cf. \cite[\S 6.1]{LMMN}). We denote the domain of discontinuity and the limit set of $G_{\mathcal{P}}$ by $\Omega(G_{\mathcal{P}})$ and $\Lambda(G_{\mathcal{P}})$, respectively.

The following proposition characterizes kissing reflection groups with connected limit sets in terms of the contact graph of the underlying circle packing.

\begin{figure}[h!]
\captionsetup{width=0.96\linewidth}
\begin{center}
\includegraphics[width=0.5\linewidth]{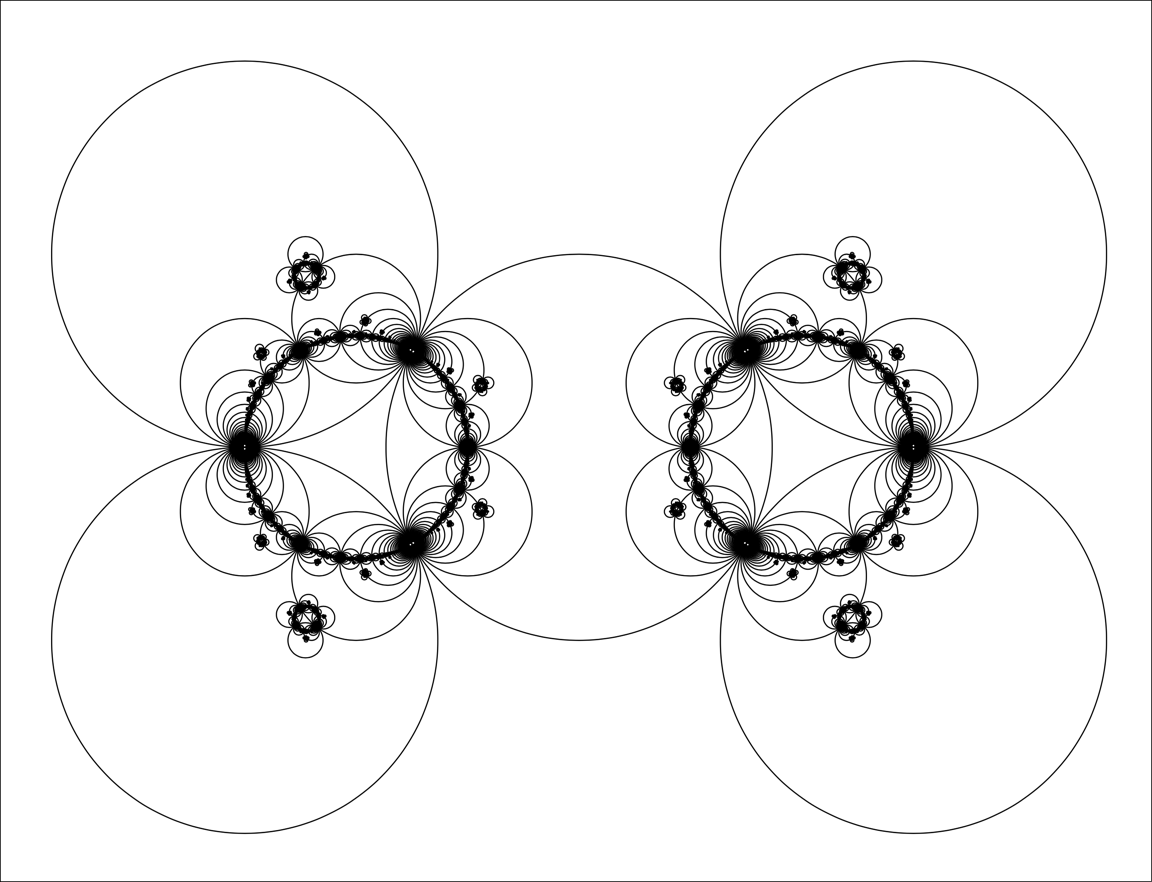}
\end{center}
\caption{A disconnected limit set for a kissing reflection group associated with a non-$2$-connected contact graph.}
\label{not_2_conn_fig}
\end{figure}

\begin{proposition}\cite[Proposition~3.4]{LLM1}\label{kissing_limit_conn_prop}
The kissing reflection group $G_\mathcal{P}$ has connected limit set if and only if the contact graph $\Gamma$ of $\mathcal{P}$ is $2$-connected (i.e., connected without cut points).
\end{proposition}
\noindent (In one direction: if there is a cut-point in $\Gamma$, then the corresponding circle of the packing intersects several components of the limit set, see Figure~\ref{not_2_conn_fig}.) 

\begin{remark}
For technical reasons, it is often important to consider \emph{marked} graphs and \emph{marked} circle packings. Such markings are particularly useful while talking about Teichm{\"u}ller spaces of kissing reflection groups (cf. \cite[\S 2, \S 3]{LLM1}).
\end{remark}

A {\em marking} of a graph $\Gamma$ is the choice of a graph isomorphism
$$
\phi: \mathscr{G}\longrightarrow \Gamma,
$$
where $\mathscr{G}$ is the underlying abstract graph of $\Gamma$.
We refer to the pair $(\Gamma, \phi)$ as a {\em marked graph}.
Similarly, a circle packing $\mathcal{P}$ is said to be marked if the associated contact graph is marked.

\subsubsection{Fundamental domain for kissing reflection groups}\label{fund_dom_subsubsec}
Let $\mathcal{P} = \{C_1,..., C_{d+1}\}$ be a circle packing, and $D_i$ be the open round disk enclosed by $C_i$.
For each $C_i$, let us consider the upper hemisphere $S_i\subset\mathbb{H}^3$ such that $\partial S_i\cap\partial\mathbb{H}^3= C_i$. The anti-M{\"o}bius reflection in $C_i$ extends naturally to the reflection in $S_i$, and defines an orientation-reversing isometry of $\mathbb{H}^3$. Let $\mathfrak{P}$ be the convex hyperbolic polyhedron (in $\mathbb{H}^3$) whose relative boundary in $\mathbb{H}^3$ is the union of the hemispheres $S_i$. Then, $\mathfrak{P}$ is a fundamental domain (called the \emph{Dirichlet fundamental polyhedron}) for the action of the group $G_{\mathcal{P}}$ on $\mathbb{H}^3$, and 
$$
\Pi(G_\mathcal{P}) :=\overline{\mathfrak{P}}\cap\Omega(G_{\mathcal{P}})
$$ 
(where the closure is taken in $\Omega(G_{\mathcal{P}})\cup\mathbb{H}^3$) is a fundamental domain for the action of $G_{\mathcal{P}}$ on $\Omega(G_{\mathcal{P}})$ (see \cite[Proposition~6.5]{LMMN}, cf. \cite[\S 3.5]{Mar16}, \cite{Vin67}). Clearly, the fundamental domain $\Pi(G_\mathcal{P})$ can also be written as
$$
\Pi(G_\mathcal{P}) =  \widehat{\C}\setminus\left(\bigcup_{i=1}^{d+1} D_i \cup \mathfrak{S}\right),
$$
where $\mathfrak{S}$ is the set consisting of points of tangency for the circle packing $\mathcal{P}$.
We remark that the fundamental domain $\Pi(G_{\mathcal{P}})$ is neither open, nor closed in $\widehat{\C}$, but is relatively closed in $\Omega(G_\mathcal{P})$. In Figure~\ref{kissing_nielsen_fig} and Figure~\ref{necklace_fig}, the fundamental domains are shaded in grey.

The above discussion shows that the action of a kissing reflection group on $\mathbb{H}^3$ admits a finite-sided polyhedron as its fundamental domain, and hence is \emph{geometrically finite}.

The index two subgroup $\widetilde{G}_{\mathcal{P}}\leqslant G_{\mathcal{P}}$ consisting of orientation-preserving elements is a Kleinian group whose domain of discontinuity coincides with $\Omega(G_\mathcal{P})$. A fundamental domain for the $\widetilde{G}_{\mathcal{P}}-$action on $\Omega(G_{\mathcal{P}})$ is given by doubling $\Pi(G_\mathcal{P})$ along $C_1$. Moreover, $\faktor{\Omega(G_{\mathcal{P}})}{\widetilde{G}_{\mathcal{P}}}$ is a finite union of punctured spheres, where each punctured sphere corresponds to the double of a component of $\Pi(G_\mathcal{P})$.

\subsubsection{Nielsen maps for kissing reflection groups}\label{nielsen_map_subsubsec}

Let $\mathcal{P}=\{C_1,\cdots, C_{d+1}\}$ be a circle packing realizing a $2$-connected simple plane graph.

\begin{definition}\label{nielsen_map_def} The \emph{Nielsen map} $\cN_{G_{\mathcal{P}}}$ is defined as: \begin{align} \cN_{G_{\mathcal{P}}} : \bigcup_{i=1}^{d+1} \overline{D_i} \rightarrow \widehat{\mathbb{C}} \nonumber \hspace{16mm} \\ \hspace{10mm} z\longmapsto \rho_i(z) \textrm{        if } z \in \overline{D_i}.  \nonumber\end{align}
\end{definition}
\noindent Note that the Nielsen map $\cN_{G_{\mathcal{P}}}$ is defined on the limit set $\Lambda(G_{\mathcal{P}})$, a fact that will be of importance later.
\begin{figure}[h!]
\captionsetup{width=0.96\linewidth}
\begin{center}
\includegraphics[width=0.42\linewidth]{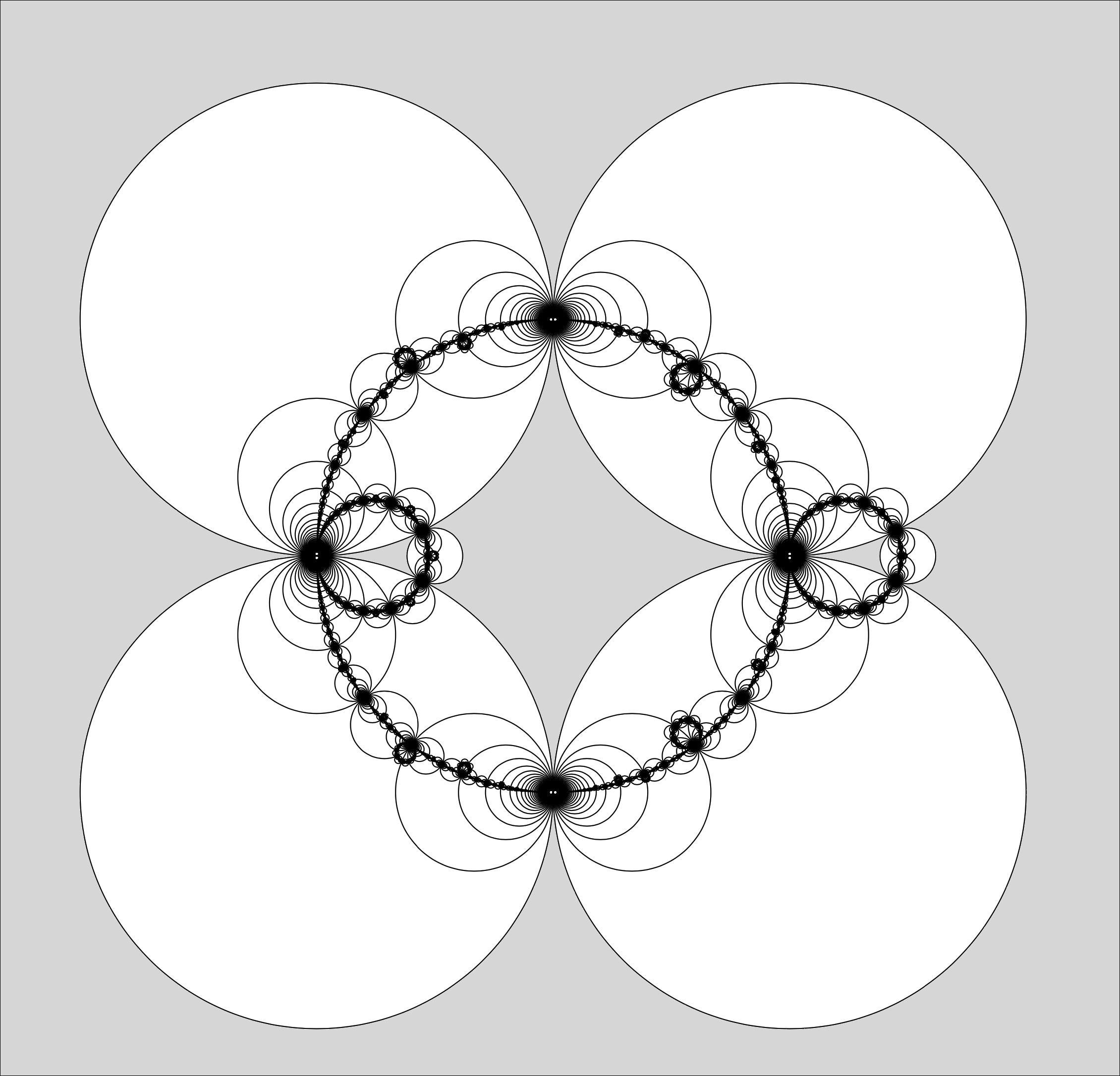}\quad \includegraphics[width=0.404\linewidth]{apollonian_gasket.png}
\includegraphics[width=0.38\linewidth]{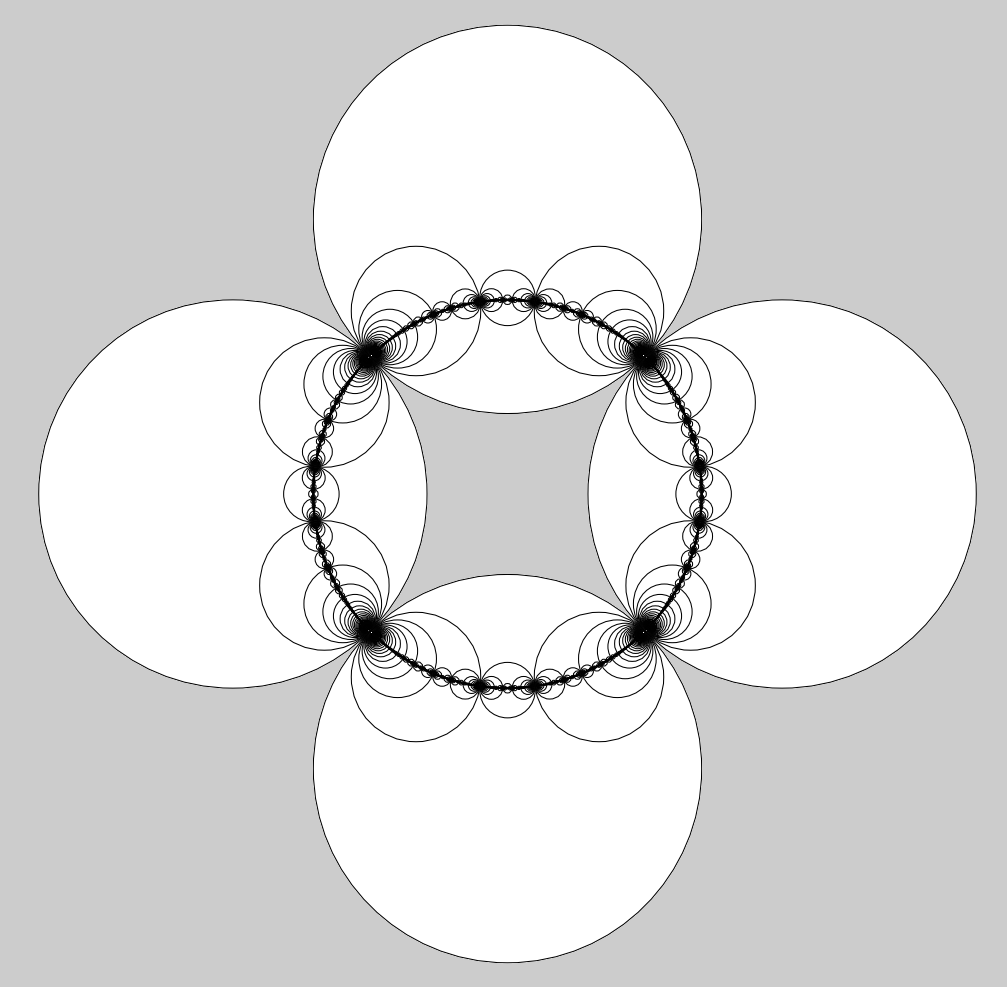}\quad \includegraphics[width=0.36\linewidth]{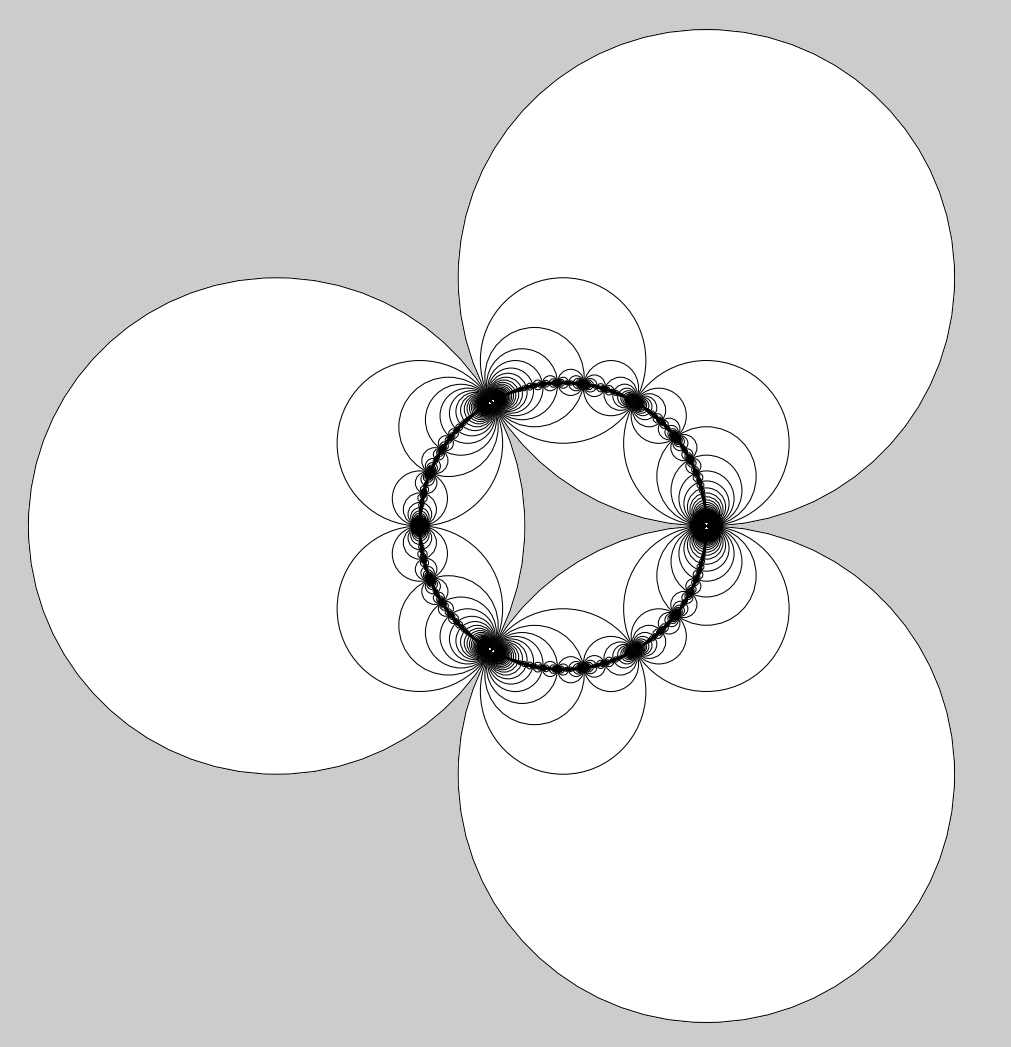}
\end{center}
\caption{Top: The fundamental domains $\Pi(G_\mathcal{P})$ of the kissing reflection groups are shaded in grey. The corresponding Nielsen maps preserve the components of the domains of discontinuity intersecting the fundamental domains. The restrictions of the Nielsen maps to these components of the domains of discontinuity are conformally conjugate to Nielsen maps of ideal pentagon and ideal triangle reflection groups. Bottom: The Nielsen maps $\pmb{\cN}_3$ and $\pmb{\cN}_2$ of regular ideal quadrilateral and ideal triangle reflection groups preserve the unit disk $\D$.}
\label{kissing_nielsen_fig}
\end{figure}

\begin{proposition}\cite[Proposition~4.1]{LLMM4}\cite[Proposition~16]{LMM2}\label{orbit_equiv_prop}
The map $\cN_{G_{\mathcal{P}}}$ is orbit equivalent to $G_{\mathcal{P}}$ on $\widehat{\C}$; i.e., for each $z\in\widehat{\C}$, the group orbit $G_{\mathcal{P}}\cdot z$ is equal to the grand orbit of $z$ under $\cN_{G_{\mathcal{P}}}$.
\end{proposition}

Let $\Pi(G_\mathcal{P})$ be as in Subsection~\ref{fund_dom_subsubsec}, and $\Pi_1,\cdots,\Pi_k$ the connected components of $\Pi(G_\mathcal{P})$. As the limit set $\Lambda(G_\mathcal{P})$ is connected, each component of the domain of discontinuity $\Omega(G_{\mathcal{P}})$ is simply connected, and each $\Pi_i$ is a closed ideal polygon in the corresponding component $\mathcal{U}_i$ of $\Omega(G_{\mathcal{P}})$ bounded by arcs of finitely many circles in the circle packing (see Figure~\ref{kissing_nielsen_fig}). Conjugating the stabilizer subgroup of $\mathcal{U}_i$ in $G_{\mathcal{P}}$ by a Riemann map of $\mathcal{U}_i$, one obtains an ideal polygon reflection group acting on $\D$. Moreover, the component $\mathcal{U}_i$ is forward invariant under the Nielsen map of $G_{\mathcal{P}}$, and the Riemann map conjugates $\cN_{G_{\mathcal{P}}}\vert_{\mathcal{U}_i}$ to the action of the Nielsen map of an ideal polygon reflection group on $\D$.
Thus, ideal polygon reflection groups play a special role while studying kissing reflection groups.

\begin{definition}\label{regular_ideal_polygon_ref_group_def} Consider the circle packing $\pmb{\mathcal{P}}_{d}:=\{\pmb{C}_1,\cdots, \pmb{C}_{d+1}\}$ where $\pmb{C}_j$ intersects $\mathbb{S}^1$ at right angles at the roots of unity $\exp{(\frac{2\pi i\cdot(j-1)}{d+1})}$, $\exp{(\frac{2\pi i\cdot j}{d+1})}$. We denote the associated kissing reflection group $G_{\pmb{\mathcal{P}}_{d}}$ by $\pmb{G}_{d}$, and call it the \emph{regular ideal polygon reflection group}.
\smallskip

We will denote the Nielsen map of $\pmb{G}_{d}$ by $\pmb{\cN}_d$.
\end{definition}

\subsubsection{Conjugation between $\pmb{\cN}_{d}$ and $\overline{z}^d$}\label{question_mark_subsubsec}

The Nielsen map $\pmb{\cN}_{d}$ restricts to a degree $d$ orientation-reversing expansive covering of $\mathbb{S}^1$. Hence, there exists a circle homeomorphism $\pmb{\mathcal{E}}_d$ that conjugates $\pmb{\cN}_{d}$ to $\overline{z}^d$, and sends $1$ to $1$. We call the homeomorphism the \emph{$d$-th Minkowski circle homeomorphism}. The existence of this conjugation between the Nielsen map of a group and an anti-polynomial lies at the heart of the connections between kissing reflection groups and anti-rational maps.

We note that the circle homeomorphism $\pmb{\mathcal{E}}_d$ conjugates an expansive circle map (with parabolic fixed points) to an expanding circle map (with only hyperbolic fixed points), and hence $\pmb{\mathcal{E}}_d$ is not a quasi-symmetric homeomorphism.

\begin{figure}[h!]
\captionsetup{width=0.96\linewidth}
\centering
\includegraphics[scale=0.42]{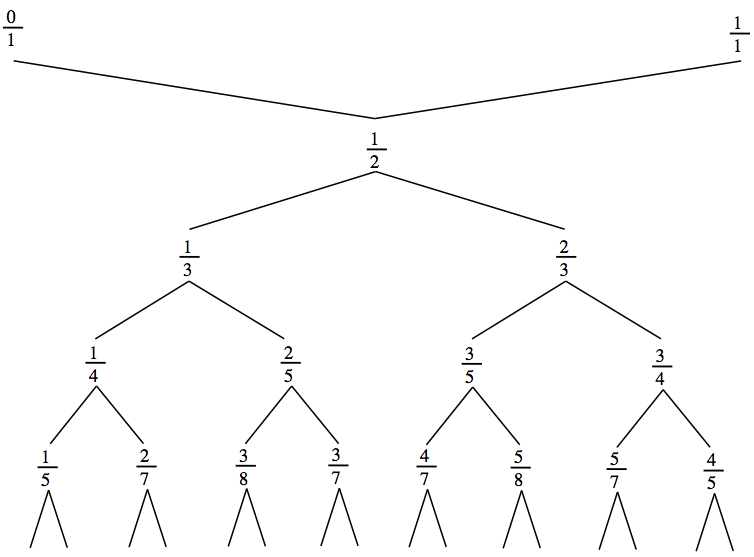}\quad \includegraphics[scale=0.42]{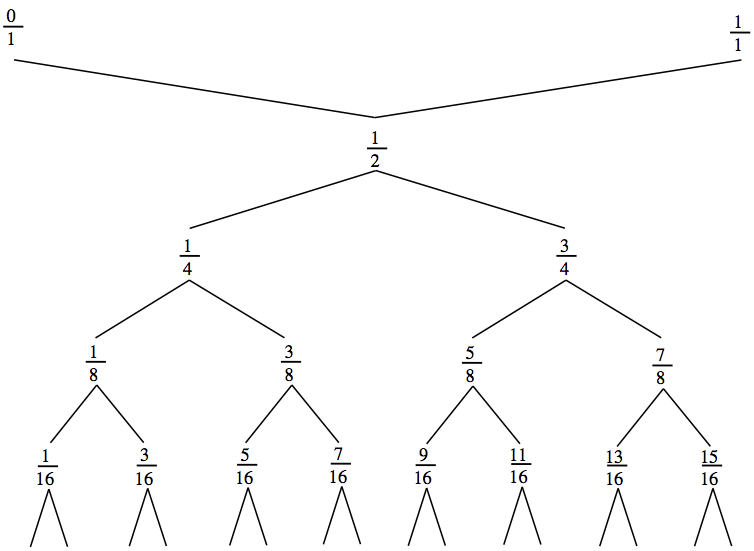} 
\caption{Left: The Farey tree. Right: The dyadic tree.}
\label{question_mark_tree_fig}
\end{figure}

We explain some connections between the circle homeomorphism $\pmb{\mathcal{E}}_2$ and classical objects in number theory and analysis.
In fact, the map $\pmb{\mathcal{E}}_2$ is a close relative of the classical Minkowski question mark function $\ciq:~ [0,1]\to [0,1].$ 
One way to define the question mark function is to set $\ciq(\frac01)=0$ and $\ciq(\frac11)=1$ and then use the recursive formula
\begin{equation}
\ciq\left(\frac{p+r}{q+s}\right)=\frac12\left (\ciq\left(\frac{p}{q}\right)+\ciq\left(\frac{r}{s}\right)\right)
\label{question_mark_recursion}
\end{equation}
which  gives us the values of $\ciq$ on all rational numbers (Farey fractions) in $[0,1]$. This defines a uniformly continuous function on $[0,1]\cap\mathbb{Q}$, whose unique continuous extension to $[0,1]$ is the Minkowski question mark function.
In particular, the map $\ciq$ is an increasing homeomorphism of $[0,1]$ that sends the vertices of level $n$ of the Farey tree to the vertices of level $n$ of the dyadic tree (see Figure~\ref{question_mark_tree_fig}). We refer the reader to \cite{Min,Den38,Sal43,Con01} for various number-theoretic and analytic properties of $\ciq$.

According to \cite[\S 4.4.2]{LLMM1}, the two maps, $\ciq$ and $\pmb{\mathcal{E}}_2$, are related by the formula
$$
\ciq^{-1}(x)=\phi\left(\pmb{\mathcal{E}}_2^{-1}\left(e^{\frac{2\pi i x}{3}}\right)\right),\ \forall\ x\in [0,1],
$$
where $\phi$ is a M{\"o}bius transformation carrying the unit disk onto the upper half-plane.
Roughly speaking, the Minkowski question mark function $\ciq$ is the restriction of the homeomorphism $\pmb{\mathcal{E}}_2$ to the arc $I:= [1,e^{2\pi i/3}]\subset\mathbb{S}^1$ written in appropriate coordinates.

\subsubsection{Deformation spaces of kissing reflection groups}\label{kissing_group_deform_space_subsubsec}

The space of all kissing reflection groups of a given rank can be organized in natural deformation spaces. As in the classical theory of Kleinian groups, the perspective of representations proves to be useful in the study of deformation spaces of reflection groups.

Let $G_0$ be a finitely generated discrete subgroup of $\textrm{Aut}^\pm(\widehat{\C})$.
A representation (i.e., a group homomorphism) $\xi: G_0 \longrightarrow \textrm{Aut}^\pm(\widehat{\C})$ is said to be \emph{weakly type-preserving} if 
\begin{enumerate}\upshape
\item $\xi(g) \in \textrm{Aut}^+(\widehat\C)$ if and only if $g\in \textrm{Aut}^+(\widehat\C)$, and
\item if $g \in \textrm{Aut}^+(\widehat\C)$, then $\xi(g)$ is parabolic whenever $g$ is parabolic.
\end{enumerate}
Note that a weakly type-preserving representation may send a loxodromic element to a parabolic one.

\begin{definition}\label{deform_space_def}
\noindent\begin{enumerate}\upshape
\item Given a kissing reflection group $G_0$, we define the \emph{algebraic deformation space}
\begin{align*}
\textrm{AH}(G_0)
&:=  \lbrace\xi: G_0\longrightarrow G\ \textrm{is a weakly type-preserving isomorphism to}
\\
&\qquad \textrm{a discrete subgroup}\ G\ \textrm{of}\ \textrm{Aut}^\pm(\widehat{\C})\rbrace / \sim,\;
\end{align*}
where $\xi_1\sim\xi_2$ if there exists a M{\"o}bius transformation $M$ such that $$\xi_2(g)=M\circ\xi_1(g)\circ M^{-1},\ \textrm{for all}\ g\in G_0.$$

\item We define the \emph{quasiconformal deformation} space  
\begin{align*}
\mathcal{QC}(G_0) 
&:= \{\xi \in \textrm{AH}(G_0): \xi(g)=\tau\circ g\circ \tau^{-1}, \text{where}\ \tau\ \textrm{is a}\\
&\qquad \textrm{quasiconformal homeomorphism of } \widehat{\C}\}.
\end{align*}

\item The \emph{Bers slice} of $\pmb{G}_{d}$ is the subspace of $\mathcal{QC}(\pmb{G}_{d})$ 
defined as 
\begin{align*}
\beta(\pmb{G}_{d}) 
&:= \{\xi\in\mathcal{QC}(\pmb{G}_{d}) :\ \tau\ \textrm{is conformal on}\ \D^*:=\widehat{\C}\setminus\overline{\D}\}.
\end{align*}
\end{enumerate}
\end{definition}

We endow $\textrm{AH}(G_0)$ with the quotient topology of algebraic convergence; more precisely, a sequence of weakly type-preserving representations $\{\xi_n\}$ converges to $\xi$ algebraically if $\{\xi_n(g_i)\}$ converges  to $\xi(g_i)$ as elements of $\textrm{Aut}^\pm(\widehat{\C})$ for (any) finite generating set $\{g_i\}$ of $G_0$.

We will now describe a natural stratification of quasiconformal deformation space closures into cells of various dimensions. Let us first recall that different realizations of a fixed marked, connected simple plane graph $\Gamma$ as circle packings $\mathcal{P}$ produce canonically isomorphic kissing reflection groups $G_{\mathcal{P}}$. Thus, the algebraic/quasiconformal deformation spaces of all such $G_{\mathcal{P}}$ can be canonically identified. Hence, it makes sense to fix a (marked) circle packing realization $\mathcal{P}$ of a (marked) graph $\Gamma$ and define $\mathcal{QC}(\Gamma):=\mathcal{QC}(G_{\mathcal{P}})$.

\begin{definition}\label{domination_def}
Let $\Gamma_0,\Gamma$ be simple plane graphs with the same number of vertices.
We say that $\Gamma$ {\em dominates} $\Gamma_0$, denoted by $\Gamma \geq \Gamma_0$, if there exists an embedding 
$\psi: \Gamma_0 \longrightarrow \Gamma$ as plane graphs (i.e., if there exists a graph isomorphism between $\Gamma_0$ and a subgraph of $\Gamma$ that extends to an orientation-preserving homeomorphism of $\widehat{\C}$).

\noindent We also define
$$
\mathrm{Emb}(\Gamma_0):=\{(\Gamma, \psi): \Gamma\geq \Gamma_0 \text{ and } \psi: \Gamma_0 \longrightarrow \Gamma \text{ is an embedding as plane graphs}\}.
$$
\end{definition}

With terminology as above, we have the following cell structure for the quasiconformal deformation space closure (where the closure is taken in $\mathrm{AH}(\Gamma_0)$). The proof of this result essentially uses the theory of pinching deformations and the Thurston Hyperbolization Theorem.
\begin{proposition}\cite[Proposition~3.17]{LLM1}\label{cell_structure_group_prop}
$$
\overline{\mathcal{QC}(\Gamma_0)} = \bigcup_{(\Gamma, \psi) \in \mathrm{Emb}(\Gamma_0)} \mathcal{QC}(\Gamma).
$$
\end{proposition}

\subsubsection{Necklace reflection groups}\label{necklace_subsubsec}

\begin{definition}\label{necklace_group_def} 
A kissing reflection group $G_{\mathcal{P}}$ is called a \emph{necklace group} if the contact graph of $\mathcal{P}$ is $2$-connected and outerplanar.
\end{definition}

This is equivalent to requiring that for the circle packing $\mathcal{P}$, each circle $C_i$ is tangent to $C_{i+1}$ (with $i+1$ taken mod $(d+1)$), and that the boundary of one of the components of $\displaystyle\widehat{\C}\setminus\bigcup_{i=1}^{d+1}\overline{D_i}$ intersects each $C_i$ (cf. \cite[Definition 6.7]{LMMN}).
\begin{figure}[h!]
\captionsetup{width=0.96\linewidth}
\centering
\includegraphics[width=0.3\textwidth]{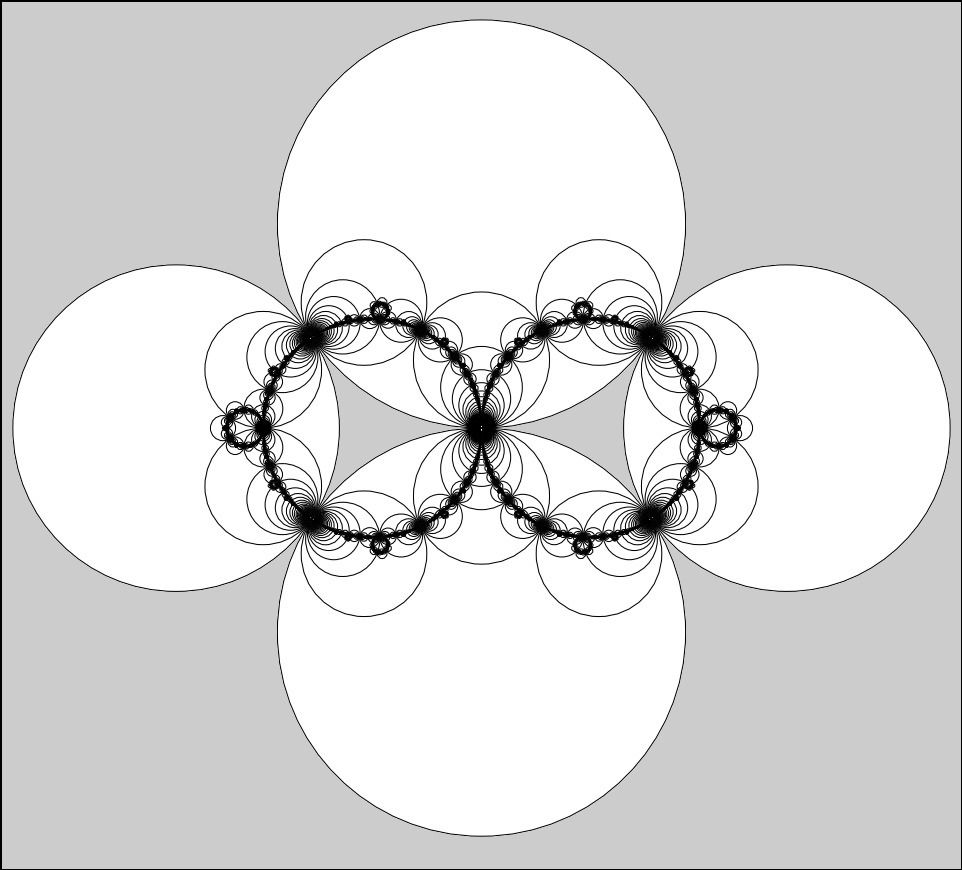}\ \includegraphics[width=0.308\textwidth]{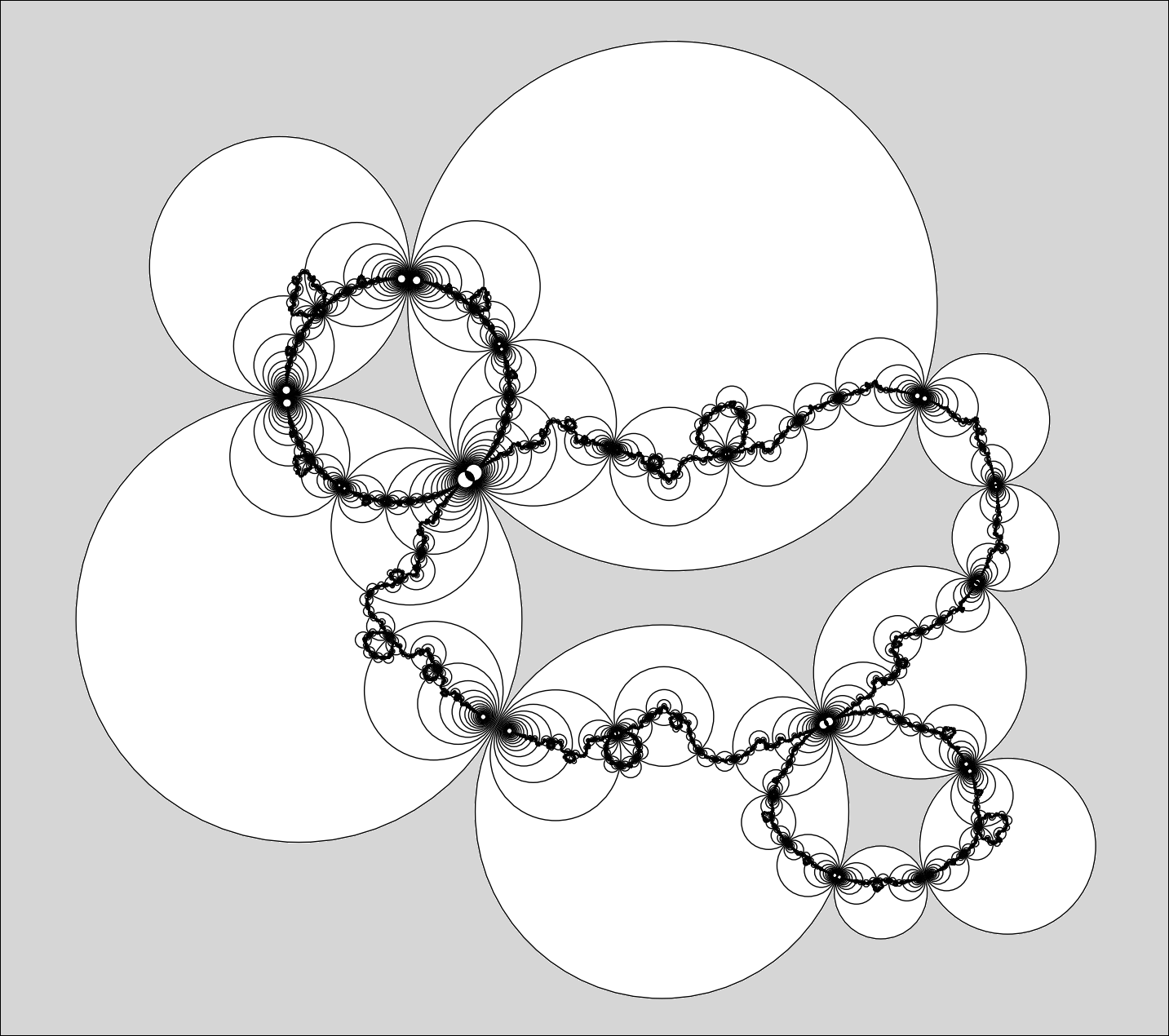}\ \includegraphics[width=0.304\textwidth]{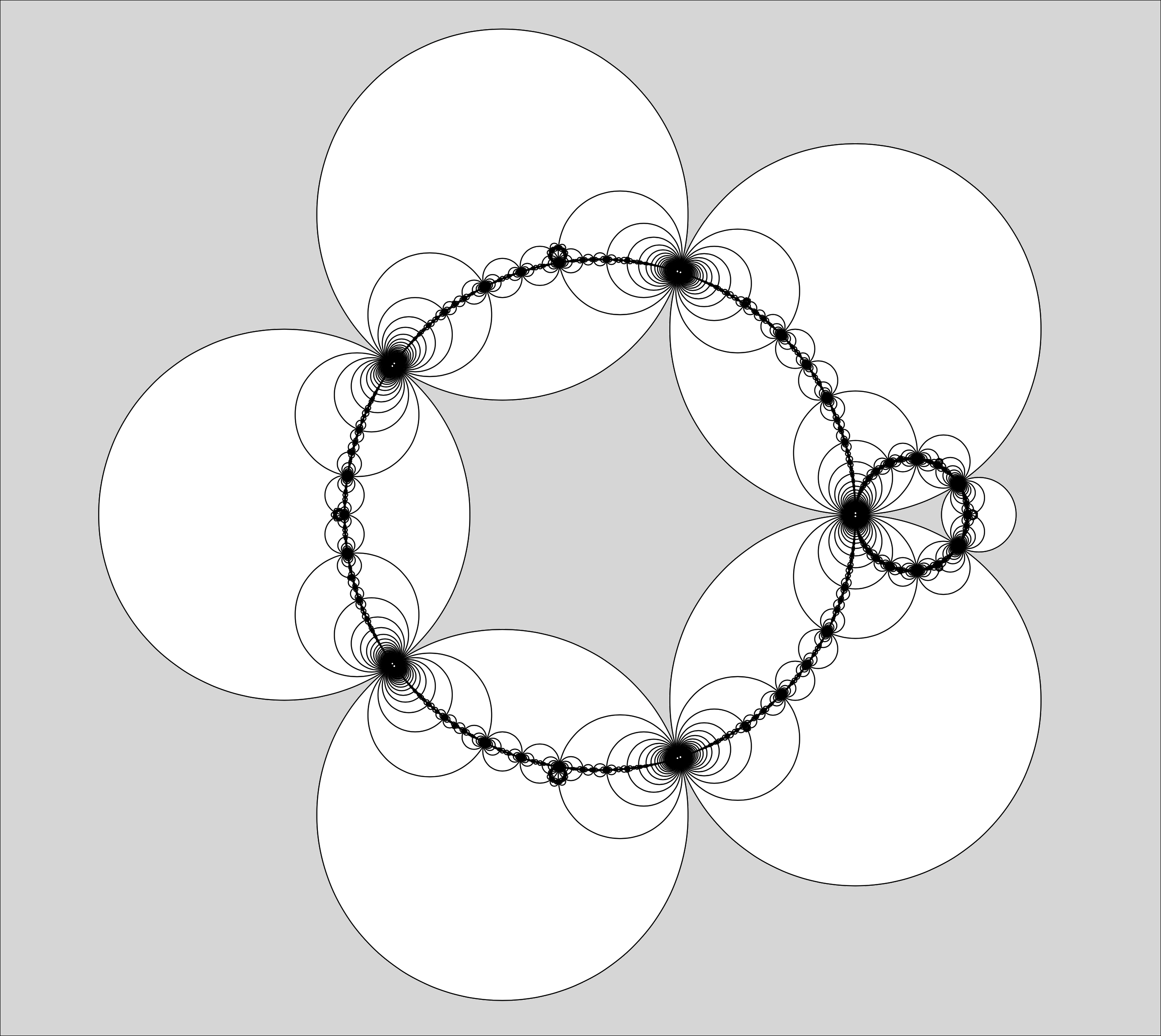}  
\caption{Limits sets of various necklace groups and their underlying circle packings are displayed. The components of the domains of discontinuity outside the limit sets are invariant under the groups.}
\label{necklace_fig}
\end{figure}

According to \cite[Proposition~3.20]{LLM1}, a kissing reflection group $G_{\mathcal{P}}$ with connected limit set is a necklace group if and only if there is an invariant component $\Omega_\infty(G_{\mathcal{P}})$ of $\Omega(G_{\mathcal{P}})$ such that the $G_{\mathcal{P}}$-action on $\Omega_\infty(G_{\mathcal{P}})$ is conformally equivalent to the action of an ideal polygon reflection group on $\D$\footnote{Note that necklace groups were termed \emph{function kissing reflection groups} in \cite{LLM1}.} (see Figure~\ref{necklace_fig}). 
\emph{After possibly quasiconformally conjugating a necklace group $G_{\mathcal{P}}$ on $\Omega_\infty(G_{\mathcal{P}})$, we can and will assume that $G_{\mathcal{P}}\vert_{\Omega_\infty(G_{\mathcal{P}})}$ is conformally equivalent to $\pmb{G}_{d}\vert_{\D^*}$.} With this convention, necklace groups are precisely those kissing reflection groups that lie on the closure of the Bers slice $\beta(\pmb{G}_{d})$.
Indeed, one can quasiconformally deform $\pmb{G}_{d}$ so that additional tangencies among the circles of the packing are introduced in the limit (cf. \cite[Proposition~11]{LMM2}).
This perspective allows one to embed the space of marked necklace groups (i.e., necklace groups associated with marked circle packings) of a given rank into the algebraic deformation space $\mathrm{AH}(\pmb{G}_{d})$ (see Definition~\ref{deform_space_def}).
In agreement with the classical theory of Kleinian groups, the Bers slice closure $\overline{\beta(\pmb{G}_{d})}$ (or equivalently, the space of necklace groups $G_{\mathcal{P}}$ with frozen conformal dynamics on $\Omega_\infty(G_{\mathcal{P}})$) is compact in $\mathrm{AH}(\pmb{G}_{d})$ (see \cite[\S 2.2]{LMM2}, \cite[\S 3.3]{LLM1} for details).

Recall that all kissing reflection groups, in particular those in $\overline{\beta(\pmb{G}_{d})}$, are geometrically finite.
This allows one to give a simple description of the dynamics of necklace groups on their limit sets.

\begin{proposition}\cite[Proposition~22]{LMM2}\label{group_lamination_prop} 
Let $G_{\mathcal{P}}$ be a necklace group associated with a marked circle packing $\mathcal{P}=\{C_1,\cdots,C_{d+1}\}$; i.e., $G_{\mathcal{P}}\in\overline{\beta(\pmb{G}_{d})}$. There exists a conformal map $\phi_{G_{\mathcal{P}}}: \D^* \rightarrow \Omega_{\infty}(G_{\mathcal{P}})$ such that 
\begin{align}\label{group_conjugacy}    
\pmb{\cN}_d (z) = \phi_{G_{\mathcal{P}}}^{-1} \circ \cN_{G_{\mathcal{P}}} \circ \phi_{G_{\mathcal{P}}}(z) \textrm{, for } z\in  \D^*\setminus \Int{\Pi(\pmb{G}_{d})}.   
\end{align}  
The map $\phi_{G_{\mathcal{P}}}$ extends continuously to a semi-conjugacy $\phi_{G_{\mathcal{P}}}: \mathbb{S}^1 \rightarrow \Lambda(G_{\mathcal{P}})$ between $\pmb{\cN}_d|_{\mathbb{S}^1}$ and $\cN_{G_{\mathcal{P}}}|_{\Lambda(G_{\mathcal{P}})}$, and for each $i$, sends the cusp of $\partial\Pi(\pmb{G}_{d})$ at $\pmb{C}_i\cap\pmb{C}_{i+1}$ to the cusp of $\partial\Pi(G_{\mathcal{P}})$ at $C_i\cap C_{i+1}$.
\end{proposition}

The continuous extension of the conformal map $\phi_{G_{\mathcal{P}}}$ to the circle is a consequence of local connectivity of $\Lambda(G_{\mathcal{P}})$ (according to \cite{AM96}, connected limit sets of geometrically finite Kleinian groups are locally connected).
Proposition~\ref{group_lamination_prop} yields a dynamically defined \emph{Carath{\'e}odory loop} for $\Lambda(G_{\mathcal{P}})$, which can be used to produce a topological model of the limit set in terms of \emph{geodesic laminations} (see Subsections~\ref{talbot_dynamics_subsubsec} and \ref{sigma_d_schwarz_dyn_subsubsec} for applications of this fact). More precisely, the fibers of $\phi_{G_{\mathcal{P}}}: \mathbb{S}^1 \rightarrow \Lambda(G_{\mathcal{P}})$ define an equivalence relation on $\mathbb{S}^1$ such that each non-trivial equivalence class consists of two points (cf. \cite[Proposition~54]{LMM2}). Connecting the points of each non-trivial equivalence class by a hyperbolic geodesic in $\D$ produces a $\pmb{G}_{d}$-invariant geodesic lamination (i.e., a closed set of mutually disjoint bi-infinite geodesics) on $\D$. This geodesic lamination can also be described as the lift to the universal cover of a collection of disjoint, simple, closed, non-peripheral geodesics on the punctured sphere $\faktor{\D}{\widetilde{\pmb{G}_d}}$ (cf. \cite[\S 4.3]{LLM1}).

For a necklace group $G_{\mathcal{P}}$, we set $K(G_{\mathcal{P}}):=\widehat{\C}\setminus\Omega_\infty(G_{\mathcal{P}})$, and call it the \emph{filled limit set} of $G_{\mathcal{P}}$. The component $\Pi(G_{\mathcal{P}})\cap\Omega_\infty(G_{\mathcal{P}})$, which is conformally equivalent to the polygon $\Pi(\pmb{G}_{d})\cap\D^*$, is denoted by $\Pi^u(G_{\mathcal{P}})$. Finally, we set $\Pi^b(G_{\mathcal{P}}) := \Pi(G_{\mathcal{P}})\setminus \Pi^u(G_{\mathcal{P}})$.

\begin{remark}
The local connectivity property holds for connected limit sets of arbitrary finitely generated Kleinian groups \cite{Mj14a}. This allows one to furnish geodesic lamination models for limit sets of Kleinian groups on boundaries of Bers slices \cite{Mj14b,Mj17}. The proofs of these results use the machinery developed for the proof of the Ending Lamination Conjecture (see \cite{Min10,BCM12}). For more detailed history of the problem, we refer the reader to \cite{Mj14a}.
\end{remark}

\subsection{Dynamics of anti-polynomials and the Tricorn}\label{tricorn_subsec}

In this subsection, we recall some known results on the dynamics of anti-polynomials, and their parameter space. Although the dynamical properties of anti-polynomials is similar to those of holomorphic polynomials, their parameter spaces have many important differences. We direct the reader to \cite[\S 2]{LLMM2} for a detailed account of the dynamics and parameter space of quadratic anti-polynomials.

\subsubsection{Some generalities}\label{anti_poly_dyn_general_subsubsec}
Let $p$ be an anti-polynomial of degree $d\geq 2$. The Fatou and Julia set of $p$ is defined to be those of the holomorphic second iterate $p^{\circ 2}$. In analogy to the holomorphic case, the set of all points that remain bounded under all iterations of $p$ is called the \emph{filled Julia set} $\mathcal{K}(p)$. The boundary of the filled Julia set equals the \emph{Julia set} $\mathcal{J}(p)$. The complement of $\mathcal{K}(p)$ is the basin of attraction of the superattracting fixed point $\infty$, and it is denoted by $\mathcal{B}_\infty(p)$.

Similar to the holomorphic case, there is a conformal map $\phi_p$ near $\infty$ that conjugates $p$ to $\overline{z}^d$. The map $\phi_p$ is unique up to multiplication by $(d+1)$-st roots of unity. If $p$ is monic (which can always be arranged by affine conjugation), then $\phi_p$ can be chosen to be tangent to the identity at $\infty$. With such normalization, the map $\phi_p$ is called the \emph{B{\"o}ttcher coordinate} of $p$ near $\infty$ \cite[Lemma~1]{Na1}. The absolute value of $\phi_p$ always extends to a continuous function on $\mathcal{B}_\infty(p)$, and the level curves of this function are called \emph{equipotential} curves of $p$. If $\mathcal{K}(p)$ is connected, then $\phi_p$ extends as a conformal conjugacy between $p\vert_{\mathcal{B}_\infty(p)}$ and $\overline{z}^d\vert_{\D^*}$. Otherwise, $\phi_p$ extends conformally to an equipotential curve containing the `fastest escaping' critical point (cf. \cite[\S 9]{Mil06}).

\begin{definition}\label{dyn_ray}
The \emph{dynamical ray} $R_p(\theta)$ of $p$ at an angle $\theta$ is defined as the pre-image of the radial line at angle $\theta$ under $\phi_p$.
\end{definition}

The dynamical ray $R_p(\theta)$ maps to the dynamical ray $R_p(-d\theta)$ under $p$. It follows that, at the level of external angles, the dynamics of $p$ can be studied by looking at the simpler map 
$$
m_{-d}:\R/\Z\to\R/\Z,\ m_{-d}(\theta)=-d\theta.
$$ 
It is well-known that if $p$ has a connected Julia set, then all rational dynamical rays of $p$ land at repelling or parabolic (pre-)periodic points. 

\begin{definition}\label{rat_lami_def}
The \emph{rational lamination} of an anti-polynomial $p$ with connected Julia set is defined as an equivalence relation on $\mathbb{Q}/\mathbb{Z}$ such that $\theta_1 \sim \theta_2$ if and only if the dynamical rays $R_p(\theta_1)$ and $R_p(\theta_2)$ land at the same point of $\mathcal{J}(p)$. The rational lamination of $p$ is denoted by $\lambda(p)$.
\end{definition}

Some of the basic properties of rational laminations are listed in the next proposition.

\begin{proposition}\cite{Kiw01}\label{rat_lami_prop}
The rational lamination $\lambda(p)$ of an anti-polynomial $p$ with connected Julia set satisfies the following properties.

\begin{enumerate}\upshape
\item $\lambda(p)$ is closed in $\Q/\Z\times\Q/\Z$.

\item Each  $\lambda(p)$-equivalence class $A$ is a finite subset of $\Q/\Z$.

\item If $A$ is a $\lambda(p)$-equivalence class, then $m_{-d}(A)$ is also a $\lambda(p)$-equivalence class.

\item If $A$ is a $\lambda(p)$-equivalence class, then $A\mapsto m_{-d}(A)$ is consecutive reversing.

\item $\lambda(p)$-equivalence classes are pairwise unlinked.
\end{enumerate}
\end{proposition}

\begin{remark}\label{lami_rmk}
For an anti-polynomial $p$ with connected Julia set, the smallest equivalence relation $\widehat{\lambda(p)}$ on $\R/\Z$ that contains the closed set $\overline{\lambda(p)}$ (in $\R/\Z\times\R/\Z$) is called the \emph{combinatorial lamination} of $p$. If $\mathcal{J}(p)$ is locally connected and $p$ has no irrationally neutral cycle, then $\mathcal{J}(p)$ is homeomorphic to the quotient of the circle by the equivalence relation~$\widehat{\lambda(p)}$ (cf. \cite[Lemma~4.17]{Kiw01}).
\end{remark}

For an antiholomorphic germ $g$ fixing a point $z_0$, the quantity $\frac{\partial g}{\partial\overline{z}}\vert_{z_0}$ is called the \emph{multiplier} of $g$ at the fixed point $z_0$. One can use this definition to define multipliers of periodic orbits of antiholomorphic maps. A cycle is called attracting (respectively, super-attracting) if the associated multiplier has absolute value between $0$ and $1$ (respectively, is equal to $0$). The dynamics of $g$ near such a point is similar to that of a holomorphic germ near a (super-)attracting fixed point \cite[\S 8, 9]{Mil06}.
On the other hand, neutral fixed points of antiholomorphic germs are special in the following sense. 
Note that for $\theta\in\R$, the map $z\mapsto e^{i\theta}\overline{z}$ is an antiholomorphic involution. Hence the second iterate of a neutral antiholomorphic germ is a tangent-to-identity holomorphic parabolic germ. In this sense, \emph{any neutral fixed point of an antiholomorphic germ $g$ is parabolic}.

\begin{proposition}\cite[Lemma~2.3]{HS}\label{normalization of fatou}
Suppose $z_0$ is a parabolic periodic point of odd period $k$ of an anti-polynomial $p$ with only one petal, and $U$ is a periodic Fatou component with $z_0 \in \partial U$. Then there is an open subset $V \subset U$ with $z_0 \in \partial V$, and $g^{\circ k} (V) \subset V$ so that for every $z \in U$, there is an $n \in \mathbb{N}$ with $g^{\circ nk}(z)\in 
V$. Moreover, there is a univalent map $\psi^{\mathrm{att}} \colon V \to \mathbb{C}$ that conjugates $g^{\circ k}$ to the glide-reflection $\zeta\mapsto\overline{\zeta}+1/2$; i.e.,
$$
\psi^{\mathrm{att}}(g^{\circ k}(z)) = \overline{\psi^{\mathrm{att}}(z)}+1/2\quad \forall\quad z\in V,
$$ and $\psi^{\mathrm{att}}(V)$ contains a right half plane. This map $\psi^{\mathrm{att}}$ is unique up to composition with a horizontal translation. 
\end{proposition}

The map $\psi^{\mathrm{att}}$ is called the \emph{attracting Fatou coordinate} for the petal $V$. The antiholomorphic iterate interchanges both ends of the {\'E}calle cylinder, so it must preserve one horizontal line around this cylinder (the \emph{equator}). The change of coordinate has been so chosen that the equator is the projection of the real axis.  We will call the vertical Fatou coordinate the \emph{{\'E}calle height}. The {\'E}calle height vanishes precisely on the equator. The existence of this distinguished real line, or equivalently an intrinsic meaning to {\'E}calle height, is specific to antiholomorphic maps and plays a crucial role in parameter space discussions (cf. Theorem~\ref{odd_hyp_bdry_tricorn_thm}).

\subsubsection{Quadratic anti-polynomials and the Tricorn}\label{tricorn_subsubsec}
Any quadratic anti-polynomial, after an affine change of coordinates, can be written in the form $f_c(z) = \overline{z}^2 + c$ for $c \in \mathbb{C}$. This leads, as in the holomorphic case, to the notion of \emph{connectedness locus} of quadratic anti-polynomials:

\begin{definition}\label{tricorn_def}
The \emph{Tricorn} is defined as $\mathcal{T} = \{ c \in \mathbb{C} : \mathcal{K}_c:=\mathcal{K}(f_c)$ is connected$\}$. 
\end{definition}

The dynamics of quadratic anti-polynomials and its connectedness locus was first studied numerically in \cite{CHRS}, where this set was called the \emph{Mandelbar set}. Their numerical experiments showed curious structural differences between the Mandelbrot set and the Tricorn; in particular, they observed that there are bifurcations from the period $1$ hyperbolic component to period $2$ hyperbolic components along arcs in the Tricorn (see Theorem~\ref{ThmBifArc} for a general statement), in contrast to the fact that bifurcations are always attached at a single point in the Mandelbrot set. The name `Tricorn' is due to Milnor, who found `copies' of the connectedness locus of quadratic anti-polynomials in parameter spaces of real cubic polynomials and other real maps \cite{Mil92, Mil00} (the name comes from the three-cornered shape of the set).
Namely, Milnor observed that the dynamics of the corresponding real maps exhibit quadratic anti-polynomial-like behavior, and this was the primary motivation to view the Tricorn as a prototypical object in the study of real slices of holomorphic maps.
Nakane proved that the Tricorn is connected, in analogy to Douady and Hubbard's classical proof of connectedness of the Mandelbrot set \cite{Na1}:

\begin{theorem}\cite{Na1}\label{RealAnalUniformization}
The map $\Phi : \mathbb{C} \setminus \mathcal{T} \rightarrow \mathbb{C} \setminus \overline{\mathbb{D}}$, defined by $c \mapsto \phi_c(c)$ (where $\phi_c$ is the B\"{o}ttcher coordinate near $\infty$ for $f_c$) is a real-analytic diffeomorphism. In particular, the Tricorn is connected.
\end{theorem}

The previous theorem also allows us to define parameter rays of the Tricorn. 
\begin{definition}
The \emph{parameter ray} at angle $\theta$ of the Tricorn $\mathcal{T}$, denoted by $\mathcal{R}_{\theta}$, is defined as $\{ \Phi^{-1}(r e^{2 \pi i \theta}) : r > 1 \}$, where $\Phi$ is the real-analytic diffeomorphism from the exterior of $\mathcal{T}$ to the exterior of the closed unit disc in the complex plane constructed in Theorem~\ref{RealAnalUniformization}.
\end{definition}

We refer the reader to \cite[\S 3]{NS}, \cite{Muk1} for details on combinatorics of landing patterns of dynamical rays for unicritical anti-polynomials.

A map $f_c$ is called hyperbolic (respectively, parabolic) if it has a (super-)attracting (respectively, parabolic) cycle. Connected components of the set of all hyperbolic parameters is called a \emph{hyperbolic component} of $\mathcal{T}$.
\begin{figure}[h!]
\captionsetup{width=0.96\linewidth}
\begin{center}
\includegraphics[width=0.32\linewidth]{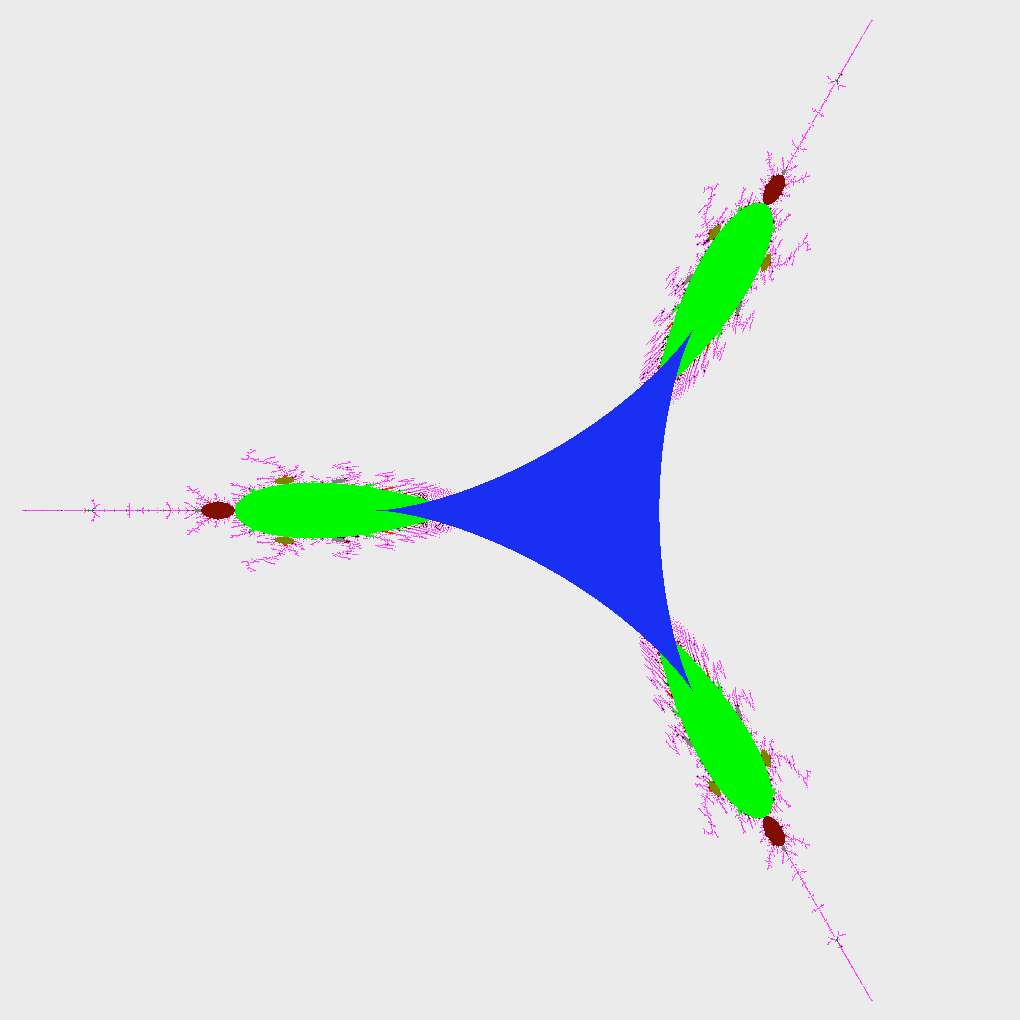}\ \includegraphics[width=0.32\linewidth]{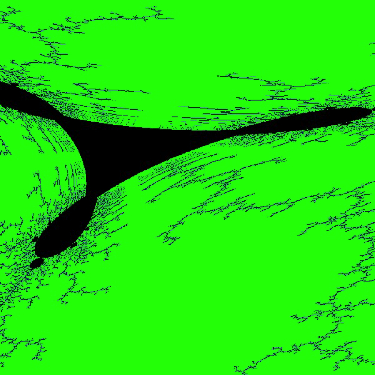}\ \includegraphics[width=0.32\linewidth]{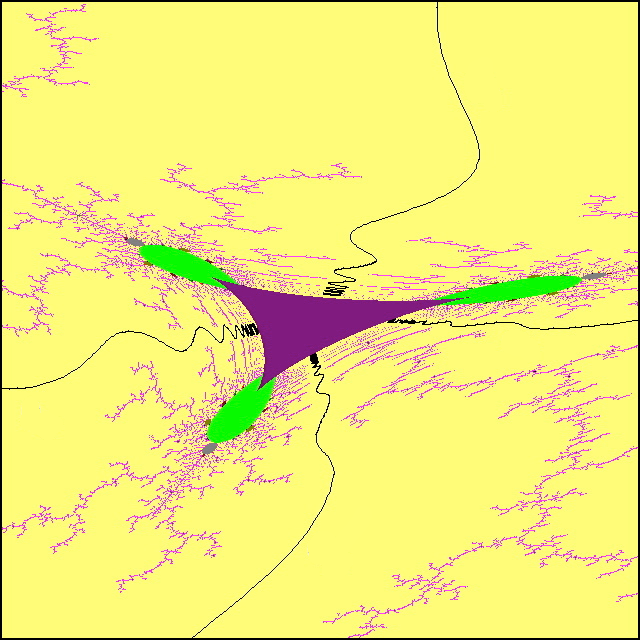}
\end{center}
\caption{Left: The Tricorn $\mathcal{T}$. Middle: Wiggling of an umbilical cord on the boundary of an odd period hyperbolic component of $\mathcal{T}$. Right: Non-trivial accumulation of parameter rays on the boundary of an odd period hyperbolic component of $\mathcal{T}$.}
\label{tricorn_fig}
\end{figure}
The even period hyperbolic components of $\cT$ are similar to the hyperbolic components of the Mandelbrot set, in the sense that they are real-analytically uniformized by the multiplier of the unique attracting cycle \cite[Theorem~5.6]{NS}, \cite[\S 2.2.4]{LLMM2}. The odd period hyperbolic components of $\cT$, however, are more delicate. 
Let $H$ be a hyperbolic component of odd period $k\neq 1$, let $c\in H$, and $z_c$ be an attracting periodic point of $f_c$. We denote the Jacobian determinant of $f_c^{\circ k}$ at $z_c$ by $\mathrm{Jac}(f_c^{\circ k},z_c)$, and note that 
$$
(f_c^{\circ 2k})'(z_c)=-\mathrm{Jac}(f_c^{\circ k},z_c)=\left|\frac{\partial f_c^{\circ k}}{\partial\overline{z}}(z_c)\right|^2
$$ 
is real and positive (cf. \cite[\S 1.1]{Muk1}). Thus, multipliers of attracting cycles fail to produce a dynamical uniformization of odd period hyperbolic components. In fact, the natural conformal invariant for maps with odd period attracting cycles, called the \emph{Koenigs ratio}, is not a purely local quantity. The Koenigs ratio captures the conformal position of the critical value in the normalized Koenigs coordinate. We refer the reader to \cite[\S 5]{NS}, \cite[\S 6]{IM2} for details.

\begin{theorem}\cite[Theorem~5.6, Theorem~5.9]{NS}\label{tricorn_hyp_unif_thm}
Let $H$ be a hyperbolic component of~$\cT$. 
\begin{enumerate}\upshape
\item If $H$ is of odd period, then the Koenigs ratio map is a real-analytic $3$-fold branched covering from $H$ onto the unit disk, ramified only over the origin.

\item If $H$ is of even period, then the multiplier map is a real-analytic diffeomorphism from $H$ onto the unit disk.
\end{enumerate}
\end{theorem}

\begin{remark}\label{koenigs_ratio_rem}
We note that for each quadratic anti-polynomial in an odd (respectively, even) period hyperbolic component, the first return map of the critical value Fatou component is a degree two proper antiholomorphic (respectively, holomorphic) map and hence has three (respectively, one) boundary fixed points. This fact is manifested in the mapping degree of the Koenigs ratio map (respectively, the multiplier map) of an odd (respectively, even) period hyperbolic component.
\end{remark}

The following results describe the boundaries of the hyperbolic components of the Tricorn and the associated bifurcation structure. Once again, the odd and even period hyperbolic components exhibit strikingly different behavior. In particular, the boundaries of odd period hyperbolic components are entirely comprised of parabolic parameters. 

\begin{theorem}\cite[Theorem~1.1]{MNS}\label{ThmEvenBif} 
\noindent\begin{enumerate}\upshape
\item The boundary of an even period hyperbolic component $H$ of $\cT$ is a topological circle. For each $c\in\partial H$, the corresponding map $f_c$ has a unique neutral cycle. 
 
\item If $f_c$ has a $2k$-periodic cycle with multiplier $e^{2\pi ip/q}$ with $\mathrm{gcd}(p,q)=1$, then $c$ sits on the boundary of a hyperbolic component of period $2kq$ of the Tricorn (and is the root thereof). 
\end{enumerate}
\end{theorem}

\begin{theorem}\cite[Lemma~2.5, Theorems~1.2,~3.2]{MNS}\label{odd_hyp_bdry_tricorn_thm}  
\noindent\begin{enumerate}\upshape
\item The boundary of a hyperbolic component $H$ of odd period $k$ consists entirely of parameters having a parabolic orbit of exact period $k$. In suitable local conformal coordinates, the $2k$-th iterate of such a map has the form $z\mapsto z+z^{q+1}+\cdots$, with $q\in\{1,2\}$. 

\item If $q=1$ for some $\widetilde{c}\in\partial H$, then $\widetilde{c}$ lies on a \emph{parabolic arc} in the  following sense: there exists a real-analytic arc of simple parabolic parameters $c(h)$ (for $h\in\mathbb{R}$) with quasiconformally equivalent but conformally distinct dynamics of which $\widetilde{c}$ is an interior point, and the {\'E}calle height of the critical value of $f_{c(h)}$ is $h$. In particular, $h\mapsto c(h)$ yields a monotone embedding of $\R$ in $\partial H$.

\item The boundary of every  odd period hyperbolic component of $\mathcal{T}$ is a topological triangle having double parabolic parameters  (i.e., $q=2$) as vertices and parabolic arcs as sides.
\end{enumerate}
\end{theorem}

For an odd period hyperbolic component $H$, as $c\in H$ approaches a simple (respectively, double) parabolic parameter on $\partial H$, an attracting periodic point merges with a repelling periodic point (respectively, two repelling periodic points) to produce a simple (respectively, double) parabolic periodic point.

When the parameter $c$ crosses a parabolic arc from the inside to the outside of $H$, then doubling bifurcation occurs. More precisely, two fixed points of $f_c^{\circ k}$ (where $k$ is the period of $H$) merge, giving rise to a period $2k$-cycle. The type of this cycle changes from attracting to neutral to repelling depending on the  crossing point on the arc.

To describe the situation precisely, we need the notion of fixed point residues, which we now recall. Let $g : U \rightarrow \mathbb{C}$ be a holomorphic function on a connected open set $U\ \left(\subset \mathbb{C}\right)$, and $\hat{z}\in U$ be an isolated fixed point of $f$. Then, the fixed point residue of $f$ at $\hat{z}$ is defined to be the complex number
$$
\displaystyle \iota(f, \hat{z}) = \frac{1}{2\pi i} \oint \frac{dz}{z-f(z)},
$$
where we integrate in a small loop in the positive direction around $\hat{z}$. It is easy to see that the fixed point residue does not depend on the choice of complex coordinates, so it is a conformal invariant (see \cite[\S 12]{Mil06} for basic properties of fixed point residues). 

By the fixed point residue of a periodic orbit of odd period of $f_c$, we will mean the fixed point residue of the second iterate $f_c^{\circ 2}$ at that periodic orbit.
Let $\mathcal{C}$ be a parabolic arc of odd period $k$ equipped with the critical {\'E}calle height parametrization $c:\R\to\mathcal{C}$  (by the above theorem). For any $h$ in $\mathbb{R}$, let us denote the fixed point residue of the unique parabolic cycle of $f_{c(h)}^{\circ 2}$ by $\ind(f_{c(h)}^{\circ 2})$. This defines a function 
$$
\ind_{\mathcal{C}}: \mathbb{R}\to\mathbb{C},\ h\mapsto \ind(f_{c(h)}^{\circ 2}).
$$

As $c$ tends to a vertex of an odd period hyperbolic component along a parabolic arc, a simple parabolic point merges with a repelling point to form a double parabolic point. In the process, the sum of the fixed point residues of the simple parabolic point and the repelling point converges to the fixed point residue of the double parabolic point. This allows one to understand the asymptotic behavior of the parabolic fixed point residue towards the ends of parabolic arcs. Finally, when a parameter on a parabolic arc is perturbed outside the odd period hyperbolic component, the simple parabolic cycle bifurcates to an attracting or a repelling cycle of twice the period, depending on whether the fixed point residue of the simple parabolic cycle is larger or smaller than $1$.

\begin{theorem}\cite[Proposition~3.7, Theorem~3.8, Corollary~3.9]{HS}, \cite[Lemma 2.12]{IM2}\label{ThmBifArc}
\begin{enumerate}\upshape
\item The function $\ind_{\mathcal{C}}$ is real-valued and real-analytic. Moreover, $$\lim_{h\to\pm\infty}\ind_{\mathcal{C}}(h)=+\infty.$$

\item Every parabolic arc of period $k$ intersects the boundary of a hyperbolic component of period $2k$ along an arc consisting of the set of parameters where the parabolic fixed point residue is at least $1$. In particular, every parabolic arc has, at both ends, an interval of positive length at which bifurcation from a hyperbolic component of odd period $k$ to a hyperbolic component of period $2k$ occurs.
\end{enumerate}
\end{theorem}

A \emph{Misiurewicz} parameter of the Tricorn is a parameter $c$ such that the critical point $0$ is strictly pre-periodic. It is well-known that for a Misiurewicz parameter, the critical point eventually maps on a repelling cycle. By classification of Fatou components, the filled Julia set of such a map has empty interior. Moreover, the Julia set of a Misiurewicz parameter is locally connected \cite[Expos{\'e} III, Proposition~4, Theorem~1]{orsay}.
These parameters play an important role in the understanding of the topology of the Tricorn.

\begin{theorem}\cite[Theorem~2.37]{LLMM2}\label{Tricorn_para_misi_ray} 
Every parameter ray of the Tricorn at a  strictly pre-periodic angle (under $m_{-2}$) lands at a Misiurewicz parameter such that in its dynamical plane, the corresponding dynamical ray lands at the critical value. Conversely, every Misiurewicz parameter $c$ of the Tricorn is the landing point of a finite (non-zero) number of parameter rays at strictly pre-periodic angles (under $m_{-2}$) such that the angles of these parameter rays are exactly the external angles of the dynamical rays that land at the critical value $c$ in the dynamical plane of $f_c$. 
\end{theorem}

We now collect some results that underscore the differences between the global topology of the Tricorn and the Mandelbrot set. Such results include non-landing of rational parameter rays, non-density of Misiurewicz parameters on the boundary, lack of local connectedness, discontinuity of straightening maps between small Tricorn-like sets and the original Tricorn, etc. The lack of quasiconformal rigidity on the boundary of the Tricorn (more precisely, the existence of arcs of quasiconformally conjugate parabolic parameters, see Theorem~\ref{odd_hyp_bdry_tricorn_thm}) lies at the heart of these topological differences.

\begin{theorem}\cite{IM1}\label{most rays wiggle}
The root (respectively, co-root) parabolic arc on the boundary of a hyperbolic component of odd period (except period one) of $\mathcal{T}$ contains the accumulation set of exactly two (respectively, one) parameter rays. The accumulation set of every such parameter ray contains an arc of positive length. On the other hand, the fixed rays at angles $0$, $1/3$ and $2/3$ land on the boundary of the period one hyperbolic component. 
\end{theorem}

\begin{remark}
Notice that the non-trivial accumulation set of a parameter ray on the boundary of an odd period hyperbolic component may be accessible from the exterior of the Tricorn. This is a manifestation of the drastic differences between conformal and real-analytic uniformizations (see Figure~\ref{tricorn_fig}).
\end{remark}

\begin{theorem}\cite{IM1}\label{thm_misi_not_dense}
Misiurewicz parameters are not dense on the boundary of $\mathcal{T}$. Indeed, there are points on the boundaries of the period $1$ and period $3$ hyperbolic components of $\mathcal{T}$ that cannot be approximated by Misiurewicz parameters. 
\end{theorem}

\begin{theorem}\cite[Theorem~6.2]{HS},\cite[Theorem~1.2]{IM2} \label{Tricorn_non_lc}
The Tricorn is not path connected. Moreover, no non-real hyperbolic component of odd period can be connected to the principal hyperbolic component by a path.
\end{theorem}

\begin{theorem}\cite{IM2}\label{Straightening_discontinuity_Tricorn}
Let $c_0$ be the center of a hyperbolic component $H$ of odd period (other than $1$) of $\mathcal{T}$, and $\mathcal{R}(c_0)$ be the corresponding $c_0$-renormalization locus (i.e. the baby Tricorn based at $H$). Then the straightening map $\chi_{c_0} : \mathcal{R}(c_0) \rightarrow \mathcal{T}$ is discontinuous at infinitely many parameters.
\end{theorem}

We conclude our discussion of the Tricorn with the definition of the real Basilica limb of the Tricorn, which will be important in Subsection~\ref{c_and_c_general_subsec}. Of course, one can give a more general definition of limbs, which can be found in \cite[\S 6]{MNS}. Let us denote the hyperbolic component of period one of $\mathcal{T}$ by $H_0$.

\begin{definition}\label{def_basilica_limb}
The connected component of $\left(\mathcal{T}\setminus\overline{H_0}\right)\cup\{-\frac{3}{4}\}$ intersecting the real line is called the \emph{real Basilica limb} of the Tricorn, and is denoted by $\mathcal{L}$.
\end{definition}

The real Basilica limb $\mathcal{L}$ is precisely the set of parameters $c$ in $\mathcal{T}$ such that in the dynamical plane of $f_c$, the rays $R_c(1/3)$ and $R_c(2/3)$ land at a common point (i.e. $1/3\sim2/3$ in $\lambda(f_c)$ for all $c\in\mathcal{L}$). The real Basilica limb of the Tricorn is depicted in Figure~\ref{c_and_c_para_fig} (right).

\subsubsection{Parabolic Tricorn and its higher degree versions}\label{para_anti_rat_gen_subsubsec}

While anti-polynomials form the simplest class of anti-rational maps, the expanding external class $\overline{z}^d$ of an anti-polynomial (with connected Julia) is a recurring source of mismatch between the dynamics of anti-polynomials and the action of groups with parabolic elements (such as kissing reflection groups). In certain situations, which will be elucidated in Subsections~\ref{chebyshev_subsec},~\ref{chebyshev_gen_subsec}, and Section~\ref{general_mating_corr_sec}, a closely related family of anti-rational maps with an expansive (but not expanding) external class enables one to overcome this difficulty. Specifically, these anti-rational maps have a completely invariant, simply connected Fatou component (like anti-polynomials with connected Julia sets), that is a parabolic basin of attraction (as opposed to the basin of attraction of $\infty$ for anti-polynomials). We now give a formal treatment of this family.
 
Note that the anti-Blaschke product 
$$
B_d(z) = \frac{(d+1)\overline z^d + (d-1)}{(d-1)\overline z^d + (d+1)}
$$ 
has a parabolic fixed point at $1$, and $\mathbb{D}$ is an invariant parabolic basin of this fixed point. Due to real-symmetry of the map $B_d$, the unique critical point $0$ of $B_d$ in $\D$ has {\'E}calle height zero.

\begin{definition}
The family $\pmb{\mathcal{B}}_d$ consists of degree $d\geq 2$ anti-rational maps $R$ with the following properties.
\begin{enumerate}\upshape
\item $\infty$ is a parabolic fixed point for $R$.
\item There is a marked parabolic basin $\mathcal{B}(R)$ of $\infty$ which is simply connected and completely invariant.
\item $R\vert_{\mathcal{B}(R)}$ is conformally conjugate to $B_d\vert_{\D}$.
\end{enumerate}
The complement $\mathcal{K}(R):=\widehat{\mathbb{C}}\setminus \mathcal{B}(R)$ of the marked parabolic basin is called the \emph{filled Julia set}  of $R$.
\end{definition}

Analogous to the connectedness locus of degree $d$ anti-polynomials, the moduli space $\left[\pmb{\mathcal{B}}_d\right]:= \faktor{\pmb{\mathcal{B}}_d}{\mathrm{Aut}(\C)}$ is also compact \cite[Proposition~4.3]{LMM3}.

The family $\pmb{\mathcal{B}}_2$ is called the \emph{parabolic Tricorn}. We direct the reader to \cite[Appendix~A]{LLMM3} for a more explicit description and the study of basic topological properties of the parabolic Tricorn.

\subsection{Schwarz reflection maps}\label{schwarz_subsec}

Every non-singular point on a real-analytic curve admits a local Schwarz reflection map.
A domain in the complex plane is called a \emph{quadrature domain} if the local Schwarz reflection maps with respect to its boundary extend anti-meromorphically to its interior.

\begin{definition}\label{qd_schwarz_def}
A domain $\Omega\subsetneq\widehat{\C}$ with $\infty\notin\partial\Omega$ and $\mathrm{int}(\overline{\Omega})=\Omega$ is called a \emph{quadrature domain} if there exists a continuous function $\sigma:\overline{\Omega}\to\widehat{\C}$ satisfying the following two properties:
\begin{enumerate}\upshape
\item $\sigma=\mathrm{id}$ on $\partial \Omega$.

\item $\sigma$ is anti-meromorphic on $\Omega$.
\end{enumerate}

The map $\sigma$ is called the \emph{Schwarz reflection map} of $\Omega$.
\end{definition}

The notion of quadrature domains first appeared in the work of Davis \cite{Dav74}, and independently in the work of Aharonov and Shapiro \cite{AS73,AS76,AS78}. It is known that except for a finite number of singular points, which are necessarily cusps and double points, the boundary of a quadrature domain consists of finitely many disjoint non-singular real-analytic curves \cite{Sak91}. 

Simplest examples of quadrature domains are given by round disks. Non-trivial examples are produced by the following characterization result.

\subsubsection{Characterization of simply connected quadrature domains and mapping properties of the associated Schwarz reflections}\label{scqd_schwarz_mapping_deg_subsubsec}

\begin{proposition}\label{simp_conn_quad_prop}\cite[Theorem~1]{AS76}
A simply connected domain $\Omega\subsetneq\widehat{\C}$ (with $\infty\notin\partial\Omega$ and $\mathrm{int}(\overline{\Omega})=\Omega$) is a quadrature domain if and only if the Riemann uniformization $f:\mathbb{D}\to\Omega$ is rational.
In this case, the Schwarz reflection map $\sigma$ of $\Omega$ is given by $f\circ\eta\circ(f\vert_{\D})^{-1}$, where $\eta(z):=1/\overline{z}$. 

Moreover, if the degree of the rational map $f$ is $d$, then $\sigma:\sigma^{-1}(\Omega)\to\Omega$ is a (branched) covering of degree $(d-1)$, and $\sigma:\sigma^{-1}(\Int{\Omega^c})\to\Int{\Omega}^c$ is a (branched) covering of degree~$d$ (where $\Omega^c:=\widehat{\C}\setminus \Omega$).
\end{proposition}
\[
  \begin{tikzcd}
    \D \arrow{d}{\eta} \arrow{r}{f} & \Omega \arrow{d}{\sigma} \\
     \widehat{\C} \arrow{r}{f}  & \widehat{\mathbb{C}}
  \end{tikzcd}
\] 

\begin{figure}[h!]
\captionsetup{width=0.96\linewidth}
\begin{tikzpicture}
\node[anchor=south west,inner sep=0] at (0,0) {\includegraphics[width=0.96\textwidth]{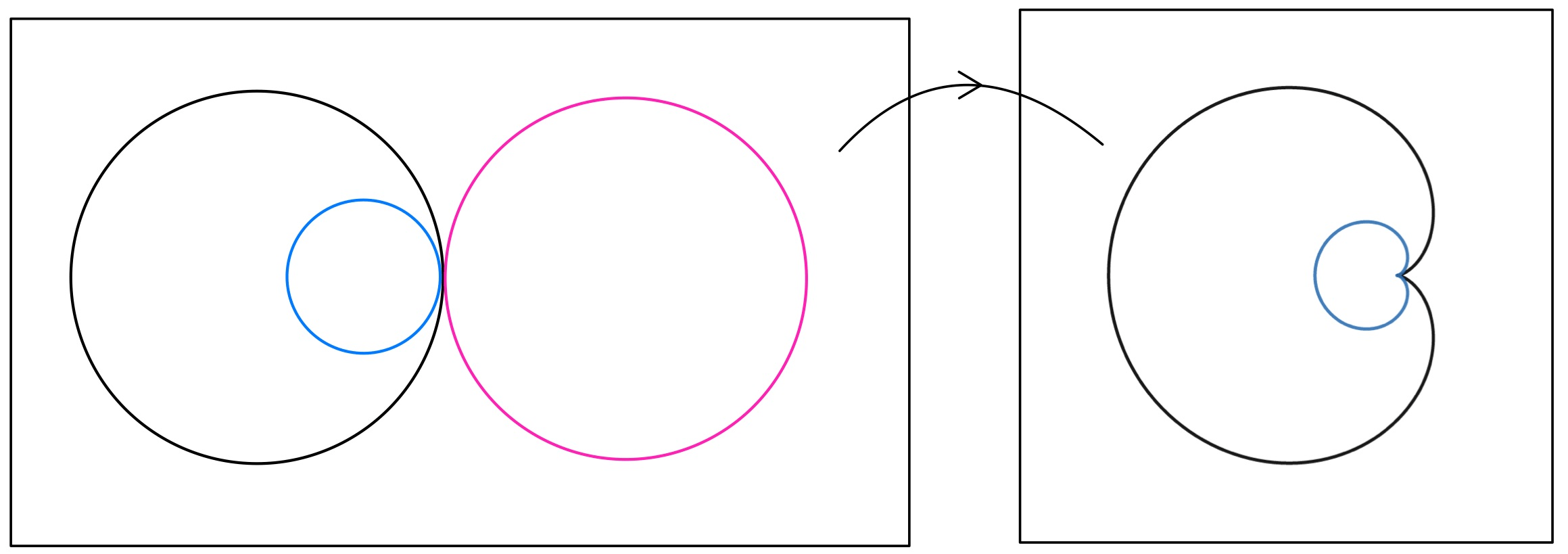}};
\node at (9.66,3) {\begin{scriptsize}$\sigma^{-1}(\Omega)$\end{scriptsize}};
\draw [->,line width=0.5pt] (9.8,2.75) to (10.5,2.2);
\node at (11.6,3.6) {\begin{scriptsize}$\Omega^c$\end{scriptsize}};
\node at (9.8,1.25) {\begin{scriptsize}$\sigma^{-1}(\Omega^c)$\end{scriptsize}};
\node at (11.56,0.55) {\begin{scriptsize}$\partial\Omega$\end{scriptsize}};
\draw [->,line width=0.5pt] (11.5,0.75) to (11.12,1.4);
\node at (7.5,3.2) {$f$};
\node at (0.8,0.9) {\begin{scriptsize}$\mathbb{S}^1$\end{scriptsize}};
\node at (3.5,3.8) {\begin{scriptsize}$f^{-1}(\Omega^c)$\end{scriptsize}};
\draw [<->,line width=0.5pt, out=300,in=210]  (2.8,1.5) to (4.2,0.8); 
\node at (1.4,2.2) {\begin{tiny}$\eta(f^{-1}(\Omega^c))$\end{tiny}};
\node at (3.5,0.5) {$\eta$};
\node at (4.8,2.1) {\begin{scriptsize}$X=f^{-1}(\Omega)\setminus\D$\end{scriptsize}};
\node at (2.8,2.15) {\begin{scriptsize}$\eta(X)$\end{scriptsize}};
\end{tikzpicture}
\caption{The map $f(z)=z/2-z^2/4$ is univalent on $\D$. The associated quadrature domain and Schwarz reflection map are denoted by $\Omega, \sigma$, respectively. The only critical point of $\sigma$ is at the origin. Left: $f^{-1}(\Omega)$ is a union of two touching disks; namely, $\D$ and $X$. Each of them maps univalently onto $\Omega$ under $f$. The complement of these two disks is $f^{-1}(\Omega^c)$. The image of $f^{-1}(\Omega^c)$ under the circular reflection $\eta$ is the croissant-shaped domain between $\mathbb{S}^1$ and the blue curve. On the other hand, the image of $X$ under $\eta$ is the small disk bounded by the blue circle. Right: The image of the disk $\eta(X)$ under $f$ is $\sigma^{-1}(\Omega)$ (the small cardioid bounded by the blue curve), which is univalently mapped by $\sigma$ onto the quadrature domain $\Omega$.
The image of $\eta(f^{-1}(\Omega^c))$ under $f$ is $\sigma^{-1}(\Omega^c)$ (the domain between the large and the small cardioids). It is mapped by $\sigma$ onto the droplet $\Omega^c$ as a $2:1$ branched cover branched only at $0$.}
\label{cardioid_disk_fig}
\end{figure}
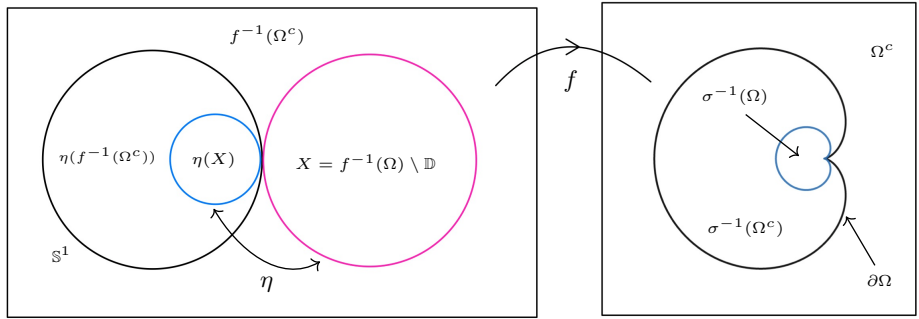

The second part of Proposition~\ref{simp_conn_quad_prop} says that the map $\sigma:\Omega\to\widehat{\C}$ is not a branched cover, but it restricts to branched covers of different degrees on two different part of its domain.
This can be seen as follows. 
Since $f$ is a degree $d$ global branched covering that maps $\D$ univalently onto $\Omega$, the map $f:f^{-1}(\Omega)\setminus\D\to\Omega$ is a branched covering of degree $d-1$. The definition of $\sigma$ implies that $\sigma^{-1}(\Omega)=f(\eta(f^{-1}(\Omega)\setminus\D))$ is mapped as a degree $d-1$ branched covering onto $\Omega$ by $\sigma$. On the other hand, the map $f:f^{-1}(\Int{\Omega^c})\to\Int{\Omega^c}$ is a branched covering of degree $d$. By definition of $\sigma$, we have that $\sigma^{-1}(\Int{\Omega^c})=f(\eta(f^{-1}(\Int{\Omega^c})))$ is mapped as a degree $d$ branched covering onto $\Int{\Omega^c}$ by $\sigma$ (see Figure~\ref{cardioid_disk_fig} for an illustration).

\begin{remark}
If $\Omega$ is a simply connected quadrature domain with associated Schwarz reflection map $\sigma$, and $M$ is a M{\"o}bius transformation, then $M(\Omega)$ is also a quadrature domain with Schwarz reflection map $M\circ\sigma\circ M^{-1}$. This fact allows one to normalize quadrature domains suitably, and will be used later.
\end{remark}

\subsubsection{Piecewise Schwarz reflections in quadrature multi-domains}\label{piecewise_schwarz_subsubsec}

Consider a finite collection of disjoint simply connected quadrature domains $\Omega_j(\subsetneq\widehat{\C})$, $j\in\{1,\cdots,k\}$, with associated Schwarz reflection maps $\sigma_j$. We define 
$$
\displaystyle\Omega:=\displaystyle\bigsqcup_{j=1}^k\Omega_j,
$$ 
and the map 
$$
\sigma:\overline{\Omega}\to\widehat{\C},\ w \mapsto \begin{array}{ll}
                    \sigma_j(w) & \mbox{if}\ w\in\overline{\Omega}_j
                                          \end{array}. 
$$
We call $\Omega$ a \emph{quadrature multi-domain} and $\sigma$ a piecewise Schwarz reflection map.

Let $f_j:\overline{\D}\to\overline{\Omega}_j$ be the Riemann uniformizations of the simply connected quadrature domains $\Omega_j$ such that each $f_j$ extends as a rational map of $\widehat{\C}$ of degree $d_j$. It follows that $\sigma_j:\sigma_j^{-1}(\Omega_j)\to\Omega_j$ is a branched covering of degree $(d_j-1)$, and $\sigma_j:\sigma_j^{-1}(\Int{\Omega_j^c})\to\Int{\Omega_j^c}$ is a branched covering of degree $d_j$. Therefore, $\sigma:\sigma^{-1}(\Omega)\to\Omega$ is a (possibly branched) covering of degree 
$$
d:=\displaystyle\sum_{j=1}^k d_j-1.
$$
We refer to $\sigma$ as a degree $d$ Schwarz reflection map.

\subsubsection{Invariant partition of the dynamical plane}\label{inv_partition_subsubsec}

With notation as in Subsection~\ref{piecewise_schwarz_subsubsec}, we~set 
$$
T(\sigma):=\widehat{\C}\setminus\Omega,\ \mathrm{and}\ T^0(\sigma):=T(\sigma)\setminus\{\mathrm{singular\ points\ on}\ \partial T\}.
$$ 
The set $T(\sigma)$ is called the \emph{droplet} of $\sigma$, and the set $T^0(\sigma)$ obtained by removing the singular points from the droplet is called the \emph{fundamental tile} of $\sigma$. We note that $T^0(\sigma)$ is neither open, nor closed. It resembles fundamental domains of kissing reflection groups (cf. Subsection~\ref{fund_dom_subsubsec}).
  
We define the \emph{escaping/tiling set} $T^\infty(\sigma)$ of $\sigma$ as the set of all points that eventually land in $T^0(\sigma)$; i.e., 
$$
T^\infty(\sigma):=\displaystyle\bigcup_{n=0}^\infty \sigma^{-n}(T^0(\sigma)).
$$ 
The \emph{non-escaping set} of $\sigma$ is defined as $K(\sigma):=\widehat{\C}\setminus T^\infty(\sigma)$.

Often, the dynamics of $\sigma$ on its non-escaping set resembles the dynamics of an anti-polynomial on its filled Julia set. On the other hand, the action of $\sigma$ on its escaping set, especially when the escaping set contains no critical point of $\sigma$, looks like the action of the Nielsen map of a necklace reflection group (for instance, see Figures~\ref{deltoid_corr_fig} (right) and~\ref{talbot_schwarz_fig}).

\subsubsection{Connections with analytic problems}\label{analysis_connect_subsubsec}
To motivate the nomenclature `quadrature domain', let us mention that a domain $\Omega\subsetneq\widehat{\C}$ (with $\infty\notin\partial\Omega$ and $\mathrm{int}(\overline{\Omega})=\Omega$) is a quadrature domain in the sense of Definition~\ref{qd_schwarz_def} if and only if there exists a rational map $R_\Omega$ with all poles inside $\Omega$ such that 
$$
\displaystyle \int_{\Omega} \phi dA=\frac{1}{2i} \oint_{\partial\Omega} \phi(z) R_{\Omega}(z) dz\quad \left(=\sum c_k \phi^{(n_k)}(a_k)\right)
$$ 
for all $\phi\in H(\Omega)\cap C(\overline{\Omega})$ (if $\infty\in\Omega$, one also requires the test function $\phi$ to vanish at $\infty$) \cite[Lemma~3.1]{LM}.
Such identities are called \emph{quadrature identities}, and they appear in various problems of complex analysis. The rational map $R_\Omega$ is called the \emph{quadrature function} of $\Omega$, and the poles of $R_\Omega$ are called the \emph{nodes} of $\Omega$. By \cite[Lemma~3.1]{LM}, the nodes of $R_\Omega$ are precisely the poles of the Schwarz reflection map~$\sigma$.

Areas of analysis where quadrature domains have found applications include quadrature identities \cite{Dav74, AS76, Sak82, Gus}, extremal problems for conformal mapping \cite{Dur83, ASS99, She00}, Hele-Shaw flows \cite{Ric72, EV92, GV06}, Richardson's moment problem \cite{Sak78, EV92, GHMP00}, free boundary problems \cite{Sha92, Sak91, CKS00}, subnormal and hyponormal operators \cite{GP17}, etc.

\subsubsection{Connections with statistical physics}\label{stat_phys_connect_subsubsec} 

Complements of quadrature domains naturally arise as accumulation sets of eigenvalues in random normal matrix models, and as accumulation sets of electrons in \emph{2D Coulomb gas models}. 

Consider $N$ electrons located at points $\lbrace z_j\rbrace_{j=1}^N$ in the complex plane, influenced by a strong ($2$-dimensional) external electrostatic field arising from a uniform non-zero charge density. Let the scalar potential of the external electrostatic field be $N\cdot Q:\C\to\R\cup\{+\infty\}$ (note that the scalar potential is rescaled so that it is proportional to the number of electrons). The combined energy of the system resulting from particle interaction and external potential is: 
$$
\mathfrak{E}_Q(z_1,\cdots,z_N)=\displaystyle\sum_{i\neq j} \ln\vert z_i-z_j\vert^{-1}+N\sum_{j=1}^N Q(z_j).
$$ 
In the equilibrium, the states of this system are distributed according to the Gibbs measure with density $$\frac{\exp(-\mathfrak{E}_Q(z_1,\cdots,z_N))}{Z_N},$$ where $Z_N$ is a normalization constant known as the \emph{partition function}. An important topic in statistical physics is to understand the limiting behavior of the `electron cloud' as the number of electrons $N$ grows to infinity. Under appropriate regularity conditions on $Q$, in the limit the electrons condensate on a compact subset $T$ of the plane, and they are distributed according to the normalized area measure of $T$ \cite{EF, HM}. Thus, the probability measure governing the distribution of the limiting electron cloud is completely determined by the shape of $T$, which is usually called a \emph{droplet}. If $Q$ is assumed to be algebraic in a suitable sense, the complementary components of the droplet $T$ turn out to be quadrature domains \cite{LM}. For example, the deltoid (the compact set bounded by the deltoid curve) is a droplet in the physically interesting case of the localized \emph{cubic external potential}, see \cite{Wie02}. 

The $2D$ Coulomb gas model described above is intimately related to logarithmic potential theory with an algebraic external field \cite{ST97} and the corresponding random normal matrix models, where the same probability measure describes the distribution of eigenvalues \cite{TBAZW, ABWZ}.

Iteration of Schwarz reflections sheds light on the topology of quadrature domains, and hence on their complementary regions. We refer the reader to \cite{LM} for details of this connection, or to \cite[\S 1.2]{LLMM2} for its brief account (cf. \cite{MR25}).

\section{Quadratic examples: deltoid, Circle-and-Cardioid, and Chebyshev cardioid}\label{quadratic_examples_sec}

In this section and the next, we will collect various explicit examples of antiholomorphic dynamical systems to illustrate the connections between anti-rational maps, Kleinian reflection groups, and hybrid dynamical systems that combine features of the former two objects in the same dynamical plane.

Let $\Omega$ be a quadrature multi-domain and $\sigma$ be the associated piecewise Schwarz reflection map as in Subsection~\ref{piecewise_schwarz_subsubsec}.
We say that $\sigma$ is \emph{quadratic} if $\sigma:\sigma^{-1}(\Omega)\to\Omega$ has degree two.
As in classical holomorphic dynamics, quadratic Schwarz reflection maps play a special role in our theory.

We use the notation of Subsection~\ref{piecewise_schwarz_subsubsec}.
When $\sigma:\sigma^{-1}(\Omega)\to\Omega$ has degree two, we have that
$$
\displaystyle\sum_{j=1}^k d_j-1=2\implies\displaystyle\sum_{j=1}^k d_j=3.
$$ 
Since each $d_j\geq 1$, it follows that $k\leq 3$. 
\smallskip

It turns out that such Schwarz reflection maps come in three interesting flavors \cite[\S 2.2]{LLMM3}. We briefly explain these possibilities below.
\smallskip

\noindent\textbf{Case 1: $k=1$.}  In this case, $\Omega=\Omega_1$ is a single quadrature domain that is the univalent image of $\D$ under some cubic rational map $f$. 
\smallskip

\noindent\textbf{Subcase 1.1.} Generically, $f_1$ has four simple critical points. A specific example of this type of quadrature domains is the exterior of a deltoid, whose associated Schwarz reflection map will be illustrated in Subsection~\ref{deltoid_subsec}. The dynamics of this map was completely described in \cite[\S 4]{LLMM1}. 
\smallskip

\noindent\textbf{Subcase 1.2.} Now suppose that the rational map $f_1$ has a unique double critical point. Then pre- and post-composing $f_1$ with M{\"o}bius maps, we can assume that $f_1(w)=w^3-3w$; and $\sigma$ is the Schwarz reflection map of $\Omega=f_1(\pmb{D})$, where $\pmb{D}$ is a round disk on which $f_1$ acts injectively. We will take a closer look at such a Schwarz reflection map in Subsection~\ref{chebyshev_subsec}. 

In fact, restricting the cubic Chebyshev map to various disks of univalence yields a natural one-parameter family of Schwarz reflections, which was studied in details in \cite{LLMM3}.
\smallskip

\noindent\textbf{Subcase 1.3.} In the final subcase, suppose that the rational map $f_1$ has two double critical points. Then pre- and post-composing $f_1$ with M{\"o}bius maps, we can assume that $f_1(w)=w^3$, and $\sigma$ is the Schwarz reflection map of $\Omega=f_1(\pmb{D})$, where $\pmb{D}$ is a round disk on which $f_1$ acts injectively. Since $\pmb{D}$ does not intersect $\{0,\infty\}$, one easily sees that $\Omega$ contains no critical value of $\sigma$. Thus in this case, $\sigma:\sigma^{-1}(\Omega)\to\Omega$ has no critical point and hence is dynamically uninteresting.
\smallskip

\noindent\textbf{Case 2: $k=2$.} We can assume that $\Deg{f_1}=2$, and $\Deg{f_2}=1$. As before, pre- and post-composing $f_1, f_2$ with M{\"o}bius maps, we can assume that $f_1(w)=w^2$, $\Omega_1$ is the univalent image of a round disk under $f_1$, and $\Omega_2$ is a round disk in $\widehat{\C}$. Of particular interest is the situation when $\Omega_1$ is a cardioid and $\Omega_2$ is the exterior of a circumcircle of $\partial\Omega_1$. 
We consider the simplest dynamically interesting map of this kind in Subsection~\ref{c_and_c_center_subsec}.

The moduli space of all Schwarz reflection maps obtained by fixing a cardioid as $\Omega_1$ and varying the center of the exterior disk $\Omega_2$ such that $\partial \Omega_2$ touches $\partial\Omega_1$ at a unique point produces the Circle-and-Cardioid family, which was the main topic of investigation in \cite[\S 5,6]{LLMM1} and \cite{LLMM2}.
\smallskip

\noindent\textbf{Case 3: $k=3$.} In this case, each $f_j$ is a M{\"o}bius map, and hence each $\Omega_j$ is a round disk. In particular, each $\sigma_j$ is the reflection in a round circle (thus, $\sigma$ has no critical point), and the resulting dynamics of $\sigma$ is completely understood.

\subsection{The deltoid reflection map}\label{deltoid_subsec}

Suppose that $f$ is a cubic rational map that is univalent on $\D^*=\widehat{\C}\setminus\overline{\D}$. Post-composing $f$ with a M{\"o}bius map, we can assume that $f(\infty)=\infty$ and $Df(\infty)=1$. Let us further assume that $f$ has a $2\pi/3$-rotation symmetry; i.e., $f$ commutes with the map $w\mapsto e^{\frac{2\pi i}{3}} w$. Then $f$ must be of the form 
$$
f_t(w)=w + \frac{t}{2w^2},\quad t\neq 0.
$$
The assumption that $f_t$ is univalent on $\D^*$ implies that critical points of $f_t$ lie in $\overline{\D}$, and hence $\vert t\vert\leq 1$.

It is easy to check that for $\vert t\vert\leq 1$, each $f_t$ is indeed univalent on $\D^*$ (cf. \cite[Proposition~B.1]{LLMM3}).
The critical points of $f_t$ lie at the origin and the third roots of $t$. 

\subsubsection{Deltoid reflection map as a limit of anti-quadratic-like maps}\label{deltoid_degeneration_subsubsec}
 Let us now try to understand the dynamics of the associated Schwarz reflection maps 
 $$
 \sigma_t:\Omega_t:=f_t(\D^*)\to\widehat{\C},\qquad t\in(0,1).
 $$ 

We first note that $\infty$ is a superattracting fixed point for each $\sigma_t$.
All critical points of $f_t$ lie in $\D$ and hence $\partial\Omega_t$ is a non-singular Jordan curve. It follows from local properties of Schwarz reflection maps that $\sigma_t^{-1}(\Omega_t)\Subset\Omega_t$. Moreover, since none of the finite co-critical points of $f_t$ (namely, $-\frac{t^{1/3}}{2}, -\frac{t^{1/3}\omega}{2}, -\frac{t^{1/3}\omega^2}{2}$) lie in $\D^*$, it follows that no finite critical value of $f_t$ lies in $\Omega_t$. The commutative diagram defining $\sigma_t$ now implies that $\sigma_t:\sigma_t^{-1}(\Omega_t)\to\Omega_t$ is a $2:1$ branched cover branched only at $\infty$. By Riemann-Hurwitz, $\sigma_t^{-1}(\Omega_t)$ must be a simply connected domain; and consequently, $\sigma_t:\sigma_t^{-1}(\Omega_t)\to\Omega_t$ is an anti-quadratic-like map (cf. \cite[\S 4]{LM}). Since $\sigma_t$ has a fixed critical point at $\infty$, we conclude that the non-escaping dynamics of $\sigma_t$ (put differently, the dynamics of the anti-quadratic-like map $\sigma_t:\sigma_t^{-1}(\Omega_t)\to\Omega_t$ on its filled Julia set) is hybrid conjugate to $\overline{z}^2$ for all $c\in(0,1)$.
\begin{figure}[h!]
\captionsetup{width=0.96\linewidth}
\begin{tikzpicture}
\node[anchor=south west,inner sep=0] at (0,0) {\includegraphics[width=0.96\textwidth]{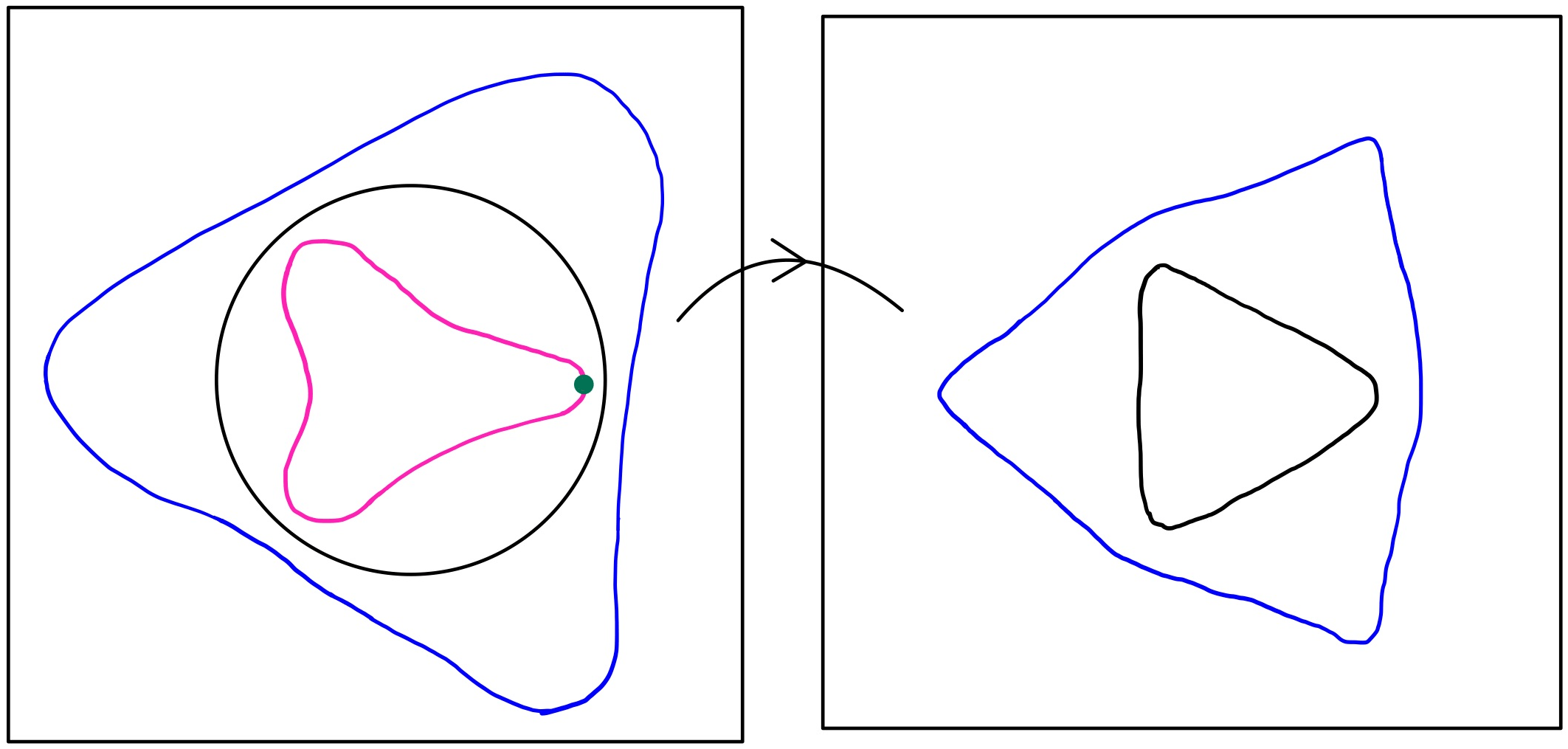}};
\node at (9.6,2.7) {$\Omega_t^c$};
\node at (10,4) {$\sigma_t^{-1}(\Omega_t^c)$};
\node at (8,4.8) {$\sigma_t^{-1}(\Omega_t)$};
\node at (7.5,1) {$\partial\Omega_t$};
\draw [->,line width=0.5pt] (7.7,1.2) to (8.75,2);
\node at (6.1,3.3) {$f_t$};
\node at (1,1) {$\mathbb{S}^1$};
\draw [->,line width=0.5pt] (1,1.2) to (1.8,2.1);
\node at (3.6,3.88) {\begin{tiny}$f_t^{-1}(\Omega_t^c)$\end{tiny}};
\node at (3.3,2.8) {\begin{tiny}$f_t^{-1}(\Omega_t)\cap\D$\end{tiny}};
\draw [<->,line width=0.5pt, out=210,in=330]  (2.16,3.3) to (0.8,3.45); 
\node at (1.2,3) {$\eta$};
\node at (4,2) {\begin{footnotesize}$\alpha_t$\end{footnotesize}};
\draw [->,line width=0.5pt, out=45,in=270]  (4.06,2.12) to (4.58,2.75);
\end{tikzpicture}
\caption{The conformal annulus $\Omega_t\setminus\overline{\sigma_t^{-1}(\Omega_t)}=\sigma_t^{-1}(\Int{\Omega_t^c})$ is anti-conformally equivalent to the conformal annulus $f_t^{-1}(\Int{\Omega_t^c})$. The green marked point $\alpha_t:=\frac{t+\sqrt{t^2+8t}}{4}$ satisfies $f_t(\alpha_t)=f_t(1)$, and hence $\alpha_t\in\partial f_t^{-1}(\Int{\Omega_t^c})$. Since $\alpha_t\to 1$ as $t\to 1^-$, it follows that the moduli of these annuli tend to $0$.}
\label{pre_deltoid_disk_fig}
\end{figure}
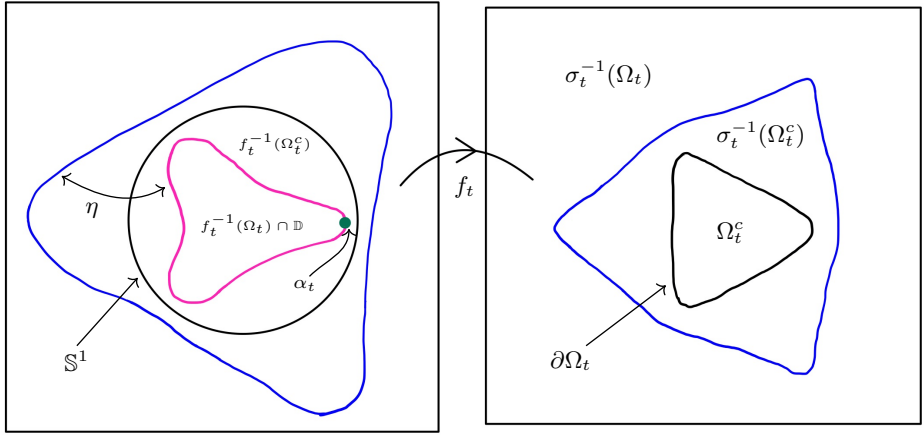

Understanding the escaping dynamics of $\sigma_t$ is equivalent to understanding the external map of the above anti-quadratic-like map. Once again, by Riemann-Hurwitz and Proposition~\ref{simp_conn_quad_prop}, $\sigma_t:\sigma_t^{-1}(\Int{\Omega_t^c})\to\Int{\Omega_t^c}$ is an annulus-to-disk branched covering of degree three. It is not hard to see that the moduli of these anti-quadratic-like maps go to $0$ (see Figure~\ref{pre_deltoid_disk_fig}).

Thus, $\{\sigma_t:t\in(0,1)\}$ defines an escaping path in the space of anti-quadratic maps. The Schwarz reflection map $\sigma_1$ naturally appears as a limit point of this path. The associated quadrature domain $f_1(\D^*)$ is the exterior of the classical deltoid curve, where the cusps of the deltoid curve arise from the three critical points of $f_1$ at the third roots of unity (see Figure~\ref{deltoid_disk_fig}).

\subsubsection{Dynamics of deltoid reflection}\label{deltoid_limit_subsubsec}

We now set $f:=f_1$, $\Omega:=f(\D^*)$, and investigate the dynamics of $\sigma:=\sigma_1:\overline{\Omega}\to\widehat{\C}$. The corresponding Schwarz reflection map $\sigma$ can be thought of as an object lying halfway between a quasi-Fuchsian group and a quasi-Blaschke product.
Indeed, $\sigma$ combines the action of a reflection group and an anti-Blaschke product in the sense that it acts like the anti-Blaschke product $\overline{z}^2\vert_{\overline{\D^*}}$ on $K(\sigma)$, and acts like the Nielsen map $\pmb{\cN}_2:\overline{\D}\setminus\Int{\Pi(\pmb{G}_2)}\to\overline{\D}$ of the ideal triangle reflection group on $\overline{T^\infty(\sigma)}$. It is the simplest example of an antiholomorphic map exhibiting such hybrid dynamics.
\begin{figure}[h!]
\captionsetup{width=0.96\linewidth}
\begin{tikzpicture}
\node[anchor=south west,inner sep=0] at (0,0) {\includegraphics[width=0.96\textwidth]{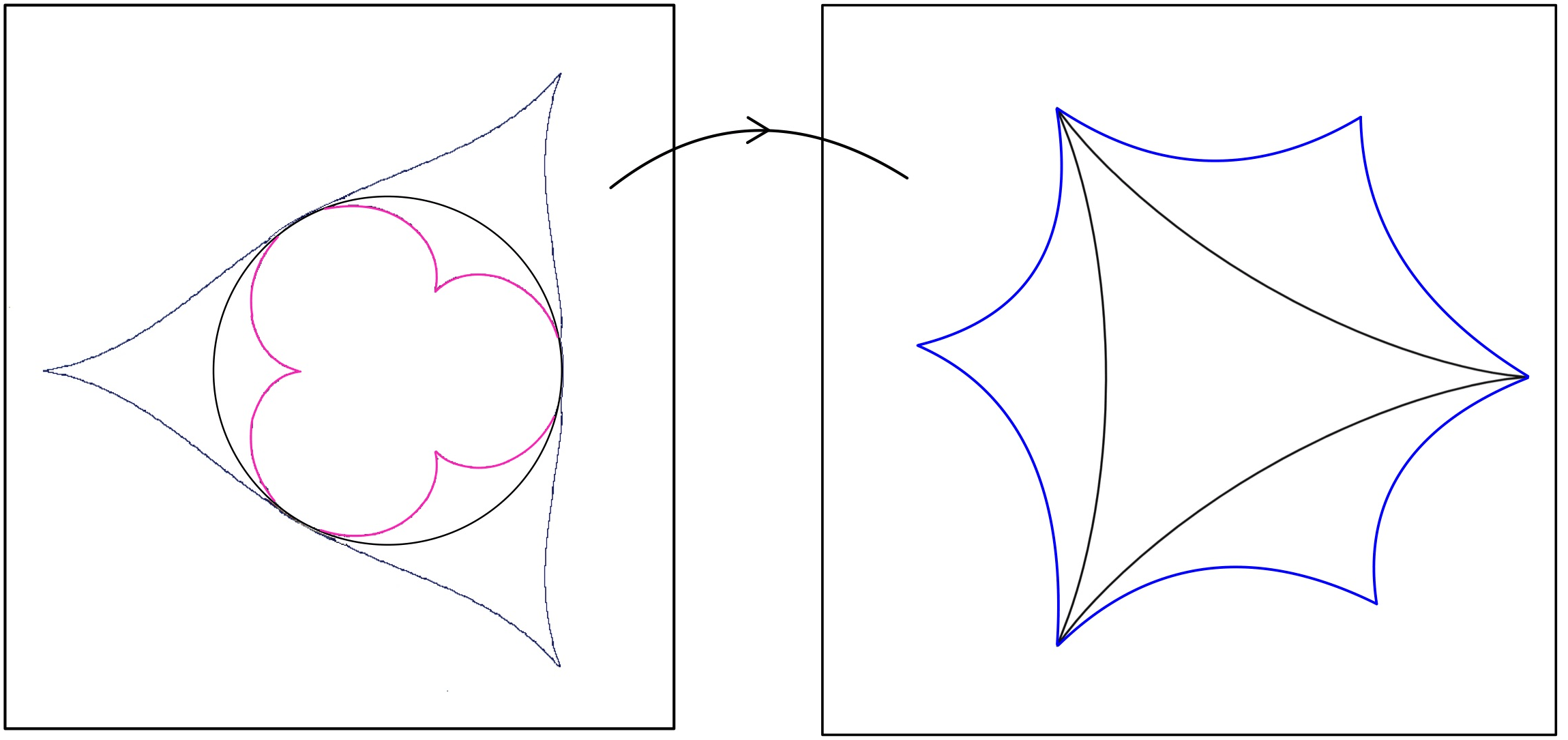}};
\node at (9.5,2.9) {\begin{scriptsize}$\Omega^c$\end{scriptsize}};
\node at (10,4) {\begin{tiny}$\sigma^{-1}(\Omega^c)$\end{tiny}};
\node at (8,3) {\begin{tiny}$\sigma^{-1}(\Omega^c)$\end{tiny}};
\node at (10,1.66) {\begin{tiny}$\sigma^{-1}(\Omega^c)$\end{tiny}};
\node at (11,5.25) {\begin{scriptsize}$\sigma^{-1}(\Omega)$\end{scriptsize}};
\node at (7.2,0.8) {\begin{scriptsize}$\partial\Omega$\end{scriptsize}};
\draw [->,line width=0.5pt] (7.4,1) to (8.4,1.9);
\node at (5.8,4) {$f$};
\node at (0.8,1) {\begin{scriptsize}$\mathbb{S}^1$\end{scriptsize}};
\draw [->,line width=0.5pt] (0.8,1.2) to (1.66,2.42);
\node at (3.6,3.88) {\begin{tiny}$X_1$\end{tiny}};
\node at (3.6,1.9) {\begin{tiny}$X_2$\end{tiny}};
\node at (1.88,2.88) {\begin{tiny}$X_3$\end{tiny}};
\node at (3.3,2.9) {\begin{tiny}$f^{-1}(\Omega)\cap\D$\end{tiny}};
\node at (1.2,2.9) {\begin{tiny}$\eta(X_3)$\end{tiny}};
\node at (3.75,4.4) {\begin{tiny}$\eta(X_1)$\end{tiny}};
\node at (3.75,1.36) {\begin{tiny}$\eta(X_2)$\end{tiny}};
\node at (1.4,5) {\begin{scriptsize}$\eta(f^{-1}(\Omega)\cap\D)$\end{scriptsize}};
\end{tikzpicture}
\caption{Left: The central simply connected domain $f^{-1}(\Omega)\cap\D$ (bounded by the pink curve) contains the critical point $0$ of $f$, and maps as a $2:1$ branched cover onto $\Omega$ under $f$. The three components $X_1, X_2, X_3$ of $f^{-1}(\Omega^c)$ and their images under the reflection map $\eta$ are marked. Right: The images of the topological triangles $\eta(X_i),\ i=1,2,3$ under $f$ comprise $\sigma^{-1}(\Omega^c)$, which is the preimage of the droplet under the deltoid Schwarz reflection. On the other hand, the image of the simply connected domain $\eta(f^{-1}(\Omega)\cap\D)$ under $f$ is $\sigma^{-1}(\Omega)$, the preimage of the quadrature domain under the Schwarz reflection map. Thus, $\sigma:\sigma^{-1}(\Omega)\to\Omega$ is a `degenerate' anti-quadratic-like map with pinching points at the three cusps of $\partial\Omega$.}
\label{deltoid_disk_fig}
\end{figure}
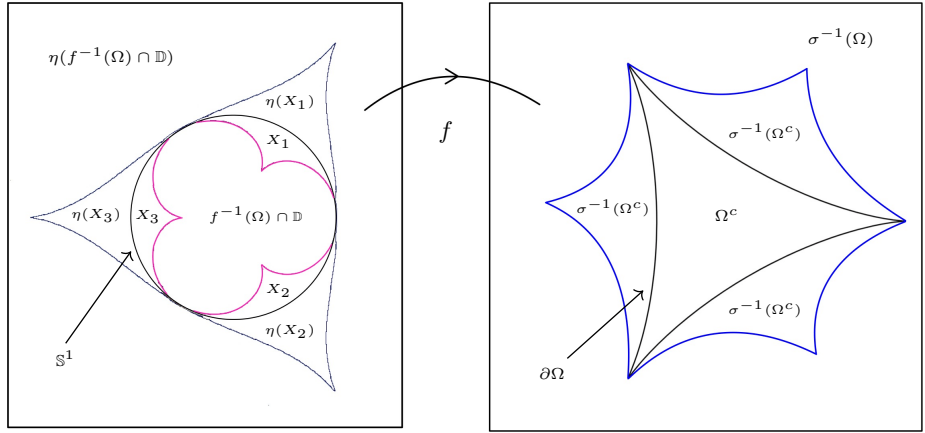

We now explain this mating phenomenon in a bit more detail.
The map $\sigma$ has a unique critical point at $\infty$ (which comes from the critical point $0$ of $f$), and this point is fixed by $\sigma$. We denote the basin of attraction of the superattracting fixed point $\infty$ by $\mathcal{B}_\infty(\sigma)$. Since $\sigma$ is unicritical with $\sigma^{-1}(\infty)=\{\infty\}$, it is natural to expect that the non-escaping set $K(\sigma)$ equals $\overline{\mathcal{B}_\infty(\sigma)}$, where $\mathcal{B}_\infty(\sigma)$ is a completely invariant simply connected domain on which $\sigma$ is conformally conjugate to $\overline{z}^2$. However, unlike the maps $\sigma_t$ studied in Subsection~\ref{deltoid_degeneration_subsubsec}, the restriction $\sigma:\sigma^{-1}(\Omega)\to\Omega$ is not an anti-quadratic-like map since $\partial\Omega$ intersects $\partial\sigma^{-1}(\Omega)$ at the three cusps of $\partial\Omega$, which are fixed points of $\sigma$ (see Figure~\ref{deltoid_disk_fig}). Hence, one cannot appeal to standard straightening theorems to analyze the non-escaping dynamics of $\sigma$.

\begin{remark}\label{deltoid_no_hybrid_rem}
In fact, the Puiseux series expansion of $\sigma$ at each of the three cusps of $\partial\Omega$ exhibit parabolic behavior, and the `attracting directions' at these three fixed points lie in the tiling set \cite[\S 4.2.1]{LLMM1}. As all fixed points of the map $\overline{z}^2$ on $\mathbb{S}^1$ are repelling, there \emph{cannot} be a hybrid conjugacy between the non-escaping dynamics of $\sigma$ and $\overline{z}^2\vert_{\overline{\D}^*}$. Another way of seeing the non-existence of such a hybrid conjugacy is to observe that the external class of $\sigma$ is $\pmb{\cN}_2$, which has parabolic fixed points, whereas the external class of an anti-polynomial is uniformly expanding.
\end{remark}

On the other hand, since the only critical point of $\sigma$ is fixed under dynamics, the dynamics on the escaping set $T^\infty(\sigma)$ is unramified (see Subsection~\ref{inv_partition_subsubsec} for the definition). Heuristically, this means that $\sigma\vert_{T^\infty(\sigma)}$ behaves like reflections in the three non-singular real-analytic arcs of $\partial\Omega$. In fact, the mapping degrees of $\sigma$ imply that the fundamental tile $T^0(\sigma)=\Omega^c\setminus\{f(1),f(\omega),f(\omega^2)\}$ (also called the rank zero tile) pulls back under $\sigma$ to three disjoint Jordan regions which are adjacent to the rank zero tile (see Figure~\ref{deltoid_disk_fig}). We call them the rank one tiles. Since the rank one tiles lie in $\Omega$, each of them has two preimages, and this pattern continues for all higher rank tiles. 
This fact can be used to deduce that the tiling set $T^\infty(\sigma)$ is a topological disk, and that $T^\infty(\sigma)\setminus T^0(\sigma)$ has three components $H_1, H_2$, and $H_3$. Let $L_i:=\partial H_i\cap T^\infty(\sigma)$ be the common edge of $H_i$ and $T^0(\sigma)$ without the cusps. Moreover, for $i=1,2,3$, the map $\sigma$ sends $H_i\cup L_i$ anti-conformally onto $T^\infty(\sigma)\setminus H_i$, and fixes $L_i$ pointwise. Hence, $\sigma_{H_i\cup L_i}$ can be extended to an anti-conformal involution $\sigma_i$ of $T^\infty(\sigma)$ pointwise fixing $L_i$. It follows that the Riemann map of $T^\infty(\sigma)$ conjugates $\sigma_i$ to an anti-M{\"o}bius involution of $\D$, and carries $T^0(\sigma)$ to an ideal hyperbolic triangle in $\D$. Since all ideal hyperbolic triangles are conformally equivalent, we can choose the Riemann map of $T^\infty(\sigma)$ so that it sends $T^0(\sigma)$ to $\Pi(\pmb{G}_2)$, and conjugates the maps $\sigma_i$ to the anti-M{\"o}bius reflections in the sides of $\Pi(\pmb{G}_2)$.

\begin{figure}[h!]
\captionsetup{width=0.96\linewidth}
\begin{tikzpicture}
\node[anchor=south west,inner sep=0] at (0,0) {\includegraphics[width=0.48\textwidth]{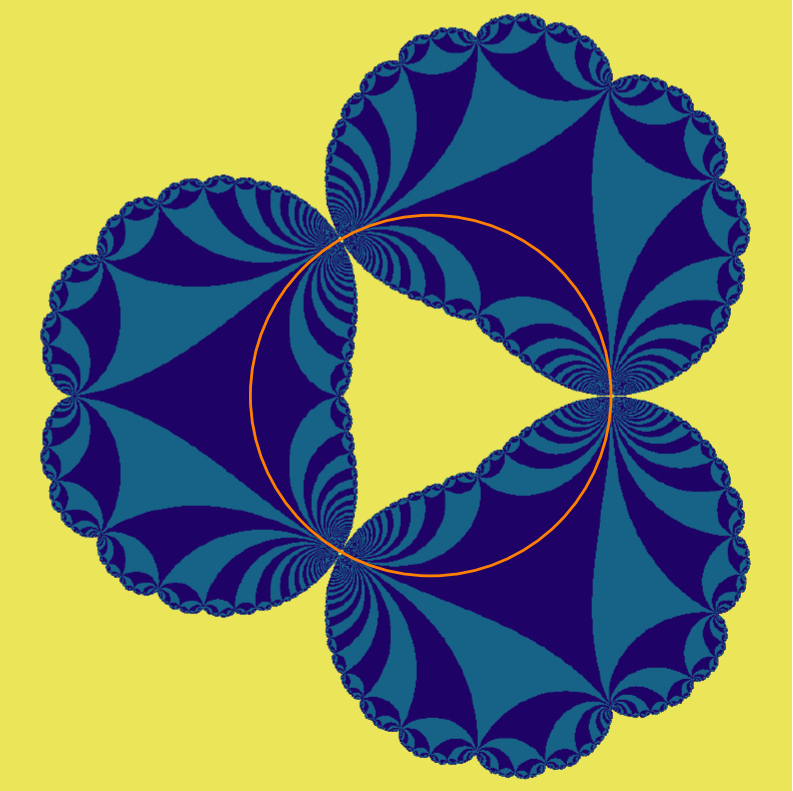}};
\node[anchor=south west,inner sep=0] at (6.6,0) {\includegraphics[width=0.48\textwidth]{deltoid_reflection_julia.png}};
\node[anchor=south west,inner sep=0] at (3.6,-6.75) {\includegraphics[width=0.48\textwidth]{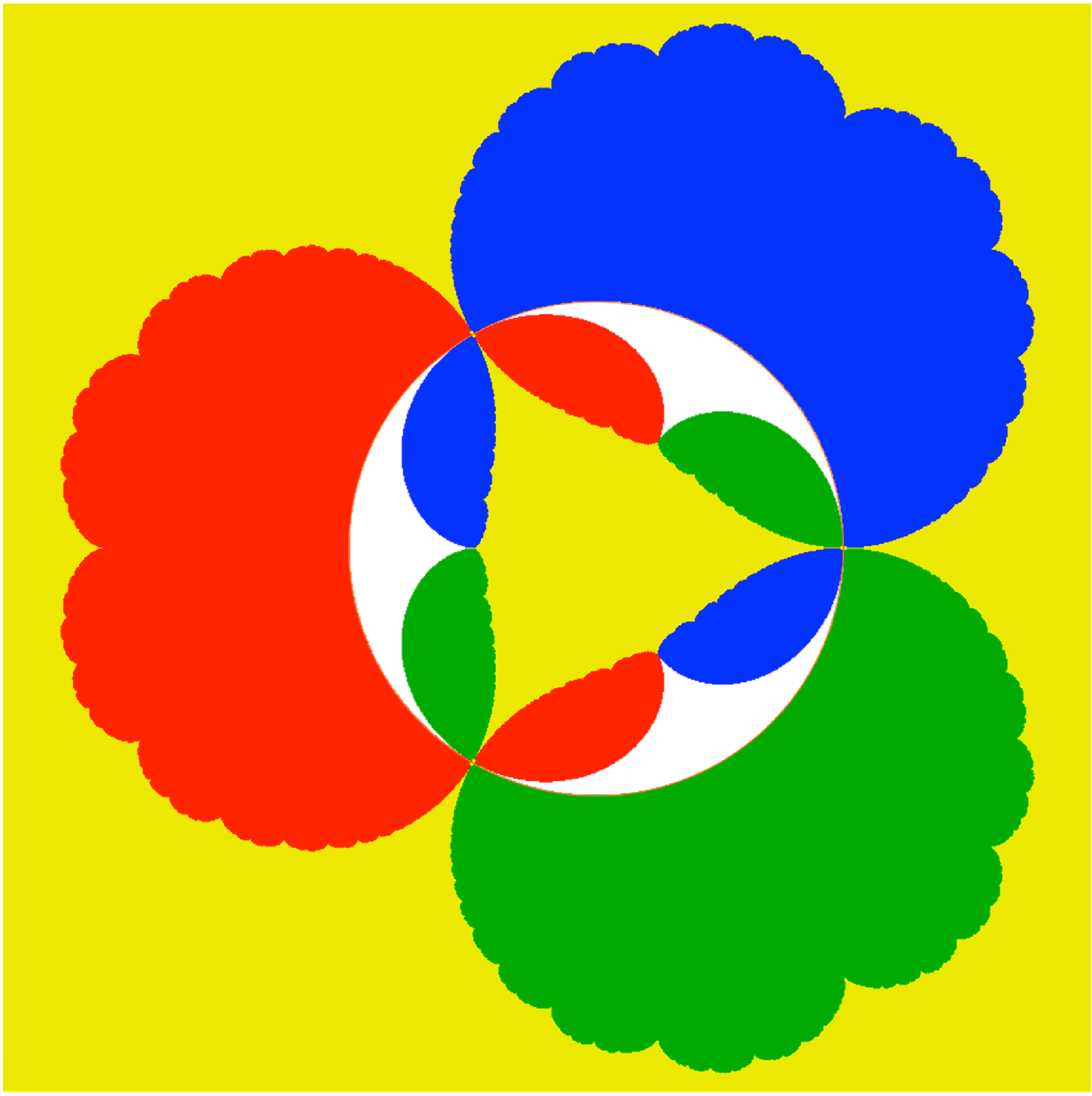}};
\node at (8.9,2.7) {\textcolor{white}{$\sigma$}};
\node at (10.6,2.4) {\textcolor{white}{$\sigma$}};
\node at (10,4) {\textcolor{white}{$\sigma$}};
\end{tikzpicture}
\caption{Top: The lifted tiling and non-escaping sets of the deltoid correspondence (left) are obtained by pulling back the tiling and non-escaping sets of the deltoid Schwarz reflection map (right) under $f$. Bottom: The lifted tiling set $\widetilde{T^\infty(\sigma)}$ of the correspondence (the complement of the yellow region) consists of three components, where each component comprises four regions shaded in white, blue, green, and red. The order three deck transformation $\tau$ of $f:\widetilde{T^\infty(\sigma)}~\to~T^\infty(\sigma)$ permutes these components cyclically preserving the coloring. Specifically, $\tau$ is given by the composition of the rigid rotation by angle $2\pi/3$ (around the origin) and an order three conformal rotation in each component of $\widetilde{T^\infty(\sigma)}$. Note the local symmetry of the picture under the reflections with respect to the third roots of unity.}
\label{deltoid_corr_fig}
\end{figure}

Adapting classical Fatou-Julia theory and \emph{puzzle piece} techniques for the setting of Schwarz reflection maps, the above picture was confirmed in \cite{LLMM1}. It was further shown that $\mathcal{B}_\infty(\sigma)$ and $T^\infty(\sigma)$ are Jordan domains with a common boundary (see Figure~\ref{deltoid_corr_fig}).  

\begin{theorem}[Dynamics of deltoid reflection]\cite[Theorem~1.1, \S 4.1, \S 4.2]{LLMM1}\label{deltoid_thm_1} 
The dynamical plane of the Schwarz reflection $\sigma$ of the deltoid can be partitioned as
$$
\widehat {\mathbb C}=T^\infty(\sigma)\sqcup \Lambda(\sigma)\sqcup \mathcal{B}_\infty(\sigma),
$$
where $T^\infty(\sigma)$ is the tiling set, $\mathcal{B}_\infty(\sigma)$ is the basin of infinity, and $\Lambda(\sigma)$ is their common boundary (which we call the \emph{limit set}). Moreover, $\sigma:\overline{T^\infty(\sigma)}\setminus\Int{T^0(\sigma)}\to \overline{T^\infty(\sigma)}$ is conformally conjugate to $\pmb{\cN}_2:\overline{\D}\setminus\Int{\Pi(\pmb{G}_2)}\to\overline{\D}$, the map $\sigma:\overline{\mathcal{B}_\infty(\sigma)}\to\overline{\mathcal{B}_\infty(\sigma)}$ is conformally conjugate to $\overline{z}^2:\overline{\D^*}\to\overline{\D^*}$, and $\Lambda(\sigma)$ is a Jordan curve.
\end{theorem}

\subsubsection{Deltoid reflection as the unique conformal mating of $\overline{z}^2$ and $\pmb{\cN}_2$}\label{deltoid_unique_mating_subsubsec}

Recall that the Minkowski circle homeomorphism $\pmb{\mathcal{E}}_2$ conjugates $\pmb{\cN}_2$ to $\overline{z}^2$, and sends $1$ to $1$. One can topologically glue $\overline{\D}$ and $\overline{\D^*}$, equipped with the two dynamical systems 
$$
\pmb{\cN}_2:\overline{\D}\setminus\Int{\Pi(\pmb{G}_2)}\to\overline{\D}\quad \mathrm{and}\quad \overline{z}^2:\overline{\D^*}\to\overline{\D^*}
$$ 
respectively, along the unit circle using the homeomorphism $\pmb{\mathcal{E}}_2$. This yields a partially defined continuous map on the topological $2$-sphere, and this map is called the \emph{topological mating} of $\overline{z}^2$ and $\pmb{\cN}_2$.

The deltoid Schwarz reflection map $\sigma$ is a \emph{conformal mating} of $\overline{z}^2$ and $\pmb{\cN}_2$ in the sense that there exists a conformal structure on the above topological $2$-sphere whose uniformizing map $\mathbb{S}^2\longrightarrow\widehat{\C}$ conjugates the topological mating to $\sigma$. 

In fact, a careful analysis of the geometry of $T^\infty(\sigma)$ reveals that it is a John domain, and hence its boundary is conformally removable (cf. \cite{Jon95}, Appendix~\ref{removable_david_subsec}). Conformal removability of $\Lambda(\sigma)$, other than being important from the complex-analytic viewpoint, implies that the conformal mating of $\overline{z}^2$ and $\pmb{\cN}_2$ is unique up to M{\"o}bius conjugation. We summarize these results below.

\begin{theorem}[Deltoid reflection as unique conformal mating]\cite[Theorem~1.1, \S 4.3, \S 4.4]{LLMM1}\label{deltoid_thm_2} 
The tiling set $T^\infty(\sigma)$ is a John domain, and hence the limit set $\Lambda(\sigma)$ is a conformally removable Jordan curve. Consequently, $\sigma$ is the unique conformal mating of the reflection map $\pmb{\cN}_2: \overline{\mathbb D}\setminus\Int{\Pi(\pmb{G}_2)}\to  \overline{\mathbb D}$ and the anti-polynomial $\overline{z}^2\vert_{\overline{\D^*}}$, up to M{\"o}bius conjugation. 
\end{theorem}

There is an alternative way of concluding uniqueness of the above conformal mating that does not appeal to conformal removability of the limit set.
In fact, any conformal mating $\widetilde{\sigma}$ of $\overline{z}^2\vert_{\overline{\D}^*}$ and $\pmb{\cN}_2:\overline{\D}\setminus\Int{\Pi(\pmb{G}_2)}\to\overline{\D}$ is an antiholomorphic map defined on the closure of a Jordan domain in the Riemann sphere. Indeed, since $\pmb{\cN}_2$ is not defined on the interior of the ideal triangle $\D\cap\Pi(\pmb{G}_2)$, any conformal mating $\widetilde{\sigma}$ must be defined on the complement of a homeomorphic copy of $\D\cap\Pi(\pmb{G}_2)$. Moreover, as $\pmb{\cN}_2$ fixes $\partial \Pi(\pmb{G}_2)$ pointwise, $\widetilde{\sigma}$ must fix the boundary of its domain of definition pointwise. It follows that the domain of definition is the closure of a simply connected quadrature domain, and $\widetilde{\sigma}$ is the associated Schwarz reflection map. Since the Nielsen map $\pmb{\cN}_2$ as well as the anti-polynomial $\overline{z}^2$ commute with the rotation $w\mapsto e^{\frac{2\pi i}{3}} w$, it follows that their mating $\widetilde{\sigma}$ also commutes with the same rotation. This in turn implies that the quadrature domain defining the Schwarz reflection $\widetilde{\sigma}$ is $2\pi/3-$rotation symmetric, and hence the uniformizing rational map of this quadrature domain commutes with $2\pi/3-$rotation. This symmetry property, combined with the fact that the Schwarz reflection map is quadratic, can be used to conclude that the uniformizing rational map can be chosen to be $f(w)=w+1/2w^2$; i.e., $\widetilde{\sigma}$ is the Schwarz reflection of the deltoid (see the discussion in the beginning of Subsection~\ref{deltoid_subsec}, cf. \cite[Appendix~A]{LLMM2}).
The above discussion also shows that Schwarz reflections are natural objects from the point of view of mating reflection groups with anti-polynomials.

\begin{remark}
If one considers cubic rational maps that are univalent on $\D^*=\widehat{\C}\setminus\overline{\D}$, have three critical points of $\mathbb{S}^1$, but are not symmetric with respect to $2\pi/3-$rotation, one obtains distorted versions of the deltoid whose Schwarz reflection maps give rise to matings of the Nielsen map $\pmb{\cN}_2$ with quadratic anti-Blaschke products with an attracting (but not superattracting) fixed point in $\D$.
\end{remark}

\subsubsection{A non-quasisymmetric welding map}\label{question_mark_welding_subsubsec}

The existence of the deltoid Schwarz reflection implies that the welding problem for the non-quasisymmetric circle homeomorphism $\pmb{\mathcal{E}}_2$ has a unique solution (see Subsection~\ref{question_mark_subsubsec}).

Specifically, according to Theorem~\ref{deltoid_thm_1}, there exist homeomorphic extensions of conformal maps $\phi^{\mathrm{in}}:\overline{\D}\to \overline{T^\infty(\sigma)}, \phi^{\mathrm{out}}:\overline{\D^*}\to \overline{\mathcal{B}_\infty(\sigma)}$, where the former conjugates $\pmb{\cN}_2$ to $\sigma$ and the latter conjugates $\overline{z}^2$ to $\sigma$. After possibly pre-composing $\phi^{\mathrm{in}}$ with a rotation, one can assume that $\phi^{\mathrm{in}}(1)=\phi^{\mathrm{out}}(1)$. Then, $\left(\phi^{\mathrm{out}}\right)^{-1}\circ\phi^{\mathrm{in}}:\mathbb{S}^1\to\mathbb{S}^1$ conjugates $\pmb{\cN}_2$ to $\overline{z}^2$, and sends $1$ to $1$. Hence, we have that $\pmb{\mathcal{E}}_2=\left(\phi^{\mathrm{out}}\right)^{-1}\circ\phi^{\mathrm{in}}$. It follows that $\pmb{\mathcal{E}}_2$ is the welding homeomorphism associated with the Jordan curve $\Lambda(\sigma)$. Moreover, conformal removability of $\Lambda(\sigma)$ implies that $\Lambda(\sigma)$ is the unique solution to the welding problem for the Minkowski circle homeomorphism $\pmb{\mathcal{E}}_2$ (up to M{\"o}bius maps).

A general result to the effect that circle homeomorphisms conjugating piecewise analytic expansive circle maps are welding maps will be discussed in Subsection~\ref{welding_subsec}.

\subsubsection{Deltoid reflection via David surgery}\label{deltoid_david_subsubsec}

Theorem~\ref{deltoid_thm_2} yields an \emph{unmating} of the deltoid reflection map into an anti-polynomial and the Nielsen map of a reflection group. Reversing the point of view, one can ask whether the topological mating of $\overline{z}^2$ and $\pmb{\cN}_2$ can be upgraded to a conformal mating without having prior knowledge of the deltoid reflection map. A direct construction of such conformal matings would require one to appeal to a suitable uniformization theorem. However, the fact that one is trying to combine an expanding dynamical system with one that has parabolic fixed points renders standard quasiconformal tools (e.g., techniques used in the proof of Bers Simultaneous Uniformization Theorem or in the construction of quasi-Blaschke products as matings of two Blaschke products) inapplicable to this setting.

In \cite[\S 9]{LLMM4}, number-theoretic properties of the map $\pmb{\mathcal{E}}_2$ (see Subsection~\ref{question_mark_subsubsec}) were used to obtain distortion estimates for $\pmb{\mathcal{E}}_2^{-1}$ and to conclude that the inverse of the Minkowski circle homeomorphism $\pmb{\mathcal{E}}_2$ continuously extends to a David homeomorphism of $\D$. This result, combined with the David Integrability Theorem (see Theorem~\ref{david_integrability_thm}) allows one to establish the existence of a unique conformal mating of $\pmb{\cN}_{2}$ and $\overline{z}^2$. 

We will return to this theme in Section~\ref{mating_anti_poly_nielsen_sec}, where a general combination theorem for anti-polynomials and Nielsen maps of reflection groups will be expounded.

\subsubsection{Antiholomorphic correspondence as lift of deltoid reflection}\label{deltoid_corr_subsubsec}

The rational map $f$ gives rise to a $2$:$2$ correspondence $\mathfrak{C}\subset\widehat{\C}\times\widehat{\C}$ whose dynamics is akin to the dynamics of a family of algebraic correspondences introduced by Bullett and Penrose in the 1990s \cite{BP}. We define the correspondence $\mathfrak{C}$ as
\begin{equation}
(z,w)\in\mathfrak{C}\quad \iff\quad \frac{f(w)-f(\eta(z))}{w-\eta(z)}\ =\ 0\quad \iff\quad \overline{z}^2 w+\overline{z}\ =\ 2w^2.
\label{corr_del_eqn}
\end{equation} 
We remark that this correspondence is called antiholomorphic because the local branches $z\mapsto w$ are antiholomorphic. Note also that the correspondence $\mathfrak{C}$ is \emph{reversible}; i.e., its forward branches $z\mapsto w$ are conjugate to its backward branches $w\mapsto z$ via the anti-conformal involution~$\eta$. The dynamics of the correspondence is generated by all possible compositions of the forward and backward branches of $\mathfrak{C}$.
The dynamical plane of $\mathfrak{C}$ splits into two invariant sets: the \emph{lifted non-escaping set} $\widetilde{K(\sigma)}:=f^{-1}(K(\sigma))$ and the \emph{lifted tiling set} $\widetilde{T^\infty(\sigma)}:=f^{-1}(T^\infty(\sigma))$ (see Figure~\ref{deltoid_corr_fig}).

The fact that the simply connected tiling set $T^\infty(\sigma)$ does not contain any critical value of $f$ implies that $f$ has an order three deck transformation $\tau$ on $\widetilde{T^\infty(\sigma)}$. By Equation~\ref{corr_del_eqn}, the branches of $\mathfrak{C}$ on the lifted tiling set $\widetilde{T^\infty(\sigma)}$ are given by $\tau\circ\eta$ and $\tau^2\circ\eta$ (i.e., compositions of $\eta$ with local deck transformations of $f$). Thanks to the relations 
$$
\tau=(\tau^{2}\circ\eta)\circ(\tau\circ\eta)^{-1}\quad \mathrm{and}\quad \eta=\tau^{-1}\circ (\tau\circ\eta),
$$ 
the grand orbits of $\mathfrak{C}$ on $\widetilde{T^\infty(\sigma)}$ are generated by $\eta$ and $\tau$.\footnote{This shows that although we divide out the map $\eta$ from the correspondence $f(w)=f(\eta(z))$, it is recovered in the grand orbit of $\mathfrak{C}$.} Moreover, $\eta$ and $\tau$ generate a subgroup of conformal and anti-conformal automorphisms of $\widetilde{T^\infty(\sigma)}$ isomorphic to the abstract modular group $\Z/2\Z\ast\Z/3\Z$ (see \cite[Proposition~B.7, Theorem~B.8]{LLMM3}).

On the other hand, the forward branch $(f\vert_{\overline{\D^*}})^{-1}\circ f\circ\eta:\widetilde{K(\sigma)}\cap\overline{\D^*}\to\widetilde{K(\sigma)}\cap\overline{\D^*}$ is conformally conjugate to $\sigma:K(\sigma)\to K(\sigma)$ via $f\vert_{\overline{\D^*}}$. Hence, this branch of $\mathfrak{C}$ is conformally conjugate to the anti-polynomial $\overline{z}^2\vert_{\overline{\D^*}}$ (see \cite[Proposition~B.6]{LLMM3}).

In light of the above discussion, the antiholomorphic correspondence $\mathfrak{C}$ can be interpreted as a \emph{mating} of $\overline{z}^2$ and $\Z/2\Z\ast\Z/3\Z$.

\subsection{Mating Basilica with ideal triangle group}\label{c_and_c_center_subsec}

We now discuss the dynamics of a specific quadratic Schwarz reflection map generated by two disjoint quadrature domains.

\subsubsection{Schwarz reflection in a cardioid and a circle}\label{c_and_c_basilica_subsubsec}
Consider the cardioid $\heartsuit:=f(\D)$, where $f(w)=w/2-w^2/4$. Since $f$ is univalent on $\D$, the cardioid $\heartsuit$ is a quadrature domain. We note that $\heartsuit$ is the principal hyperbolic component of the Mandelbrot set, and $f$ is the usual \emph{multiplier map} of this hyperbolic component (i.e., the uniformization of $\heartsuit$ by the multiplier of the unique finite attracting fixed point of the corresponding quadratic polynomial).

Here and in the rest of the paper, we will use the notation $B(a,r), \overline{B}(a,r)$ to denote the open, closed (respectively) disks centered at $a\in\C$ with radius $r>0$.

The circle $\{\vert z\vert=3/4\}$ is a circumcircle to $\partial\Omega$. We define the quadrature multi-domain
$$
\Omega := \heartsuit\cup\overline{B}(0,3/4)^c,
$$ 
and the piecewise Schwarz reflection map 
$$
F:\overline{\Omega}\to\widehat{\C},\qquad 
z \mapsto \left\{\begin{array}{ll}
                    \sigma(z) & \mbox{if}\ z\in\overline{\heartsuit}, \\
                    \sigma_0(z) & \mbox{if}\ z\in B(0,3/4)^c, 
                                          \end{array}\right. 
$$
where $\sigma$ is the Schwarz reflection of $\heartsuit$, and the map $\sigma_0$ is reflection with respect to the circle $\{\vert z\vert=3/4\}$. 
\begin{figure}[h!]
\captionsetup{width=0.96\linewidth}
\begin{tikzpicture}
\node[anchor=south west,inner sep=0] at (0.66,0) {\includegraphics[width=0.32\textwidth]{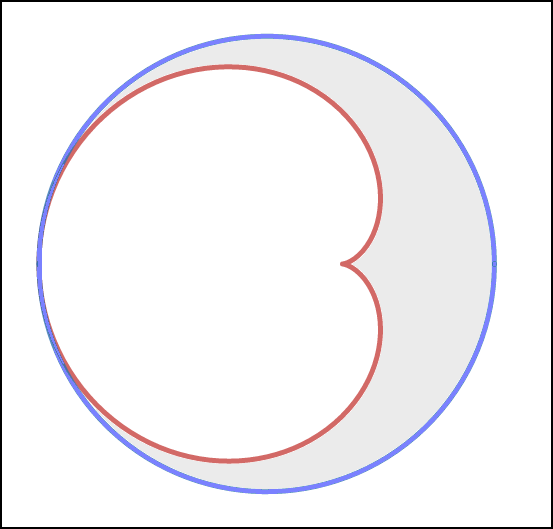}};
\node[anchor=south west,inner sep=0] at (5.6,0) {\includegraphics[width=0.5\textwidth]{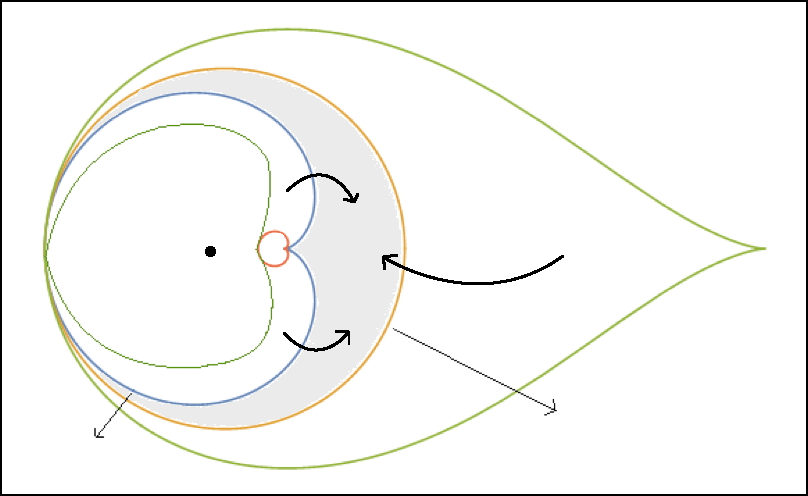}};
\node[anchor=south west,inner sep=0] at (2,-5) {\includegraphics[width=0.66\textwidth]{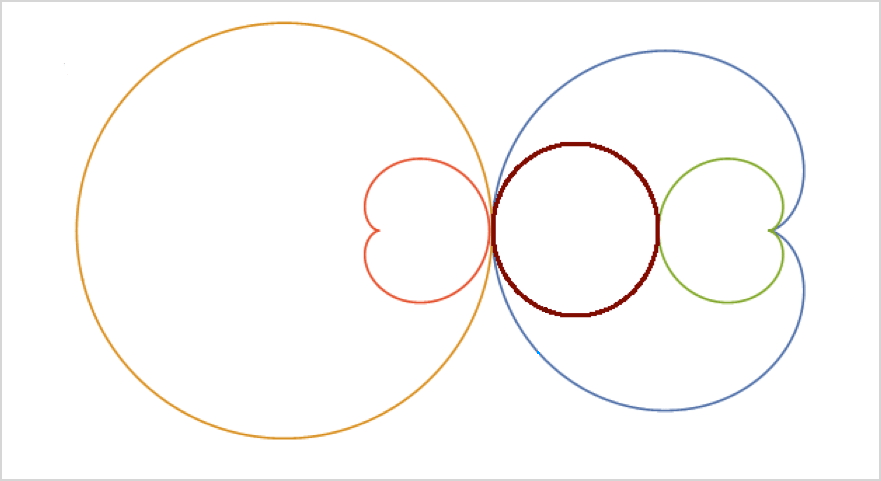}};
\node at (4.04,0.28) {\begin{scriptsize}$\overline{B}(0,3/4)^c$\end{scriptsize}};
\node at (2,2) {$\heartsuit$};
\node at (3.8,1.9) {\begin{scriptsize}$T^0(F)$\end{scriptsize}};
\node at (6.2,0.28) {\begin{scriptsize}$\partial\heartsuit$\end{scriptsize}};
\node at (10.2,0.36) {\begin{scriptsize}$\partial B(0,3/4)$\end{scriptsize}};
\node at (6.2,0.28) {\begin{scriptsize}$\partial\heartsuit$\end{scriptsize}};
\node at (6.2,0.28) {\begin{scriptsize}$\partial\heartsuit$\end{scriptsize}};
\node at (7.25,1.75) {\begin{scriptsize}$0$\end{scriptsize}};
\node at (2.75,-1.4) {$V_2$};
\node at (9.75,-1.4) {$V_1$};
\node at (8.8,-2.66) {$U_1$};
\node at (7.4,-2.66) {$U_2$};
\node at (6.06,-2.66) {$U_1'$};
\end{tikzpicture}
\caption{Top left: The cardioid $\heartsuit$ and the exterior disk $\overline{B}(0,3/4)^c$ are shown. The shaded region is the fundamental tile $T^0(F)$, which is obtained by removing the singular points $1/4,-3/4$ from the complement of $\Omega= \heartsuit\cup\overline{B}(0,3/4)^c$. Top right: The three preimages of $T^0(F)$ under $F$ are shown. Bottom: A topological model for the map $F$ is shown. In these coordinates,  $V_1$ is the cardioid $\heartsuit$, $V_2$ is the disk $\overline{B}(0,3/4)^c$, and the domains $U_1$, $U_1'$ and $U_2$ are $\sigma^{-1}(\heartsuit)$, $\sigma_0^{-1}(\heartsuit)$ and $\sigma^{-1}(\overline{B}(0,3/4)^c)$, respectively. The map $F:U_1\cup U_1'\cup U_2\longrightarrow V_1\cup V_2$ is a degenerate anti-quadratic-like map with a simple critical point at $0$.}
\label{basilica_schwarz_fig}
\end{figure}
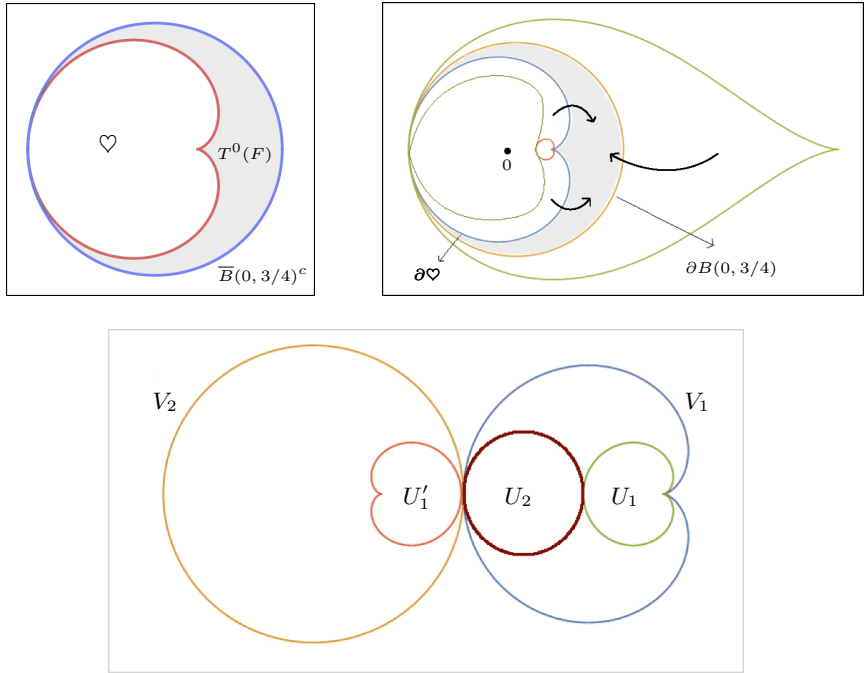

Recall from Figure~\ref{cardioid_disk_fig} that the only critical point of $\sigma$ is at the origin. Hence, the map $F$ has a unique critical point at $0$. Moreover, by construction, $\{0,\infty\}$ is a superattracting $2$-cycle for $F$.

\subsubsection{Triangle group structure}\label{c_and_c_reflection_group_subsubsection}

As the fundamental tile 
$$
T^0(F)=\widehat{\C}\setminus\left(\Omega\cup\{1/4,-3/4\}\right)
$$ 
is a triangle (see Figure~\ref{basilica_schwarz_fig}), and the unique critical point of $F$ does not escape to the tiling set $T^\infty(F)$, it is not hard to see that the action of $F$ on its tiling set is conformally conjugate to $\pmb{\cN}_2$. More precisely, there exists a conformal isomorphism $\psi:\D\to T^\infty(F)$ that conjugates $\pmb{\cN}_2:\D\setminus\Int{\Pi(\pmb{G}_2)}\to\D$ to $F:T^\infty(F)\setminus\Int{T^0(F)}\to T^\infty(F)$ (see \cite[Proposition~5.38]{LLMM1}).

\subsubsection{A degenerate anti-quadratic-like structure}\label{c_and_c_basilica_pinched_anti_quad_subsubsec}

By Proposition~\ref{simp_conn_quad_prop}, the map $F:F^{-1}(\Omega)\to\Omega$ is a $2:1$ branched cover branched only at the origin. As depicted in Figure~\ref{basilica_schwarz_fig}, $F^{-1}(\Omega)\subset \Omega$, and $\partial\Omega$ intersects $\partial F^{-1}(\Omega)$ at two points. Hence, $F:F^{-1}(\Omega)\to\Omega$ can be thought of as a degenerate anti-quadratic-like map.
 
\begin{figure}[h!]
\captionsetup{width=0.96\linewidth}
\begin{center}
\includegraphics[width=0.36\linewidth]{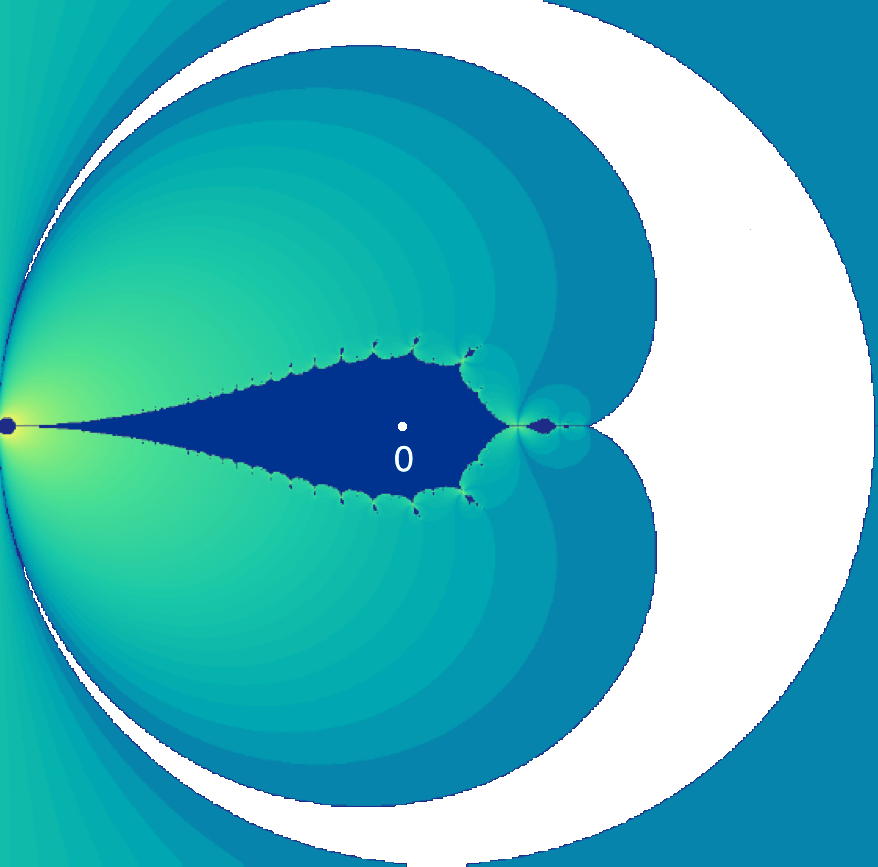}\quad \includegraphics[width=0.5\linewidth]{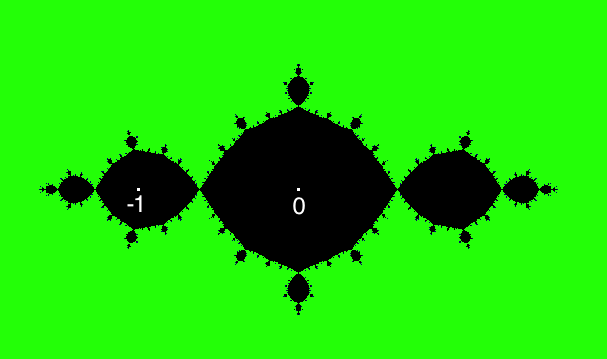}
\end{center}
\caption{Left: Part of the dynamical plane of $F$. Right: The filled Julia set of the Basilica anti-polynomial $\overline{z}^2-1$.}
\label{basilica_anti_poly_fig}
\end{figure}

Since $\{0,\infty\}$ is a superattracting $2$-cycle for $F$, it is natural to expect that the non-escaping set dynamics of $F$ is topologically conjugate to the filled Julia set dynamics of the \emph{Basilica} anti-polynomial $p(z):=\overline{z}^2-1$ (which is the unique quadratic anti-polynomial with a superattracting $2$-cycle) such that the conjugacy is conformal on the interior (see Figure~\ref{basilica_anti_poly_fig}). However, since $\partial\Omega\cap\partial F^{-1}(\Omega)$ consists of two neutral fixed points and the attracting directions at these neutral fixed points lie in the tiling set $T^\infty(F)$, there cannot be a hybrid conjugacy between the non-escaping dynamics of $F$ and the filled Julia set dynamics of $p$ (this is analogous to the deltoid case, see Remark~\ref{deltoid_no_hybrid_rem}).

One is thus forced to take a more combinatorial/topological route to prove the existence of such a conjugacy. The strategy of studying the non-escaping set of $F$ from outside (i.e., from the tiling set) turns out to be fruitful.
Since $F$ is postcritically finite, one can adapt standard arguments from polynomial dynamics to show that $\Lambda(F):=\partial T^\infty(F)$ is locally connected \cite[Proposition~6.4]{LLMM1}. Hence, the conformal conjugacy $\psi$ between $\pmb{\cN}_2$ and $F$ (see Subsection~\ref{c_and_c_reflection_group_subsubsection}) extends continuously to yield a semi-conjugacy $\psi:\mathbb{S}^1\to\Lambda(F)$ between $\pmb{\cN}_2$ and $F$. Thus, $\Lambda(F)$ can be topologically modeled as the quotient of $\mathbb{S}^1$ by an $\pmb{\cN}_2$-invariant equivalence relation, which we call the \emph{lamination} associated with $F$ \cite[\S 5.3.5]{LLMM1}. 
\begin{definition}\label{push_lami_def}
The \emph{push-forward} $\phi_{\ast}(\lambda)$ of an equivalence relation $\lambda$ under a circle homeomorphism $\phi$ is defined as the image of $\lambda\subset\mathbb{S}^1\times\mathbb{S}^1$ under $\phi\times\phi$. Clearly, $\phi_{\ast}(\lambda)$ is an equivalence relation on $\mathbb{S}^1$. 
\end{definition}
Applying combinatorial tools from the study of polynomial Julia sets, one can give a complete description of the lamination of $F$ and conclude that the push-forward of this lamination under the Minkowski circle homeomorphism $\pmb{\mathcal{E}}_2$ is precisely the combinatorial lamination of $p$ (cf. Definition~\ref{rat_lami_def} and Remark~\ref{lami_rmk}). This provides a topological conjugacy between $F\vert_{\Lambda(F)}$ and $p\vert_{\mathcal{J}(p)}$. Finally, a detailed study of the Fatou components of $F$ (i.e., the components of the interior of $K(F)$) allows one to conformally extend this topological conjugacy between limit and Julia sets to a conjugacy between $F\vert_{K(F)}$ and $p\vert_{\mathcal{K}(p)}$ (see \cite[Proposition 5.30, Corollary 5.33]{LLMM1} and \cite[Proposition 11.1]{LLMM2}).

In light of the above statements, it can be justified that $F$ is the conformal mating of $\overline{z}^2-1$ and $\pmb{\cN}_2$ in a precise sense (see \cite[\S 7]{LLMM1}).

\subsubsection{Uniqueness of the conformal mating of $\overline{z}^2-1$ and $\pmb{\cN}_2$}\label{c_and_c_basilica_unique_mating_subsubsec}

As the invariant external dynamical rays of $\overline{z}^2-1$ at the angles $1/3$ and $2/3$ land at the same point, it follows (from the definition of topological mating of $\overline{z}^2-1$ and $\pmb{\cN}_2$) that two ideal vertices of $\D\cap\Pi(\pmb{G}_2)$ are identified in the topological mating. Hence, the domain of definition of the topological mating is the closure of the union of two Jordan domains touching at a single point. Furthermore, as the Nielsen map $\pmb{\cN}_2$ fixes the boundary of its domain of definition pointwise, one concludes that any conformal mating of $\overline{z}^2-1$ and $\pmb{\cN}_2$ is a piecewise Schwarz reflection map associated with two disjoint simply connected quadrature domains touching at a single point.
Using mapping degrees of $\overline{z}^2-1$ and $\pmb{\cN}_2$, one can now argue in light of Proposition~\ref{simp_conn_quad_prop} that one of these quadrature domains is uniformized by a degree one rational map, while the other is uniformized by a degree two rational map. Finally, the critical orbit relation of $\overline{z}^2-1$ can be exploited to deduce that this piecewise Schwarz reflection map is given by $F$, up to M{\"o}bius conjugacy (see \cite[Appendix~A]{LLMM2}). Thus, $F$ is the unique conformal mating between the maps $\overline{z}^2-1$ and $\pmb{\cN}_2$. 

\begin{theorem}\cite[\S 7]{LLMM1}, \cite[Appendix~A]{LLMM2}\label{c_and_c_basilica_dynamics_thm}
The map $F\vert_{K(F)}$ is topologically conjugate to $p\vert_{\mathcal{K}(p)}$ such that the conjugacy is conformal on the interior, where $p(z)=\overline{z}^2-1$. On the other hand, the map $F:\overline{T^\infty(F)}\setminus\Int{T^0(F)}\to \overline{T^\infty(F)}$ is topologically semi-conjugate to $\pmb{\cN}_2:\overline{\D}\setminus\Int{\Pi(\pmb{G}_2)}\to\overline{\D}$ such that the semi-conjugacy restricts to a conformal conjugacy on $T^\infty(F)$. Moreover, up to M{\"o}bius conjugacy, $F$ is the unique conformal mating of $\pmb{\cN}_2$ and $p$.
\end{theorem}

To conclude this subsection, let us mention that as in the deltoid case, the existence and uniqueness of the above conformal mating can also be proved using David surgery (see Section~\ref{mating_anti_poly_nielsen_sec}).

\subsection{An antiholomorphic correspondence from cubic Chebyshev polynomial}\label{chebyshev_subsec}

Consider the cubic Chebyshev polynomial $f(w)=w^3-3w$. By \cite[Proposition~3.2]{LLMM3}, $f$ is injective on the closed disk $\overline{B}(3,2)=\{w\in\C:\vert w-3\vert\leq 2\}$. The image $\Omega:= f(B(3,2))$ is a simply connected quadrature domain. As $f$ has a critical point at $1$, the boundary $\partial\Omega$ has a cusp at $f(1)=-2$ and is non-singular otherwise. Thus, $T^0(\sigma)=\widehat{\C}\setminus\cup\left(\Omega\cup\{-2\}\right)$.

We denote the reflection in the circle $\partial B(3,2)$ by $\widehat{\eta}$. Specifically, $\widehat{\eta}(w)=3+\frac{4}{\overline{w}-3}$.
As the only critical points of $f$ outside $\overline{B}(3,2)$ are at $-1$ and $\infty$, the associated Schwarz reflection map $\sigma=f\circ\widehat{\eta}\circ(f\vert_{B(3,2)})^{-1}$ has two critical points: a simple critical point at $f(\widehat{\eta}(-1))=f(2)=2$ and a double critical point at $f(\widehat{\eta}(\infty))=f(3)$.
Moreover, $\sigma$ fixes the critical point at $2$ and sends the double critical point at $f(3)$ to $\infty$.

\subsubsection{The external map of $\sigma$}\label{nielsen_first_return_external_map_subsubsec}

Since $\sigma$ sends the double critical point at $f(3)$ to $\infty\in T^0(\sigma)$, it follows that the rank one tile of $T^\infty(\sigma)$ contains a critical point. However, since the other critical point $2$ of $\sigma$ is fixed, it does not lie in the tiling set, and hence $\sigma$ acts like reflection in $\partial\Omega\setminus\{-2\}$ on tiles of higher ranks. It turns out that a conformal model for the tiling set dynamics of $\sigma$ (also called the external map of $\sigma$) arises from a discrete subgroup of $\mathrm{Aut}^\pm(\D)$ generated by a circular reflection and a torsion element.
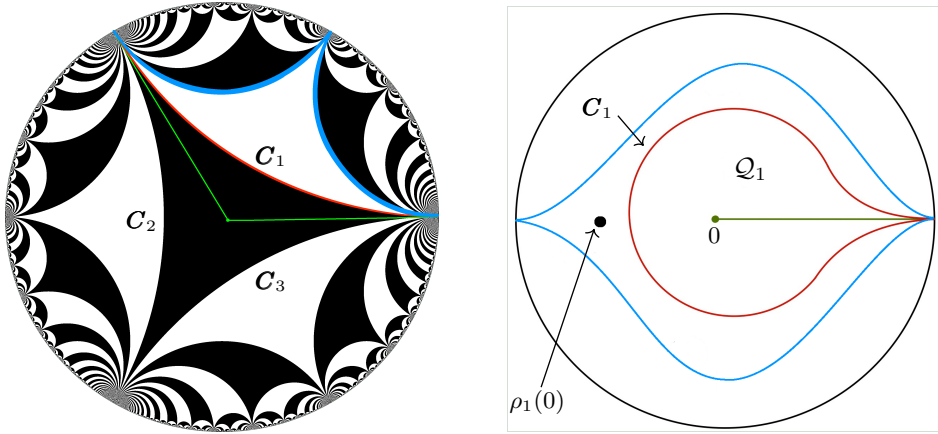
\begin{figure}[h!]
\captionsetup{width=0.96\linewidth}
\begin{tikzpicture}
\node[anchor=south west,inner sep=0] at (0,0) {\includegraphics[width=0.45\linewidth]{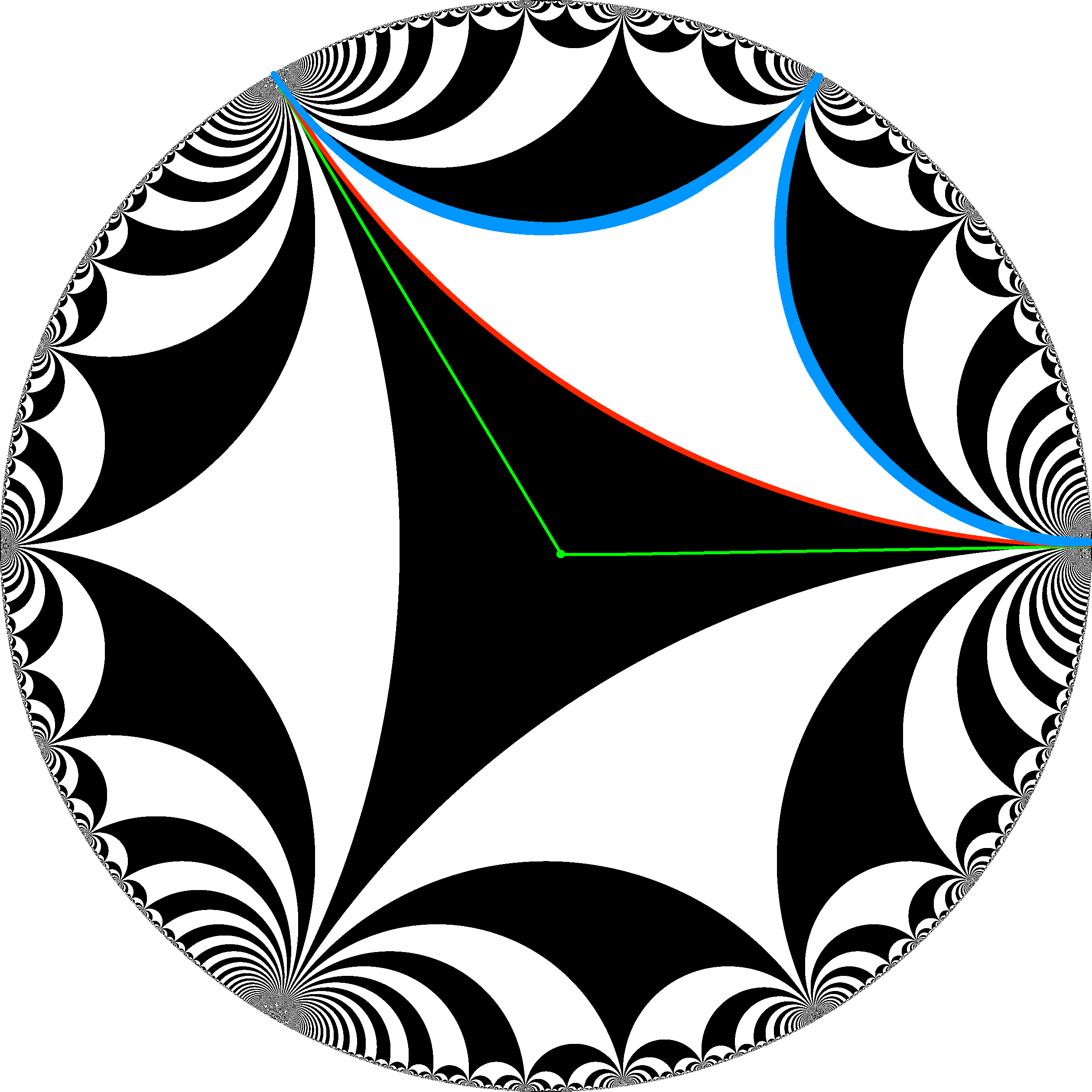}};
\node[anchor=south west,inner sep=0] at (6.6,0) {\includegraphics[width=0.45\linewidth]{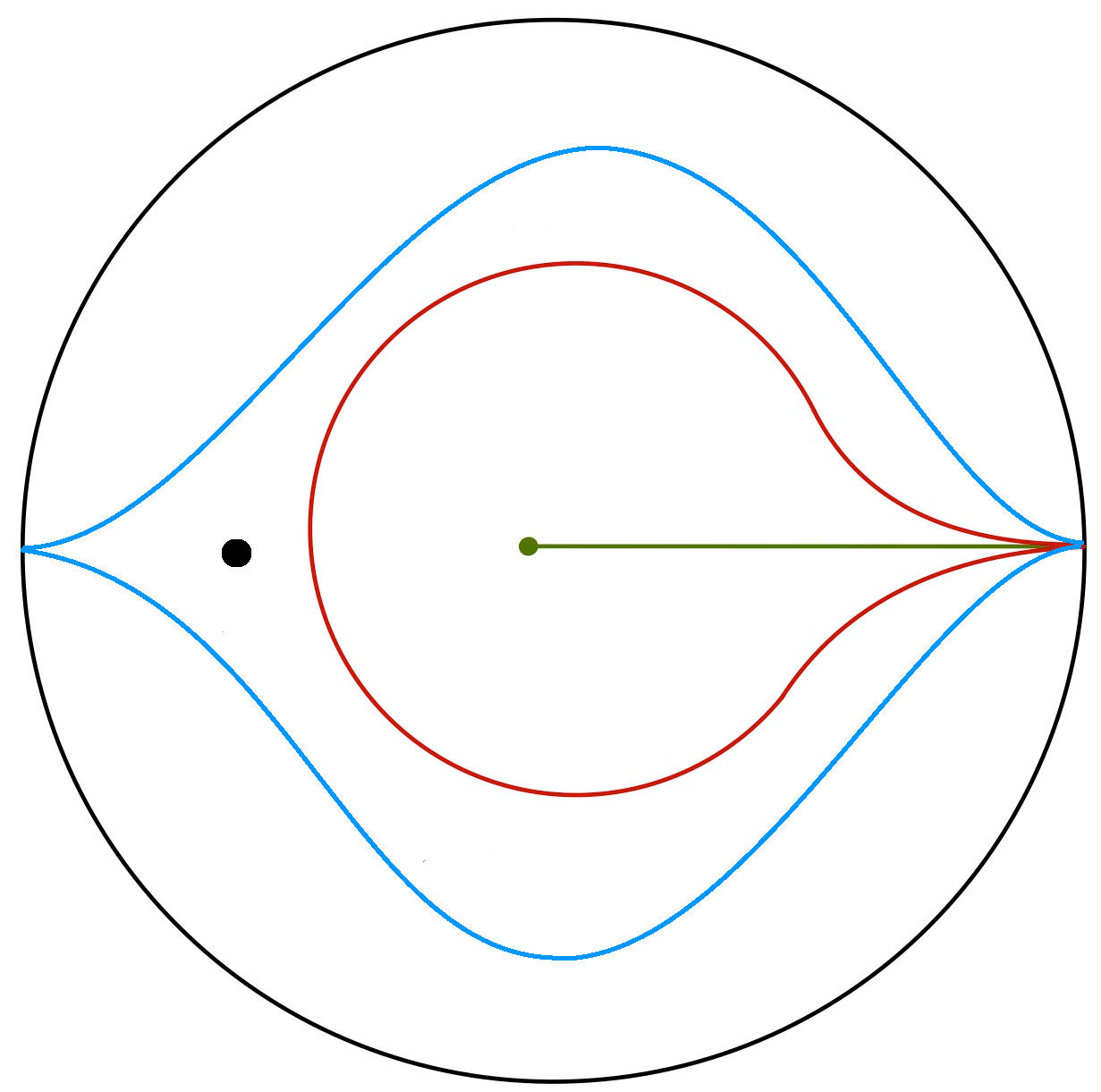}};
\node at (3.5,3.66) {$\pmb{C}_1$};
\node at (1.8,2.8) {$\pmb{C}_2$};
\node at (3.5,2) {$\pmb{C}_3$};
\node at (7.8,4.3) {$\pmb{C}_1$};
\draw [->,line width=0.5pt] (8.06,4.2) to (8.4,3.8);
\node at (9.8,3.5) {\begin{large}$\mathcal{Q}_0$\end{large}};
\node at (7,0.4) {$\rho_1(0)$};
\draw [->,line width=0.5pt] (7.05,0.6) to (7.75,2.7);
\node at (9.33,2.62) {$0$};
\end{tikzpicture}
\caption{Left: The group $\mathbbm{G}_2$ is generated by reflection in $\pmb{C}_1$ and rotation by $2\pi/3$. The central black ideal triangle is $\Pi^{\D}(\pmb{G}_2)$, and the part of it bounded by the two green radial lines is the fundamental domain $\Pi(\mathbbm{G}_2)$ for the $\mathbbm{G}_2$-action on $\D$. The surface $\mathcal{Q}$ (respectively, $\mathcal{Q}_0$) is obtained from the region of $\D$ bounded by the green lines (respectively, from $\Pi(\mathbbm{G}_2)$) by identifying the green lines under $M_\omega$. The tessellation of $\D$ under $\pmb{G}_2$ induces a tessellation of $\mathcal{Q}$. Right: The surface $\mathcal{Q}$, after being uniformized to the disk $\D$, is depicted.}
\label{itg_nielsen_first_return_fig}
\end{figure}

Specifically, consider the group $\mathbbm{G}_2(\geqslant\pmb{G}_2)$ generated by $M_\omega(z):=\omega z$ (where $\omega=e^{\frac{2\pi i}{3}}$) and the reflection $\rho$ in the hyperbolic geodesic $\pmb{C}_1$ of $\D$ connecting $1$ and $\omega$ (see Definition~\ref{regular_ideal_polygon_ref_group_def}). A fundamental domain $\Pi(\mathbbm{G}_2)$ for the $\mathbbm{G}_2$-action on $\D$ is given by one-third of the fundamental domain $\Pi^{\D}(\pmb{G}_2):=\Pi(\pmb{G}_2)\cap\D$ for the $\pmb{G}_2$-action on $\D$ (see Figure~\ref{itg_nielsen_first_return_fig}). Note that the Riemann surface $\mathcal{Q}:=\faktor{\D}{\langle M_\omega\rangle}$ is biholomorphic to $\D$, and $\mathcal{Q}_0:=\faktor{\Pi^{\D}(\pmb{G}_2)}{\langle M_\omega\rangle}$ is a simply connected region embedded in $\mathcal{Q}$. We will use the region in $\D$ (respectively, in $\Pi^{\D}\left(\pmb{G}_2\right)$) bounded by the radial lines at angles $0$ and $2\pi/3$ as coordinates on $\mathcal{Q}$ (respectively, on $\mathcal{Q}_0$). 
We define the map
$$
\pmb{\cF}_2:\mathcal{Q}\setminus\Int{\mathcal{Q}_0}\longrightarrow\mathcal{Q}
$$
as the map $\rho$ post-composed with the quotient map from $\D$ to $\mathcal{Q}$. Note that the surface $\mathcal{Q}$ has a natural tessellation structure with $\mathcal{Q}_0$ as the rank $0$ tile and components of $\pmb{\cF}_2^{-n}(\mathcal{Q}_0)$ as tiles of rank $n$. It is worth pointing out that there is a unique rank one tile for this tessellation given by $\rho\left(\Pi^{\D}(\pmb{G}_2)\right)$ (see Figure~\ref{itg_nielsen_first_return_fig}).
The salient features of the map $\pmb{\cF}_2$ are that 
$$
\pmb{\cF}_2:\left(\pmb{\cF}_2^{-1}(\mathcal{Q}_0),\rho(0)\right) \longrightarrow \left(\mathcal{Q}_0,0\right)
$$ 
is a $3:1$ branched cover between two pointed disks with a double critical point at $\rho(0)$, while $\pmb{\cF}_2$ is injective on tiles of higher ranks.

According to \cite[Proposition~4.15]{LLMM3}, there exists a conformal isomorphism
$$
\psi:\left(\mathcal{Q},0\right)\longrightarrow\left(T^\infty(\sigma),\infty\right)
$$
that conjugates $\pmb{\cF}_2:\mathcal{Q}\setminus\Int{\mathcal{Q}_0}\to\mathcal{Q}$ to $\sigma:T^\infty(\sigma)\setminus\Int{T^0(\sigma)}\to T^\infty(\sigma)$ . The map $\psi$ can be constructed by lifting the conformal isomorphism $(\mathcal{Q}_0,0)\to (T^0(\sigma),\infty)$ (whose homeomorphic boundary extension sends $1\in\partial\mathcal{Q}_0$ to the cusp $-2\in\partial\Omega$) by iterates of $\pmb{\cF}_2$ and $\sigma$. Such a lifting procedure can be performed since 
\smallskip

\noindent i) the maps $\pmb{\cF}_2:\pmb{\cF}_2^{-1}(\mathcal{Q}_0) \longrightarrow \mathcal{Q}_0$ and $\sigma:\sigma^{-1}(T^0(\sigma))\to T^0(\sigma)$ are $3:1$ branched coverings having a unique (double) critical point with associated critical value at $0, \infty$, respectively (equivalently, they are doubly ramified over $0, \infty$, respectively, and are unramified otherwise),
\smallskip

\noindent ii) the maps $\pmb{\cF}_2, \sigma$ have no other critical point in $\mathcal{Q},  T^\infty(\sigma)$, respectively, and
\smallskip

\noindent iii) the maps $\pmb{\cF}_2, \sigma$ act as the identity map on $\partial\mathcal{Q}_0, \partial T^0(\sigma)$.

In other words, the map $\pmb{\cF}_2$ is the external map for $\sigma$. We refer the reader to \cite[\S 4.4]{LLMM3} for details of this construction and more properties of the map $\pmb{\cF}_2$\footnote{The map $\pmb{\cF}_2$ was termed $\rho$ in \cite{LLMM3}.}.

\subsubsection{Hybrid conjugacy between $\sigma$ and a quadratic parabolic anti-rational map}\label{chebyshev_center_hybrid_conj_subsubsec}

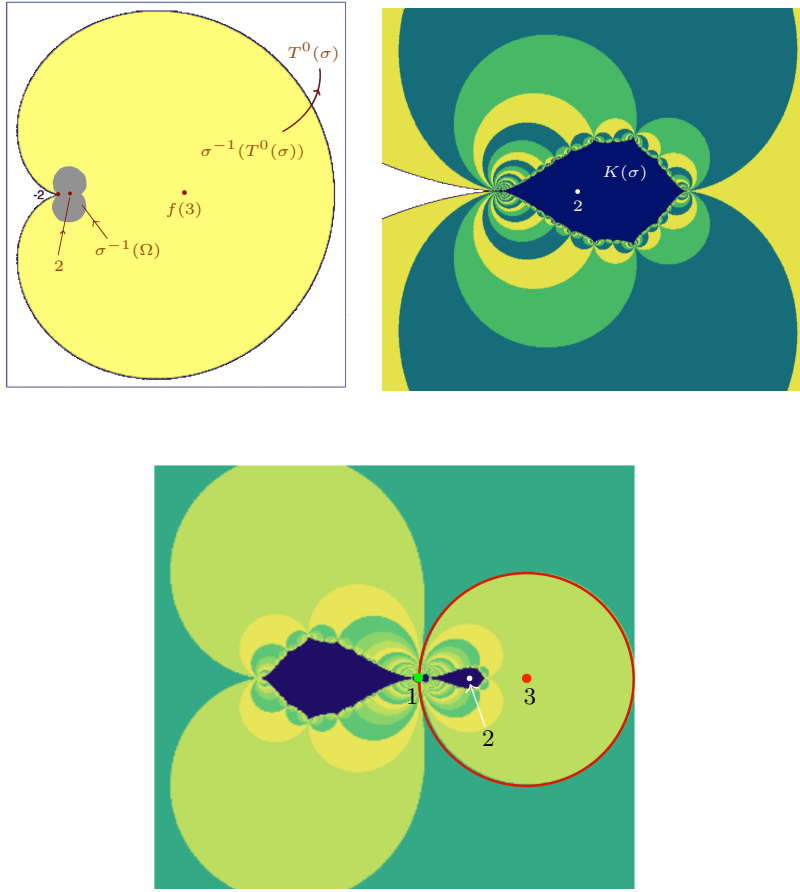
\begin{figure}[h!]
\captionsetup{width=0.96\linewidth}
\begin{tikzpicture}
\node[anchor=south west,inner sep=0] at (0,0) {\includegraphics[width=0.36\linewidth]{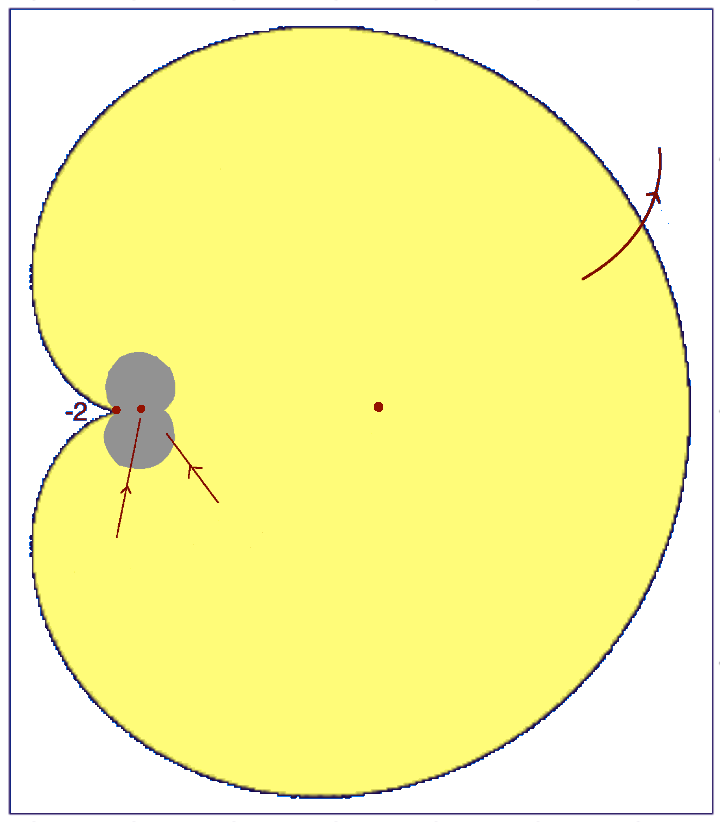}};
\node[anchor=south west,inner sep=0] at (5,0) {\includegraphics[width=0.44\linewidth]{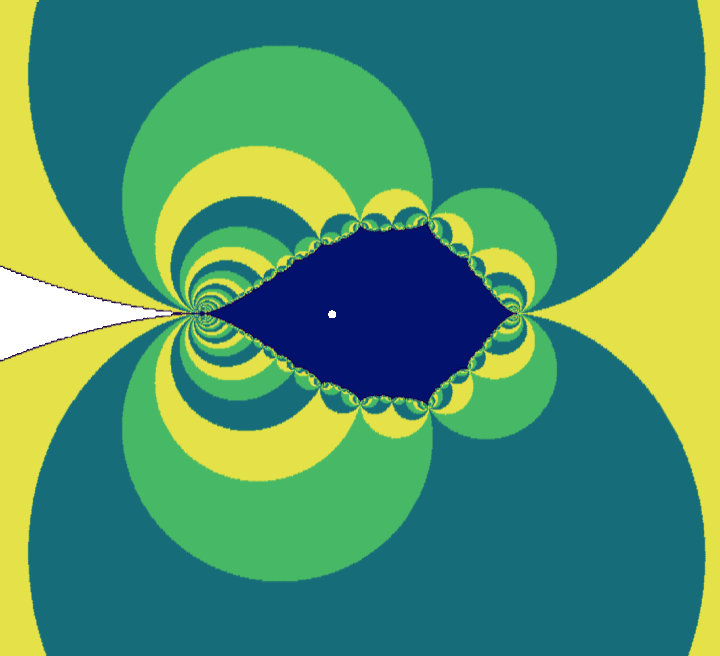}};
\node[anchor=south west,inner sep=0] at (2,-6.6) {\includegraphics[width=0.5\linewidth]{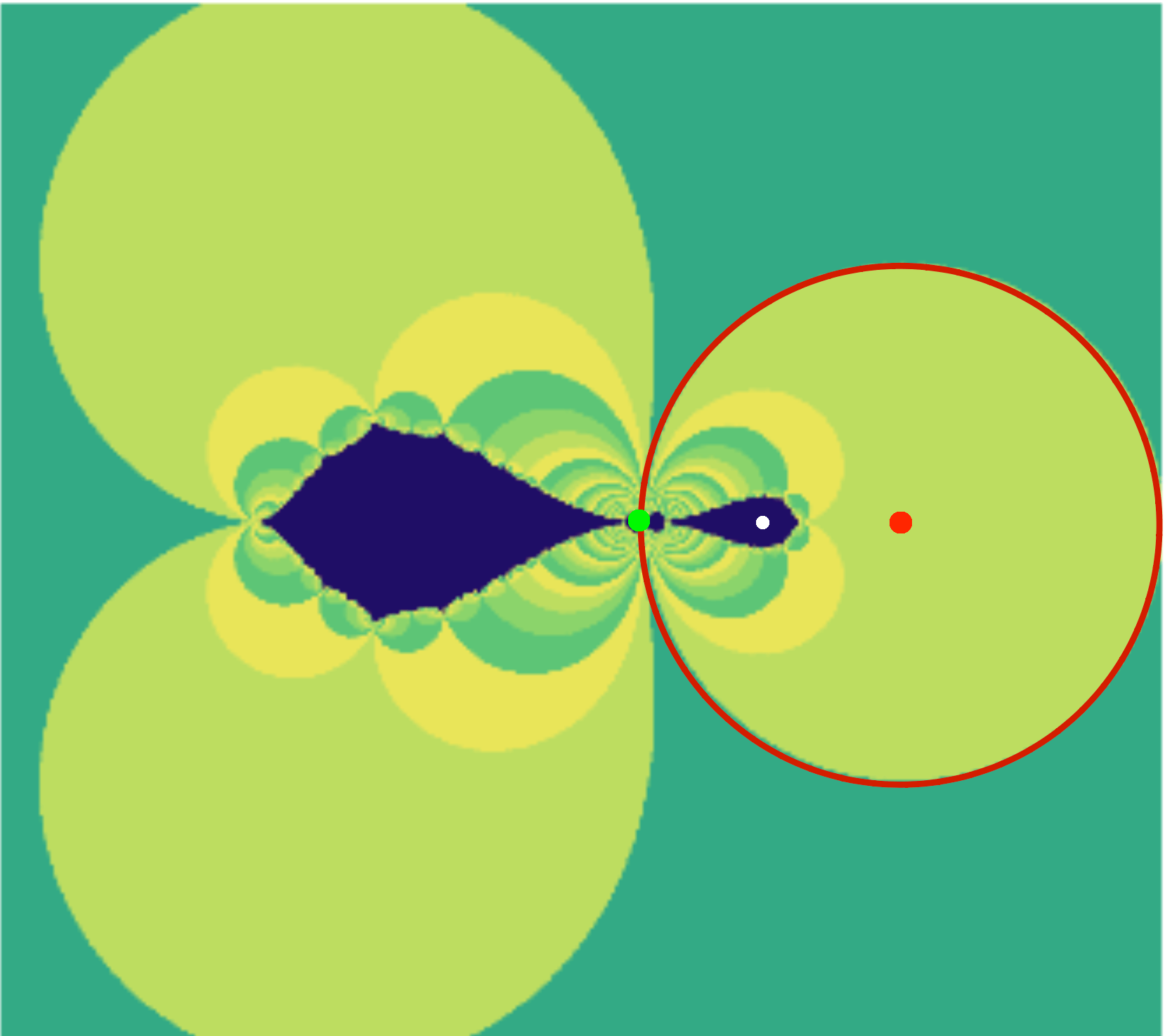}};
\node at (3.3,3.2) {\begin{scriptsize}\textcolor{RawSienna}{$\sigma^{-1}(T^0(\sigma))$}\end{scriptsize}};
\node at (1.66,1.88) {\begin{scriptsize}\textcolor{RawSienna}{$\sigma^{-1}(\Omega)$}\end{scriptsize}};
\node at (0.75,1.66) {\begin{scriptsize}\textcolor{RawSienna}{$2$}\end{scriptsize}};
\node at (2.4,2.4) {\begin{scriptsize}\textcolor{RawSienna}{$f(3)$}\end{scriptsize}};
\node at (4.1,4.5) {\begin{scriptsize}\textcolor{RawSienna}{$T^0(\sigma)$}\end{scriptsize}};
\node at (7.58,2.46) {\begin{scriptsize}\textcolor{White}{$2$}\end{scriptsize}};
\node at (8.2,2.9) {\begin{scriptsize}\textcolor{White}{$K(\sigma)$}\end{scriptsize}};
\node at (6.95,-4.05) {$3$};
\node at (6.4,-4.6) {$2$};
\draw [->,line width=0.5pt,white] (6.36,-4.45) to (6.16,-3.88);
\node at (5.4,-4.05) {$1$};
\end{tikzpicture}
\caption{Top left: The large cardioid is the quadrature domain $\Omega$, the grey `double cardioid' region is $\sigma^{-1}(\Omega)$, and the yellow region is the rank one tile $\sigma^{-1}(T^0(\sigma))$ which contains the double critical point $f(3)$. Top right: The non-escaping set of $\sigma$ is a closed topological disk (in dark blue) with the cardioid cusp $-2$ on its boundary. Bottom: Pulling back the tiling and non-escaping sets of $\sigma$ under $f$ gives the lifted tiling and lifted non-escaping sets of the correspondence $\mathfrak{C}$.}
\label{chebyshev_center_fig}
\end{figure}

By Proposition~\ref{simp_conn_quad_prop} and the discussion of the critical points of $\sigma$, the map $\sigma:\sigma^{-1}(\Omega)\to\Omega$ is a $2:1$ branched covering with a superattracting fixed point at $2$. Moreover, $\sigma^{-1}(\Omega)\subset\Omega$ and $\partial\sigma^{-1}(\Omega)\cap\partial\Omega=\{-2\}$ (note that $\partial\Omega$ has a unique singular point, see Figure~\ref{chebyshev_center_fig}). Thus, $\sigma:\sigma^{-1}(\Omega)\to\Omega$ exhibits a degenerate anti-quadratic-like structure with a superattracting fixed point. In some sense, the situation is analogous to the deltoid reflection map, and it is not hard to employ similar techniques to prove that the non-escaping dynamics of $\sigma$ is topologically conjugate to $\overline{z}^2\vert_{\overline{\D}}$ with the conjugacy being conformal on the interior (the non-escaping set of $\sigma$ is shown in Figure~\ref{chebyshev_center_fig}). Moreover, for reasons similar to the ones mentioned in Remark~\ref{deltoid_no_hybrid_rem}, there is no hybrid conjugacy between this degenerate anti-quadratic-like map and~$\overline{z}^2$.
\begin{figure}[ht!]
\captionsetup{width=0.96\linewidth}
\begin{tikzpicture}
\node[anchor=south west,inner sep=0] at (0,0) {\includegraphics[width=0.47\textwidth]{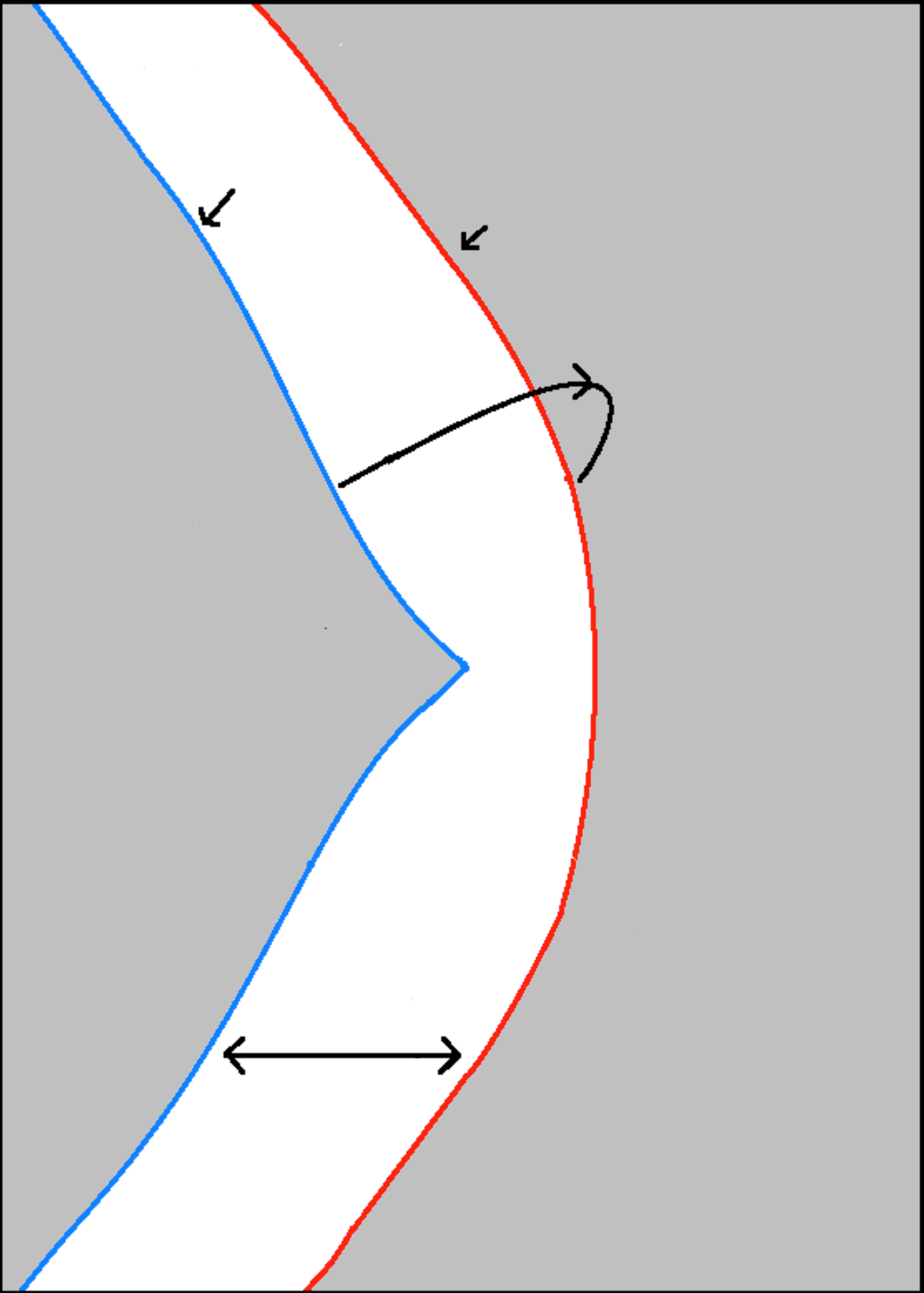}};
\node[anchor=south west,inner sep=0] at (6.3,0) {\includegraphics[width=0.5\textwidth]{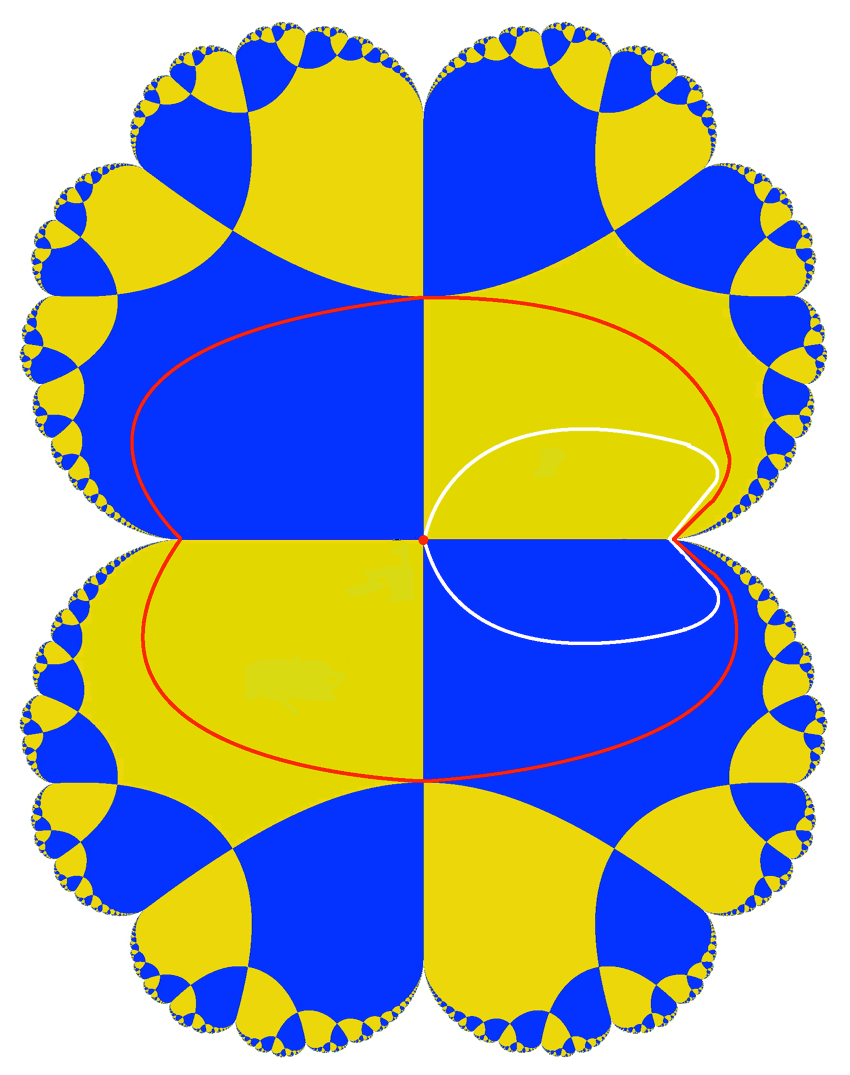}};
\node at (10.76,4.5) {$P$};
\node at (9.09,3.66) {$-\frac12$};
\node at (8.4,2.84) {$q^{-1}(P)$};
\node at (1.64,7.36) {$\mathbf{\partial\mathbf{U}}$};
\node at (3.45,7.1) {$\mathbf{\partial\mathbf{V}}$};
\node at (4,6.18) {$\mathbf{F}$};
\node at (1.4,4.6) {$\mathbf{U}$};
\node at (5.08,2.7) {$\widehat{\C}\setminus\overline{\mathbf{V}}$};
\node at (2.45,1.9) {$\frac12$};
\end{tikzpicture}
\caption{Left: Depicted is a pinched anti-quadratic-like map $\mathbf{F}:(\overline{\mathbf{U}},\infty)\to(\overline{\mathbf{V}},\infty)$. Right: The filled Julia set of the parabolic anti-polynomial $q(z)=\overline{z}+\overline{z}^2$ is shown. The dynamics of $q$ on its filled Julia set is conformally conjugate to the parabolic anti-Blaschke product $B_2$. The attracting petal $P$ of $q$ subtends a positive angle at the parabolic fixed point~$0$. The critical point $-\frac12$ (of $q$) lies on the boundary of the petal $P$. The pre-image of $P$ (under $q$) is a simply connected domain, which maps $2:1$ onto $P$ branched only at $-\frac12$. The pinched anti-quadratic-like map $\mathbf{F}$ is straightened to a parabolic anti-rational map by quasiconformally gluing $(q^{-1}(P),P)$ into $(\widehat{\C}\setminus\overline{U},\widehat{\C}\setminus\overline{V})$ in a boundary equivariant way.}
\label{anti_quad_like_fig}
\end{figure}

However, there is a key difference between the current setting and the examples described in the previous two subsections. Here, $\partial\sigma^{-1}(\Omega)\cap\partial\Omega$ is a singleton; while in both the deltoid and Circle-and-Cardioid examples the domain of the degenerate anti-quadratic-like map touches the range at more than one points. Equivalently, the external map $\pmb{\cF}_2$ of $\sigma$ has a \emph{unique} parabolic fixed point, while the external map $\pmb{\cN}_2$ in the previous two examples has three parabolic fixed points.
It is natural to expect that  $\pmb{\cF}_2\vert_{\partial\mathcal{Q}}$ is quasisymmetrically conjugate to the action of a Blaschke product with a parabolic fixed point on $\mathbb{S}^1$. On the other hand, there is no such candidate Blaschke product for the external map $\pmb{\cN}_2$ as a Blaschke product cannot have multiple parabolic points on the circle.

The above observation leads one to the parabolic Tricorn $\pmb{\mathcal{B}}_2$ consisting of quadratic anti-rational maps with a completely invariant parabolic basin (see Subsection~\ref{para_anti_rat_gen_subsubsec}). In order to straighten the non-escaping dynamics of $\sigma$ to parabolic anti-rational maps, a class of maps called \emph{pinched anti-quadratic-like maps} was introduced in \cite[\S 5]{LLMM3} (see Figure~\ref{anti_quad_like_fig}) and it was shown that pinched anti-quadratic-like maps are hybrid conjugate to maps in the parabolic Tricorn (cf. \cite{Lom15}). This is done via a quasiconformal surgery that glues an attracting petal of the parabolic anti-Blaschke product $B_2$ outside a pinched anti-quadratic-like map; or equivalently, replaces the external dynamics of a pinched anti-quadratic-like map with the map $B_2$ (see Figure~\ref{anti_quad_like_fig}). However, the existence of the pinching point makes this straightening theorem subtler than the classical straightening for polynomial-like maps. In particular, since the fundamental domain of a pinched anti-polynomial-like map is a pinched annulus, one needs to perform quasiconformal interpolation in a topological strip. This necessitates one to control the asymptotics of conformal maps between topological strips, which can be achieved by a result of Warschawski (cf. \cite{War42}).

As a consequence of the aforementioned straightening theorem, it can be concluded that the non-escaping dynamics of $\sigma$ is hybrid conjugate to the dynamics of the quadratic parabolic anti-rational map $z\mapsto\overline{z}+1/\overline{z}+1$ on its filled Julia set (which has a superattracting fixed point).

\begin{theorem}\cite[Proposition~4.15, Theorem~5.4]{LLMM3}\label{chebyshev_center_dynamics_thm}
\noindent\begin{enumerate}\upshape
\item The map $\sigma\vert_{K(\sigma)}$ is hybrid conjugate to $R\vert_{\mathcal{K}(R)}$, where $R(z)=\overline{z}+1/\overline{z}+1$. 
\item The maps $\sigma:T^\infty(\sigma)\setminus\Int{T^0(\sigma)}\to T^\infty(\sigma)$ and $\pmb{\cF}_2:\mathcal{Q}\setminus\Int{\mathcal{Q}_0}\to\mathcal{Q}$ are conformally conjugate.
\end{enumerate}
\end{theorem}

\subsubsection{The associated antiholomorphic correspondence}\label{chebyshev_center_corr_subsubsec}

As in the deltoid setting (see Subsection~\ref{deltoid_corr_subsubsec}), the polynomial $f(w)=w^3-3w$ and the reflection map $\widehat{\eta}$ define a $2$:$2$ correspondence $\mathfrak{C}\subset\widehat{\C}\times\widehat{\C}$ by the formula:
\begin{equation}
(z,w)\in\mathfrak{C}\iff \frac{f(w)-f(\widehat{\eta}(z))}{w-\widehat{\eta}(z)}=0,
\label{cheby_center_corr_eqn}
\end{equation} 
and the dynamical plane of $\mathfrak{C}$ admits an invariant partition into the \emph{lifted non-escaping set} $\widetilde{K(\sigma)}:=f^{-1}(K(\sigma))$ and the \emph{lifted tiling set} $\widetilde{T^\infty(\sigma)}:=f^{-1}(T^\infty(\sigma))$ (see Figure~\ref{chebyshev_center_fig} and \cite[\S 10]{LLMM3}). Similar to the deltoid situation, the correspondence $\mathfrak{C}$ defined by the cubic Chebyshev polynomial $f$ is reversible; i.e., $\widehat{\eta}$ conjugates the forward branches of $\mathfrak{C}$ to its backward branches.

The forward branch $(f\vert_{\overline{B}(3,2)})^{-1}\circ f\circ\widehat{\eta}:\widetilde{K(\sigma)}\cap\overline{B}(3,2)\to\widetilde{K(\sigma)}\cap\overline{B}(3,2)$ is conformally conjugate to $\sigma:K(\sigma)\longrightarrow K(\sigma)$ via $f\vert_{\overline{B}(3,2)}$. By Theorem~\ref{chebyshev_center_dynamics_thm}, this branch of $\mathfrak{C}$ is hybrid conjugate to $R\vert_{\mathcal{K}(R)}$, where $R(z)=\overline{z}+1/\overline{z}+1$ (see \cite[Proposition~10.6]{LLMM3}).

Although the Schwarz reflection $\sigma$ has a critical point in its tiling set, the dynamics of the correspondence $\mathfrak{C}$ on its lifted tiling set has a group structure. This is due to the following reason. The point at $\infty$ is a fully ramified critical point for $f$ and hence $f: \widetilde{T^\infty(\sigma)}\setminus\{\infty\}\longrightarrow T^\infty(\sigma)\setminus\{\infty\}$ is a degree $3$ covering map between two topological annuli. Thus, it is a Galois covering with deck transformation group isomorphic to $\Z/3\Z$. We choose a generator $\widehat{\tau}$ for this deck group, and observe that the branches of $\mathfrak{C}$ on $\widetilde{T^\infty(\sigma)}$ are given by $\widehat{\tau}\circ\widehat{\eta}$ and $\widehat{\tau}^{\circ 2}\circ\widehat{\eta}$.
It can be verified using a \emph{ping-pong} argument that the grand orbits of $\mathfrak{C}$ on the lifted tiling set are generated by $\widehat{\eta}$ and $\widehat{\tau}$, and the group generated by these two maps is the free product of $\langle\widehat{\eta}\rangle$ and $\langle\widehat{\tau}\rangle$ (see \cite[Proposition~10.5]{LLMM3}).

In fact, according to \cite[Proposition~2.18, Remark~3.11]{LMM3}, there exists a uniformizing map $\widetilde{\psi}:\D\longrightarrow\widetilde{T^\infty(\sigma)}$ that conjugates the group $\mathbbm{G}_2\vert_{\D}$ to the group $\langle\widehat{\eta}\rangle\ast\langle\widehat{\tau}\rangle$ . This can be seen as follows. One can lift the conformal map $\psi: \left(\mathcal{Q},0\right) \longrightarrow \left(T^\infty(\sigma),\infty\right)$ via the two branched coverings appearing in the vertical arrows of the commutative diagram below to construct a conformal isomorphism $\widetilde{\psi}:\left(\D,0\right)    \longrightarrow \left(\widetilde{T^\infty(\sigma)},\infty\right)$. Since $\psi$ conjugates $\pmb{\cF}_2$ to $\sigma$, its lift $\widetilde{\psi}$ can be chosen so that it conjugates the circular reflection $\rho$ in $\pmb{C}_1$ to $\widehat{\eta}$. On the other hand, since $M_\omega$ is a generator of the deck transformation group for the projection map $\D\rightarrow\mathcal{Q}$, the lifted map $\widetilde{\psi}$ conjugates $M_\omega$ to the deck transformation $\widehat{\tau}$ (after possibly replacing $\widehat{\tau}$ with some iterate). Hence, $\widetilde{\psi}$ conjugates $\mathbbm{G}_2=\langle\rho\rangle\ast\langle M_\omega\rangle$ to $\langle\widehat{\eta}\rangle\ast\langle\widehat{\tau}\rangle$.
\[
  \begin{tikzcd}
  \left(\D,0\right)    \arrow{d}{\mathrm{proj}} \arrow{r}{\widetilde{\psi}} &  \left(\widetilde{T^\infty(\sigma)},\infty\right) \arrow{d}{f} \\
   \left(\mathcal{Q},0\right)   \arrow{r}{\psi}  & \left(T^\infty(\sigma),\infty\right)
  \end{tikzcd}
\]

\begin{theorem}\cite[Theorem~10.7]{LLMM3}  \cite[Proposition~2.18]{LMM3}\label{chebyshev_center_corr_thm}
\noindent\begin{enumerate}\upshape
\item Each of the sets $\widetilde{T^\infty(\sigma)}$ and $\widetilde{K(\sigma)}$ is completely invariant under the correspondence~$\mathfrak{C}$.

\item On $\widetilde{T^\infty(\sigma)}$, the grand orbits of the correspondence $\mathfrak{C}$ are generated by $\widehat{\eta}$ and $\widehat{\tau}$. Moreover, the group $\langle\widehat{\eta},\widehat{\tau}\rangle=\langle\widehat{\eta}\rangle\ast\langle\widehat{\tau}\rangle$ is conformally conjugate to $\mathbbm{G}_2\cong\Z/2\Z\ast\Z/3\Z$.

\item On $\widetilde{K(\sigma)}\cap\overline{B}(3,2)$, one branch of the forward correspondence is hybrid conjugate to $R\vert_{\mathcal{K}(R)}$, where $R(z)=\overline{z}+1/\overline{z}+1$. The other branch maps $\widetilde{K_a}\cap\overline{B}(3,2)$ onto $\widetilde{K_a}\setminus B(3,2)$. On the other hand, the forward correspondence preserves $\widetilde{K(\sigma)}\setminus B(3,2)$. 

\noindent The backward branches of the correspondence are conjugate to the forward branches via~$\widehat{\eta}$.
\end{enumerate}
\end{theorem}

\section{Cubic examples: Talbot curve, Deltoid-and-Circle, and Apollonian gasket}\label{cubic_examples_sec}

Up to M{\"o}bius conjugacy, there is a unique kissing reflection group generated by circle packings consisting of three circles; namely, the ideal triangle reflection group. Circle packings consisting of four circles give rise to a more interesting collection of kissing reflection groups. This collection contains the \emph{maximal cusp} necklace group of Figure~\ref{necklace_fig} (left) and the classical Apollonian gasket reflection group shown in Figure~\ref{kissing_nielsen_fig} (top right).

In this section, we will describe the dynamics of two explicit examples of piecewise Schwarz reflection maps $\sigma$ such that $\sigma:\sigma^{-1}(\Omega)\to\Omega$ has degree three (in the sense of Subsection~\ref{piecewise_schwarz_subsubsec}). Just like the ideal triangle reflection group shows up in the dynamical study of quadratic Schwarz reflection maps, the two kissing reflection groups mentioned above will naturally appear in the examples of this section.

\subsection{Schwarz reflection in a Talbot curve}\label{talbot_subsec}

Among all simply connected unbounded quadrature domains uniformized by rational maps of global degree $d+1$, the simplest ones have a unique node at $\infty$ (equivalently, their Schwarz reflections have a unique pole at $\infty$, see Subsection~\ref{analysis_connect_subsubsec}). It is easy to see in light of Proposition~\ref{simp_conn_quad_prop} that such a quadrature domain admits a uniformization $f:\D^*\to\Omega$, where $f(z)=z+\frac{a_1}{z}+\cdots+\frac{a_d}{z^d}$, after possibly replacing the quadrature domain by an affine image of it (\cite[Proposition~2.13]{LMM1}). A simply connected unbounded quadrature domain is called \emph{extremal} if it has a unique node at $\infty$ and the boundary $\partial\Omega$ has $d+1$ cusps and $d-2$ double points. The term `extremal' is justified by the fact that these are the maximal possible numbers for a given degree \cite[Lemma~2.4]{LM1}. 
The deltoid is the unique example of such an unbounded quadrature domain for $d=2$. 

\subsubsection{Extremal unbounded quadrature domain in degree three, and associated Schwarz reflection}\label{talbot_dynamics_subsubsec}

\begin{figure}[h!]
\captionsetup{width=0.96\linewidth}
\begin{tikzpicture}
\node[anchor=south west,inner sep=0] at (0,0) {\includegraphics[width=0.5\linewidth]{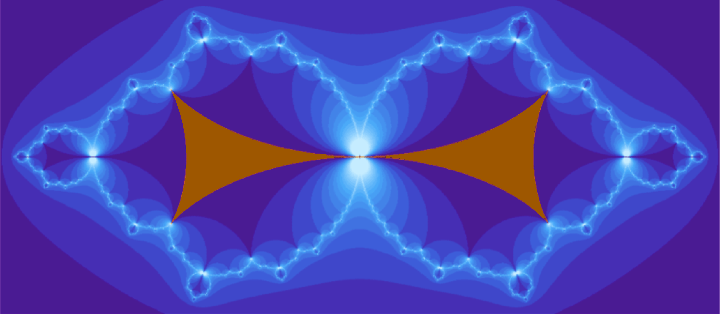}}; 
\node[anchor=south west,inner sep=0] at (7,0) {\includegraphics[width=0.44\linewidth]{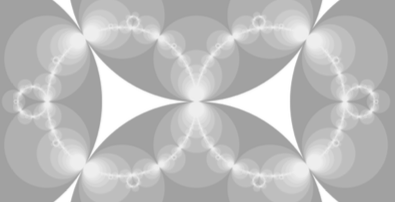}}; 
\node[anchor=south west,inner sep=0] at (3.6,-2.8) {\includegraphics[width=0.45\linewidth]{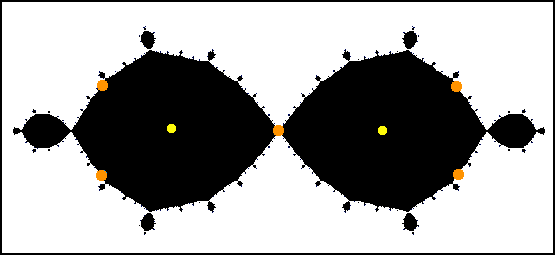}};
\end{tikzpicture}
\caption{Top left: The boundary of the brown pair of triangles is the Talbot curve $\partial\Omega_{\mathrm{ext}}$. The limit set $\Lambda(\sigma_{\textrm{ext}})$ is the glowing blue fractal. Top right: The limit set of the \emph{Julia necklace group} $G$ is homeomorphic to $\Lambda(\sigma_{\textrm{ext}})$. The dynamics of $\sigma_{\mathrm{ext}}$ on its tiling set closure is conformally conjugate to the dynamical of $\cN_G$ on its filled limit set. Bottom: The Julia set of the critically fixed \emph{Julia anti-polynomial} $p(z)=(3\overline{z}-\overline{z}^3)/2$ is homeomorphic to $\Lambda(\sigma_{\textrm{ext}})$, and $p\vert_{\mathcal{J}(p)}$ is topologically conjugate to $\sigma_{\mathrm{ext}}\vert_{\Lambda(\sigma_{\textrm{ext}})}$.}
\label{talbot_schwarz_fig}
\end{figure}

When $d=3$, there is a unique extremal quadrature domain $\Omega_{\textrm{ext}}$, up to the action of $\mathrm{Aut}(\C)$, where $\Omega_{\textrm{ext}}$ is the univalent image of $\D^*$ under $f_{\textrm{ext}}(z)=z+\frac{2}{3z}-\frac{1}{3z^3}$ (cf. \cite[Table~1, Theorem~5.1]{LMM1}). In Figure~\ref{talbot_schwarz_fig} (top left), $\Omega_{\textrm{ext}}$ is the complement of the brown region whose boundary is a so-called \emph{Talbot curve} (cf. \cite[p. 157]{Loc61}). 

As $f_{\textrm{ext}}$ has a triple pole at the origin, the associated Schwarz reflection map $\sigma_{\textrm{ext}}$ has a superattracting fixed point at $\infty$ of local degree three. The basin of attraction of the superattracting fixed point $\infty$ is the simply connected domain given by the exterior of the limit set $\Lambda(\sigma_{\textrm{ext}})$, which is the glowing blue curve in Figure~\ref{talbot_schwarz_fig} (top left). Thus, the action of $\sigma_{\textrm{ext}}$ on its basin of infinity is conformally conjugate to the action of $\overline{z}^3$ on the disk.

The above structure allows one to study $\Lambda(\sigma_{\textrm{ext}})$ from outside using external dynamical rays (as in the case of polynomials). More precisely, one can give a topological model of the limit set as a quotient $\faktor{\mathbb{S}^1}{\sim}$, where $\sim$ is the $m_{-3}$-invariant equivalence relation generated by the angles of the two $2$-periodic rays (under $m_{-3}$) landing at the unique double point of $\partial\Omega_{\textrm{ext}}$ (see \cite[Proposition~48]{LMM2}.

This description allows one to show that the limit set of $\sigma_{\mathrm{ext}}$ is homeomorphic to the limit set of the \emph{Julia necklace group}\footnote{The nomenclature `Julia necklace group' will be justified in Subsection~\ref{talbot_crit_fixed_anti_poly_subsubsec}.} $G$ displayed in Figure~\ref{talbot_schwarz_fig} (top right), and this homeomorphism conjugates the Schwarz reflection $\sigma_{\mathrm{ext}}$ to the Nielsen map $\cN_G$ (see \cite[Proposition~56, Remark~57]{LMM2}). Roughly speaking, this can be seen from the following facts:
\smallskip

\noindent$\bullet$ the limit set of $G$ is homeomorphic to $\faktor{\mathbb{S}^1}{\lambda}$, where $\lambda$ is the $\pmb{\mathcal{N}}_3$-invariant equivalence relation generated by the angles of the two $2$-periodic rays (under $\pmb{\mathcal{N}}_3$) landing at the unique accidental parabolic of $G$ lying on $\partial\Pi(G)$ (cf. Subsection~\ref{necklace_subsubsec}), and
\smallskip

\noindent$\bullet$ the equivalence relation $\sim$ (which gives a lamination model of $\Lambda(\sigma_{\textrm{ext}})$) is the push-forward of $\lambda$ under the Minkowski circle homeomorphism $\pmb{\mathcal{E}}_3$, that conjugates $\pmb{\mathcal{N}}_3$ to $m_{-3}$ (see Definition~\ref{push_lami_def}).

In fact, the above homeomorphism between the limit sets of $\sigma_{\mathrm{ext}}$ and $G$ can be conformally extended to the tiling set producing a conjugacy between the dynamics of $\sigma_{\textrm{ext}}$ on its tiling set closure and the dynamics of $\cN_G$ on its filled limit set \cite[Theorem~A, \S 4.6]{LMM2}.
In this sense, $\sigma_{\textrm{ext}}$ is a mating of $\overline{z}^3$ with the Nielsen map of the necklace reflection group $G$.

\subsubsection{Talbot reflection as a limit of pinching deformation}\label{talbot_pinching_subsubsec}

The group $G$ can be constructed as the limit of a quasiconformal deformation of the regular ideal quadrilateral reflection group $\pmb{G}_3$. Such a quasiconformal deformation pinches a pair of opposite sides of the fundamental domain $\Pi(\pmb{G}_3)\cap\D$ so that the corresponding circles of the packing touch in the limit. Thus, it is natural to ask whether the Schwarz reflection map $\sigma_{\textrm{ext}}$ (equivalently, the quadrature domain $\Omega_{\textrm{ext}}$) can be obtained as the limit of a quasiconformal deformation of a base Schwarz reflection map $\sigma_0$ such that
\smallskip

\noindent$\bullet$ the droplet associated with $\sigma_0$ is a quadrilateral (akin to $\Pi(\pmb{G}_3)\cap\D$), and
\smallskip

\noindent$\bullet$ the quasiconformal deformations pinch a pair of opposite sides of the droplet.

An affirmative answer to this question was given in \cite{LMM1}. Specifically, one can take the quadrature domain $\Omega_0:=f_0(\D^*)$, where $f_0(z)=z-\frac{1}{3z^3}$, as a starting point of the desired pinching deformation. The boundary of the quadrature domain $\Omega_0$ is a classical \emph{astroid curve}. The corresponding Schwarz reflection $\sigma_0$ behaves much like the deltoid reflection map, and techniques mentioned in Subsections~\ref{deltoid_limit_subsubsec},~\ref{deltoid_unique_mating_subsubsec} can be used to justify that $\sigma_0$ is a mating of $\overline{z}^3$ and the Nielsen map $\pmb{\cN}_3$ associated with $\pmb{G}_3$ (cf. \cite[Appendix~B]{LLMM3}). 

One can now deform the Schwarz reflection map $\sigma_0$ quasiconformally so that the non-escaping dynamics remains conformally equivalent to $\overline{z}^3$, while the moduli of the deformed droplets tend to $\infty$. Since one side of the dynamics is `frozen', standard compactness arguments show that there exists a limiting quadrature domain which is extremal; i.e., it has a unique node at $\infty$, and its boundary has four cusps and one double point (see \cite[\S 4.1]{LMM1} for details). That this extremal quadrature domain is affinely equivalent to $\Omega_{\mathrm{ext}}$ follows from the rigidity theorem \cite[Theorem~5.1]{LMM1}.

In Subsection~\ref{talbot_david_subsubsec}, we will outline a completely different recipe for constructing such extremal quadrature domains.

\subsubsection{Relation with a critically fixed anti-polynomial}\label{talbot_crit_fixed_anti_poly_subsubsec}

Consider the critically fixed cubic anti-polynomial $p(z)=(3\overline{z}-\overline{z}^3)/2$, which has two fixed critical points in $\C$. We call the map $p$ the \emph{Julia anti-polynomial} since the dynamics of the cubic polynomial $\overline{p(z)}$ was originally studied by Julia (cf. \cite[p. 51]{julia-1918}). The filled Julia set of $p$ is displayed in Figure~\ref{talbot_schwarz_fig}. Standard arguments from polynomial dynamics show that the lamination of $p$ is precisely $\sim$, and hence $\mathcal{J}(p)\cong\faktor{\mathbb{S}^1}{\sim}$ is equivariantly homeomorphic to the limit set of $\sigma_{\mathrm{ext}}$ \cite[Proposition~64, Remark~66]{LMM2}.

An interesting consequence of the above fact is that $\Lambda(G)$ is homeomorphic to $\mathcal{J}(p)$ and the homeomorphic conjugates $\mathcal{N}_G\vert_{\Lambda(G)}$ to the action of the Julia anti-polynomial $p$ on its Julia set $\mathcal{J}(p)$. (This is the reason why the group $G$ is called the `Julia necklace group'.)

Like the Schwarz reflection map $\sigma_{\mathrm{ext}}$, the anti-polynomial $p$ also enjoys an extremality property. It is a cubic anti-polynomial with the maximal possible number of planar fixed points (these fixed points are marked in yellow/orange in Figure~\ref{talbot_schwarz_fig}.

\subsubsection{Constructing $\sigma_{\mathrm{ext}}$ via David surgery}\label{talbot_david_subsubsec}

The critically fixed anti-polynomial $p$ of Subsection~\ref{talbot_crit_fixed_anti_poly_subsubsec} has two invariant bounded Fatou components $\mathcal{U}_i$, and the restriction $p\vert_{\overline{\mathcal{U}_i}}$ is conformally conjugate to $\overline{z}^2\vert_{\D}$, for $i\in\{1,2\}$. Since the inverse of the Minkowski circle homeomorphism $\pmb{\mathcal{E}}_2$ (that conjugates $\overline{z}^2$ to $\pmb{\cN}_2$) admits a David extension to $\D$, one can replace $p\vert_{\overline{\mathcal{U}_i}}$ with $\pmb{\cN}_2:\overline{\D}\setminus\Int{\Pi(\pmb{G}_2)}\to\overline{\D}$, and uniformize this partially defined topological map of $\mathbb{S}^2$ to a Schwarz reflection map $\sigma$. 

Since $p$ was not altered on its basin of infinity, one concludes that the Schwarz reflection map $\sigma$ has a superattracting fixed point of local degree three. Moreover, since the fixed points of $p$ on $\partial\mathcal{U}_1\cup\partial\mathcal{U}_2$ (the orange points in Figure~\ref{talbot_schwarz_fig}) are identified with the ideal boundary points of $\Pi(\pmb{G}_2)$, it is easy to see that the droplet $T(\sigma)$ is the union of two topological triangles touching at a common vertex. Moreover, $W^{1,1}$-removability of $\mathcal{J}(p)$ and analytic properties of David homeomorphisms imply that the limit set of $\sigma$ is conformally removable (cf. Subsection~\ref{conf_removable_subsec}). Using these facts, one can argue that the Schwarz reflection map $\sigma$ is M{\"o}bius conjugate to $\sigma_{\mathrm{ext}}$ studied above (see \cite[Theorem~12.8]{LMMN} for details).

As a fallout of this construction, one concludes that $\Lambda(\sigma_{\mathrm{ext}})$ is conformally removable, and hence $\sigma_{\mathrm{ext}}$ is the unique conformal mating of $\overline{z}^3$ and $\cN_G$.

We summarize some key points from the above discussion in the following theorem.

\begin{theorem}\label{talbot_schwarz_group_anti_poly_thm}
Let $\sigma_{\mathrm{ext}}$ be the Schwarz reflection map associated with the quadrature domain $f_{\mathrm{ext}}(\D^*)$, where $f_{\textrm{ext}}(z)=z+\frac{2}{3z}-\frac{1}{3z^3}$. Then the following hold.
\noindent\begin{enumerate}\upshape
\item $\sigma_{\mathrm{ext}}$ is the unique conformal mating of $\overline{z}^3\vert_{\overline{\D}}$ and the Nielsen map $\mathcal{N}_G\vert_{K(G)}$, where $G$ is the Julia necklace group shown in Figure~\ref{talbot_schwarz_fig}.

\item The dynamical systems $\sigma_{\mathrm{ext}}\vert_{\Lambda(\sigma_{\mathrm{ext}})},\ \mathcal{N}_G\vert_{\Lambda(G)}$, and $p\vert_{\mathcal{J}(p)}$ are topologically conjugate (where $p$ is the Julia anti-polynomial).
\end{enumerate}
\end{theorem}

\subsection{Apollonian gasket and its cousins}\label{apollo_group_map_schwarz_subsec}

As in Subsection~\ref{talbot_subsec}, we will look at three homeomorphic fractals in this subsection: the limit set of a kissing reflection group, the Julia set of a critically fixed anti-rational map, and the limit set of a cubic Schwarz reflection~map.

\subsubsection{From the Apollonian reflection group $G$ to a cubic anti-rational map $R$ via the Thurston Realization Theorem}\label{apollo_group_rat_map_subsubsec}

The Apollonian gasket is the limit set of the reflection group $G$ generated by reflections in the four red circles displayed in Figure~\ref{kissing_nielsen_fig} (top right). Let us call the four components of $\Omega(G)$ that intersect the fundamental domain $\Pi(G)$ the \emph{principal components} of $\Omega(G)$, and denote them by $\mathcal{U}_i$, $i\in\{1,2,3,4\}$. The restriction of the Nielsen map $\mathcal{N}_G$ to each $\overline{\mathcal{U}_i}$ is  M{\"o}bius-conjugate to $\pmb{\mathcal{N}}_2$. We extend $\pmb{\mathcal{E}}_2:\mathbb{S}^1\to\mathbb{S}^1$ to a self-homeomorphism of $\overline{\D}$, and define a global orientation-reversing critically fixed branched cover of degree three as:
$$
\widetilde{R}:\widehat{\C}\to\widehat{\C},\quad 
z \mapsto \left\{\begin{array}{ll}
                    \mathcal{N}_G(z) & \mbox{if}\ z\in \widehat{\C}\setminus\bigcup_{i=1}^4 \mathcal{U}_i, \\
                   \phi_i^{-1}\left(\pmb{\mathcal{E}}_2^{-1} \left(\overline{\pmb{\mathcal{E}}_2(\phi_i(z))}^2\right)\right) & \mbox{if}\ z\in \mathcal{U}_i,
                            \end{array}\right.
$$
where $\phi_i:\mathcal{U}_i\to\D$ is a M{\"o}bius map, $i\in\{1,2,3,4\}$.

\begin{figure}[h!]
\captionsetup{width=0.96\linewidth}
\centering
\includegraphics[width=0.4\linewidth]{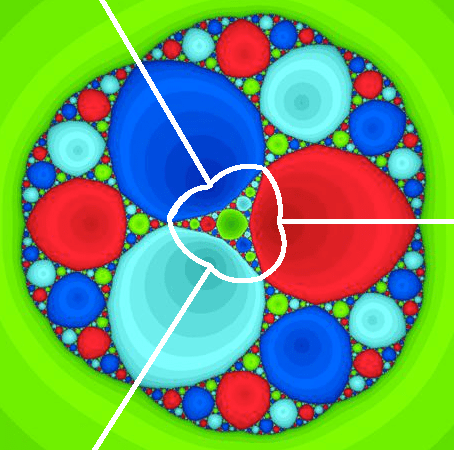}\quad \includegraphics[width=0.4\linewidth]{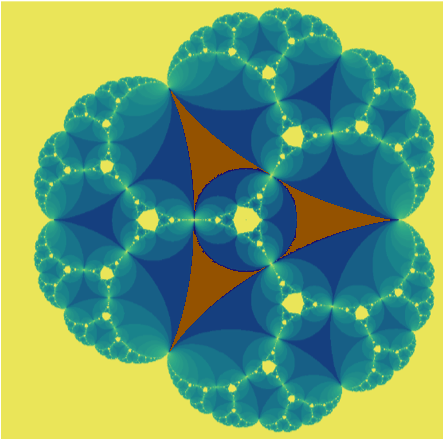}
\caption{Left: The Julia set of $R$ and its Tischler graph is shown (in the chosen coordinates, the point at $\infty$ is a superattracting fixed point). Right: The dynamical plane of the Deltoid-and-Circle Schwarz reflection map is displayed. The droplet consists of the three brown triangles, while the non-escaping set and the tiling set are shown in yellow and blue/green, respectively.}
\label{apollo_cousins_fig}
\end{figure}

It was shown in \cite{LLMM4} that $\widetilde{R}$ has no Thurston obstruction, and hence by an antiholomorphic version of the Thurston Realization Theorem (see \cite[Proposition~6.1]{LLMM4}, \cite[Theorem~3.9]{Gey20}), it is equivalent to a critically fixed cubic anti-rational map~$R$. In fact, the expansiveness property of $\widetilde{R}$ on $\Lambda(G)$ (coming from circular reflections) was used to show that $\widetilde{R}\vert_{\Lambda(G)}\equiv \mathcal{N}_G\vert_{\Lambda(G)}$ is topologically conjugate to $R\vert_{\mathcal{J}(R)}$. Moreover, $R$ can be chosen to be 
$$
R(z)=\frac{3\overline{z}^2}{2\overline{z}^3+1}
$$ 
(see Figure~\ref{apollo_cousins_fig}) \cite[Corollary~8.2]{LLMM4}. In particular, $\mathcal{J}(R)$ is homeomorphic to~$\Lambda(G)$.

There is a natural forward-invariant graph in the dynamical plane of $R$, that contains all the fixed critical points. This graph, which captures the touching structure of the fixed Fatou components of $R$, is called the \emph{Tischler graph} of $R$ (see Figure~\ref{apollo_cousins_fig}). It is worth pointing out that the planar dual of the Tischler graph of $R$ is isomorphic as a plane graph to the contact graph of the circle packing giving rise to $G$ (in fact, both graphs are isomorphic to the $1$-skeleton of a tetrahedron).
We will discuss this duality in a general framework in Subsection~\ref{new_line_dict_dyn_subsec}.

\subsubsection{Pinching geodesics and matings}\label{qf_bdry_mating_subsubsec}

The Apollonian group $G$ lies on the boundary of the quasi-Fuchsian deformation space of the regular ideal quadrilateral reflection group $\pmb{G}_3$. Indeed, the group $G$ can be constructed as a limit of a quasiconformal deformation $\{G_n\}$ of $\pmb{G}_3$ such that the moduli of the two (quadrilateral) components of $\Pi(G_n)$ go to zero and hence the non-adjacent circles of the circle packings defining $G_n$ touch pairwise in the limit. Equivalently, such a deformation pinches suitable simple closed non-peripheral geodesics on the punctured spheres $\D/\widetilde{\pmb{G}}_3$ and $(\widehat{\C}\setminus\overline{\D})/\widetilde{\pmb{G}}_3$, where $\widetilde{\pmb{G}}_3$ is the index two Fuchsian subgroup of $\pmb{G}_3$ (cf. \cite[\S 3.3, \S 3.4]{LLM1}). According to \cite[Proposition~3.21]{LLM1}, the group $G$ can also be interpreted as the mating of two copies of the Julia necklace group considered in Subsection~\ref{talbot_dynamics_subsubsec} (and displayed in Figure~\ref{talbot_schwarz_fig}).

The Apollonian anti-rational map $R$ also enjoys a mating description akin to that of $G$. Recall that the limit set of the Julia necklace group is equivariantly homeomorphic to the Julia set of the Julia anti-polynomial $p$ introduced in Subsection~\ref{talbot_crit_fixed_anti_poly_subsubsec}. It turns out that $R$ is a mating of two copies of $p$
(cf. \cite[\S 4.2, Remark~6.14]{LLMM4}, \cite[Corollary~4.17, \S 4.3]{LLM1}). 

Thus, $G$ is a mating of two necklace groups and $R$ is a mating of two critically fixed anti-polynomials in a compatible manner.

\subsubsection{A nearly affine model for the Apollonian anti-rational map}\label{affine_model_subsubsec}

The anti-rational map $R$ constructed in Subsection~\ref{apollo_group_rat_map_subsubsec} admits an anti-quasiregular model $\mathfrak{R}$ on a tetrahedron \cite[\S 8]{LLMM4}. In fact, this model has the added advantage of being piecewise affine outside its `invariant Fatou components'.
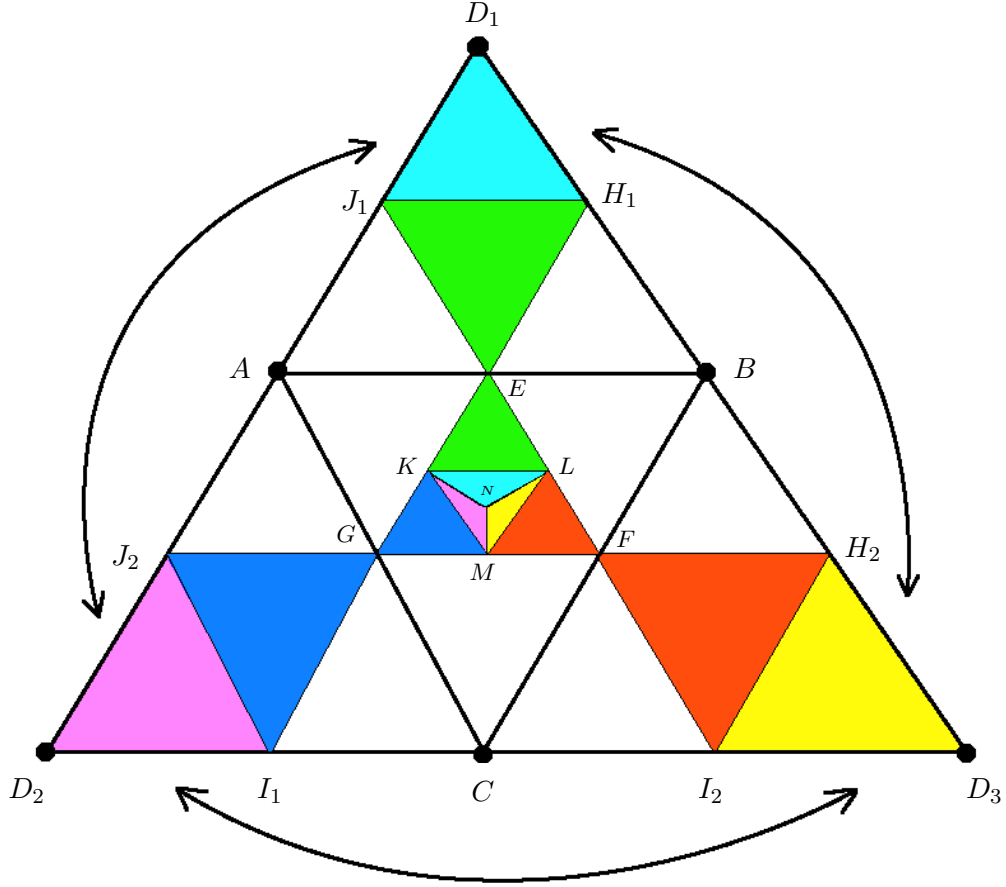
\begin{figure}[h!]
\captionsetup{width=0.96\linewidth}
\begin{tikzpicture}
\node[anchor=south west,inner sep=0] at (-4,0) {\includegraphics[width=1\textwidth]{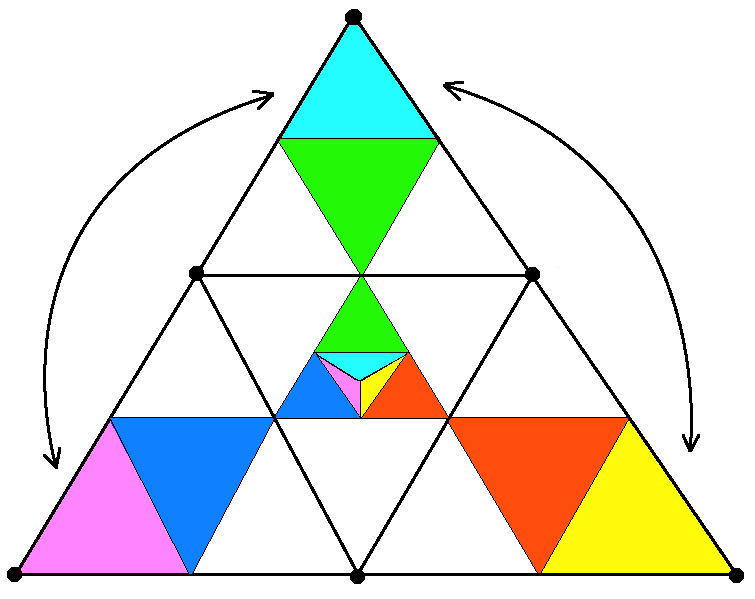}};
\node[anchor=south west,inner sep=0] at (-2.4,-1.5) {\includegraphics[width=0.75\textwidth]{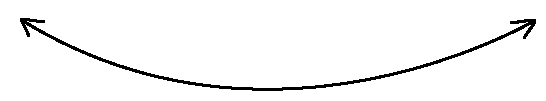}};
\node at (-1.2,5.5) {\begin{Large}$A$\end{Large}};
\node at (5.45,5.5) {\begin{Large}$B$\end{Large}};
\node at (2,-0.1) {\begin{Large}$C$\end{Large}};
\node at (2,10.2) {\begin{Large}$D_1$\end{Large}};
\node at (-4,-0.1) {\begin{Large}$D_2$\end{Large}};
\node at (8.6,-0.1) {\begin{Large}$D_3$\end{Large}};
\node at (-2.7,3) {\begin{Large}$J_2$\end{Large}};
\node at (0.32,7.7) {\begin{Large}$J_1$\end{Large}};
\node at (0.2,3.3) {$G$};
\node at (2.45,5.24) {$E$};
\node at (3.88,3.24) {$F$};
\node at (5,-0.1) {\begin{Large}$I_2$\end{Large}};
\node at (-0.8,-0.1) {\begin{Large}$I_1$\end{Large}};
\node at (7,3.1) {\begin{Large}$H_2$\end{Large}};
\node at (3.8,7.8) {\begin{Large}$H_1$\end{Large}};
\node at (1,4.2) {$K$};
\node at (3.1,4.2) {$L$};
\node at (2,2.8) {$M$};
\node at (2.06,3.88) {\begin{tiny}$N$\end{tiny}};
\end{tikzpicture}
\caption{Pictured is a triangle that can be folded to construct a tetrahedron. The arrows indicate the required foldings.}
\label{affine_model_fig}
\end{figure}
We describe the map $\mathfrak{R}$ using Figure~\ref{affine_model_fig}. The tetrahedron can be obtained from the union of the four equilateral triangles $ABC, ABD_1, ACD_2$, and $BCD_3$, where one identifies the sides $AD_1, BD_1$, and $CD_2$ with $AD_2, BD_3$, and $CD_3$, respectively, Thus, the vertices $D_1, D_2, D_3$ correspond to the same vertex on the tetrahedron. 

The union of the triangles $D_1J_1H_1, D_2 I_1J_2,$ and $D_3 I_2H_2$ corresponds to a Jordan domain in the tetrahedron. The map $\mathfrak{R}$ is quasiconformally conjugate to $\overline{z}^2$ on this Jordan domain such that it restricts to a piecewise affine orientation-reversing double covering map on the boundary. For instance, it maps the edge $J_1H_1$ affinely to the union of the edges $J_2I_1$ and $H_2I_2$. We call this Jordan domain an \emph{invariant Fatou component}.
The map $\mathfrak{R}$ is similarly defined on the other three invariant Fatou components containing $A, B, C$ (these components cover the white region).

It remains to specify $\mathfrak{R}$ on the equilateral triangles $EFG, EH_1J_1, GJ_2I_1$, and $FH_2I_2$. 
For definiteness, let us consider the triangle $EFG$. One can map the quadrilaterals $EKNL$, $GMNK$, and $FLNM$ onto the quadrilaterals $EJ_1D_1H_1$, $GI_1D_2J_2$, and $FH_2D_3I_2$ in a color-preserving way such that the maps are anti-conformal (Euclidean) reflections on the triangles $EKL$, $GMK$, $FLM$, and affine on the triangles $KNL$, $MNK$, $LNM$. The definition of $\mathfrak{R}$ on the other three triangles  $EH_1J_1$, $GJ_2I_1$, and $FH_2I_2$ is symmetric.

In fact, the affine nature of the construction guarantees that $\mathfrak{R}$ is a degree three anti-quasiregular map of the tetrahedron. By construction, the only critical points of $\mathfrak{R}$ are fixed and lie in the four invariant Fatou components. Finally, the fact that $\mathfrak{R}$ is anti-conformal outside the first preimages of its invariant Fatou components implies that it can be straightened to an anti-rational map. Since this cubic anti-rational map has four distinct fixed critical points, it is easily seen to be M{\"o}bius conjugate to the map $R$ of Subsection~\ref{apollo_group_rat_map_subsubsec} (cf. \cite[Proposition~8.1, Corollary~8.2]{LLMM4}).

\subsubsection{From the cubic anti-rational map $R$ to the Apollonian reflection group $G$ via David surgery}\label{rat_map_apollo_group_subsubsec}

One can also go backwards from the anti-rational map $R$ to the Apollonian reflection group $G$ using David surgery. More precisely, since $\pmb{\mathcal{E}}_2^{-1}$ admits a David extension to $\overline{\D}$, one can replace the dynamics of $R$ on its critically fixed Fatou components with $\pmb{\cN}_2:\overline{\D}\setminus\Int{\Pi(\pmb{G}_2)}\to\overline{\D}$ (cf. Subsections~\ref{deltoid_david_subsubsec} and~\ref{talbot_david_subsubsec}). Finally, one can invoke the David Integrability Theorem to uniformize this partially defined topological map of $\mathbb{S}^2$ to a Schwarz reflection map, and use the mapping properties of the resulting Schwarz reflection to justify that it is indeed the Nielsen map of $G$ (see \cite[\S 10.2]{LLMM4} for details of this construction).

\subsubsection{Schwarz reflection in the deltoid and a circle}\label{d_and_c_schwarz_subsubsec}

One can modify the construction of Subsection~\ref{rat_map_apollo_group_subsubsec} to obtain a hybrid dynamical system. Specifically, if one replaces the dynamics of $R$ on three of its invariant Fatou components with the Nielsen map $\pmb{\cN}_2:\overline{\D}\setminus\Int{\Pi(\pmb{G}_2)}\to\overline{\D}$ using the above David surgery (for instance, on the three $2\pi/3-$rotation symmetric invariant Fatou components in Figure~\ref{apollo_cousins_fig} (left)) and leaves the action of $R$ on the fourth invariant Fatou component unaltered, then one obtains a piecewise Schwarz reflection map $\sigma$. Using an explicit characterization of the deltoid Schwarz reflection map, it was shown in \cite[Proposition~10.5]{LLMM4} that $\sigma$ arises from the deltoid (of Subsection~\ref{deltoid_subsec}) and an inscribed circle. Roughly speaking, the appearance of the Deltoid-and-Circle Schwarz reflection map in the current setting can be seen as follows. The $2\pi/3-$rotational symmetry of the construction and the fact that $\overline{z}^2$ commutes with rotation by $2\pi/3$ together guarantee that $\sigma$ also admits a $2\pi/3-$rotational symmetry. One can combine this symmetry property of $\sigma$ with precise information about the critical/singular points of $\sigma$ to conclude that $\sigma$ is the Schwarz reflection map of the deltoid and an inscribed circle (see the first paragraph of Subsection~\ref{deltoid_subsec})).

In a suitable sense, the Deltoid-and-Circle Schwarz reflection $\sigma$ combines the action of the Nielsen map $\cN_G$ of the Apollonian reflection group and the critically fixed cubic anti-rational map $R$ (see \cite[Theorem~10.8]{LLMM4} for a precise statement). Alternatively, one can explore the dynamics of the Schwarz reflection of the Deltoid-and-Circle (using the dynamics of the deltoid reflection map discussed in Subsection~\ref{deltoid_subsec}) to directly recognize  this map as a mating of $R$ and $\cN_G$.
The dynamical plane of $\sigma$ is depicted in Figure~\ref{apollo_cousins_fig}~(right).

\subsubsection{Quasisymmetry groups}\label{qs_grp_subsubsec}

Note that the boundaries of the fixed Fatou components of $R$ intersect at positive angles, while the boundaries of the principal components of $\Omega(G)$ intersect tangentially. This implies that the dynamically meaningful homeomorphism between $\mathcal{J}(R)$ and $\Lambda(G)$ described in Subsection~\ref{apollo_group_rat_map_subsubsec} does not admit a quasiconformal extension to the sphere. In fact, it has been proved in \cite{LZ23b} that there does not exist any global quasiconformal map that induces a homeomorphism between the Julia set of $R$ and the limit set of $G$.

However, the quasiconformal non-equivalence of $\mathcal{J}(R)$ and $\Lambda(G)$ is \emph{not} captured by their groups of quasisymmetries. According to \cite[Theorems~3.8,~7.2]{LLMM4}, these fractal have isomorphic quasisymmetry groups, and these quasisymmetry groups coincide with the respective groups of self-homeomorphisms (see \cite{LZ23a} for a general quasisymmetric rigidity result which implies that the homeomorphism group of $\Lambda(G)$ coincides with its conformal symmetry group).

On the other hand, the quasisymmetry group of the limit set of the Deltoid-and-Circle reflection map is a strict subgroup of its homeomorphism group. This is essentially a consequence of the fact that the green complementary components of this limit set have inward pointing cusps, but all cusps on the boundaries of the yellow components are outward pointing cusps, implying that there is no quasisymmetry carrying the boundary of a green component to the boundary of a yellow component. 

The main results discussed in this subsection are encapsulated below.

\begin{theorem}\label{apollo_schwarz_group_anti_rat_thm}
Let $G$ be the Apollonian reflection group, $R$ be the critically fixed anti-rational map $\frac{3\overline{z}^2}{2\overline{z}^3+1}$, and $\sigma$ be the Deltoid-and-Circle Schwarz reflection map.
Then the following hold.
\begin{itemize}
\item There exists a global David homeomorphism which restricts to a topological conjugacy between $R\vert_{\mathcal{J}(R)}$ and $\cN_G\vert_{\Lambda(G)}$.

\item There exists a global David homeomorphism which restricts to a topological conjugacy between $R\vert_{\mathcal{J}(R)}$ and $\sigma\vert_{\Lambda(\sigma)}$.

\item The quasisymmetry groups of $\mathcal{J}(R)$ and $\Lambda(G)$ coincide with the respective self-homeomorphism groups, and hence they are isomorphic. On the other hand, the quasisymmetry group of $\Lambda(\sigma)$ is strictly smaller than its group of self-homeomorphisms.

\item $G$ is the mating of two copies of a necklace group $G_1$ and $R$ is the mating of two copies of a critically fixed anti-polynomial $p$, where $\mathcal{N}_{G_1}\vert_{\Lambda(G_1)}$ and $p\vert_{\mathcal{J}(p)}$ are topologically conjugate.
\end{itemize}
\end{theorem}

\section{Parameter spaces of special families of Schwarz reflections}\label{schwarz_para_space_sec}

After analyzing interconnections between various examples of Schwarz reflection maps, anti-rational maps and Kleinian reflection groups, we will  now explicate some general results about families of Schwarz reflection maps. Specifically, we will focus on families of Schwarz reflections containing the examples from Subsections~\ref{c_and_c_center_subsec},~\ref{chebyshev_subsec}, and~\ref{talbot_subsec}. In all three cases, mating descriptions of the Schwarz reflection maps can be given via suitable straightening theorems, and this in turn can be used to relate the topology of parameter spaces of Schwarz reflections to that of parameter spaces of appropriate families of anti-rational maps or reflection groups.

\subsection{The Circle-and-Cardioid family}\label{c_and_c_general_subsec}

Recall from Subsection~\ref{c_and_c_center_subsec} that $\heartsuit=f(\D)$, where $f(w)=w/2-w^2/4$ is a simply connected quadrature domain whose Schwarz reflection map is denoted by $\sigma$. Moreover, the only critical point of $\sigma$ is at the origin. 

For $a\in\C$, let $B(a,r_a)$ be the smallest disk containing $\heartsuit$ and centered at $a$; i.e., $\partial B(a,r_a)$ is the circumcircle to $\heartsuit$.
By \cite[Proposition~5.9]{LLMM1}, we have the following dichotomy:
\begin{itemize}
\item for $a\in\left(-\infty,-1/12\right)$, the circumcircle $\partial B(a,r_a)$ touches $\partial\heartsuit$ at exactly two points, and
\item for any $a\in\C\setminus\left(-\infty,-1/12\right)$, the circumcircle $\partial B(a,r_a)$ touches $\partial\heartsuit$ at exactly one point.
\end{itemize}
In order to extract triangle group structure from our family of Schwarz reflection maps, we will restrict to parameters $a\in\C\setminus\left(-\infty,-1/12\right)$.
We will denote the anti-M{\"o}bius reflection in $\partial B(a,r_a)$ by $\sigma_a$. For $a\in\C\setminus\left(-\infty,-1/12\right)$, the unique touching point of $\partial B(a,r_a)$ and $\partial\heartsuit$ is denoted by $\alpha_a$.
The union $\Omega_a$ of the disjoint quadrature domains $\heartsuit$ and $\overline{B}(a,r_a)^c$ 
gives rise to a piecewise Schwarz reflection map
$$
F_a:\overline{\Omega_a}\to\widehat{\C},\quad 
z\mapsto
\begin{cases}
\sigma(z)\quad \textrm{if}\quad z\in\overline{\heartsuit},\\
\sigma_a(z)\quad \textrm{if}\quad z\in B(a,r_a)^c,
\end{cases}
$$
(cf. Subsection~\ref{piecewise_schwarz_subsubsec}).
We denote this family of maps by $\mathcal{S}$; i.e., 
$$
\mathcal{S}:=\lbrace F_a:\overline{\Omega}_a\to\widehat{\C}:a\in\C\setminus (-\infty,-1/12)\rbrace,
$$
and call it the \emph{Circle-and-Cardioid} or \emph{C\&C} family.
Note that $1/4$ and $\alpha_a$ are the only singular points on the boundary of $\Omega_a$. Following Subsection~\ref{inv_partition_subsubsec}, we define the fundamental tile, the tiling set, and the non-escaping set as
$$
T^0_a\equiv T^0(F_a):=\widehat{\C}\setminus\left(\Omega_a\cup\{1/4,\alpha_a\}\right),\quad T^\infty_a\equiv T^\infty(F_a):=\bigcup_{n\geq 0} F_a^{-n}(T^0_a),\quad \mathrm{and}
$$
$$
K_a\equiv K(F_a):=\widehat{\C}\setminus T^\infty(F_a).
$$

\begin{figure}[h!]
\captionsetup{width=0.96\linewidth}
\includegraphics[width=0.32\linewidth]{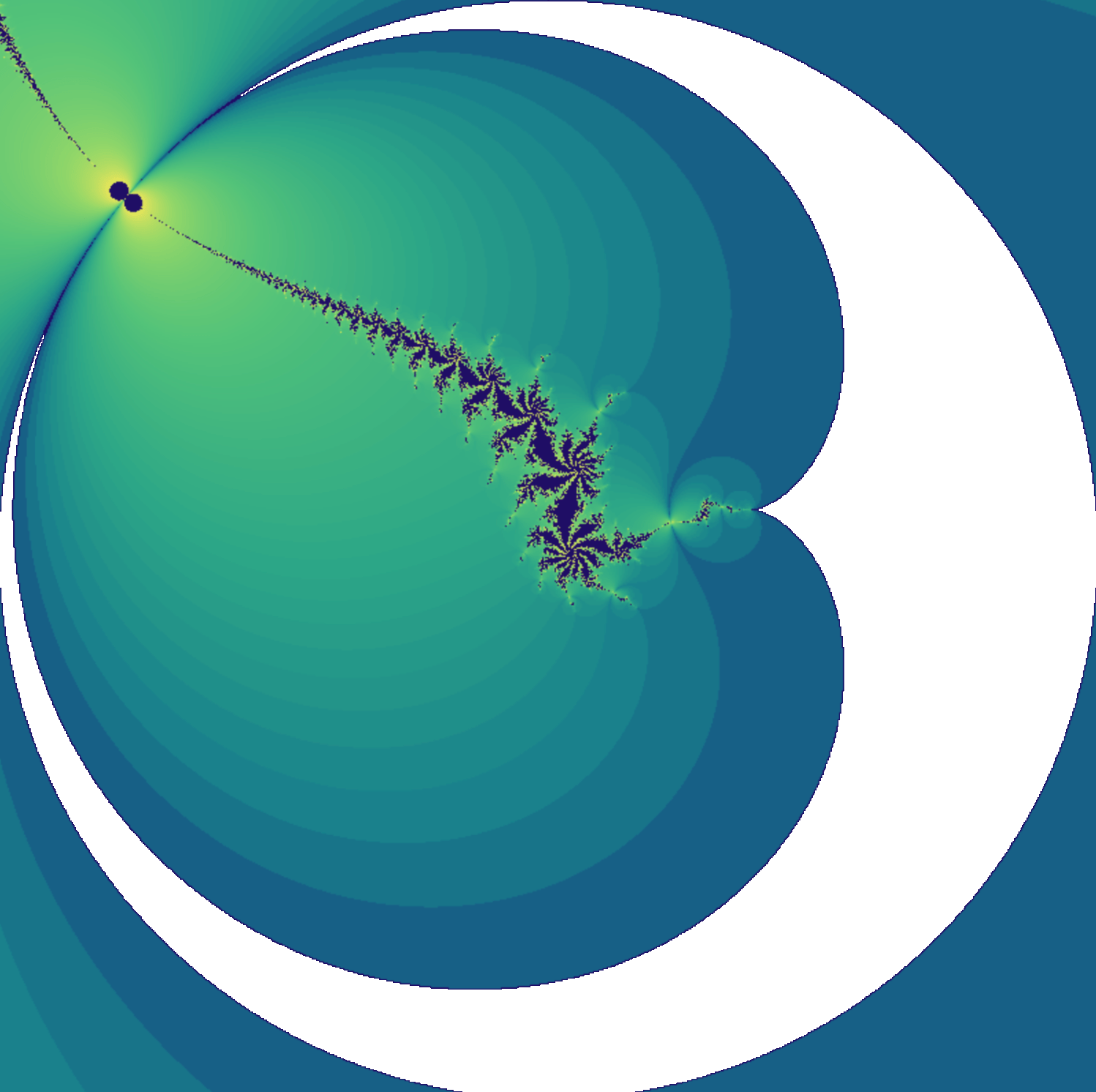}\ \includegraphics[width=0.32\linewidth]{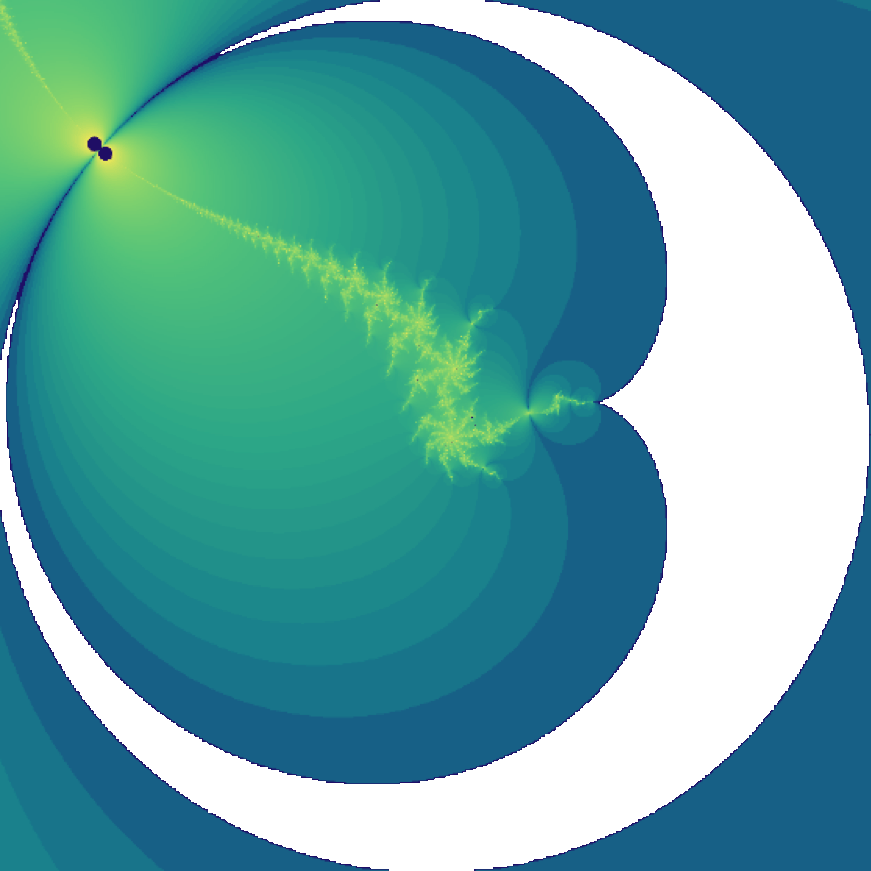}\ \includegraphics[width=0.32\linewidth]{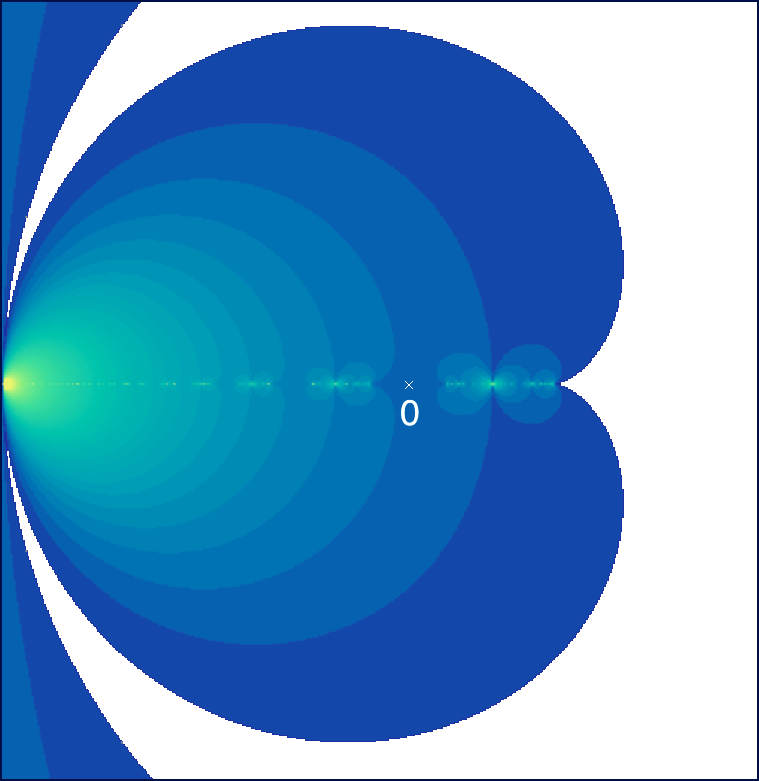}
\caption{Parts of the non-escaping sets of various maps in the family $\cS$ are shown.}
\label{c_and_c_conn_cantor_limit_fig}
\end{figure}

\subsubsection{The connectedness locus of $\cS$}\label{c_and_c_conn_locus_subsubsec}

As in the case of (anti-)polynomial-like maps, it turns out that connectedness of the non-escaping set $K(F_a)$ is equivalent to the requirement that the orbit of the unique critical point $0$ of $F_a$ does not escape. Indeed, if the critical point $0$ does not lie in the tiling set of $F_a$, then the vertex-preserving conformal isomorphism $\psi_a:T_a^0\to\Pi(\pmb{G}_2)$ can be extended via iterated lifting (or Schwarz reflections) to a conformal isomorphism $\psi_a:T^\infty_a\to\D$ conjugating $F_a$ to the Nielsen map $\pmb{\mathcal{N}}_2$.

\begin{theorem}[Basic dichotomy]\label{group_dyn_c_and_c_thm}\cite[Theorem 1.2]{LLMM1}
\noindent\noindent\begin{enumerate}\upshape 
\item If the critical point of $F_a$ does not escape to the fundamental tile $T_a^0$, then the conformal map $\psi_a$ from $T_a^0$ onto $\Pi(\pmb{G}_2)$ extends to a biholomorphism between the tiling set $T_a^\infty$ and the unit disk $\D$. Moreover, the extended map $\psi_a$ conjugates $F_a$ to the Nielsen map $\pmb{\mathcal{N}}_2$. In particular, $K_a$ is connected.

\item If the critical point of $F_a$ escapes to the fundamental tile, then the corresponding non-escaping set $K_a$ is a Cantor set.
\end{enumerate}
\end{theorem}

The above theorem leads to the notion of the \emph{connectedness locus} of the family~$\mathcal{S}$.

\begin{definition}\label{c_and_c_conn_locus_def}
The connectedness locus of the family $\mathcal{S}$ is defined as 
$$
\cC(\mathcal{S})=\{a\in\C\setminus(-\infty,-1/12): 0\notin T_a^\infty\}=\{a\in\C\setminus(-\infty,-1/12): K_a\ \textrm{is\ connected}\}.
$$ 
The complement of the connectedness locus in the parameter space is called the \emph{escape~locus}.
\end{definition}

\subsubsection{Combinatorial straightening and mating description for geometrically finite maps}\label{c_and_c_straightening_mating_subsubsec}

According to Theorem~\ref{group_dyn_c_and_c_thm}, one observes the full tessellation structure of the ideal triangle reflection group in the tiling set of $F_a$ if and only if $a\in\cC(\cS)$ (see Figure~\ref{c_and_c_conn_cantor_limit_fig}). In light of the mating description for the map $F_a$ with $a=0$ given in Subsection~\ref{c_and_c_center_subsec}, it is natural to expect that the non-escaping set dynamics of any postcritically finite map $F_a$ in $\cC(\mathcal{S})$ is conjugate to the filled Julia set dynamics of a quadratic anti-polynomial. 

However, as explained in Subsection~\ref{c_and_c_basilica_pinched_anti_quad_subsubsec}, one cannot \emph{straighten} Schwarz reflections in $\cC(\cS)$ to quadratic anti-polynomials in the Tricorn using quasiconformal tools. In \cite{LLMM2}, a \emph{combinatorial straightening} theory for `nice' maps in $\cC(\cS)$ was developed. In fact, for maps in $\mathcal{C}(\mathcal{S})$, the conformal map $\psi_a:T^\infty_a\to\D$ plays a role akin to that of B{\"o}ttcher coordinates in polynomial dynamics. For postcritically finite maps in $\cC(\cS)$, the \emph{limit set} (i.e., the common boundary of $K_a$ and $T_a^\infty$) turns out to be locally connected (see \cite[Theorem~1.4]{LLMM1}), and hence the conformal map $\psi_a^{-1}$ extends continuously to produce a topological semi-conjugacy between $\pmb{\cN}_2\vert_{\mathbb{S}^1}$ and $F_a\vert_{\partial T_a^\infty}$. This enables one to construct a topological model for the limit set of $F_a$ as a quotient of $\mathbb{S}^1$ under an $\pmb{\cN}_2$-invariant lamination. This lamination can then be turned into an $m_{-2}$-invariant lamination using the topological conjugacy $\pmb{\mathcal{E}}_2$ between $\pmb{\cN}_2\vert_{\mathbb{S}^1}$ and $m_{-2}\vert_{\mathbb{S}^1}$. Finally, standard results in polynomial dynamics (which use the Thurston Realization Theorem or landing properties of external parameter rays of the Tricorn) imply that such $m_{-2}$-invariant laminations are indeed realized by Julia sets of postcritically finite maps in the Basilica limb $\mathcal{L}$ of the Tricorn (see Definition~\ref{def_basilica_limb}). 

A map in $\cC(\cS)$ is called \emph{geometrically finite} if it has an attracting/parabolic cycle or if its unique critical point (at the origin) is non-escaping and strictly pre-periodic.
A good understanding of the dynamics of geometrically finite maps in $\cC(\cS)$ (see \cite[\S 5,6]{LLMM1}) allows one to push the above arguments to all geometrically finite maps and prove the following mating statement.

\begin{theorem}\label{c_and_c_general_mating_thm}\cite[Theorem~1.1]{LLMM2}
Every geometrically finite map $F_a\in\cC(\cS)$ is a conformal mating of a unique geometrically finite quadratic anti-polynomial and the Nielsen map $\pmb{\cN}_2$. More precisely, for each geometrically finite map $F_a\in\cC(\cS)$, the following hold.
\noindent\begin{enumerate}\upshape
\item There exists a topological semi-conjugacy between
$$
F_a:\overline{T_a^\infty}\setminus\Int{T_a^0}\to \overline{T_a^\infty}\quad \textrm{and}\quad \pmb{\cN}_2:\overline{\D}\setminus \Int{\Pi(\pmb{G}_2)}\to\overline{\D}
$$ 
such that the semi-conjugacy restricts to a conformal conjugacy on $T_a^\infty$.

\item There exists a unique geometrically finite quadratic anti-polynomial $f_c\in\mathcal{L}$ such that $F_a\vert_{K_a}$ is topologically conjugate to $f_c\vert_{\mathcal{K}_c}$ and the conjugacy is conformal on the interior of $K_a$. Moreover, the Minkowski circle homeomorphism $\pmb{\mathcal{E}}_2$ pushes forward the lamination associated with the limit set of $F_a$ to the lamination associated with the Julia set of $f_c$.
\end{enumerate}
\end{theorem}
\noindent (See Definition~\ref{push_lami_def} for the definition of push-forward of a lamination.)

The uniqueness part of the above theorem follows from rigidity of geometrically finite polynomials with prescribed lamination and prescribed conformal class of Fatou set~dynamics.

\subsubsection{Bijection between geometrically finite maps of $\cC(\cS)$ and $\mathcal{L}$}\label{c_and_c_geom_fin_bijection_subsubsec}

The combinatorial straightening map from geometrically finite maps in $\cC(\cS)$ to those in $\mathcal{L}$ given by Theorem~\ref{c_and_c_general_mating_thm} is in fact a bijection. To prove injectivity of this map, one needs to establish appropriate rigidity theorems for geometrically finite maps in $\cC(\cS)$ (namely, that they are uniquely determined by their lamination and conformal class of the Fatou set dynamics). The proof of this fact uses a pullback argument and involves an analysis of the boundary behavior of conformal maps near cusps and double points (see \cite[\S 8]{LLMM2}). In the postcritically finite case, one can also appeal to conformal removability of the limit sets (which follows from the John property of the tiling sets) to prove combinatorial rigidity.

On the other hand, surjectivity of the above map amounts to constructing geometrically finite maps in $\cC(\cS)$ with prescribed lamination and conformal data. To achieve this, one needs a dynamically natural uniformization of the escape locus of $\mathcal{S}$. For $a\notin\mathcal{C}(\mathcal{S})\cup(-\infty,-1/12)$, let $n(a)$ be the smallest integer such that $F_a^{\circ n(a)}(\infty)$ lands on the fundamental tile $T_a^0$. Then the action of $F_a$ on the union of preimages of $T_a^0$ up to time $n(a)$ is unramified and hence conformally conjugate to the action of $\pmb{\cN}_2$ (see \cite[Proposition~5.38]{LLMM1}). This conjugacy, which we call $\psi_a$, can be thought of as an analog of B{\"o}ttcher coordinates for polynomials with disconnected Julia sets. Analogous to the uniformization of the exterior of the Mandelbrot set or the Tricorn, one can uniformize the escape locus of $\mathcal{S}$ by the conformal position of the critical value $\infty$; i.e., by the map $\pmb{\Psi}:a\mapsto\psi_a(\infty)$ \cite[Theorem~1.3, \S 6]{LLMM2}. The map $\pmb{\Psi}$ gives rise to parameter tiles in the escape locus of $\cS$, which provides a phase-parameter duality. One can exploit this phase-parameter duality to study landing/accumulation properties of these parameter tiles, and construct the desired geometrically finite maps in $\cC(\cS)$ as limit points of suitable sequences of parameter tiles (cf. \cite[\S 9]{LLMM2}). This is akin to constructing geometrically finite maps of the Mandelbrot set or the Tricorn as limit points of suitable parameter rays.

\begin{theorem}\cite[Theorem~1.2]{LLMM2}
There exists a natural bijection $\chi$ between the geometrically finite parameters in $\cC(\cS)$ and those in the Basilica limb $\mathcal{L}$ of the Tricorn such that the laminations of the corresponding maps are related by the Minkowski circle homeomorphism $\pmb{\mathcal{E}}_2$ and the dynamics on the respective periodic Fatou components are conformally conjugate.
\end{theorem}

In Section~\ref{mating_para_space_sec}, we will discuss higher degree generalizations of the C\&C family and extensions of the above result to generic maps in such families.

\subsubsection{Combinatorial model of $\cC(\cS)$}\label{c_and_c_conn_locus_model_subsubsec}  

The landing/accumulation properties of the parameter tiles of the escape locus of $\cS$ also allows one to study the topology of $\cC(\cS)$ from outside. It turns out that the landing/accumulation patterns of parameter rays at parabolic/Misiurewicz points of $\cC(\cS)$ and $\mathcal{L}$ are compatible.
\begin{figure}[h!]
\captionsetup{width=0.96\linewidth}
\includegraphics[width=0.4\linewidth]{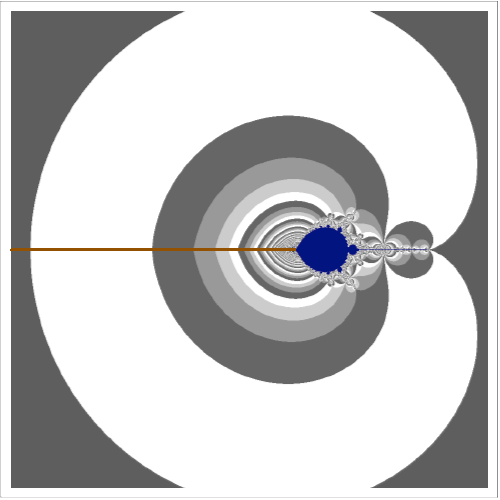}\ \includegraphics[width=0.56\linewidth]{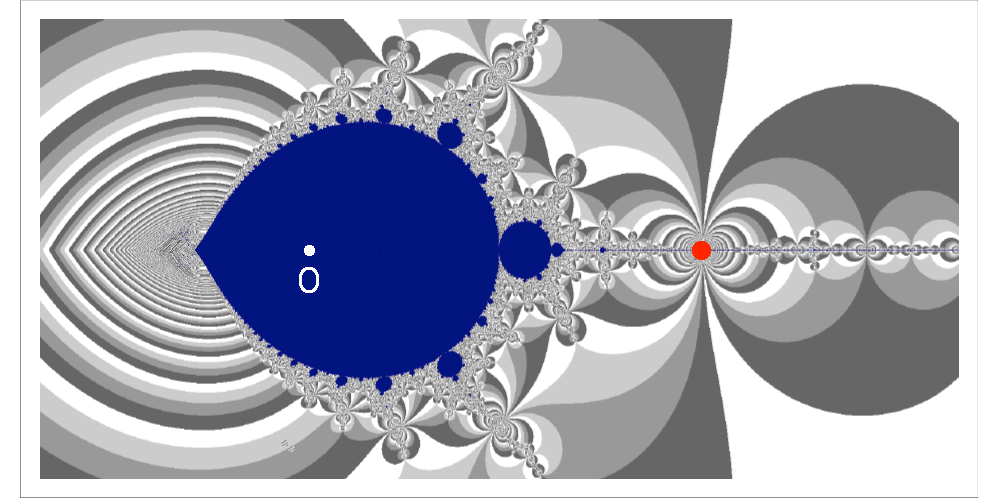}\\ 
\includegraphics[width=0.6\linewidth]{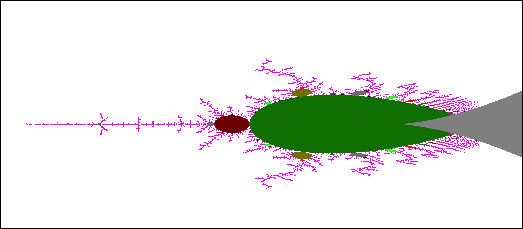}
\caption{Left: The parameter tiles of the exterior of $\cC(\cS)$ are displayed. The brown horizontal line is the slit $(-\infty,-1/12)$. The parameter $a=-1/12$ corresponds to the unique map in $\cC(\cS)$ for which $\alpha_a$ has an attracting direction in $K_a$. Some parameter tiles meet the slit $(-\infty,-1/12)$ along non-degenerate intervals; such intervals are collapsed to the points $\omega$ or $\omega^2$ by the uniformization $\pmb{\Psi}$ of the escape locus. Right: A blow-up of $\cC(\cS)$ around the principal hyperbolic component (having its center at $a=0$) is shown. The marked red parameter corresponds to a Misiurewicz map for which the critical point $0$ lands on the fixed point $\alpha_a$ in three iterates. The boundary points of parameter tiles that are cut-points (respectively, tips) of $\cC(\cS)$ are parameters for which the critical point $0$ eventually maps to $\alpha_a$ (respectively, to $1/4$). Bottom: The region to the left of the grey region (which is a part of the principal hyperbolic component of the Tricorn) is the real Basilica limb of the Tricorn.}
\label{c_and_c_conn_locus_fig}
\end{figure}
More precisely, the angles of the parameter rays of $\cS$ landing/accumulating at a parabolic/Misiurewicz point $a\in\cC(\cS)$ are mapped by $\pmb{\mathcal{E}}_2$ to the angles of the parameter rays of the Tricorn landing/accumulating at the parabolic/Misiurewicz point $\chi(a)\in\mathcal{L}$. This property can be utilized to construct a homeomorphism between suitable pinched disk models of the two connectedness loci.

\begin{theorem}\cite[Theorem~1.4]{LLMM2}
The pinched disk model of $\mathcal{C}(\mathcal{S})$ is homeomorphic to that of the Basilica limb $\mathcal{L}$ of the Tricorn.
\end{theorem}

\subsubsection{A different perspective: David surgery}\label{c_and_c_conn_david_surjery_subsubsec}

One can use David surgery techniques to give an alternative proof of surjectivity of the combinatorial straightening map of Theorem~\ref{c_and_c_general_mating_thm} (which is essentially different from the proof sketched in Subsection~\ref{c_and_c_geom_fin_bijection_subsubsec}). We sketch the key steps below. 
\smallskip

i) Let $f_c\in\mathcal{L}$ be a hyperbolic/Misiurewicz anti-polynomial. One can apply \cite[Lemma~7.1]{LMMN} to replace $f_c\vert_{\mathcal{B}_\infty(f_c)}$ (where $\mathcal{B}_\infty(f_c)$ is the basin of infinity of $f_c$) with the Nielsen map $\pmb{\cN}_2\vert_{\D}$. This produces a piecewise Schwarz reflection map $F$ associated with two disjoint quadrature domains.
\smallskip

ii) The dynamics of the piecewise Schwarz reflection map $F$ on its non-escaping set is topologically conjugate to the dynamics of $f_c$ on its filled Julia set such that the conjugacy is conformal on the interior. Finally, one repeats the arguments of Subsection~\ref{c_and_c_basilica_unique_mating_subsubsec} to deduce that $F$ is a piecewise Schwarz reflection map associated with the disjoint union of two touching simply connected quadrature domains one of which is a round disk and the other is a cardioid. This identifies $F$ as a member of the C\&C family (cf.~\cite[\S A.2]{LLMM2}).

\subsection{Cubic Chebyshev family and associated correspondences}\label{chebyshev_gen_subsec}

The Schwarz reflection map described in Subsection~\ref{chebyshev_subsec} sits in a natural one-parameter family of maps arising from univalent restrictions of the cubic Chebyshev polynomial $f(w)=w^3-3w$ on a varying collection of round disks. The family of $2$:$2$ reversible correspondences associated with these Schwarz reflections can be viewed as an antiholomorphic analog of the family of algebraic correspondences introduced and studied extensively by Bullett, Penrose, and Lomonaco \cite{BP,BuLo1,BuLo2,BuLo3}. A different motivation behind studying this family of quadratic Schwarz reflection maps comes from the discussion in the beginning of Subsection~\ref{quadratic_examples_sec}.

The goal of the current subsection is to give an overview of the dynamics of these Schwarz reflections (respectively, correspondences) and the structure of their parameter space. Our exposition is based on \cite{LLMM3} and its subsequent improvement in \cite{LMM3}.

\subsubsection{The family of Schwarz reflections}\label{cubic_cheby_schwarz_def_subsec}

We set $\HP_R:=\{a\in\C: \re(a)>0\}$. For $a\in\HP_R\setminus\{1\}$, we will denote the disk $B(a,\vert a-1\vert)$ by $\Delta_a$. Note that $\Delta_a$ is centered at $a$ and has the critical point $1$ of $f$ on its boundary. In order to define Schwarz reflection maps via univalent restrictions $f\vert_{\Delta_a}$, we are naturally led to the set
$$
\widehat{S}:=\{a\in\HP_R\setminus\{1\}: f(\partial \Delta_a)\ \mathrm{is\ a\ Jordan\ curve}\}.
$$ 
Elementary arguments show that $\widehat{S}\neq\emptyset$, and $\widehat{S}\subset \{a\in\C: 0<\re(a)\leq4\}$. An explicit description of the set $\widehat{S}$ can be found in \cite[\S 3]{LLMM3}, where it is shown that $\widehat{S}$ is a topological quadrilateral. 

Note that since the critical points of $f$ lie outside $\Delta_a$ for $a\in\HP_R$, univalence of $f\vert_{\overline{\Delta_a}}$ follows from the condition that $f(\partial\Delta_a)$ is a Jordan curve. Hence, $\Omega_a:=f(\Delta_a)$ is a simply connected quadrature domain for $a\in\widehat{S}$. 
Moreover, the presence of the critical point $1$ (of $f$) on $\partial\Delta_a$ implies that $\partial\Omega_a$ has a conformal cusp at $f(1)=-2$ and is non-singular away from $-2$.

We denote the reflection in the circle $\partial\Delta_a$ by $\eta_a$. Then the Schwarz reflection map $\sigma_a$ of $\Omega_a$ is given by $f\circ\eta_a\circ \left(f\vert_{\overline{\Delta}_a}\right)^{-1}$. 
It is easily checked that $\sigma_a$ has two distinct critical points; namely, $c_a:=f(\eta_a(-1))$ and $c_a^*:=f(\eta_a(\infty))$. Furthermore, $\sigma_a$ maps $c_a$ (respectively, $c_a^*$) to $2$ (respectively, to $\infty$) with local degree two (respectively, three). By Proposition~\ref{simp_conn_quad_prop}, $\sigma_a:\sigma_a^{-1}(\Omega_a)\to\Omega_a$ is a two-to-one (possibly branched) covering, and $\sigma_a:\sigma_a^{-1}(\Int{\Omega_a^c})\to \Int{\Omega_a^c}$ is a branched covering of degree three. 

It turns out that for $a\in\widehat{S}$ with $\re(a)\leq 3/2$, the set $\sigma_a^{-1}(\Omega_a)$ has two connected components and $\sigma_a:\sigma_a^{-1}(\Omega_a)\to\Omega_a$ is an unbranched covering map. On the other hand, for $a\in\widehat{S}$ with $\re(a)> 3/2$, we have the following properties:
\begin{itemize}
\item $\Omega_a, \sigma_a^{-1}(\Omega_a)$ are topological disks,
\item $\sigma_a^{-1}(\Omega_a)\subset\Omega_a$ with $\partial\Omega_a\cap\partial\sigma_a^{-1}(\Omega_a)=\{-2\}$, and
\item $\sigma_a:\sigma_a^{-1}(\Omega_a)\to\Omega_a$ is a proper antiholomorphic map of degree two with a unique critical point.
\end{itemize}
(cf. Subsection~\ref{chebyshev_center_hybrid_conj_subsubsec} and \cite[\S 3.2]{LLMM3}.)
Since we are interested in maps $\sigma_a$ exhibiting anti-polynomial-like behavior on its non-escaping set, it is natural to work with the space of Schwarz reflections
$$
\cS:=\{\sigma_a:\overline{\Omega_a}\to\widehat{\C}: a\in\widehat{S},\ \re(a)\in(3/2,4]\}.
$$
Note that $T^0_a\equiv T^0(\sigma_a)=\widehat{\C}\setminus\left(\Omega_a \cup\{-2\}\right)$. We denote the tiling and non-escaping sets of $\sigma_a$ by $T_a^\infty$ and $K_a$, respectively.

We also remark that for all $\sigma_a\in\cS$ with $\re(a)<4$, the cusp $-2$ of $\partial\Omega_a$ is a $(3,2)$-cusp, and $\sigma_a^{\circ 2}$ has a repelling direction in $K_a$ at $-2$. On the other hand, for $\sigma_a\in\cS$ with $\re(a)=4$, the cusp $-2$ of $\partial\Omega_a$ is a $(\nu,2)$-cusp with $\nu\in\{5,7\}$, and $\sigma_a^{\circ 2}$ has at least one attracting direction in $K_a$ at $-2$ (cf. \cite[\S 4.2]{LLMM3}).

\subsubsection{Connectedness locus: coexistence of anti-rational map and reflection group structure}\label{cubic_cheby_qc_straightening_subsubsec}

\begin{figure}[h!]
\captionsetup{width=0.96\linewidth}
\begin{tikzpicture}
\node[anchor=south west,inner sep=0] at (0,0) {\includegraphics[width=0.5\textwidth]{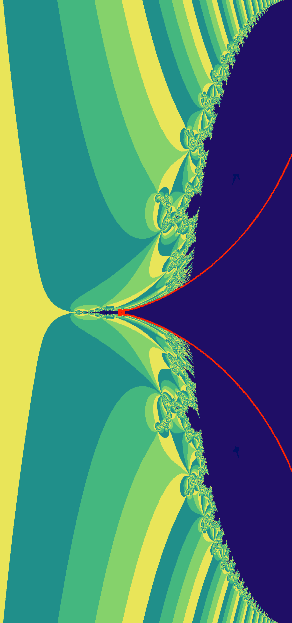}}; 
\node at (5.4,10.8) {\textcolor{white}{Period}};
\node at (5.45,10.4) {\textcolor{white}{two}};
\node at (5.4,6.6) {\textcolor{white}{Period}};
\node at (5.45,6.2) {\textcolor{white}{one}};
\node at (5.4,3.2) {\textcolor{white}{Period}};
\node at (5.45,2.8) {\textcolor{white}{two}};
\draw [->,line width=0.5pt] (-1,4) to (2.28,6.66); 
\node at (-1.2,3.7) {Period};
\node at (-1.15,3.3) {two};
\node at (0,-0.2) {2};
\node at (1.6,-0.2) {2.5};
\node at (3.3,-0.2) {3};
\node at (5,-0.2) {3.5};
\node at (6.3,-0.2) {3.8};
\end{tikzpicture}
\caption{A part of the connectedness locus $\cC(\cS)$ is shown in blue. The period one and two hyperbolic components are labeled. Unlike the Tricorn, the connectedness locus $\cC(\cS)$ does not have a $2\pi/3$-rotation symmetry. This is because the external map $\overline{z}^2$ (of anti-polynomials in the Tricorn) commutes with rotation by $2\pi/3$, while the external map $\pmb{\cF}_2$ (of Schwarz reflections in $\cC(\cS)$) does not have such a symmetry. The exterior of $\cC(\cS)$ is tessellated by parameter tiles. The strips emanating from the top and bottom edges of the depicted parameter rectangle represent tiles that meet the top and bottom sides of the parameter quadrilateral $\widehat{S}$ in intervals of positive length.}
\label{cheby_conn_locus_fig}
\end{figure}
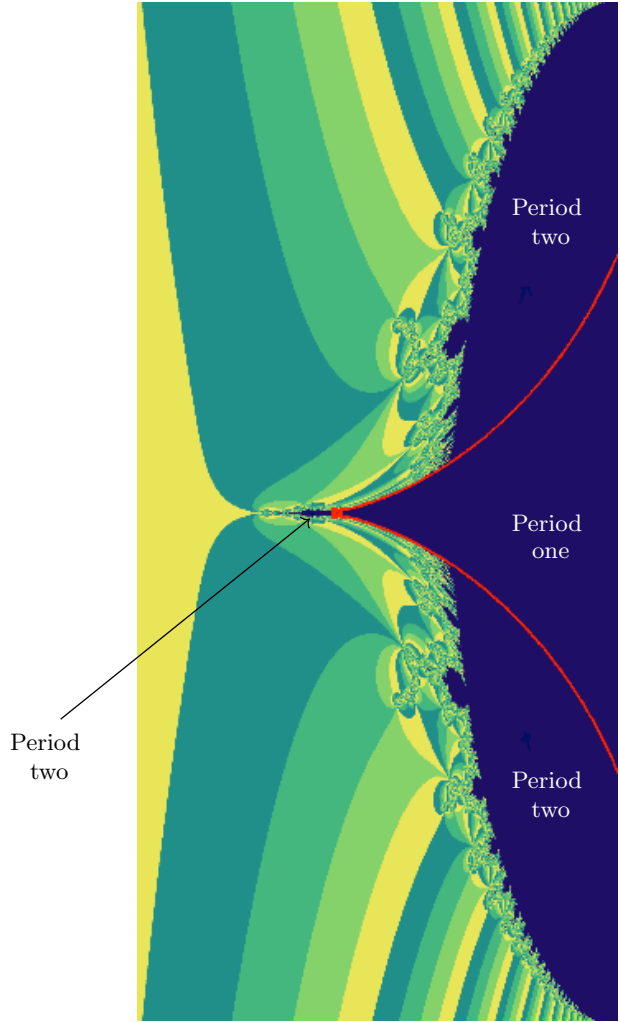

Recall that each $\sigma_a$ has a \emph{passive} critical point at $c_a^*$ that escapes under one iterate of $\sigma_a$. Thus, the map $\sigma_a$ has a unique \emph{active/free} critical point, namely $c_a$. As in the case for unicritical anti-polynomials, it is easy to see that the non-escaping set $K_a$ (of $\sigma_a$) is connected if and only if the unique free critical point $c_a$ does not escape. More precisely, the \emph{connectedness locus} of $\cS$ has the following description 
$$
\cC(\cS)=\{\sigma_a\in\cS: K_a\ \textrm{is\ connected}\}=\{\sigma_a\in\cS: 2\notin T_a^\infty\}
$$
(see Figure~\ref{cheby_conn_locus_fig}, cf. \cite[\S 4.3]{LLMM3}). 

By definition, the free critical value $2$ of $\sigma_a$ never hits the fundamental tile $T^0_a$ for $\sigma_a\in\cC(\cS)$. As a consequence, the conformal isomorphism $\mathcal{Q}_0\to T^0_a$ that carries $0$ to to $\infty$ and sends $1$ to $-2$ can be extended by iterated lifting to produce a conformal conjugacy between $\pmb{\cF}_2:\mathcal{Q}\setminus\Int{\mathcal{Q}_0}\to\mathcal{Q}$ and $\sigma_a:T^\infty_a\setminus\Int{T^0_a}\to T^\infty_a$ (see Subsection~\ref{nielsen_first_return_external_map_subsubsec} for the definitions of $\mathcal{Q}_0$, the map $\pmb{\cF}_2$, and comments on the construction of this conjugacy).

For $\sigma_a\in\cC(\cS)$ with $\re(a)<4$, the restriction $\sigma_a:\sigma_a^{-1}(\Omega_a)\to\Omega_a$ gives rise to a pinched anti-quadratic-like map (in the sense of \cite[Definition~5.1]{LLMM3}) with connected filled Julia set. By \cite[Lemma~5.3, Theorem~5.4]{LLMM3}, such a pinched anti-quadratic-like restriction is hybrid conjugate to a unique anti-rational map in the parabolic Tricorn $\pmb{\mathcal{B}}_2$ (with a simple parabolic fixed point at $\infty$). We refer the reader to Subsection~\ref{chebyshev_center_hybrid_conj_subsubsec} for the main ideas in the proof of this theorem and to \cite[Theorem~4.12]{LMM3} for an alternative method of straightening. 

For $\sigma_a\in\cC(\cS)$ with $\re(a)=4$, the existence of an attracting direction at $-2$ in $K_a$ turns out to be an obstruction to a pinched anti-quadratic-like restriction of $\sigma_a$ (in the sense of \cite[Definition~5.1]{LLMM3}). However, the alternative straightening surgery of \cite[Theorem~4.12]{LMM3} allows one to straighten such a map $\sigma_a$ to a unique anti-rational map in $\pmb{\mathcal{B}}_2$ (with a higher order parabolic fixed point at $\infty$). This surgery construction is carried out in the following two steps.
\smallskip

i) One first shows using a quasiconformal interpolation argument (similar to the one used in \cite[Lemma~5.3]{LLMM3}) that the restriction of $\pmb{\cF}_2$ on the closure of a neighborhood of $\partial\mathcal{Q}\setminus\pmb{\cF}_2^{-1}(1)$ is quasiconformally conjugate to the restriction of the parabolic anti-Blaschke product $B_2(z)=\frac{3\overline{z}^2+1}{\overline{z}^2+3}$ on the closure of a neighborhood of $\mathbb{S}^1\setminus B_2^{-1}(1)$ (see \cite[Lemma~4.10]{LMM3}).

\smallskip

ii) Next, one can use the Riemann map of $T^\infty_a$ (which conjugates $\sigma_a$ to $\pmb{\cF}_2$) and the above quasiconformal conjugacy between $\pmb{\cF}_2$ and $B_2$ to glue the action of $B_2$ outside the non-escaping set of $\sigma_a$. This produces a quasiregular map, which can be straightened to an anti-rational map in $\pmb{\mathcal{B}}_2$ (see \cite[Theorem~4.12]{LMM3} for details).

\begin{theorem}\label{cheby_mating_1_thm}
Each $\sigma_a\in\cC(\cS)$ is a mating of $\pmb{\cF}_2$ and a unique anti-rational map in the parabolic Tricorn $\pmb{\mathcal{B}}_2$.
\end{theorem}

The above theorem defines a \emph{straightening map} $\chi:\cC(\cS)\to\pmb{\mathcal{B}}_2$.

\begin{remark}\label{straightening_regularity_rem}
1) The main difference between the straightening of pinched anti-quadratic-like maps carried out in \cite[\S 5]{LLMM3} and the alternative straightening surgery of \cite[\S 4.4.1]{LMM3} is that in the former surgery, the quasiconformal interpolation is carried out directly in the Schwarz dynamical plane, while in the latter, the interpolation takes place in the dynamical plane of the model map $\pmb{\cF}_2$.

2) The straightening surgery of \cite[\S 4.4.1]{LMM3} applies to all maps in $\cC(\cS)$. However, due to the appearance of the Riemann map of the tiling set in the proof, it is harder to control parameter dependence of this surgery.
\end{remark}

\subsubsection{Properties of the straightening map and homeomorphism between models of parameter spaces}\label{cheby_chi_prop_subsubsec}

The straightening map $\chi$ turns out to be a bijection.
Injectivity of $\chi$ is a consequence of the fact that all maps in $\cC(\cS)$ have the same external dynamics $\pmb{\cF}_2$. Thus, if two maps in $\cC(\cS)$ are hybrid conjugate to the same anti-rational map, then these two Schwarz reflections have the same hybrid class and the same external class; which in turn implies that they are affinely conjugate (see \cite[Proposition~5.9]{LLMM3}).
 
Surjectivity of $\chi$ can be proved using an inverse construction to straightening. While a potentially weaker version of this statement (to the effect that the image of $\chi$ contains the closure of all geometrically finite maps in $\pmb{\mathcal{B}}_2$) was proved in \cite{LLMM3}, the full surjectivity of $\chi$ was established in \cite[\S 5]{LMM3}. As in the construction of $\chi$, one can give two different proofs of surjectivity: one using pinched anti-quadratic-like restrictions of maps in $\pmb{\mathcal{B}}_2$ (with simple parabolic fixed point at $\infty$), and the other using the quasiconformal conjugation between $\pmb{\cF}_2$ and $B_2$ constructed in \cite[Lemma~4.10]{LMM3} (see \cite[Theorems~5.1, 5.2]{LMM3}).

The parameter dependence of the quasiconformal surgeries giving rise to the straightening map and its inverse (that involve pinched anti-quadratic-like maps) were investigated in \cite[\S 6]{LMM3} (cf. \cite[\S 8.1]{LLMM3}), and it was proved that $\chi$ is continuous at hyperbolic and quasiconformally rigid parameters of $\cC(\cS)$. On the other hand, $\chi$ is not necessarily continuous at quasiconformally non-rigid parabolic parameters \cite[\S 8.1]{LLMM3}. 

Using the above properties, it was shown in \cite[\S 9]{LLMM3} that $\chi$ induces a homeomorphism between topological models of $\cC(\cS)$ and $\pmb{\mathcal{B}}_2$. Roughly speaking, these models are constructed by pinching appropriate regions of (possible) discontinuity of $\chi$ to points.

On the other hand, the exterior of $\cC(\cS)$ in the parameter space is simply connected, and it admits a natural uniformization via the conformal position of the escaping critical value (see \cite[\S 7]{LLMM3} and Figure~\ref{cheby_conn_locus_fig}).

\subsubsection{Antiholomorphic analog of Bullett-Penrose correspondences}\label{cheby_corr_gen_mating_subsubsec}

The construction of the $2$:$2$ antiholomorphic correspondence carried out for the map $\sigma_a$ with $a=3$ in Subsection~\ref{chebyshev_center_corr_subsubsec} generalizes verbatim to all maps $\sigma_a\in\cC(\cS)$, with the reflection map $\widehat{\eta}$ in Equation~\eqref{cheby_center_corr_eqn} replaced with the anti-M{\"o}bius reflection $\eta_a$ in the circle $\partial \Delta_a$. Furthermore, the proof of Theorem~\ref{chebyshev_center_corr_thm} also holds for these correspondences $\mathfrak{C}_a$ and shows that each $\mathfrak{C}_a$ is a mating of the group $\mathbbm{G}_2\cong\Z/2\Z\ast\Z/3\Z$ and the map $\chi(\sigma_a)\in\pmb{\mathcal{B}}_2$. Combining this fact with bijectivity of $\chi$, one concludes the following result.

\begin{theorem}\cite{LLMM3} \cite{LMM3}\label{cheby_corr_gen_mating_thm}
\noindent\begin{enumerate}\upshape
\item The straightening map $\chi:\cC(\cS)\to\pmb{\mathcal{B}}_2$ is a bijection.

\item Each antiholomorphic correspondence $\mathfrak{C}_a$ is a mating of the group $\mathbbm{G}_2\cong\Z/2\Z\ast\Z/3\Z$ and the map $\chi(\sigma_a)\in\pmb{\mathcal{B}}_2$.

\item For each $R\in\pmb{\mathcal{B}}_2$, there exists a unique $\sigma_a\in\cC(\cS)$ such that the associated correspondence $\mathfrak{C}_a$ is a mating of $\mathbbm{G}_2$ and $R$. 
\end{enumerate}
\end{theorem}

\subsection{Matings of $\overline{z}^d$ with necklace reflection groups}\label{sigma_d_subsec}

We recall from Subsection~\ref{talbot_subsec} that a simply connected unbounded quadrature domain with a unique node at $\infty$ admits a uniformization $f:\D^*\to\Omega$, where $f(z)=z+\frac{a_1}{z}+\cdots+\frac{a_d}{z^d}$, after possibly replacing the quadrature domain by an affine image of it. By \cite[Propositions~2.7,2.13]{LMM1}, the corresponding quadrature domain $\Omega_f:=f(\D^*)$ has precisely $d+1$ cusps on its boundary (or equivalently, $f$ has $d+1$ critical points on $\mathbb{S}^1$) if and only if $a_d=-1/d$. The dynamics of Schwarz reflection maps associated with two specific quadrature domains of this type were explored in Subsections~~\ref{deltoid_subsec} and~\ref{talbot_subsec} (namely, the deltoid and Talbot Schwarz reflections). In this subsection, we will explicate the general situation by describing the dynamics of Schwarz reflections arising from maps in the family
$$ 
\Sigma_d^* := \left\{ f(z)= z+\frac{a_1}{z} + \cdots -\frac{1}{d z^d} : f\vert_{\D^*} \textrm{ is univalent}\right\}.
$$
The exposition will mostly follow \cite{LMM1,LMM2}.
\begin{figure}[h!]
\captionsetup{width=0.96\linewidth}
\begin{tikzpicture}
\node[anchor=south west,inner sep=0] at (0,0) {\includegraphics[width=0.28\textwidth]{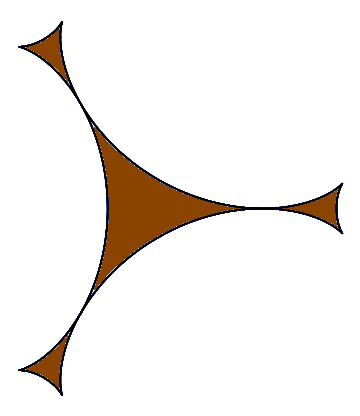}}; 
\node[anchor=south west,inner sep=0] at (6,0) {\includegraphics[width=0.28\textwidth]{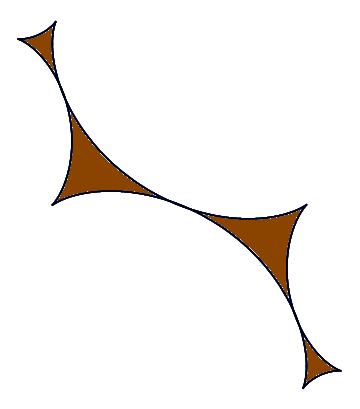}}; 
\end{tikzpicture}
\caption{Depicted are the images of $\mathbb{S}^1$ under two maps $f$ in $\Sigma_5^*$. The curves $f(\mathbb{S}^1)$ have six cusps and three double points. The corresponding desingularized droplets have four components, which are shaded in brown. The exterior white regions are the quadrature domains $\Omega_f=f(\D^*)$.}
\label{sigma_extremal_fig}
\end{figure}

For each $f\in\Sigma_d^*$, we set $\Omega_f:=f(\D^*)$, and denote the Schwarz reflection map of the simply connected quadrature domain $\Omega_f$ by $\sigma_f=f\circ\eta\circ\left(f\vert_{\D^*}\right)^{-1}$. We recall that there are $d+1$ cusps and at most $d-2$ double points on $\partial\Omega_f$ (cf. \cite[Lemma~2.4]{LM1}), and removing these singular points from $T(\sigma_f)=\widehat{\C}\setminus\Omega_f$, we obtain the fundamental tile $T^0(\sigma_f)$. Note also that $\Int{T^0(\sigma_f)}$ has at most $d-1$ connected components (see Figure~\ref{sigma_extremal_fig}).

\subsubsection{Basin of infinity and non-escaping set for Schwarz reflections arising from $\Sigma_d^*$}\label{basin_non_escaping_sigma_d_subsubsec}

By definition, each $\sigma_f$ has a $d$-fold pole at $\infty$; i.e., the point $\infty$ is a superattracting fixed point of $\sigma_f$ with local degree $d$. Since $\deg{f}=d+1$, it has $2d$ critical points in $\widehat{\C}$, of which $d+1$ lie on $\mathbb{S}^1$ and the remaining $d-1$ lie at $0$. Consequently, $f$ has no critical values in $\Omega_f\setminus\{\infty\}$. It follows that $\sigma_f$ has no critical points in $\Omega_f\setminus\{\infty\}$ either. This implies that the basin of attraction $\mathcal{B}_\infty(\sigma_f)$ of the superattracting fixed point $\infty$ is a simply connected, completely invariant domain in $\widehat{\C}$, and $\sigma_f\vert_{\mathcal{B}_\infty(\sigma_f)}$ is conformally conjugate to $\overline{z}^d\vert_{\D}$ \cite[Proposition~3.2]{LMM1} (cf. \cite[\S 9]{Mil06}). We normalize the conformal conjugacy between $\overline{z}^d\vert_{\D}$ and $\sigma_f\vert_{\mathcal{B}_\infty(\sigma_f)}$ such that its derivative at $\infty$ has argument $\frac{\pi}{d+1}$ (see \cite[Remark~2]{LMM2}), and call this normalized conjugacy the \emph{B{\"o}ttcher coordinate} of $\sigma_f$.

According to \cite[Corollary~41]{LMM2}, one has that 
$$
\widehat{\C}=\mathcal{B}_\infty(\sigma_f)\sqcup\Lambda(\sigma_f)\sqcup T^\infty(\sigma_f),
$$ 
where $\Lambda(\sigma_f)=\partial \mathcal{B}_\infty(\sigma_f)=\partial T^\infty(\sigma_f)$ is the \emph{limit set} of $\sigma_f$. The proof of existence of this invariant partition of the dynamical plane of $\sigma_f$ uses a combination of classical arguments of Fatou adapted for the setting of partially defined maps \cite[Proposition~3.4]{LMM1}, local connectivity of the boundary of $\mathcal{B}_\infty(\sigma_f)$ (which essentially follows from expansiveness of $\sigma_f$ near $\partial\mathcal{B}_\infty(\sigma_f)$) \cite[Lemma~32]{LMM2}, and visibility of the cusps of $\partial T^0(\sigma_f)$ (which also lie on $\partial T^\infty(\sigma_f)$) from the basin of infinity \cite[Proposition~35]{LMM2}. It follows that the non-escaping set $K(\sigma_f)$ is the closure of $\mathcal{B}_\infty(\sigma_f)$.

\subsubsection{Mating structure for Schwarz reflections arising from $\Sigma_d^*$}\label{sigma_d_schwarz_dyn_subsubsec}

The above decomposition of $\widehat{\C}$ and local connectivity of $\Lambda(\sigma_f)$ show that $\overline{T^\infty(\sigma_f)}$ is homeomorphic to the quotient of $\overline{\D}$ under an $m_{-d}$-invariant equivalence relation defined by co-landing of dynamical rays in $\mathcal{B}_\infty(\sigma_f)$ (i.e., images of radial lines under B{\"o}ttcher coordinates). As in the special case described in Subsection~\ref{talbot_subsec}, this lamination $\lambda(\sigma_f)$ is generated by the angles of the pairs of $2$-periodic rays (under $m_{-d}$) landing at the double points of $\partial\Omega_f$.

On the other hand, according to \cite[Proposition~27]{LMM2}, there is a unique marked necklace group $G_f$ (up to M{\"o}bius conjugacy) such that $T(\sigma_f)$ is conformally isomorphic to $\overline{\Pi^b(G_f)}$ in a cusp-preserving manner (see Subsection~\ref{necklace_subsubsec} for the definition of $\Pi^b(G_f)$). The uniqueness of $G_f$ is a consequence of standard rigidity theorems for geometrically finite Kleinian groups (cf. \cite[Theorem~4.2]{Tuk85}). The existence of the necklace group $G_f$ is demonstrated in the following two steps.

i) The removal of the cusp points from the droplet boundary $\partial T(\sigma_f)$ yields $d+1$ (open) non-singular real-analytic arcs. These arcs may only touch at the double points of $\partial T(\sigma_f)$. One can extend each of these arcs to a Jordan curve in $\overline{\Omega_f}$ such that no further tangency/intersection among the curves is introduced in the process. This results in a \emph{Jordan curve packing}; i.e., a connected finite collection of oriented Jordan curves in the plane with disjoint interiors. We treat the landing point of the $0$-ray of $\sigma_f$ as a marked cusp on $\partial T(\sigma_f)$, and this yields a marking for this Jordan curve packing (i.e., a labeling of the Jordan curves). The Circle Packing Theorem now produces a (marked) circle packing that is homeomorphic to the Jordan curve packing. This circle packing defines a (marked) necklace group $G'$ such that $T(\sigma_f)$ is quasiconformally isomorphic to $\overline{\Pi^b(G')}$ in a cusp-preserving manner.

\begin{figure}[h!]
\captionsetup{width=0.96\linewidth}
\begin{tikzpicture}
\node[anchor=south west,inner sep=0] at (0,4) {\includegraphics[width=0.48\textwidth]{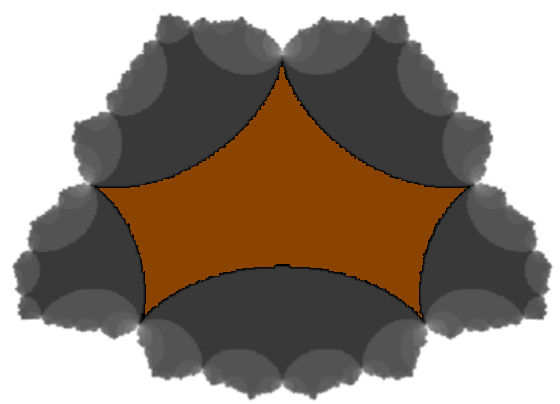}}; 
\node[anchor=south west,inner sep=0] at (6.6,4) {\includegraphics[width=0.44\textwidth]{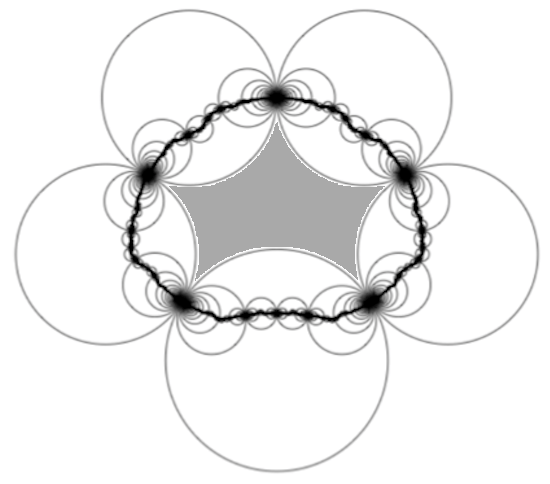}}; 
\node[anchor=south west,inner sep=0] at (0,0) {\includegraphics[width=0.48\textwidth]{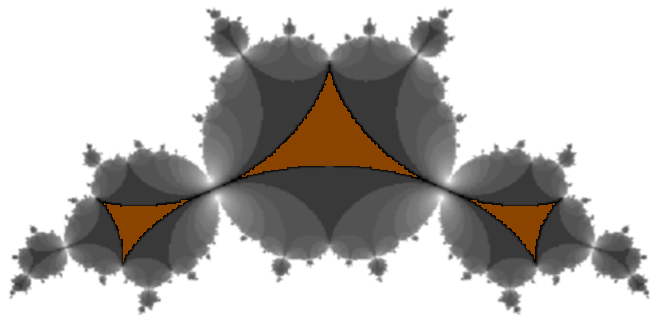}}; 
\node[anchor=south west,inner sep=0] at (6.6,0) {\includegraphics[width=0.44\textwidth]{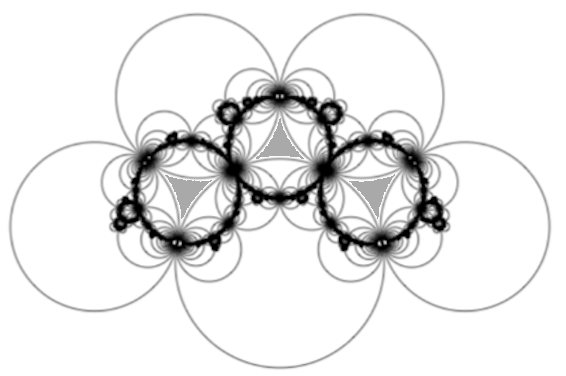}}; 
\end{tikzpicture}
\caption{The dynamical planes of two Schwarz reflection maps arising from $\Sigma_4^*$ (left) and the corresponding necklace groups (right) are displayed.}
\label{sigma_d_mating_fig}
\end{figure}

ii) Next, one quasiconformally deforms the group $G'$ to obtain the required group $G_f$ ensuring that the (marked) components of $\Pi^b(G_f)$ are conformally equivalent to those of $T^0(\sigma_f)$ in a cusp-preserving way.

Moreover, one can use Proposition~\ref{group_lamination_prop} to give an explicit topological model for the limit set $\Lambda(G_f)$ as a quotient of $\mathbb{S}^1$ under an $\pmb{\cN}_d$-invariant lamination $\lambda(G_f)$. Specifically, this lamination is generated by the angles of the pairs of $2$-periodic rays (under $\pmb{\cN}_d$) landing at the double points of $\partial\Pi^b(G_f)$ \cite[\S 4.4]{LMM2}. Alternatively, one can arrive at the same lamination using the perspective of pinching deformation of necklace groups as explained in \cite[\S 3.3, \S 4.3]{LLM1}.

The fact that $T(\sigma_f)$ and $\overline{\Pi^b(G_f)}$ have the same topology translates to a combinatorial equivalence between the laminations $\lambda(\sigma_f)$ and $\lambda(G_f)$. More precisely, the topological conjugacy $\pmb{\mathcal{E}}_d$ between $\pmb{\cN}_d\vert_{\mathbb{S}^1}$ and $\overline{z}^d\vert_{\mathbb{S}^1}$ carries the lamination $\lambda(G_f)$ to $\lambda(\sigma_f)$ \cite[Proposition~56]{LMM2}. Thus, the limit set dynamics of $\sigma_f$ and $\mathcal{N}_{G_f}$ are topologically conjugate. Finally, the conformal equivalence of $T(\sigma_f)$ and $\overline{\Pi^b(G_f)}$ enables one to extend this topological conjugacy to a topological conjugacy between $\sigma_f: \overline{T^\infty(\sigma_f)}\setminus\Int{T^0(\sigma_f)}\to \overline{T^\infty(\sigma_f)}$ and $\mathcal{N}_{G_f}: K(G_f)\setminus\Int{\Pi^b(G_f)}\to K(G_f)$ in such a way that the conjugacy is conformal on the interior.

\begin{theorem}\cite[Theorem~A]{LMM2}\label{sigma_d_mating_thm}
For each $f\in\Sigma_d^*$, the Schwarz reflection $\sigma_f$ is a conformal mating of $\overline{z}^d$ and the Nielsen map of a unique marked necklace group $G_f$. More precisely, 
\noindent\begin{enumerate}\upshape
\item the maps 
$$
\sigma_f:\overline{T^\infty(\sigma_f)}\setminus\Int{T^0(\sigma_f)}\to \overline{T^\infty(\sigma_f)}\ \mathrm{and}\ \mathcal{N}_{G_f}: K(G_f)\setminus\Int{\Pi^b(G_f)}\to K(G_f)
$$
are topologically conjugate in such a way that the conjugacy is conformal on $\Int{T^\infty(\sigma_f)}$, and

\item $\sigma_f\vert_{K(\sigma_f)}$ is topologically conjugate to $\overline{z}^d\vert_{\overline{\D}}$ such that the conjugacy is conformal on $\Int{K(\sigma_f)}$. Moreover, the Minkowski circle homeomorphism $\pmb{\mathcal{E}}_d$ maps the lamination associated with the limit set of $G_f$ to the lamination associated with the limit set of $\sigma_f$.
\end{enumerate}
\end{theorem}

\subsubsection{Homeomorphism between $\Sigma_d^*$ and the Bers slice closure $\overline{\beta(\pmb{G}_{d})}$}\label{sigma_d_bers_homeo_subsubsec}

The map $f\mapsto G_f$ connects two disparate objects: a space of univalent maps and a space of Kleinian reflection groups. Interestingly, this map turns out to be a homeomorphism between $\Sigma_d^*$ and the Bers slice closure of $\pmb{G}_{d}$ (see Definition~\ref{deform_space_def}).

We first sketch a proof of surjectivity, whose main idea is similar to the one described in Subsection~\ref{talbot_pinching_subsubsec}. One starts with the Schwarz reflection associated with the Jordan quadrature domain $\Omega_0:=f_0(\D^*)$, where $f_0(z)=z-\frac{1}{dz^d}$, and quasiconformally deforms the droplet to introduce additional tangencies on the quadrature domain boundary in the limit. We remark that the quadrature domain $\Omega_0$ is the exterior of a classical \emph{hypocycloid} (or a \emph{($d+1$)-deltoid}) curve. For any given necklace group $G$, this pinching procedure can be used to create a map in $f\in\Sigma_d^*$ such that the associated droplet $T(\sigma_f)$ has the topology of $\overline{\Pi^b(G)}$ (cf. \cite[Theorem~4.14]{LMM1}, \cite[Proposition~29]{LMM2}).
 
Alternatively, one can construct such an $f\in\Sigma_d^*$ using David surgery. To accomplish this, one first invokes the bijection between necklace reflection groups and critically fixed anti-polynomials (which will be discussed in Subsection~\ref{dyn_conseq_bijection_subsubsec}) to obtain a critically fixed anti-polynomial $p$ whose Julia set dynamics is topologically conjugate to $\cN_G\vert_{\Lambda(G)}$, where $G$ is a given necklace group. Subsequently, one glues in Nielsen maps of appropriate ideal polygon reflection groups in the bounded fixed Fatou components of $p$, where the choices of the ideal polygons are dictated by the components of $\Pi^b(G)$. This produces a Schwarz reflection map in a simply connected quadrature domain, which is easily seen to be uniformized by some $f\in\Sigma_d^*$. An added advantage of this approach is that the limit set of the Schwarz reflection thus produced is conformally removable: indeed, the limit set of such a Schwarz reflection is the image of a $W^{1,1}$-removable compact set (namely, the connected Julia set of a hyperbolic polynomial) under a global David homeomorphism (see \cite[Theorem~2.7]{LMMN}). We refer the reader to \cite[\S 12]{LMMN} for details of this construction. 

Recall that all Schwarz reflections arising from $\Sigma_d^*$ have conformally equivalent dynamics (conjugate to $\overline{z}^d$) on their basin of attraction of $\infty$. Additionally, if two of them have conformally isomorphic droplets, then their dynamics on the respective tiling set closures are also conformally equivalent. Injectivity of the map $f\mapsto G_f$ now follows from conformal removability of the limit sets of Schwarz reflections constructed above (cf. \cite[Theorem~12.8]{LMMN}). An alternative argument for injectivity (that uses the pullback argument) can be found in \cite[Theorem~5.1]{LMM1}.

Finally, continuity of the map $f\mapsto G_f$ is deduced from a continuity property of moduli of quadrilaterals \cite[Theorems~25, 31]{LMM2}.

\begin{theorem}\cite[Theorems~A, C]{LMM2}\label{sigma_d_bers_homeo_thm}
\noindent\begin{enumerate}\upshape
\item The map $f\mapsto G_f$ of Subsection~\ref{sigma_d_schwarz_dyn_subsubsec} is a homeomorphism between $\Sigma_d^*$ and the Bers slice closure $\overline{\beta(\pmb{G}_{d})}$.
\item For each $f\in\Sigma_d^*$, there exists a necklace group $G_f$ and a critically fixed anti-polynomial $p_f$ such that $\sigma_f\vert_{\Lambda(\sigma_f)}$, $\cN_{G_f}\vert_{\Lambda(G_f)}$, and $p_f\vert_{\mathcal{J}(p_f)}$ are topologically conjugate.
\end{enumerate}
\end{theorem}

\subsubsection{Cell complex structure of $\Sigma_d^*$}\label{sigma_d_cell_structure_subsec}

The space $\Sigma_d^*$ of univalent rational maps admits a natural cell complex structure. For $0\leq k\leq d-2$, let us denote by $\Sigma_{d,k}^*$ the collection of maps $f\in\Sigma_d^*$ such that the boundary $\partial T(\sigma_f)$ has exactly $k$ double points. Then, the components of $\Sigma_{d,k}^*$ are the $(d-2-k)$-dimensional cells of $\Sigma_d^*$, and each such cell is homeomorphic to a quasiconformal deformation space of Schwarz reflections. There is a unique cell of dimension $(d-2)$ in $\Sigma_d^*$, which contains the map $f_0(z)=z-\frac{1}{dz^d}$. We refer to the $0$-dimensional cells as \emph{vertices} of $\Sigma_d^*$. A vertex $f$ of $\Sigma_d^*$ is also called a \emph{Suffridge map} (cf. \cite{LM1,Suf72}); the corresponding desingularized droplet comprises $(d-1)$ triangles. 

The Bers slice closure $\overline{\beta(\pmb{G}_{d})}$ also has a similar cell complex structure, where the cells of various dimensions are defined by the number (of conjugacy classes) of accidental parabolics of the necklace groups (see Proposition~\ref{cell_structure_group_prop}). By construction, the homeomorphism of Theorem~\ref{sigma_d_bers_homeo_thm} respects these cell complex structures.

\section{David surgery}\label{david_surgery_sec}

The fact that quasisymmetric circle homeomorphisms admit quasiconformal extensions to $\D$ plays an important role in the theory of Kleinian groups as well as in rational dynamics; such as in the Bers Simultaneous Uniformization Theorem, in the construction of quasi-Blaschke products, and in the Douady-Ghys surgery to construct maps having Siegel disks with controlled geometry. However, as we saw in various examples described in the previous sections (see Subsections~\ref{deltoid_david_subsubsec},~\ref{talbot_david_subsubsec},~\ref{rat_map_apollo_group_subsubsec},~\ref{c_and_c_conn_david_surjery_subsubsec},~\ref{sigma_d_bers_homeo_subsubsec}), there is a number of mating/surgery frameworks where one needs similar extension theorems for naturally arising non-quasisymmetric circle homeomorphisms (such as the Minkowski circle homeomorphisms and circle homeomorphisms conjugating expanding Blaschke products to parabolic ones).
In this section, we will formulate such an extension theorem for topological conjugacies between piecewise real-analytic, expansive, covering maps of $\mathbb{S}^1$ satisfying certain regularity conditions such that the conjugacy carries parabolics to parabolics, but can send hyperbolics to parabolics as well. While the most general version of this result that was proved in \cite{LMMN} does not require the circle coverings to be $C^1$, we will only state a special case assuming $C^1$-regularity as this will suffice for all the applications discussed in this article.

\subsection{A David extension theorem for circle homeomorphisms}\label{david_ext_subsec}
\begin{definition}\label{expansive_def}
1) A continuous map $f\colon \mathbb{S}^1\to \mathbb{S}^1$ is called \emph{expansive} if there exists a constant $\delta>0$ such that for any $a,b\in \mathbb{S}^1$ with $a\neq b$ we have $|f^{\circ n}(a)-f^{\circ n}(b)|>\delta$ for some $n\in \N$.

2) We say that a periodic point $a$ of an expansive $C^1$ map $f\colon \mathbb{S}^1\to \mathbb{S}^1$ is \emph{parabolic} (respectively, \emph{hyperbolic}) if the derivative of the first orientation-preserving return map of $f$ to $a$ has absolute value equal to (respectively, larger than) $1$. 
\end{definition}

Let $f\colon \mathbb{S}^1\to \mathbb{S}^1$ be a $C^1$, expansive, covering map of degree $d\geq 2$ admitting a Markov partition $\mathcal P(f;\{a_0,\dots,a_r\})$.
We define $A_k=\arc{[a_k,a_{k+1}]}$ for $k\in \{0,\dots,r\}$ (with the convention $a_{r+1}=a_0$), and recall that $f_k:=f|_{\Int{A_k}}$ is injective by the definition of a Markov partition. We assume, that $f_k$ is analytic and that there exist open neighborhoods $U_k$ of $\Int{A_k}$ and $V_k$ of $f_k(\Int{A_k})$ in the plane such that $f_k$ has a conformal extension from $U_k$ onto $V_k$. We still denote the extension by $f_k$. We impose the condition that 
\begin{align}\label{condition:uv}
V_k=f_k(U_k)\supset U_j \quad \textrm{whenever}\quad  f_k(A_k)\supset A_j,
\end{align} 
for $j,k\in \{0,\dots,r\}$. We also require that 
\begin{align}\label{condition:holomorphic}
\textrm{$f_k$ extends holomorphically to neighborhoods of $a_k$ and $a_{k+1}$}
\end{align}
for each $k\in \{0,\dots,r\}$.

Now suppose that $a\in \{a_0,\dots,a_r\}$ is a parabolic periodic point of $f$. 
By condition~\eqref{condition:holomorphic}, the left and right branches of the first orientation-preserving return map of $f$ to $a$ have holomorphic extensions valid in complex neighborhoods of $a$. These extensions define a pair of parabolic germs fixing $a$. We say that $a$ is \emph{symmetrically parabolic} if these two parabolic germs have the same parabolic multiplicity (see \cite[\S 10]{Mil06} for background on parabolic germs).

\begin{remark}\label{symm_parabolic_rem}
If $f$ is \emph{orientation-reversing} and $a\in \{a_0,\dots,a_r\}$ is a parabolic fixed point, then it is automatically symmetrically parabolic. This is because the map $f$ itself defines a topological conjugacy between the corresponding left and right branches of the first orientation-preserving return maps (see \cite[Remark~4.7]{LMMN} for a similar assertion in a more general setup).
\end{remark}

With these preliminary definitions at our disposal, let us now turn to the main extension theorem of this section.

\begin{theorem}\cite[Theorem~4.9 (special case)]{LMMN}\label{david_extension_general_thm}
Let $f,g\colon \mathbb{S}^1\to \mathbb{S}^1$ be $C^1$, expansive, covering maps with the same degree and orientation, and $\mathcal P(f;\{a_0,\dots,a_r\})$, $\mathcal P(g;\{b_0,\dots,b_r\})$ be Markov partitions satisfying conditions \eqref{condition:uv} and \eqref{condition:holomorphic}.  Suppose that the map $h\colon \{a_0,\dots,a_r\} \to \{b_0,\dots,b_r\}$ defined by $h(a_k)=b_k$, $k\in \{0,\dots,r\}$, conjugates $f$ to $g$ on the set $\{a_0,\dots,a_r\}$ and assume that for each periodic point $a\in \{a_0,\dots,a_r\}$ of $f$ and for $b=h(a)$ one of the following alternatives occur.
\smallskip

\begin{enumerate}\upshape
\item[\scriptsize(\textbf{H${\rightarrow}$H})]\label{HH} Both $a,b$ are hyperbolic.

\item[\scriptsize(\textbf{P${\rightarrow}$P})]\label{PP} Both $a,b$ are symmetrically parabolic

\item[\scriptsize(\textbf{H$\to$P})]\label{HP} $a$ is hyperbolic and $b$ is symmetrically parabolic.
\end{enumerate}

\noindent
Then the map $h$ extends to a homeomorphism $\widetilde h$ of $\overline{\D}$ such that $\widetilde h|_{\mathbb S^1}$ conjugates $f$ to $g$ and $\widetilde h|_{\D}$ is a David homeomorphism. Moreover, if the alternative \begin{scriptsize}$\mathbf{(H\to P)}$\end{scriptsize} does not occur, then $\widetilde h|_{\D}$ is a quasiconformal map and $\widetilde h|_{\mathbb S^1}$ is a quasisymmetry.
\end{theorem}

The facts that 
\begin{itemize}
\item the map $f(z)=z^d$ or $\overline{z}^d$ satisfies conditions \eqref{condition:uv} and \eqref{condition:holomorphic} for every Markov partition of $\mathbb{S}^1$, and
\item each periodic point of $f(z)=z^d$ or $\overline{z}^d$ is hyperbolic
\end{itemize}
immediately yield the following corollary which is often useful in practice.

\begin{corollary}\cite[Theorem~4.12]{LMMN}\label{power_map_cor}
Let $g\colon \mathbb{S}^1\to \mathbb{S}^1$ be a $C^1$, expansive, covering map of degree $d\geq 2$ and let $\mathcal P(g;\{b_0,\dots,b_r\})$ be a Markov partition satisfying conditions \eqref{condition:uv} and \eqref{condition:holomorphic}, and with the property that $b_k$ is either hyperbolic or symmetrically parabolic for each $k\in \{0,\dots,r\}$. Then there exists an orientation-preserving homeomorphism $h\colon \mathbb{S}^1\to \mathbb{S}^1$ that conjugates the map $z\mapsto z^d$ or $z\mapsto \overline{z}^d$ to $g$ and has a David extension in~$\D$.
\end{corollary}

The proof of Theorem~\ref{david_extension_general_thm} involves careful distortion estimates for the circle homeomorphisms in question using local normal forms for hyperbolic and parabolic periodic points. Specifically, one verifies that the scalewise distortion function of the topological conjugacy $h$ grows at most as $\log(1/t)$ as the scale $t$ goes to $0$ and then applies the David extension criterion of Chen--Chen--He and Zakeri (see Theorem~\ref{david_extension_criterion_thm} and the preceding discussion).

We now illustrate the above extension theorems with two important examples.

\begin{example}[Minkowski circle homeomorphism.]\upshape\label{example_1}
Recall that the inverse $\pmb{\mathcal{E}}_d^{-1}$ of the $d$-th Minkowski circle homeomorphism  conjugates $\overline{z}^d\vert_{\mathbb{S}^1}$ to the Nielsen map $\pmb{\cN}_d\vert_{\mathbb{S}^1}$. It is trivial to check from the explicit formula of $\pmb{\cN}_d$ that the Markov partition of $\pmb{\cN}_d\vert_{\mathbb{S}^1}$ defined by the $(d+1)$-st roots of unity satisfies conditions~\eqref{condition:uv},~\eqref{condition:holomorphic} as well as all the conditions of Corollary~\ref{power_map_cor} (see \cite[Example~4.3]{LMMN} for details). Hence, $\pmb{\mathcal{E}}_d^{-1}$ admits a David extension to $\D$. 
 
More generally, since any degree $d$ anti-Blaschke product with an attracting fixed point in $\D$ is quasisymmetrically conjugate to $\overline{z}^d$ on $\mathbb{S}^1$ and the Nielsen map of any polygonal reflection group generated by $d+1$ reflections is quasisymmetrically conjugate to $\pmb{\cN}_d$ on $\mathbb{S}^1$, one concludes (for instance, using the Ahlfors-Beurling Extension Theorem) that a circle homeomorphism topologically conjugating a degree $d$ anti-Blaschke product with an attracting fixed point in $\D$ to the Nielsen map of a polygonal reflection group generated by $d+1$ reflections admits a David extension to $\D$ (see \cite[Theorem 4.13]{LMMN} for details).
\end{example}

\begin{example}[Circle homeomorphisms conjugating expanding Blaschke products to parabolic ones.]\upshape\label{example_2} 
The Blaschke product $B_d(z)=\frac{(d+1)z^d+(d-1)}{(d-1)z^d+(d+1)}$ has a parabolic fixed point at $1$ and is an expansive covering of degree $d$. The Markov partition $\mathcal P(B_d;\{b_0, \dots, b_{2d-1}\})$, where $b_0, \dots, b_{2d-1}$ are $d$-th roots of $1$ and $-1$ with $b_0=1$, satisfies both conditions~\eqref{condition:uv} and \eqref{condition:holomorphic}. Here, in condition~\eqref{condition:uv}, the set $U_k$, $k\in \{0,1,\dots 2d-1\}$, is an open sector with vertex 0 and angle $\pi/d$, whose boundary contains the points $b_k$ and $b_{k+1}$, indices taken modulo $2d$; the sets $V_k$ are the upper half-plane for even $k$ and the lower half-plane for odd $k$. Moreover, as $B_d$ is a global holomorphic map of $\widehat{\C}$, each $b_k$ is either hyperbolic or symmetrically parabolic. Hence, Corollary~\ref{power_map_cor} guarantees the existence of an orientation-preserving homeomorphism $h\colon \mathbb{S}^1\to \mathbb{S}^1$ that conjugates the map $z\mapsto z^d$ to $B_d$ and has a David extension in $\D$.

More generally, since any degree $d$ Blaschke product with an attracting fixed point in $\D$ is quasisymmetrically conjugate to $z^d$ on $\mathbb{S}^1$ and any degree $d$ Blaschke product with a double parabolic fixed point on $\mathbb{S}^1$ is quasisymmetrically conjugate to $B_d$ on $\mathbb{S}^1$ (cf. \cite[Theorem~6.1, Proposition~6.8]{McM98}), one concludes that a circle homeomorphism topologically conjugating a degree $d$ Blaschke product with an attracting fixed point in $\D$ to a degree $d$ Blaschke product with a double parabolic fixed point on $\mathbb{S}^1$ admits a David extension to $\D$.
\end{example}

\subsection{David surgery to pass from hyperbolic to parabolic dynamics}\label{david_surgery_hyp_para_subsec}

The David Extension Theorem~\ref{david_extension_general_thm}, combined with the David Integrability Theorem~\ref{david_integrability_thm} gives a unified approach to turn hyperbolic anti-rational maps to parabolic anti-rational maps, Kleinian reflection groups, and Schwarz reflection maps that are matings of anti-polynomials and Nielsen maps of kissing reflection groups. This is done by replacing the attracting dynamics of an anti-rational map on suitable invariant Fatou components by Nielsen maps of ideal polygon groups or parabolic anti-Blaschke products. 

Recall that a rational map is called \emph{subhyperbolic} if every critical orbit is either finite or converges to an attracting periodic orbit.

\begin{lemma}\cite[Lemma~7.1]{LMMN}\label{david_surgery_lemma}
Let $R$ be a subhyperbolic anti-rational map with connected Julia set, and let $U_i,\ i=1,\cdots, k$, be  invariant Fatou components of $R$ such that $R\vert_{\partial U_i}$ have degree $d_i$. Then, there is a global David surgery that replaces the dynamics of $R$ on each $U_i$ by the dynamics of $\pmb{\cN}_{d_i}\vert_{\D}$ (respectively, $\overline{B_d(z)}\vert_{\D}$), transferred to $U_i$ via a Riemann map. More precisely, there exists a global David homeomorphism $\Psi$, and an antiholomorphic map $F$, defined on a subset of $\widehat{\C}$, such that $F\vert_{\Psi(U_i)}$ is conformally conjugate to $\pmb{\cN}_{d_i}\vert_{\D}$ (respectively, $\overline{B_d(z)}\vert_{\D}$), and $F$ is conformally conjugate to $R$ outside the grand orbit of $\displaystyle\bigcup_{i=1}^k \Psi(U_i)$.
\end{lemma}

The proof of this result can be split into the following main steps.
\smallskip

\noindent $\bullet$ As $R$ is subhyperbolic, $R\vert_{U_i}$ is conformally conjugate to the action of an anti-Blaschke product with an attracting fixed point in $\D$, and hence $R$ induces the action of an expanding anti-Blaschke product on the ideal boundary $\mathbb{S}^1$ of $U_i$. As a topological conjugacy $h_i$ between this anti-Blaschke product and $\pmb{\cN}_{d_i}\vert_{\mathbb{S}^1}$ (or $B_d\vert_{\mathbb{S}^1}$) admits a David extension to $\D$ (also denoted by $h_i$), one can glue the dynamics of $\pmb{\cN}_{d_i}\vert_{\D}$ (or $B_d\vert_{\D}$) in $U_i$ via the composition of $h_i$ with a Riemann uniformization $\phi_i:U_i\to\D$. Recall that the existence of such David extensions follow from the discussion in Examples~\ref{example_1} and~\ref{example_2}. This produces an orientation-reversing map $\widetilde{R}$ on a subset of $\widehat{\C}$ with the desired topological dynamics.
\smallskip

\noindent $\bullet$ To straighten the above orientation-reversing map to an antiholomorphic map, one needs to apply to the David Integrability Theorem, which necessitates the construction of an $\widetilde{R}$-invariant David coefficient on $\widehat{\C}$. This is done by first pulling back the standard complex structure on $\D$ under $h_i\circ\phi_i:U_i\to\D$, and then spreading it throughout the grand orbit of $U_i$ by iterates of $R$. Outside the grand orbits of the various $U_i$, one uses the standard complex structure. That this Beltrami coefficient $\mu$ satisfies the David condition on each $U_i$ follows from the John property of the Fatou components $U_i$, which is a consequence of subhyperbolicity of $R$ and connectedness of $\mathcal{J}(R)$ (see \cite{Mih11} and Appendix~\ref{integrable_thms_subsec}). To show that $\mu$ satisfies the David condition on all of $\widehat{\C}$, one again appeals to subhyperbolicity of $R$. Specifically, subhyperbolicity of $R$ implies that there is a neighborhood of the closures of all but finitely many Fatou components that is disjoint from the postcritical set of $R$ and that
the Fatou components of $R$ are uniform John domains \cite{Mih11}. These facts allow one to employ the Koebe Distortion Theorem to control the area distortion under the inverse branches of $R$, which yields the global David property of~$\mu$.
\smallskip

\noindent $\bullet$ By the David Integrability Theorem, there exists a global David homeomorphism $\Psi$ that solves the Beltrami equation with coefficient $\mu$. The final step of the proof is to demonstrate that $\Psi\circ\widetilde{R}\circ\Psi^{-1}$ is antiholomorphic. This follows essentially from the local uniqueness of David homeomorphisms integrating a David coefficient (see Theorem~\ref{stoilow_thm}) and $W^{1,1}$-removability of $\displaystyle\bigcup_{i=1}^k \partial U_i$. 
\smallskip

In the next two sections, we will present various applications of Lemma~\ref{david_surgery_lemma} which convert hyperbolic anti-rational maps to kissing reflection groups and Schwarz reflection maps exhibiting hybrid dynamics.

\section{Kissing reflection groups vs critically fixed anti-rational maps}\label{new_line_dict_sec}

In Subsections~\ref{apollo_group_map_schwarz_subsec} and~\ref{sigma_d_subsec}, we discussed the existence of dynamically meaningful homeomorphisms between limit sets of certain kissing reflection groups and Julia sets of certain critically fixed anti-rational maps (see Theorem~\ref{apollo_schwarz_group_anti_rat_thm} for the Apollonian gasket reflection group and Theorem~\ref{sigma_d_bers_homeo_thm} for necklace reflection groups). These examples were generalized in \cite{LLM1}, where a precise dynamical relation between kissing reflection groups and critically fixed anti-rational maps was established. Subsequently, parameter space ramifications of this relation were studied in \cite{LLM2}. It transpired that there are stark similarities between the global topological properties of the deformation spaces of kissing reflection groups and critically fixed anti-rational maps.
The goal of the current section is to expound these general results.

\subsection{A bijection between the two classes of conformal dynamical systems}\label{new_line_dict_dyn_subsec}

\subsubsection{Construction of the bijection: Minkowski surgery}\label{group_map_bijection_subsubsec}

Recall from Subsection~\ref{kissing_group_subsubsec} that each connected simple plane graph $\Gamma$ gives rise to a kissing reflection group, whose limit set is connected if and only if $\Gamma$ is $2$-connected. It is a straightforward consequence of rigidity of geometrically finite Kleinian groups that if the circle packings $\mathcal{P}_1, \mathcal{P}_2$ defining two kissing reflection groups $G_{\mathcal{P}_1}, G_{\mathcal{P}_2}$ have isomorphic contact graphs (as plane graphs), then the groups $G_{\mathcal{P}_1}, G_{\mathcal{P}_2}$ are quasiconformally conjugate. Such a quasiconformal conjugacy only deforms the conformal classes of the polygons in the fundamental domains $\Pi(G_{\mathcal{P}_1}), \Pi(G_{\mathcal{P}_2})$ (see Subsection~\ref{fund_dom_subsubsec}). Thus, the space of planar isomorphism classes of $2$-connected, simple, plane graphs with $d+1$ vertices is in bijective correspondence with the space of quasiconformal conjugacy classes of kissing reflection groups of rank $d+1$ with connected limit set.

It was observed in \cite{LLM1} that the surgery construction of Subsection~\ref{apollo_group_rat_map_subsubsec}, which turns the Apollonian gasket reflection group to a critically fixed cubic anti-rational map, can be applied to general kissing reflection groups of rank $d+$~$1$ producing  critically fixed anti-rational maps of degree $d$. Roughly speaking, this surgery procedure modifies the Nielsen map of a kissing reflection group on appropriate components of the domain of discontinuity by gluing in the dynamics of suitable power maps.
 Specifically, given a kissing reflection group $G_{\mathcal{P}}$ arising from a circle packing $\mathcal{P}$, one looks at the so-called \emph{principal components} of $\Omega(G_{\mathcal{P}})$; i.e., the components $\mathcal{U}_1,\cdots,\mathcal{U}_k$ of $\Omega(G_{\mathcal{P}})$ that intersect the fundamental domain $\Pi(G_{\mathcal{P}})$ non-trivially. If the components of $\Pi(G_{\mathcal{P}})$ are $(r_i+1)$-gons (where $i\in\{1,\cdots,k\}$), then the restriction of the Nielsen map $\cN_{G_{\mathcal{P}}}$ to the principal component $\mathcal{U}_i$ is quasiconformally conjugate to the Nielsen map $\pmb{\cN}_{r_i}$ of the regular $(r_i+1)$-gon reflection group. The Minkowski circle homeomorphism $\pmb{\mathcal{E}}_{r_i}$, which conjugates $\pmb{\cN}_{r_i}\vert_{\mathbb{S}^1}$ to $\overline{z}^{r_i}\vert_{\mathbb{S}^1}$, allows one to construct a degree $d$ orientation-reversing branched cover of $\mathbb{S}^2$ by replacing the action of $\cN_{G_{\mathcal{P}}}$ on $\mathcal{U}_i$ with the power map $\overline{z}^{r_i}\vert_{\D}$ (see Subsection~\ref{apollo_group_rat_map_subsubsec}). Note that all the critical points of this branched covering come from the introduction of power maps and hence they are all fixed. To justify the absence of Thurston obstructions for this critically fixed branched cover, one can use a result of Pilgrim and Tan \cite{PT} that reduces, under these circumstances, the collection of possible Thurston obstructions to a small checkable list. Moreover, the expansiveness of this branched cover on $\Lambda(G_{\mathcal{P}})$ (which comes from expansivity of circular reflections) implies that it is topologically conjugate to some critically fixed anti-rational map $R_\Gamma$, where $\Gamma$ is the contact graph of the circle packing $\mathcal{P}$. (A more explicit relation between $\Gamma$ and $R_\Gamma$ will be articulated below.) In particular, the Nielsen map $\cN_{G_{\mathcal{P}}}\vert_{\Lambda(G_{\mathcal{P}})}$ is topologically conjugate to $R_\Gamma\vert_{\mathcal{J}(R_\Gamma)}$.
\begin{figure}[h!]
\captionsetup{width=0.96\linewidth}
\includegraphics[width=0.4\linewidth]{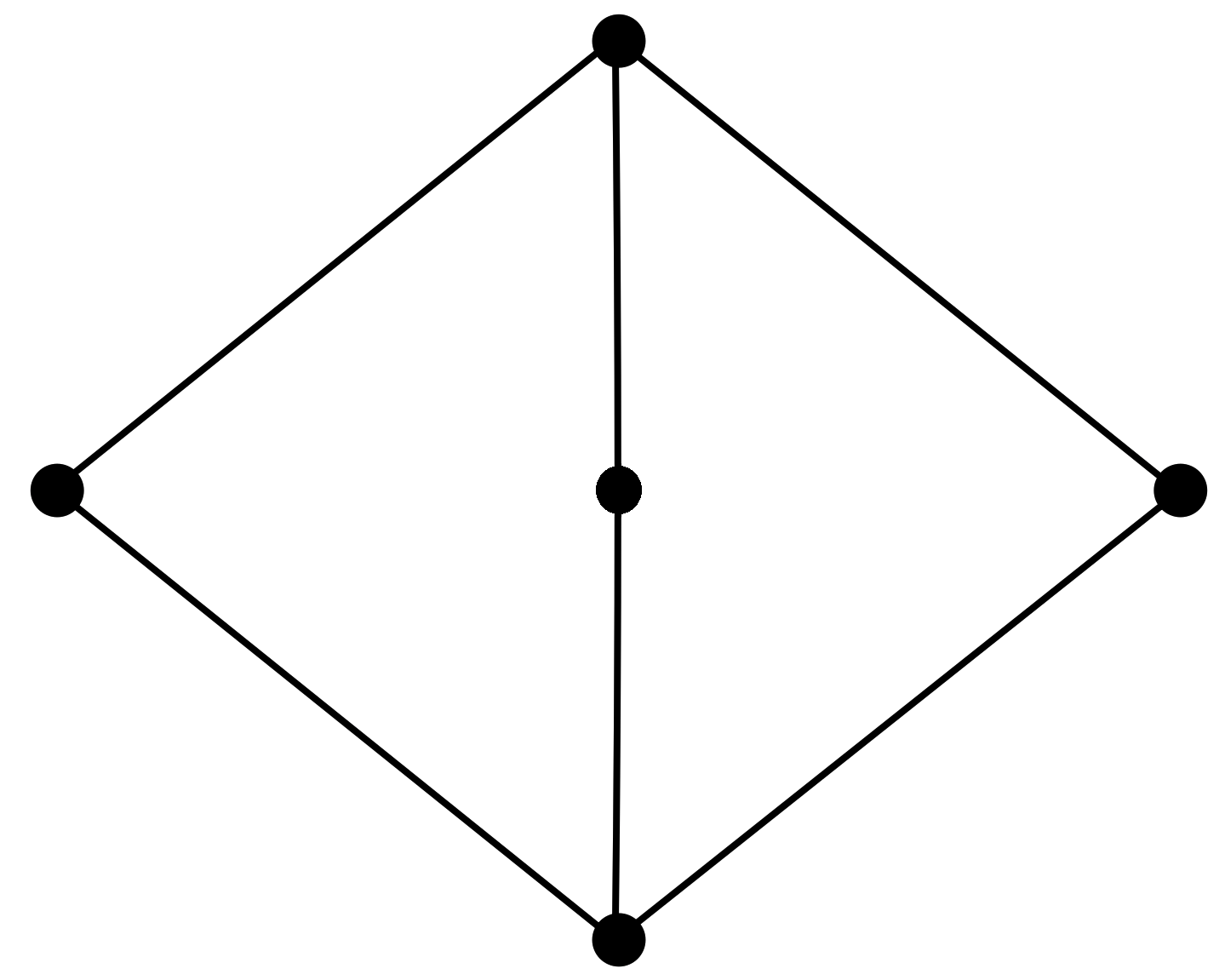}\\
\includegraphics[width=0.53\linewidth]{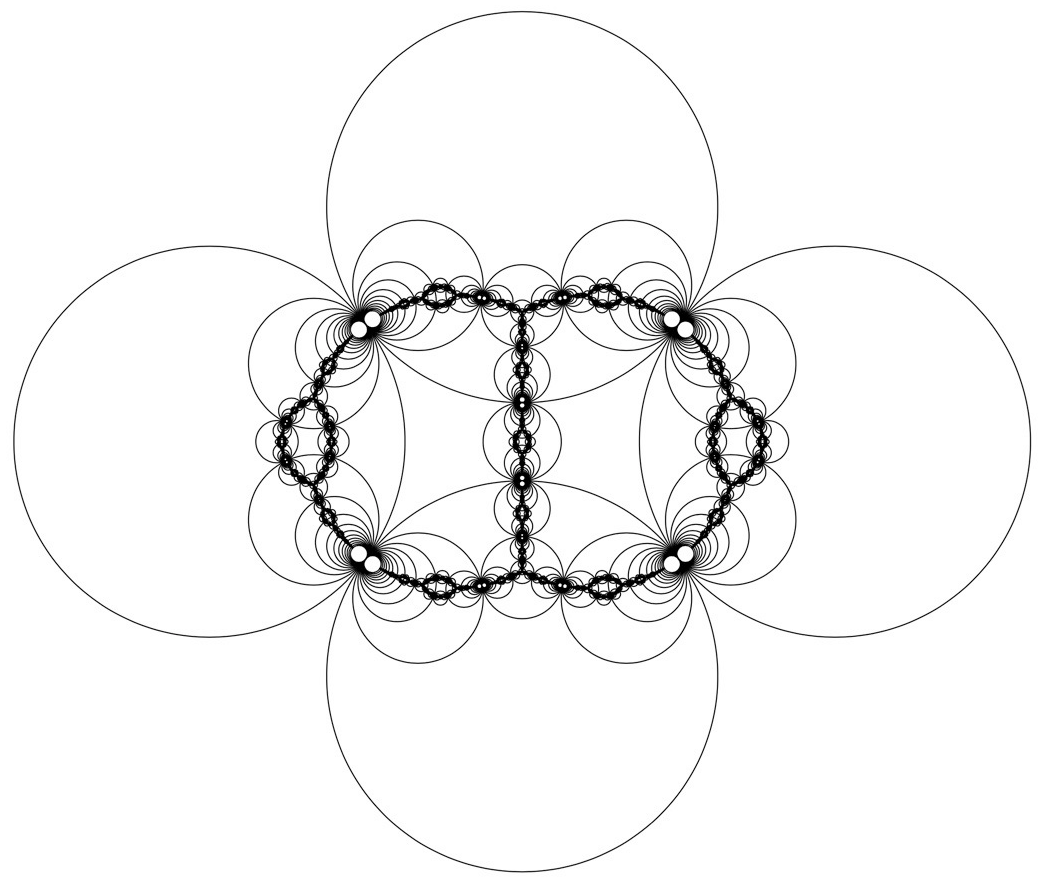}\ \includegraphics[width=0.46\linewidth]{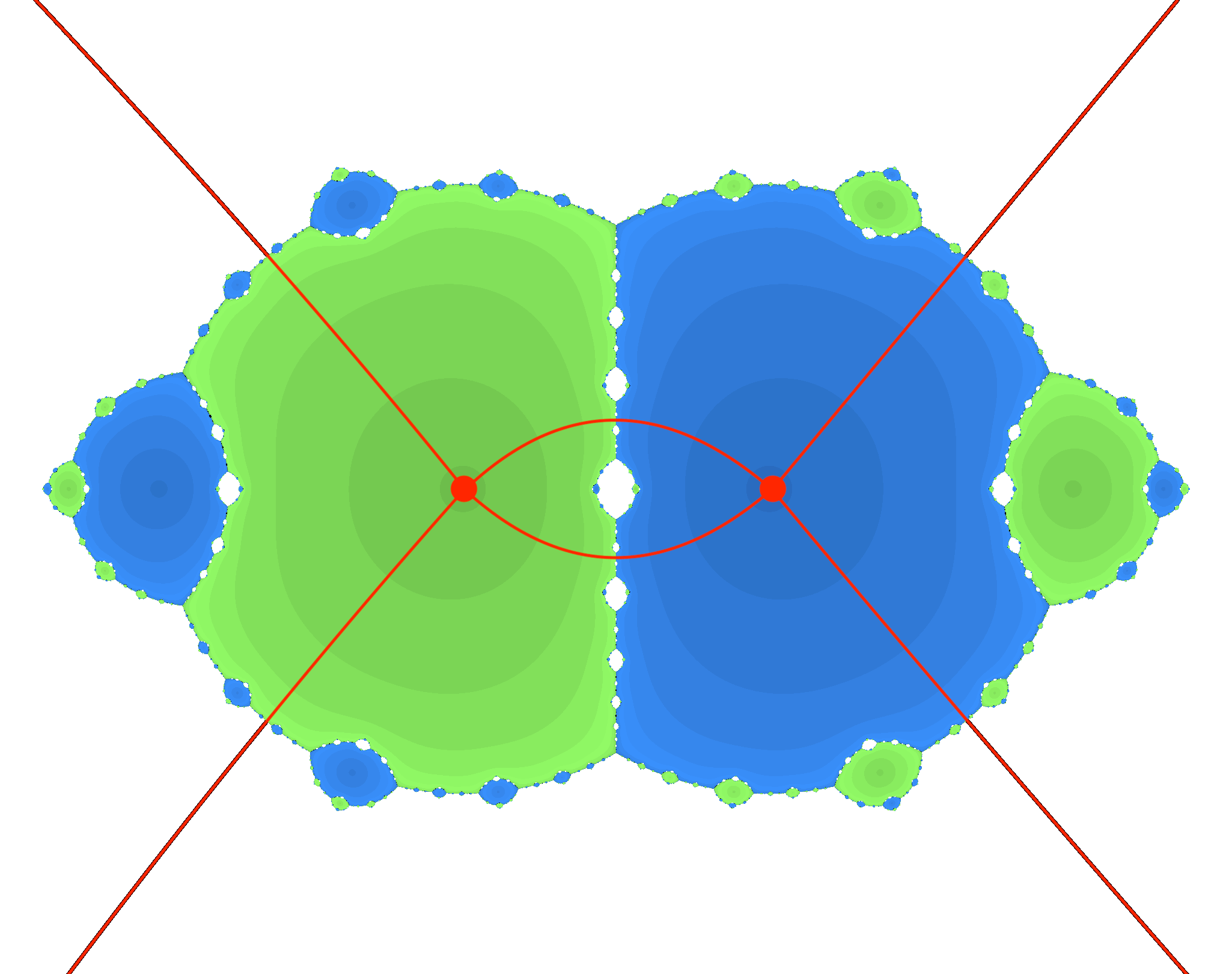}
\caption{Top: a simple plane $2$-connected graph $\Gamma$ with $5$ vertices. Bottom left: the corresponding circle packing $\mathcal{P}$ and the limit set of the associated kissing reflection group $G_{\mathcal{P}}$ of rank $5$. Bottom right: The Julia set of the corresponding degree $4$ critically fixed anti-rational map $R_\Gamma$ with its Tischler graph $\mathscr{T}(R_\Gamma)$ drawn in red.}
\label{non_gasket_fig}
\end{figure}

To see that the critically fixed anti-rational map $R_\Gamma$ is unique (up to M{\"o}bius conjugacy), one needs to look at a combinatorial invariant called \emph{Tischler graph}. Roughly speaking, Tischler graphs are to critically fixed anti-rational maps what Hubbard trees are to polynomials. More precisely, the Tischler graph $\mathscr{T}(R)$ of a critically fixed anti-rational map $R$ is defined as the union (of the closures) of all invariant rays in various fixed Fatou components (recall that the restriction of $R$ to a fixed Fatou component is conformally conjugate to $\overline{z}^r\vert_{\D}$, for some $r\geq 2$). Thus, it is a forward invariant graph containing all the critical points of $R$. It is easily checked using Thurston rigidity that a critically fixed anti-rational map is completely determined by its Tischler graph.
It turns out from the construction of $R_\Gamma$ that its Tischler graph $\mathscr{T}(R_\Gamma)$ is the planar dual of $\Gamma$ (see Figure~\ref{non_gasket_fig}), and hence $R_\Gamma$ is unique (up to M{\"o}bius conjugacy). Thus, the surgery procedure described above yields an injective map from the space of quasiconformal conjugacy classes of kissing reflection groups of rank $d+1$ with connected limit set to the space of M{\"o}bius conjugacy classes of critically fixed anti-rational maps of degree $d$.

The fact that this map is a surjection follows from combinatorial properties of the Tischler graphs of critically fixed anti-rational maps. Indeed, according to \cite[Lemma~4.9]{LLM1}, the planar dual of such a Tischler graph is simple and $2$-connected. Then we apply the above construction to a kissing reflection group $G_{\mathcal{P}}$ to construct the desired anti-rational map $R$. 

Let us list the main ingredients of the proof of the above properties of Tischler graphs:
\begin{enumerate}\upshape
\item landing patterns of the invariant rays at repelling fixed points (due to orientation-reversal, at most two invariant rays can land at the same point),
\item each face of a Tischler graph is a Jordan domain which is mapped to its complement as an orientation-reversing homeomorphism (which is a consequence of the fact that the edges of the Tischler graphs are invariant under the maps and the faces do not contain critical points), and
\item the boundaries of two faces of a Tischler graph meet at most along one edge (which follows from the absence of Levy cycles for rational maps).
\end{enumerate} 

We summarize the upshot of the above discussion in the following theorem.

\begin{theorem}\cite[Theorem~1.1, Proposition 4.10]{LLM1}\label{new_line_dict_thm}
The following three sets are in natural bijective correspondence:
\begin{itemize}
\item $\{2$-connected, simple, plane graphs $\Gamma$ with $d+1$ vertices up to isomorphism of plane graphs$\}$,\\

\item $\{$Kissing reflection groups $G$ of rank $d+1$ with connected limit set up to QC conjugacy$\}$,\\

\item $\{$Critically fixed anti-rational maps $R$ of degree $d$ up to M{\"o}bius conjugacy$\}$.
\end{itemize}

\noindent Moreover, 
\begin{enumerate}\upshape
\item the Tischler graph of $R_\Gamma$ is the planar dual of $\Gamma$, and
\item there exists a homeomorphism between the limit set $\Lambda(G_\Gamma)$ and the Julia set $\mathcal{J}(R_\Gamma)$ that conjugates $\cN_{G_{\mathcal{P}}}$ to $R_\Gamma$.
\end{enumerate}
\end{theorem}

As a by-product of the above discussion, we have a full classification of critically fixed anti-rational maps. In particular, the above result shows that the Julia set of a critically fixed anti-rational map is connected.
(A slightly different, but ultimately equivalent, classification of critically fixed anti-rational maps was given independently by Geyer~\cite{Gey20}.)

\subsubsection{From critically fixed anti-rational maps to kissing reflection groups by David surgery}\label{david_regularity_subsubsec}

Using the results of Section~\ref{david_surgery_sec}, one can invert the topological surgery of Subsection~\ref{group_map_bijection_subsubsec}.
Specifically, suppose that the critically fixed anti-rational map $R$ corresponds to the kissing reflection group $G$ in the bijection of Theorem~\ref{new_line_dict_thm}. 
One can use Lemma~\ref{david_surgery_lemma} to replace the dynamics of $R$ on its invariant Fatou components with Nielsen maps $\pmb{\cN}_r$ of appropriate degree thus producing a Schwarz reflection map $\sigma$ (cf. Subsection~\ref{rat_map_apollo_group_subsubsec}). The touching structure of the faces of the Tischler graph $\mathscr{T}(R)$ (recall that each face is a Jordan domain and two faces meet at most along one edge) can be used to conclude that $\mathrm{Dom}(\sigma)$ is the closure of the union of finitely many disjoint Jordan domains whose touching structure is given by the planar dual of $\mathscr{T}(R)$. Since the Nielsen map $\pmb{\cN}_r$ fixes the boundary of its domain of definition pointwise, one also deduces that each component of $\Int{\mathrm{Dom}(\sigma)}$ is a Jordan quadrature domain and $\sigma$ is the piecewise Schwarz reflection map associated with a quadrature multi-domain. Finally, the fact that $R$ maps each face of $\mathscr{T}(R)$ homeomorphically to its complement implies that each of the above quadrature domains is a round disk. Consequently, $\sigma$ is the piecewise circular reflection map on the disks of a circle packing with associated contact graph given by the planar dual of $\mathscr{T}(R)$. In order words, $\sigma$ is the Nielsen map of a kissing reflection group that is quasiconformally conjugate to $G$.

\begin{figure}[h!]
\captionsetup{width=0.96\linewidth}
\includegraphics[width=0.22\linewidth]{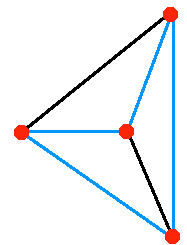}\ \hspace{4mm}\  \includegraphics[width=0.3\linewidth]{apollonian_gasket.png}\ \hspace{4mm}\  \includegraphics[width=0.3\linewidth]{apollonian_julia.png}
\smallskip

\includegraphics[width=0.22\linewidth]{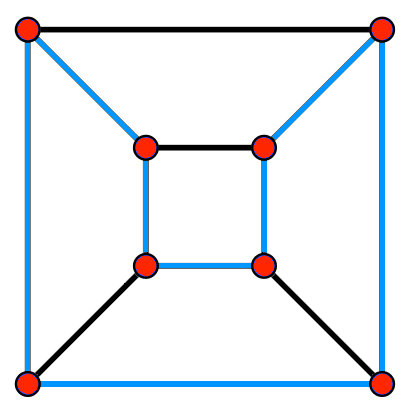}\ \hspace{4mm}\ \includegraphics[width=0.3\linewidth]{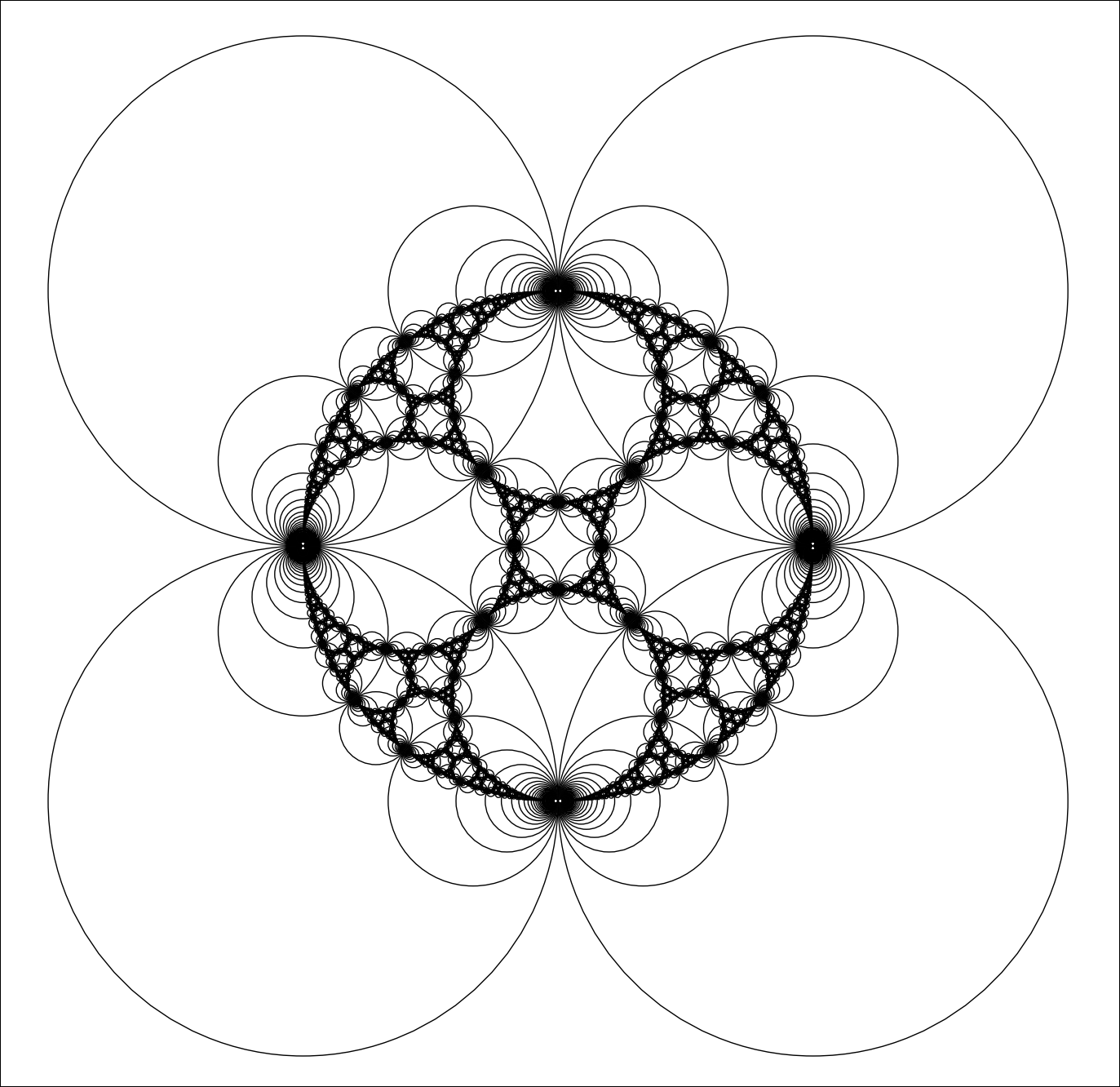}\ \hspace{4mm}\  \includegraphics[width=0.3\linewidth]{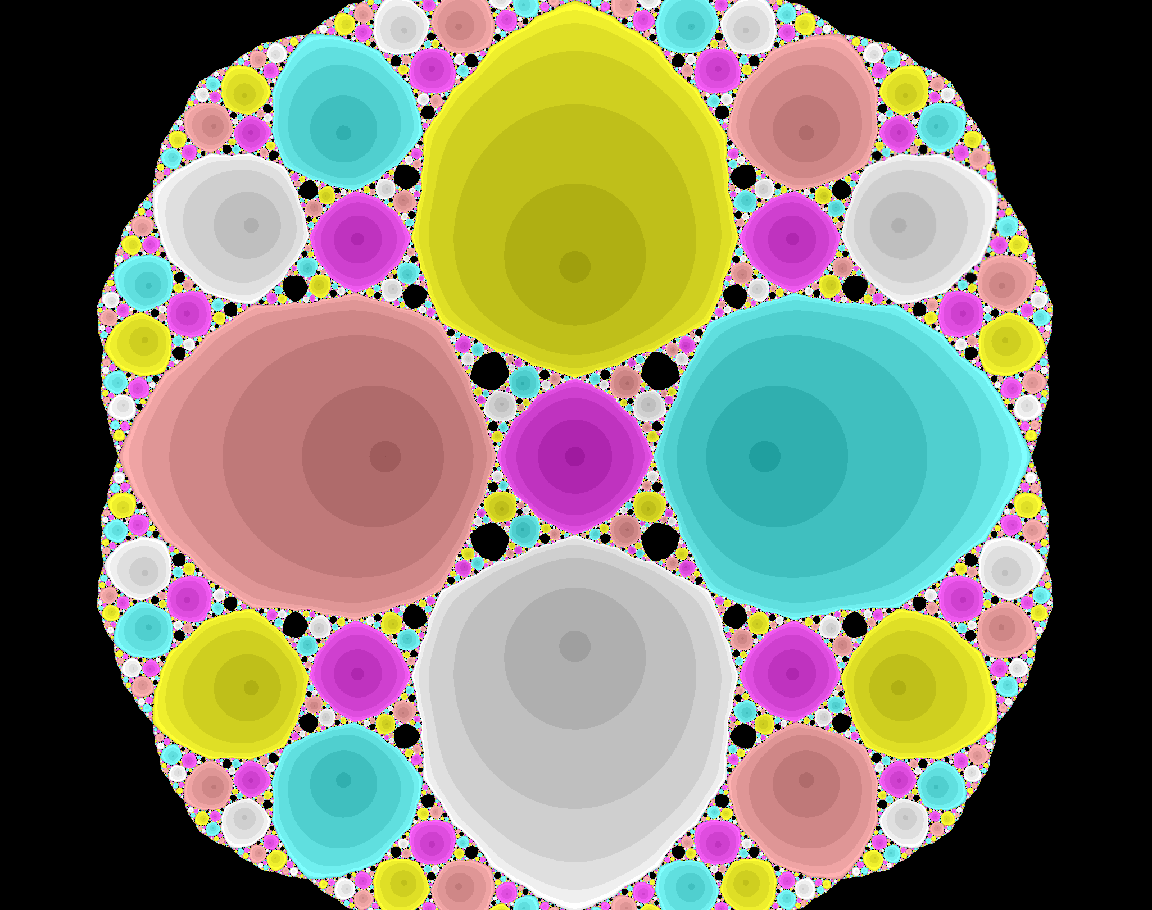}  
\caption{In each row, a polyhedral and Hamiltonian graph along with a choice of a Hamiltonian cycle (left), the corresponding circle packing and the limit set of the associated kissing reflection group (middle), and the gasket Julia set of the corresponding critically fixed anti-rational map (right) are displayed. }
\label{hamiltonian_polyhedral_fig}
\end{figure}

The duality between the Tischler graph $\mathscr{T}(R)$ and the contact graph $\Gamma$ of the circle packing defining $G$ can be seen directly as follows. The Tischler graph decomposes the sphere into finitely many faces/cells. Each such face is a polygon $\Delta_m$, and two faces can share at most one edge of $\mathscr{T}(R)$ on their boundary. Further, each face $\Delta_m$ is mapped by $R$ to its complement as an orientation-reversing homeomorphism (informally speaking, $R$ acts on $\Delta_m$ as  a `reflection'  in $\partial\Delta_m$). 
To pass from $R$ to the reflection group $G$, let us first replace each face $\Delta_m$ with an inscribed Jordan disk $D_m$ touching $\mathscr{T}(R)$ at the repelling fixed points $\alpha_j$ (one on each edge of the face), which are touching points between appropriate immediate basins $U_i$ of $R$.  More specifically, we can connect each pair $\alpha_j$ and $\alpha_k$ that belongs to  the same $\overline U_i$ by the hyperbolic geodesic  in $U_i$.  (Inside each  basin $U_i$, it amounts to  `blowing up' the star of internal rays to an ideal  hyperbolic  $(r_i+1)-$gon $\Pi_i$  meeting  the boundary $\partial U_i$ at the repelling fixed points.) Moreover, the touching pattern between the disks $D_m$ coincides with  the adjacency pattern between the faces $\Delta_m$, so they have the same dual graph $\Gamma=\mathscr{T}(R)^\vee$. Realizing the topological circle packing $(D_m)$ geometrically (as the round circle packing $(\D_m)$), we obtain the desired kissing reflection group $G$ generated by reflections in the circles $\partial \D_m$ (cf. Figures~\ref{non_gasket_fig} and~\ref{hamiltonian_polyhedral_fig}).
Note finally that on each basin $U_i$, this surgery amounts to replacing the $\bar z^{r_i}-$action with the action of the Nielsen map corresponding to the ideal polygon $\Pi_i$.

It follows from the above David surgery construction that the topological conjugacy from the action of a critically fixed anti-rational map on its Julia set to the action of the Nielsen map of a kissing reflection group on its limit set is the restriction of a David homeomorphism of $\widehat{\C}$. We refer the reader to \cite[\S 8]{LMMN} for details of the above~construction.

\subsubsection{Dynamical ramifications of the bijection}\label{dyn_conseq_bijection_subsubsec}

Various properties of the graph $\Gamma$ translate to interesting dynamical properties of the corresponding kissing reflection group and anti-rational map. In particular, it provides bijections between interesting subclasses of kissing reflection groups and critically fixed anti-rational maps, which we will discuss next. (For further applications of these ideas, see Subsection~\ref{bdd_thm_parallels_subsubsect}.)
\smallskip

\noindent \emph{Class 1: Necklace groups vs critically fixed anti-polynomials.} Recall from Subsection~\ref{necklace_subsubsec} that a kissing reflection group is called necklace if the associated graph $\Gamma$ is $2$-connected and outerplanar. Since an outerplanar graph $\Gamma$ has a face that contains all the $d+1$ vertices (of $\Gamma$) on its boundary, it follows from Theorem~\ref{new_line_dict_thm} that the Tischler graph $\mathscr{T}(R_\Gamma)$ of the corresponding anti-rational map $R_\Gamma$ has a vertex of valence $d+1$. Hence, the degree $d$ anti-rational map $R_\Gamma$ has a critical point of local degree $d$; i.e., $R_\Gamma$ is an anti-polynomial (cf. \cite[Corollary~4.17]{LLM1}). This can be alternatively seen from the fact that a necklace group $G$ has a marked component $\Omega_\infty(G)$ in its domain of discontinuity where the action of $\cN_G$ is quasiconformally conjugate to $\pmb{\cN}_d$, and hence the surgery construction of Subsection~\ref{group_map_bijection_subsubsec} replaces $\cN_G\vert_{\Omega_\infty(G)}$ with $\overline{z}^d$. Thus, the bijection of Theorem~\ref{new_line_dict_thm} restricts to a bijection between necklace groups and critically fixed anti-polynomials (cf. Theorem~\ref{sigma_d_bers_homeo_thm}). Moreover, the geodesic lamination of a necklace group is carried to the lamination of the Julia set of the corresponding critically fixed anti-polynomial by the Minkowski circle homeomorphism $\pmb{\mathcal{E}}_d$ (cf. \cite[\S 4.3]{LLM1}, \cite[Propositions~56,~64]{LMM2}).
\smallskip

\noindent \emph{Class 2: Quasi-Fuchsian closure vs mating.} Let us suppose that $\Gamma$ is a Hamiltonian graph. Then $\Gamma$ can be obtained by drawing additional edges connecting the vertices of a polygonal graph. Thus, a circle packing $\mathcal{P}$ corresponding to $\Gamma$ can be thought of as a limit of deformations of the circle packing $\pmb{\mathcal{P}}_{d}$ where the deformations introduce additional tangencies among the circles in the packing (see Definition~\ref{regular_ideal_polygon_ref_group_def}). The kissing reflection group $G_{\mathcal{P}}$ can be realized as a limit of a quasiconformal deformation of the base group $\pmb{G}_{d}$, where the quasiconformal maps deform the conformal class of the polygons of the fundamental domain $\Pi(\pmb{G}_{d})$ (cf. Subsection~\ref{qf_bdry_mating_subsubsec}). Hence, for a circle packing $\mathcal{P}$ with a Hamiltonian contact graph, the group $G_{\mathcal{P}}$ lies in the closure of the quasi-Fuchsian deformation space of $\pmb{G}_{d}$. We refer the reader to \cite[Propositions~3.17, 3.18]{LLM1} for a complete proof of this statement that relates the above quasiconformal deformations to pinching appropriate simple closed geodesics on the surfaces $\D/\widetilde{\pmb{G}}_{d}$ and $(\widehat{\C}\setminus\overline{\D})/\widetilde{\pmb{G}}_{d}$, where $\widetilde{\pmb{G}}_{d}$ is the index two Fuchsian subgroup of $\pmb{G}_{d}$.

A Hamiltonian cycle also induces a splitting of the Hamiltonian graph $\Gamma$ into a pair of $2$-connected outerplanar graphs $\Gamma^+$ and $\Gamma^-$. It was shown in \cite[Proposition~3.21]{LLM1} that an associated kissing reflection group $G_{\mathcal{P}}$ (where the circle packing $\mathcal{P}$ has $\Gamma$ as its contact graph) can be interpreted as a \emph{mating} of two necklace groups $G_{\mathcal{P}^+}$ and $G_{\mathcal{P}^-}$, where the circle packings $\mathcal{P}^{\pm}$ have the outerplanar graphs $\Gamma^{\pm}$ as their contact graphs. Roughly speaking, this means that the filled limit sets $K(G_{\mathcal{P}^{\pm}})$ can be pasted together along their limit sets to produce the Riemann sphere in such a way that the action of $G_{\mathcal{P}^{\pm}}$ on $K(G_{\mathcal{P}^{\pm}})$ is semi-conjugate (in fact, conformally conjugate away from the limit set) to the action of $G_{\mathcal{P}}$ on a suitable subset of $\widehat{\C}$. Using the relation between necklace groups and critically fixed anti-polynomials, this mating statement can be transported to the rational map side proving that the anti-rational map $R_\Gamma$ is a conformal mating of the two anti-polynomials $R_{\Gamma^{\pm}}$ that correspond to the groups $G_{\mathcal{P}^{\pm}}$ (see \cite[Corollary~4.17]{LLM1} for details). Hence, Theorem~\ref{new_line_dict_thm} yields a bijection between groups in the closure of the quasi-Fuchsian deformation space of $\pmb{G}_{d}$ and critically fixed anti-rational maps that are matings of two anti-polynomials such that the bijection commutes with the operation of mating in respective categories. (For a specific example, see the Apollonian reflection group and anti-rational map discussed in Subsection~\ref{qf_bdry_mating_subsubsec}.) We also remark that the unmating of $R_{\Gamma}$ into a pair of anti-polynomials depends on the choice of a Hamiltonian cycle in $\Gamma$, and hence distinct Hamiltonian cycles in $\Gamma$ may lead to different unmatings of $R_{\Gamma}$. This gives rise to many examples of shared matings in the world of critically fixed anti-rational maps (cf. \cite[Appendix~B.1]{LLM2}). For earlier examples of shared matings, see \cite{BEKMPRT12}. 

Note that a more intrinsic way of unmating critically fixed anti-rational maps associated with Hamiltonian graphs can be derived from Meyer's results \cite[Theorem~4.2]{Meyer14}. Indeed, the existence of a Hamiltonian cycle for $\Gamma$ translates to the existence of an equator (in the sense of \cite[Definition~4.1]{Meyer14}) for the corresponding anti-rational map $R_\Gamma$.

We conclude this subsection with the following question, which is inspired by the discussion in Subsection~\ref{qs_grp_subsubsec}.
\begin{problem}
Describe the quasisymmetry groups of the Julia sets of critically fixed anti-rational maps and the limit sets of the corresponding kissing reflection groups. Can such Julia sets be quasiconformally equivalent to the limit sets of the corresponding kissing reflection groups?
\end{problem}

\subsection{Parameter space ramifications of the bijection}\label{new_line_dict_para_subsec}

The dictionary between kissing reflection groups and critically fixed anti-rational maps has surprising parameter space consequences. Recall from Subsection~\ref{kissing_group_deform_space_subsubsec} that each kissing reflection group $G$ sits in a quasiconformal deformation space $\mathcal{QC}(G)$. To discuss some of the analogies between the parameter spaces of reflection groups and anti-rational maps, we will first introduce an appropriate counterpart of quasiconformal deformation spaces in the anti-rational map world.

\subsubsection{Pared deformation space of anti-rational maps}\label{pared_def_space_anti_rat_subsubsec}

Let us set $\mathrm{Rat}_d^-(\C)$ to be the space of degree $d$ anti-rational maps, and
$\mathcal{M}_d^-:= \mathrm{Rat}_d^-(\C)/\mathrm{PSL}_2(\C)$ to be the corresponding moduli space. For a critically fixed anti-rational map $R_\Gamma$ (where the Tischler graph of $R_\Gamma$ is dual to $\Gamma$), we denote the component of hyperbolic maps in $\mathcal{M}_d^-$ containing $[R_\Gamma]$ by $\mathcal{H}_\Gamma$. Although $\mathcal{H}_\Gamma$ consists of anti-rational maps whose Julia set dynamics are quasiconformally conjugate to that of $R_\Gamma$, it turns out that $\mathcal{H}_\Gamma$ is not quite the correct analog of $\mathcal{QC}(\Gamma)$ due to the following reason.

If a kissing reflection group $G_\mathcal{P}$ corresponds to the anti-rational map $R_\Gamma$ in the bijection of Theorem~\ref{new_line_dict_thm}, then the cusps of the $3$-manifold with boundary $\mathcal{M}(G_\mathcal{P}) := \left(\mathbb{H}^3 \cup \Omega(G_\mathcal{P})\right)/\widetilde{G}_\mathcal{P}$ (where $\widetilde{G}_\mathcal{P}$ is the index two Kleinian subgroup of $G_\mathcal{P}$) correspond bijectively to the repelling fixed points of $R_\Gamma$ (see \cite[\S 4.3]{LLM1}). While the parabolics of $G_\mathcal{P}$ remain parabolic throughout $\mathcal{QC}(G_\mathcal{P})$, the multipliers of the repelling fixed points can grow arbitrarily large in $\mathcal{H}_\Gamma$. Thus, preserving the parabolics of $G_\mathcal{P}$ on the group side is analogous to controlling the multipliers of all repelling fixed points of maps in $\mathcal{H}_\Gamma$. 

We now mention the necessary adjustment that leads to a more natural deformation space of $R_\Gamma$ from the perspective of our dictionary.  Note that throughout $\mathcal{H}_\Gamma$, the dynamics of the maps on any invariant Fatou component is conformally conjugate to the dynamics of an anti-Blaschke product on $\D$. We refer to such an anti-Blaschke product as the uniformizing model. For $K>0$, we define the \emph{pared deformation space} $\mathcal{H}_\Gamma(K)\subset\mathcal{H}_\Gamma$ of $R_\Gamma$ to be the connected component containing $[R_\Gamma]$ of the set of anti-rational maps $[R]\in\mathcal{H}_\Gamma$ where the multiplier of any repelling fixed point in any uniformizing model is bounded by $K$. We refer the reader to \cite[\S 4.1]{LLM2} for a detailed discussion of pared deformation spaces of critically fixed anti-rational maps.

\begin{remark}
The terminology `pared deformation space' is borrowed from Kleinian group literature. While studying the deformation space of a Kleinian group or a $3$-manifold with cusps, it is often important to restrict attention to deformations preserving parabolic elements. Such parabolic-preserving deformation spaces are called \emph{pared deformation space} in the $3$-manifold literature (cf. \cite{Thu86}). Since the imposition of a bound on the multipliers of the repelling fixed points of critically fixed anti-rational maps (in uniformizing models) is analogous to parabolic-preserving deformation of kissing reflection groups, we use the same term `pared deformation space' to denote the subset $\mathcal{H}_\Gamma(K)$.
\end{remark}

\subsubsection{Parallel results: boundedness theorems}\label{bdd_thm_parallels_subsubsect}

According to \cite[Proposition~3.10]{LLM1}, the graph $\Gamma$ is polyhedral/$3$-connected if and only if the limit set $\Lambda(G_\mathcal{P})$ is a gasket (where the contact graph of the circle packing $\mathcal{P}$ is isomorphic to $\Gamma$ as a plane graph); i.e., it is homeomorphic to a set $\Lambda$ that is the closure of some infinite circle packing such that the complement of $\Lambda$ is a union of round disks which is dense in $\widehat{\C}$. In fact, the requirement that $\Gamma$ is polyhedral is equivalent to the conditions that each component of $\Omega(G_\mathcal{P})$ is a Jordan domain and the closures of any two different components of $\Omega(G_\mathcal{P})$  intersect at most at a single point which is necessarily a cusp. The latter property is, in turn, equivalent to the so-called \emph{acylindricity} property for the $3$-manifold $\mathcal{M}(G_\mathcal{P})$ \cite[Proposition~3.6]{LLM1}. Indeed, the condition that the closures of any two different components of $\Omega(G_\mathcal{P})$ may intersect only at a cusp means that there are no essential cylinder other than the cusp pairing cylinders in the $3$-manifold $\mathcal{M}(G_\mathcal{P})$ (see \cite[\S 3.2]{LLM1} for a formal definition of acylindricity). For us, the importance of the acylindricity property stems from the Thurston Compactness Theorem \cite{Thu86}, which implies the following:
$$
\mathcal{QC}(G_\mathcal{P})\ \textrm{is bounded (i.e., pre-compact in}\ \mathrm{AH}(G_\mathcal{P})\textrm{)}\ \iff \mathcal{M}(G_\mathcal{P})\ \mathrm{is\ acylindrical.}
$$
Combining this with \cite[Proposition~3.6]{LLM1}, we conclude that
$$
\mathcal{QC}(G_\mathcal{P})\ \textrm{is bounded}\ \iff \Gamma\ \textrm{is polyhedral.}
$$

It is natural to ask whether an analog of the above statement holds in the anti-rational map setting. In fact, inspired by the Thurston Compactness Theorem in the convex-cocompact setting, McMullen conjectured in \cite[Question~5.3]{McM95} that hyperbolic components of rational maps with Sierpinski carpet Julia sets are bounded. It turns out that a counterpart of the Thurston Compactness Theorem and McMullen's conjecture indeed holds in the anti-rational map setting.

\begin{theorem}\cite[Theorem 1.1]{LLM2}\label{compactness_parallel_thm}
If $\Gamma$ is polyhedral, then for any $K>0$, the pared deformation space $\mathcal{H}_\Gamma(K)$ is bounded in $\mathcal{M}_d^-$. Conversely, if $\Gamma$ is not polyhedral, then there exists some $K>0$ such that $\mathcal{H}_\Gamma(K)$ is unbounded in $\mathcal{M}_d^-$.
\end{theorem}

We first note that pared deformation spaces play an essential role in the above theorem since the full hyperbolic component $\mathcal{H}_\Gamma$ is never bounded (see \cite[Proposition~4.16]{LLM1}).

\begin{figure}[h!]
\captionsetup{width=0.96\linewidth}
\includegraphics[width=0.6\linewidth]{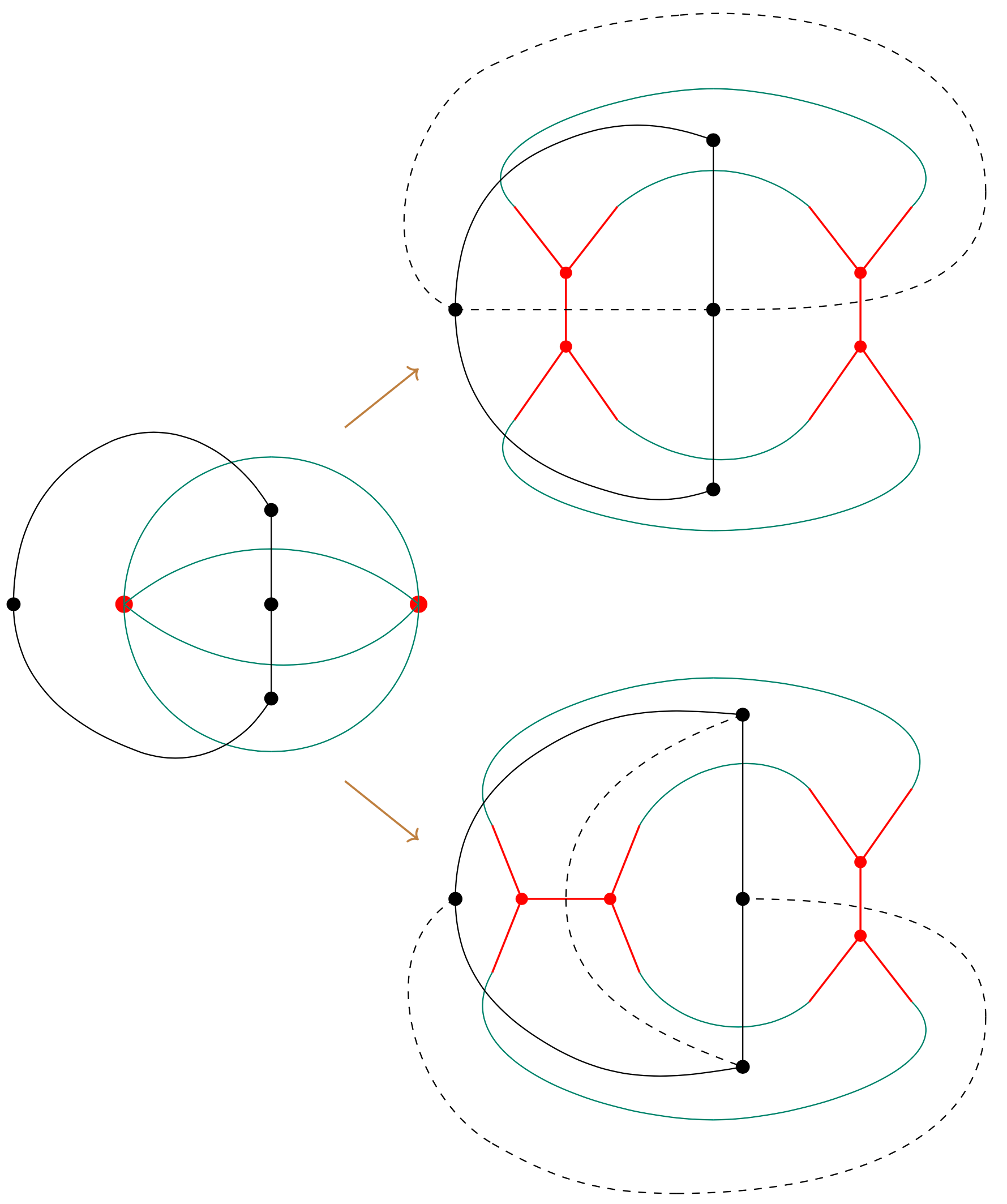}
\caption{Left: The Tischler graph of $\overline{z}^3$ is shown in green and its dual is depicted in black. Right: Two enrichments of the Tischler graph of $\overline{z}^3$ are shown in green/red, where the vertices of the original Tischler graph are blown up to trees. In the top figure, the dual of the enrichment has a bigon; while in the bottom figure, the dual graph is simple. The duals of both enrichments dominate the dual of the original Tischler graph.}
\label{enrichments_fig}
\end{figure}

Let us now spend a few words on the proof of Theorem~\ref{compactness_parallel_thm}. The uniform upper bound on the derivatives of repelling fixed points for maps in $\mathcal{H}_\Gamma(K)$ gives uniform control on the displacement of the critical points under dynamics. This enables to study the limiting dynamics of an escaping sequence in $\mathcal{H}_\Gamma(K)$ on each invariant Fatou component in terms of a \emph{quasi-fixed tree} (see \cite[\S 2]{LLM2} for details of the construction). This, in turn, allows one to record the global limiting dynamics of degenerations in $\mathcal{H}_\Gamma(K)$ by blowing up or tuning the Tischler graph of $R_\Gamma$ by the above quasi-fixed trees (see \cite[\S 4]{LLM2}). The graphs obtained by this blowing up procedure are called \emph{enriched Tischler graphs}, and they capture the combinatorics of degenerations in $\mathcal{H}_\Gamma(K)$. An important property of the planar duals of enriched Tischler graphs is that they have no self-loop or topologically trivial bigon, and that they dominate the graph $\Gamma$, which is the planar dual of the Tischler graph of the center of $\mathcal{H}_\Gamma$ (see Definition~\ref{domination_def}, cf. \cite[Definition~4.5, Proposition~4.6]{LLM2}). Examples of such domination are displayed in Figure~\ref{enrichments_fig}.

Using rescaling limit arguments, it was proved in \cite[Theorem~4.11]{LLM2} that an escaping sequence in $\mathcal{H}_\Gamma(K)$ converges in $\mathcal{M}_d^-$ if and only if the planar dual $\Gamma'$ of the enriched Tischler graph is simple; i.e., has no non-trivial bigon (see Figure~\ref{enrichments_fig}). Theorem~\ref{compactness_parallel_thm} now follows from the graph-theoretic fact that $\Gamma$ is polyhedral if and only if all planar graphs $\Gamma'$ dominating $\Gamma$ are simple (see \cite[Lemma~4.10]{LLM2}).

\subsubsection{Parallel results: interaction between deformation spaces}\label{def_space_touching_parallels_subsubsect}

Recall from Proposition~\ref{cell_structure_group_prop} that for a pair of simple, $2$-connected, plane graphs $\Gamma, \Gamma'$ with the same number of vertices, we have that
$$
\mathcal{QC}(G_{\Gamma'})\subset \overline{\mathcal{QC}(G_\Gamma)}\ \iff\ \Gamma' > \Gamma.
$$
Remarkably, the various pared deformation spaces of critically fixed anti-rational maps also have the same pattern of interaction.

\begin{theorem}\cite[Theorem~1.2]{LLM2}\label{def_space_interaction_parallel_thm}
Let $\Gamma,\Gamma'$ be two distinct $2$-connected, simple, plane graphs with the same number of vertices. For all large $K>0$, the pared deformation space $\mathcal{H}_\Gamma(K)$ parabolic bifurcates to $\mathcal{H}_{\Gamma'}(K)$ if and only if $\Gamma'\geq\Gamma$.
\end{theorem}
\noindent In the above theorem, $\mathcal{H}_\Gamma(K)$ is said to parabolic bifurcate to $\mathcal{H}_{\Gamma'}(K)$ if $\Gamma\neq\Gamma'$ and the intersection $\partial\mathcal{H}_\Gamma(K)\cap\partial\mathcal{H}_{\Gamma'}(K)$ contains a parabolic map whose Julia dynamics is topologically conjugate to that of maps in $\mathcal{H}_{\Gamma'}$ (an anti-rational map $R$ is parabolic if every critical point of $R$ lies in the Fatou set and $R$ has at least one parabolic cycle). Such a parabolic map is called a \emph{root} of $\mathcal{H}_{\Gamma'}(K)$.

The proof of this theorem also employs the notion of enriched Tischler graphs associated with degenerations in $\mathcal{H}_\Gamma(K)$. The necessity of the condition $\Gamma'\geq\Gamma$ in the parabolic bifurcation statement comes from the fact that the dual of an enriched Tischler graph dominates the dual of the original Tischler graph (see \cite[Proposition~4.6]{LLM2}). On the other hand, the demonstration of the fact that $\mathcal{H}_\Gamma(K)$ indeed parabolic bifurcates to $\mathcal{H}_{\Gamma'}(K)$ for $\Gamma'\geq\Gamma$ involve the following key steps. Here we assume that $\Gamma,\Gamma'$ are two distinct $2$-connected, simple, plane graphs with $\Gamma'\geq\Gamma$. 
\smallskip

\noindent$\bullet$ The vertices of the Tischler graph of $R_\Gamma$ (which is dual to $\Gamma$) can be blown up to trees in such as way that the  planar dual of the resulting graph is isomorphic to the graph $\Gamma'$ \cite[Lemma~4.8]{LLM2}.

\noindent$\bullet$ Each of the above trees arises as quasi-fixed trees associated with degenerating sequences of Blaschke products, and hence the planar dual of $\Gamma'$ is realized as an enriched Tischler graph for some escaping sequence in $\mathcal{H}_\Gamma(K)$ \cite[Theorem~3.1, Proposition~4.4]{LLM2}.

\noindent$\bullet$ As the enriched Tischler graph $\Gamma'$ is simple, such an escaping sequence in $\mathcal{H}_\Gamma(K)$ converges in $\mathcal{M}_d^-$, and the limiting map is parabolic whose Julia dynamics is conjugate to that of maps in $\mathcal{H}_{\Gamma'}(K)$ \cite[Theorem~4.11]{LLM2}.
\smallskip

Since the collection of all $2$-connected simple plane graphs is connected under the domination relation, we have the following simple consequence of Theorem~\ref{def_space_interaction_parallel_thm}.

\begin{corollary}\label{crit_fixed_hyp_comps_closure_conn_cor}
The union of the closures of hyperbolic components (respectively, pared deformation spaces) of degree $d$ critically fixed anti-rational maps is connected.
\end{corollary}

We refer the reader to \cite[Appendix~A.4]{LLM2} for a refinement of Theorem~\ref{def_space_interaction_parallel_thm} that counts the number of accesses from $\mathcal{H}_{\Gamma}(K)$ to $\mathcal{H}_{\Gamma'}(K)$ in terms of the number of nonequivalent embeddings of $\Gamma$ into $\Gamma'$ and to \cite[Appendix~B]{LLM2} for applications of this count to the phenomena of shared matings and existence of self-bumps on boundaries of hyperbolic components (for earlier examples of self-bumps on boundaries of hyperbolic components, see \cite{Luo23}).

\begin{remark}
The boundaries of two hyperbolic components $\mathcal{H}_\Gamma, \mathcal{H}_{\Gamma'}$ may have wild intersection. One works with pared deformation spaces to circumvent this difficulty; indeed, the boundary of a pared deformation space $\mathcal{H}_\Gamma(K)$ only consists of parabolic maps and hence bifurcations between pared deformation spaces are tame.
\end{remark}

\subsubsection{Parallel results: global topology}\label{global_top_parallels_subsubsec}

In light of Corollary~\ref{crit_fixed_hyp_comps_closure_conn_cor}, it is natural to seek finer information about the topology of the union of the closures of all hyperbolic components (respectively, pared deformation spaces) of degree $d$ critically fixed anti-rational maps. On the group side, the moduli space $\mathfrak{M}_n$  of marked circle packings with $n$ circles in $\C$ was studied by Hatcher and Thurston in \cite{HT}, where they showed that the natural map $\Psi: \mathfrak{M}_n \longrightarrow \mathfrak{S}_n$ (where $\mathfrak{S}_n$ denotes the configuration space of $n$ marked points in $\C$) that associates the centers of the circles to a marked circle packing induces a homotopy equivalence between the spaces.

In order to state a counterpart of the Hatcher--Thurston result in the setting of anti-rational maps, let us define $\mathcal{H}_{\Gamma, Rat}(K) \subset \mathcal{H}_{\Gamma, Rat} \subset \mathrm{Rat}_d^-(\C)$ to be the lifts (i.e., preimages under the projection map $\mathrm{Rat}_d^-(\C) \longrightarrow \mathcal{M}_d^-$) of $\mathcal{H}_\Gamma(K) \subset \mathcal{H}_\Gamma \subset \mathcal{M}_d^-$. 
We further set
$$
\mathcal{X}_d(K) := \displaystyle\bigcup_{\Gamma} \overline{\mathcal{H}_{\Gamma, Rat}(K)} \subset \mathrm{Rat}_d^-(\C).
$$

\begin{theorem}\cite[Theorem~1.4]{LLM2}\label{monodromy_parallel_thm}
Let $d\geq 3$. For all large $K$, there exists a surjective monodromy representation
$$
\rho: \pi_1(\mathcal{X}_d(K))  \twoheadrightarrow \mathrm{Mod}(S_{0,d+1}),
$$
where $\mathrm{Mod}(S_{0,d+1})$ is the mapping class group of the $(d+1)$-punctured sphere
\end{theorem}

For maps in $\mathcal{X}_d(K)$, the analog of the circles in the circle packing are suitable Markov partition pieces of the Julia set determined by Tischler graphs. The proof of Theorem~\ref{monodromy_parallel_thm} is carried out by showing that these Markov pieces move `continuously' (away from the grand orbits of fixed points or $2$-cycles) and that their braiding patterns along curves in $\mathcal{X}_d(K)$ can be followed to construct a monodromy representation of $\mathcal{X}_d(K)$ into $\mathrm{Mod}(S_{0,d+1})$. Surjectivity of this representation is demonstrated by exhibiting that suitable half Dehn twists (that generate the mapping class group) are images of certain explicit paths in $\mathcal{X}_d(K)$ (see \cite[\S 5]{LLM2} for details).

\section{Matings of anti-polynomials and necklace groups}\label{mating_anti_poly_nielsen_sec}

In the earlier sections, we discussed the dynamics of various Schwarz reflection maps that arise as matings of anti-polynomials and Nielsen maps of kissing reflection groups. This provided us with matings of geometrically finite quadratic anti-polynomials with the Nielsen map of the ideal triangle group (see Theorems~\ref{deltoid_thm_2},~\ref{c_and_c_general_mating_thm}) and of the specific anti-polynomial $\overline{z}^d$ with Nielsen maps of arbitrary necklace groups (see Theorem~\ref{sigma_d_mating_thm}). These examples lead to the following questions.
\begin{enumerate}\upshape
\item Do Schwarz reflection maps provide a general framework for mating anti-polynomials with Nielsen maps of necklace groups?

\item What is the topological structure of the parameter space of such matings?
\end{enumerate}

The first question above was addressed in \cite{LMMN}, where an existence theorem for matings of large classes of anti-polynomials and necklace groups was established. 

\subsection{Definition of conformal mating}\label{conf_mating_def_subsec}

We now explicate the notion of conformal matings of Nielsen maps of necklace groups with anti-polynomials. The idea is analogous to the classical definition of conformal matings of two (anti-)polynomials.

Let $G \in \overline{\beta(\pmb{G}_{d})}$ be a necklace group associated with a circle packing $\mathcal{P}=\{C_1,\cdots, C_{d+1}\}$. By Proposition~\ref{group_lamination_prop}, there is a canonical semiconjugacy $\phi_{G}: \mathbb{S}^1 \rightarrow \Lambda(G)$ between $\pmb{\cN}_d\vert_{\mathbb{S}^1}$ and $\cN_G\vert_{\Lambda(G)}$, and hence $\cN_G\vert_{\Lambda(G)}$ is a factor or $\pmb{\cN}_d\vert_{\mathbb{S}^1}$. On the other hand, if $P$ is a monic, centered, anti-polynomial of degree $d$ such that $\mathcal{J}(P)$ is connected and locally connected, then $P\vert_{\mathcal{J}(P)}$ is a factor of $\overline{z}^d\vert_{\mathbb{S}^1}$ via the continuous boundary extension $\phi_P:\mathbb{S}^1\to\mathcal{J}(P)$ (of the inverse) of the normalized B{\"o}ttcher coordinate of $\mathcal{B}_\infty(P)$.

Recall also that $\pmb{\mathcal{E}}_d: \mathbb{S}^1 \rightarrow \mathbb{S}^1$ is a topological conjugacy between $\pmb{\cN}_d\vert_{\mathbb{S}^1}$ and $z\mapsto\overline{z}^{d}\vert_{\mathbb{S}^1}$. In other words, the circle coverings induced by the action of the anti-polynomial $P$ on its Julia set and the Nielsen map $\cN_G$ on its limit set are topologically conjugate. This compatibility allows one to glue $\mathcal{K}(P)$ with $K(G)$ along their boundaries and obtain a partially defined continuous map on the resulting topological space.

\begin{definition}\label{conf_mating_equiv_reltn} 
We define the equivalence relation $\sim$ on $K(G) \sqcup \mathcal{K}(P)$ generated by $\phi_G(t)\sim\phi_P(\overline{\pmb{\mathcal{E}}_d(t)})$ for all $t\in\mathbb{S}^1$.
\end{definition}

Clearly, the equivalence relation $\sim$ is preserved by the map 
\begin{center}
	$P\sqcup \cN_{G}: \mathcal{K}(P)\sqcup (K(G)\setminus\Int{\Pi(G)})\to \mathcal{K}(P)\sqcup K(G),$\\
	$P\sqcup \cN_{G}\vert_{\mathcal{K}(P)}=P,\quad  P\sqcup \cN_{G}\vert_{K(G)\setminus\Int{\Pi(G)}}=\cN_G,$
\end{center}
and hence $P\sqcup\ \cN_G$ descends to a continuous map $P\mate \cN_G$ to the quotient 
$$
\mathcal{K}(P)\mate K(G):=\left(\mathcal{K}(P)\sqcup K(G)\right)/\sim.
$$ 
The map $P\mate\cN_G$ is the \emph{topological mating} of $P$ and $\cN_G$. If $\mathcal{K}(P)\mate K(G)$ is homeomorphic to a $2$-sphere, the topological mating is said to be \emph{Moore-unobstructed}.
Finally, one says that $P$ and $\cN_G$ are \emph{conformally mateable} if their topological mating is Moore-unobstructed, and if the topological $2$-sphere $\mathcal{K}(P)\mate K(G)$ admits a complex structure that turns the topological mating $P\mate \cN_G$ into an antiholomorphic map.

Alternatively, an antiholomorphic map $F$ (defined on a subset of the Riemann sphere) is a conformal mating of $P$ and $\cN_G$ if there exist continuous semi-conjugacies from $\mathcal{K}(P), K(G)$ (equipped with the actions of $P, \cN_G$, respectively) into the dynamical plane of $F$ such that the semi-conjugacies are conformal on the interiors, the images of the semi-conjugacies fill up the whole sphere and intersect only along their boundaries as prescribed by the equivalence relation $\sim$. We refer the reader to \cite[\S 10.2]{LMMN} for precise definitions and further details.

\subsection{Existence of conformal matings}\label{conf_mating_gen_thm_subsec}

We will now state a general result that guarantees the existence of conformal matings of necklace groups and anti-polynomials. By definition, conformal mateability of $P$ and $\cN_G$ (as in the previous subsection) requires their topological mating to be Moore-unobstructed. It turns out that for hyperbolic anti-polynomials, this is the only obstruction to conformal mating; i.e., whenever $\mathcal{K}(P)\mate K(G)$ is homeomorphic to $\mathbb{S}^2$, one can upgrade the topological mating of $P$ and $\cN_G$ to a conformal mating.

\begin{theorem}\cite[Lemma~10.17, Theorem~10.20]{LMMN}\label{antipoly_nielsen_mating_thm}
	Let $P$ be a monic hyperbolic anti-polynomial of degree $d$, and let $G\in \overline{\beta(\pmb{G}_{d})}$ be a necklace group. Then, $P$ and $\cN_G$ are conformally mateable if and only if their topological mating is Moore-unobstructed. 
	
	Moreover, if $F:\mathrm{Dom}(F)\to\widehat{\C}$ is a conformal mating of $P$ and $\cN_G$, then each component of the interior of $\mathrm{Dom}(F)$ is a simply connected quadrature domain, and $F$ is the piecewise Schwarz reflection map associated with a quadrature multi-domain.
\end{theorem}

As expected, promoting the topological mating of $P$ and $\cN_G$ to an antiholomorphic map lies at the heart of the difficulty. To achieve this, one uses a combination of the Thurston Realization Theorem and the David surgery procedure of Section~\ref{david_surgery_sec}, which we now outline.

\noindent\textbf{Step I:} According to Theorem~\ref{new_line_dict_thm} (also see the discussion in Subsection~\ref{dyn_conseq_bijection_subsubsec}), there exists a degree $d$ critically fixed anti-polynomial $P_G$ whose Julia set dynamics is topologically conjugate to the dynamics of $\cN_G\vert_{\Lambda(G)}$. The first step in the construction of a conformal mating of $P$ and $\cN_G$ is to mate the two anti-polynomial $P$ and $P_G$. In fact, one readily verifies that the absence of Moore obstruction for the topological mating of $P$ and $G$ is equivalent to the absence of Moore obstruction for the topological mating of $P$ and $P_G$, and then invokes the conformal mateability criterion \cite[Proposition~4.23]{LLM1} (which is a statement about classical polynomial matings, and is an application of the Thurston Realization Theorem) to conclude that if $P\mate\cN_G$ is Moore-unobstructed, then the anti-polynomial $P$ and $P_G$ are conformally mateable.

\noindent\textbf{Step II:} The next step is to turn the conformal mating $R$ of the anti-polynomials $P$ and $P_G$ into a conformal mating of $P$ and $\cN_G$. To this end, one only needs to glue Nielsen maps of ideal polygon reflection groups in suitable critically fixed Fatou components of $R$ (these Fatou components correspond to the critically fixed Fatou components of $P_G$). This is precisely where the David Surgery Lemma~\ref{david_surgery_lemma} comes into play.

The statement that a conformal mating of $P$ and $\cN_G$ is necessarily a piecewise Schwarz reflection map follows from the observation that $\cN_G$ fixes the $\partial \Pi(G)$ pointwise. This fact establishes the naturality of Schwarz reflection maps in the world of combination theorems on firm footing. 

\begin{remark}
A simple algorithm to check whether the topological mating of $P$ and $P_G$ has a Moore obstruction was given in \cite[Lemma~4.22]{LLM1}, which makes Theorem~\ref{antipoly_nielsen_mating_thm} effective in concrete situations. We state this condition here for completeness. Note that each repelling fixed point of the critically fixed anti-polynomial $P_G$ which is a cut-point of $\mathcal{J}(P_G)$ is the landing point of a $2$-cycle of external dynamical rays. We call these finitely many external rays the \emph{principal external rays of $P_G$}. Then, the only possible Moore obstructions are $2$-cycles and $4$-cycles of \emph{extended external rays} (for the formal mating of $P$ and $P_G$) containing principal external rays of $P_G$. We refer the reader to \cite[Section~4.4]{LLM1} for details.
\end{remark}

\subsection{Recognizing conformal matings}\label{conf_mating_identify_subsec}

Suppose that the piecewise Schwarz reflection map 
$$
F:\mathrm{Dom}(F):=\bigcup_{i=1}^k\overline{\Omega_i}\longrightarrow\widehat{\C},\ z\mapsto \sigma_i(z)\ \mathrm{for}\ z\in\overline{\Omega_i}
$$ 
is a conformal mating of a marked degree $d$ hyperbolic anti-polynomial $P$ (with connected Julia set) and the Nielsen map $\cN_G$ of a necklace group $G$ of rank $d+1$. One can determine the topology of $\mathrm{Dom}(F)$ in terms of the following finite combinatorial data associated with $P$ and the structure of accidental parabolics of the group $G$.

Recall that an anti-polynomial $P$ acts on the angles of its external dynamical rays by $m_{-d}:\mathbb{S}^1\to\mathbb{S}^1,\ \theta\mapsto -d\theta$. 
\begin{definition}[Fixed ray lamination]\label{fixed_ray_lamination_def}
Set $\mathrm{Fix}(m_{-d}):= \{0,\frac{1}{d+1},\cdots,\frac{d}{d+1}\}$, the set of angles that are fixed by $m_{-d}$. We define the equivalence relation $\mathfrak{L}_P$ on $\mathrm{Fix}(m_{-d})$ as: $\theta_1\sim\theta_2$ if and only if the external dynamical rays of $P$ at angles $\theta_1, \theta_2$ land at the same point of $\mathcal{J}(P)$. The equivalence relation $\mathfrak{L}_P$ is called the \emph{fixed ray lamination} of $P$.
\end{definition}
We remark that each equivalence class contains at most two elements, and we refer to an equivalence class $\mathfrak{L}_P$ consisting of two elements as a \emph{leaf}. 

The closure of $\Pi^b(G)=\Pi(G)\cap K(G)$ is a tree of polygons, and its boundary meets the limit set $\Lambda(G)$ at finitely many points. We call the set of these points $S_G$. Let $S_G^{\mathrm{cusp}}\subset S_G$ be those $d+1$ points in $S_G$ that do not separate $\Lambda(G)$ (these are the images of the $(d+1)$-st roots of unity under $\phi_G$). Note that in the topological mating of $P$ and $\cN_G$, the points on $\mathcal{J}(P)$ with external address in $\mathrm{Fix}(m_{-d})$ are glued with points in $S_G^{\mathrm{cusp}}$.

The next result, which counts the number of connected components of $\Int{\mathrm{Dom}(F)}$, follows from the fact that each leaf of $\mathfrak{L}_P$ forces a pair of points of $S_G^{\mathrm{cusp}}$ to be identified in the topological mating, thereby creating a disconnection of the interior of $\mathrm{Dom}(F)$. 

\begin{proposition}\cite{LLM23}\label{count_qd_comps_prop}
The number of connected components of $\mathrm{Dom}(F)$ is equal to the number of gaps of the lamination $\mathfrak{L}_P$ (equivalently, one more than the number of leaves in $\mathfrak{L}_P$).
\end{proposition}

We enumerate the gaps of $\mathfrak{L}_P$ cyclically as $\mathcal{G}_1,\cdots,\mathcal{G}_k$ such that the arc $\arc{(0,\frac{1}{d+1})}\subset\partial\mathcal{G}_1$ and $\mathcal{G}_i$ corresponds to the quadrature domain $\Omega_i$ (after possibly renumbering the quadrature domains). We now state a formula for the degrees of rational maps uniformizing the simply connected quadrature domains.  

\begin{proposition}\cite{LLM23}\label{gap_qd_order_prop}
Let $f_i$ be a rational map of degree $d_i$ that carries $\D$ univalently onto $\Omega_i$. Then, $\partial\mathcal{G}_i\cap\mathbb{S}^1$ contains exactly $d_i$ arcs of $\mathbb{S}^1\setminus\mathrm{Fix}(m_{-d})$, for each $i\in\{1,\cdots, k\}$. In particular, $\sum_{i=1}^k d_i=d+1$.
\end{proposition}

For a proof of Proposition~\ref{gap_qd_order_prop}, note that if $\partial\mathcal{G}_i\cap\mathbb{S}^1$ contains exactly $q_i$ arcs of $\mathbb{S}^1\setminus\mathrm{Fix}(m_{-d})$, then the map $m_{-d}\vert_{\partial\mathcal{G}_i\cap\mathbb{S}^1}$ covers $\partial\mathcal{G}_i\cap\mathbb{S}^1$ exactly $(q_i-1)$ times.
Hence, the part of the limit set of $F$ that lies in $\Omega_i$ covers itself $(q_i-1)$ times under the map $F$, which implies that $\sigma_i:\sigma_i^{-1}(\Omega_i)\to\Omega_i$ is a degree $(q_i-1)$ branched covering whence the result follows from Proposition~\ref{simp_conn_quad_prop}.
\begin{figure}[h!]
\captionsetup{width=0.96\linewidth}
\begin{tikzpicture}
\node[anchor=south west,inner sep=0] at (0,0) {\includegraphics[width=0.56\linewidth]{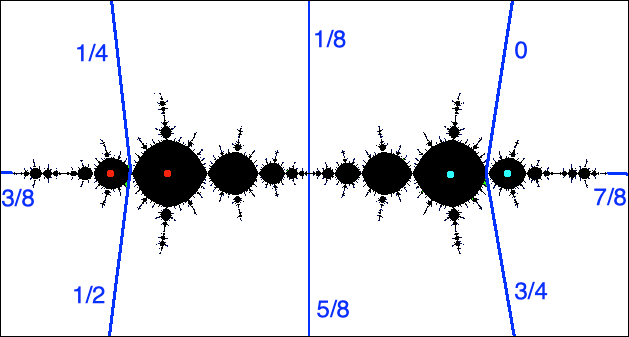}};
\node[anchor=south west,inner sep=0] at (7.4,0) {\includegraphics[width=0.33\linewidth]{bers_slice_cusp.png}};
\node[anchor=south west,inner sep=0] at (0.6,-10.5) {\includegraphics[width=0.8\linewidth]{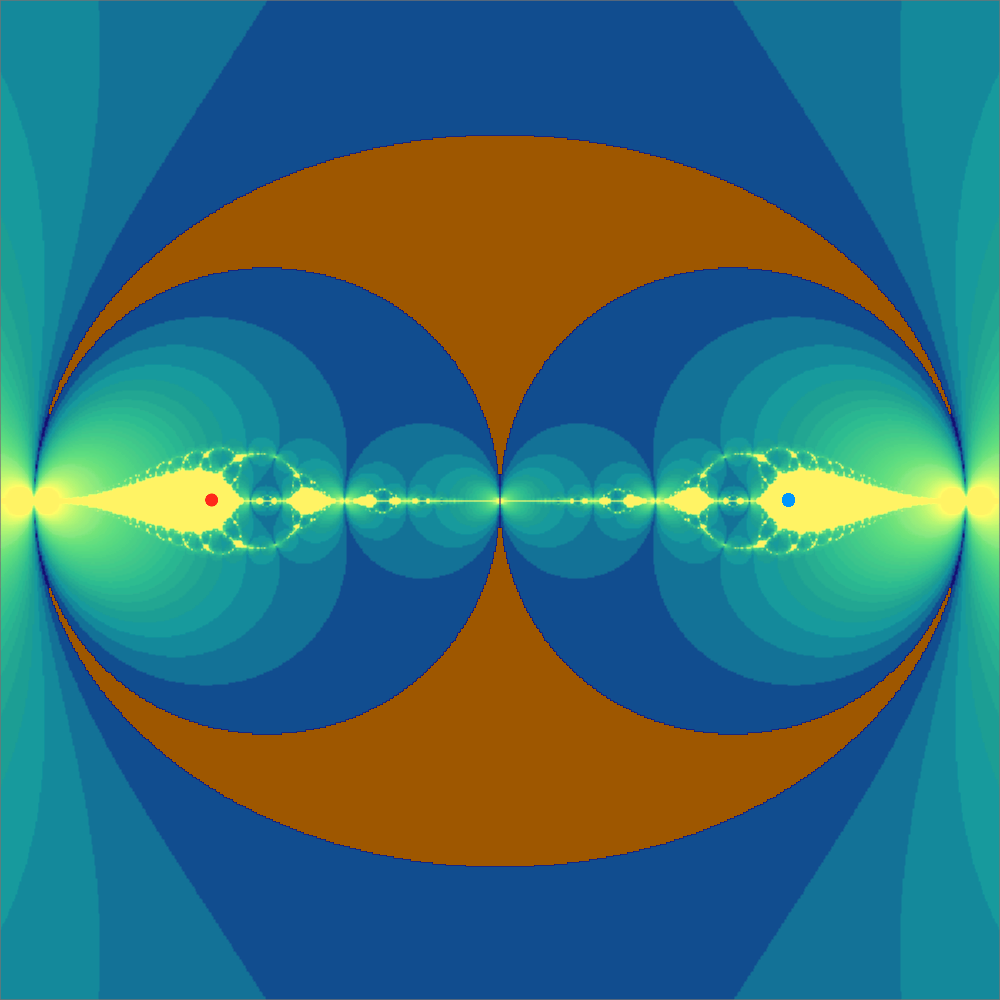}};
\node at (10.5,3.5) {$C_1$};
\node at (7.8,1.05) {$C_2$};
\node at (10.5,0.32) {$C_3$};
\node at (11.32,1.05) {$C_4$};
\end{tikzpicture}
\caption{Top left: The filled Julia set of a cubic anti-polynomial $P$ with two $2$-periodic critical points is displayed. The external dynamical rays of period $1$ and $2$ are also drawn in. Top right: The filled limit set of the Julia necklace group $G$ and the corresponding marked circle packing are depicted. Bottom: Pictured is a part of the dynamical plane of the conformal mating of $P$ and $\cN_G$, which is realized as the piecewise Schwarz reflection map in the exterior of an ellipse and a pair of touching inscribed disks. The droplet is the union of the two brown triangles. The red and blue dots represent the $2$-periodic critical values of the mating. The corresponding $2$-periodic critical points lie in Fatou components contained in the exterior of the ellipse.}
\label{ellipse_disk_fig}
\end{figure}

Note that under the gluing of $\mathcal{K}(P)$ and $K(G)$, the $\phi_P$-image of each arc of $\mathbb{S}^1\setminus\mathrm{Fix}(m_{-d})$ is identified with the $\phi_G$-image of a unique arc of $\mathbb{S}^1\setminus \sqrt[d+1]{1}$. As the $\phi_G$-image of any arc of $\mathbb{S}^1\setminus \sqrt[d+1]{1}$ is enclosed by a unique circle $C_i$ of the packing $\mathcal{P}$, this defines a bijective correspondence between the components of $\mathbb{S}^1\setminus\mathrm{Fix}(m_{-d})$ and the circles $C_1,\cdots, C_{d+1}$. 

The facts that $\cN_G$ has no anti-conformal extension to a relative neighborhood in $K(G)$ of any point of $S_G$ and that the only possible singularities on boundaries of quadrature domains are double points and conformal cusps \cite{Sak91} imply the following structure of double points and cusps on $\partial\Omega$.
\begin{proposition}\cite{LLM23}\label{sing_structure_qd_prop}
\noindent 1) The number of cusps on $\partial\Omega_i$ is equal to the number of points of $\mathrm{Fix}(m_{-d})\cap \partial\mathcal{G}_i$ that are not endpoints of leaves of $\mathfrak{L}_P$.
\smallskip

\noindent 2) The quadrature domains corresponding to the pair of gaps bordering on a leaf of $\mathfrak{L}_P$ have a tangential intersection.
\smallskip

\noindent 3) Suppose that $C_r\cap C_s\neq\emptyset$ with $ r-s\neq \pm 1\ (\mathrm{mod}\ d+1)$. If the arcs of $\mathbb{S}^1\setminus\mathrm{Fix}(m_{-d})$ corresponding to $C_r$ and $C_s$ lie on the boundary of the same gap $\mathcal{G}_i$, then the quadrature domain boundary $\partial\Omega_i$ has a double point corresponding to the pair $(r,s)$. On the other hand, if the arcs of $\mathbb{S}^1\setminus\mathrm{Fix}(m_{-d})$ corresponding to $C_r$ and $C_s$ lie on the boundaries of two distinct gaps $\mathcal{G}_i, \mathcal{G}_j$, then the quadrature domain boundaries $\partial\Omega_i$ and $\partial\Omega_j$ have a tangential intersection (so that $\partial\Omega$ has a double point) corresponding to the pair $(r,s)$.
\end{proposition}

We conclude this subsection with an explicit example. Let $P(z)=\overline{z}^3-\frac{3i}{\sqrt{2}}\overline{z}$ be the anti-polynomial with a pair of $2$-periodic critical points and $G$ be the Julia necklace group introduced in Subsection~\ref{talbot_dynamics_subsubsec}. The anti-polynomial $P$ is obtained by tuning each of the two (bounded) invariant Fatou components of the Julia anti-polynomial by the Basilica anti-polynomial $\overline{z}^2-1$.
Recall from Theorem~\ref{talbot_schwarz_group_anti_poly_thm} that the critically fixed anti-polynomial $P_G$ corresponding to $G$ is given by the Julia anti-polynomial $z\mapsto(3\overline{z}-\overline{z}^3)/2$. It is easily checked using \cite[Lemma~4.22]{LLM1} that the topological mating of $P$ and $P_G$ is Moore-unobstructed, and hence $P$ and $\cN_G$ are conformally mateable. 

The fixed ray lamination of $P$ is given by $\{\{0,3/4\},\{1/4,1/2\}\}$. It thus follows from Propositions~\ref{count_qd_comps_prop} and~\ref{gap_qd_order_prop} that $\Int{\mathrm{Dom}(F)}$ (where $F$ is the conformal mating of $P$ and $\cN_G$) has three components $\Omega_1, \Omega_2, \Omega_3$, two of which (say $\Omega_1,\Omega_2$) are round disks and the other one is the univalent image of $\D$ under a quadratic rational map. Moreover, Proposition~\ref{sing_structure_qd_prop} implies that the quadrature domains $\Omega_1, \Omega_2, \Omega_3$ have non-singular boundaries and they touch each other pairwise. These facts can be used to show that up to M{\"o}bius conjugacy, $\Omega_3$ is the exterior of an ellipse and $\Omega_1,\Omega_2$ are touching round disks inscribed in the ellipse $\partial\Omega_3$ (see  Figure~\ref{ellipse_disk_fig}). For a detailed discussion of this specific example, we refer the reader to \cite[\S 11.2]{LMMN}.

\section{Polygonal Schwarz reflections and connectedness loci of anti-polynomials}\label{mating_para_space_sec}

In this section, which is based on \cite{LLM23}, we introduce a family of degree $d$ piecewise Schwarz reflections that generalizes the C\&C family and the deltoid reflection. We will expound how the mating operation between generic anti-polynomials (with connected Julia set) and the Nielsen map of the ideal $(d+1)$-gon reflection group yields dynamical relations between this family of Schwarz reflections and the connectedness locus $\mathscr{C}_d$ of monic, centered anti-polynomials of degree $d$. 

\subsection{Polygonal Schwarz reflections}\label{poly_schwarz_subsec}

Let $\sigma:\overline{\Omega}\to\widehat{\C}$ be a piecewise Schwarz reflection map associated with a quadrature multi-domain $\Omega=\sqcup_{j=1}^k\Omega_j$ (see Subsection~\ref{piecewise_schwarz_subsubsec}). We further assume that $\overline{\Omega}$ is connected and simply connected. This is equivalent to requiring that each $\Omega_j$ is a Jordan domain, and the \emph{contact graph} of the quadrature multi-domain (i.e., a graph having a vertex for each $\Omega_j$ and an edge connecting two vertices if the corresponding quadrature domains touch) is a tree. We refer to such an $\Omega$ as a \emph{tree-like} quadrature multi-domain.
For a tree-like quadrature multi-domain, the desingularized droplet $T^0(\sigma)$ is homeomorphic to an ideal polygon in the hyperbolic plane.

\begin{definition}\label{poly_schwarz_def}
A degree $d$ piecewise Schwarz reflection map $\sigma:\overline{\Omega}\to\widehat{\C}$ associated with a tree-like quadrature multi-domain is said to be \emph{polygonal} if $T^0(\sigma)$ is homeomorphic to an ideal $(d+1)$-gon in $\D$ such that the homeomorphism is conformal on the interior and if $\sigma$ has no critical values in $T^0(\sigma)$. It is called \emph{regular polygonal} if $T^0(\sigma)$ is conformally equivalent to the regular ideal $(d+1)$-gon $\Pi(\pmb{G}_d)$ respecting the vertices. 
\end{definition}

We denote the space of degree $d$ polygonal (respectively, regular polygonal) piecewise Schwarz reflection maps by $\mathscr{S}_d$ (respectively, $\mathscr{S}_d^{\mathrm{reg}}$).

\begin{remark}
1) For $\sigma\in\mathscr{S}_d$, the condition that $T^0(\sigma)$ is homeomorphic to an ideal $(d+1)$-gon is equivalent to saying that the sum of the number of cusps and twice the number of double points on $\partial\Omega$ is equal to $d+1$.

2) The nomenclature `polygonal Schwarz reflection maps' is justified by the observation that such maps are characterized by having the Nielsen map of a polygonal reflection group as the conformal model of their dynamics on the rank zero and rank one tiles (cf.~\cite[\S4.1]{LLM23}).

\end{remark}

\begin{proposition}\cite{LLM23}\label{polygonal_conn_locus_prop}
Let $(\sigma:\overline{\Omega}\to\widehat{\C})\in \mathscr{S}_d$.
Then the following are equivalent.
\noindent\begin{enumerate}\upshape
\item $K(\sigma)$ is connected.
\item $\sigma^{\circ n}(c)\notin T^\infty(\sigma)$, $\forall\ c\in \mathrm{crit}(\sigma)$.
\item $\sigma:T^\infty(\sigma)\setminus\Int{T^0(\sigma)}\longrightarrow T^\infty(\sigma)$ is conformally conjugate to the Nielsen map $\cN_G:\D\setminus \Int{\Pi(G)}\to \D$, where $G$ is an ideal $(d+1)$-gon reflection group.
\end{enumerate}
\end{proposition}

\begin{definition}
We define the {\em connectedness locus} of $\mathscr{S}_d^{\mathrm{reg}}$ as
$$
\mathcal{S}_{\pmb{\cN}_d}:= \{\sigma \in \mathscr{S}_d^{\mathrm{reg}}: K(\sigma) \text{ is connected}\}.
$$ 
\end{definition}

By Proposition~\ref{polygonal_conn_locus_prop}, the connectedness locus $\mathcal{S}_{\pmb{\cN}_d}$ of $\mathscr{S}_d^{\mathrm{reg}}$ consists precisely of maps having the Nielsen map $\pmb{\cN}_d$ as the conformal model of their tiling set dynamics.

The maps in $\cS_{\pmb{\cN}_d}$ are normalized so that the conformal conjugacy $\psi_\sigma:\D\to T^\infty(\sigma)$ between $\pmb{\cN}_d$ and $\sigma$ sends the origin to $\infty$ and has asymptotics $z\mapsto 1/z+O(z)$ as $z\to 0$.

\subsection{Dynamical rays}\label{polygonal_comb_model_subsec}

To obtain a combinatorial model for the limit set of a Schwarz reflection $\sigma\in\cS_{\pmb{\cN}_d}$, one needs the following notion of dynamical rays. These are defined in terms of the Cayley graph of $\pmb{G}_d$ with respect to the generating set $\pmb{\rho}_1,\cdots, \pmb{\rho}_{d+1}$, where $\pmb{\rho}_i$ is the anti-M{\"o}bius reflection in the circle $\pmb{C}_i$ (see Definition~\ref{regular_ideal_polygon_ref_group_def}). 

\begin{definition}\label{Gd_rays}
Let $(i_1, i_2, \cdots)\in\{1,\cdots, d+1\}^{\N}$ with $i_j\neq i_{j+1}$ for all $j$. The corresponding infinite sequence of tiles $\{\Pi(\pmb{G}_d),\rho_{i_1}(\Pi(\pmb{G}_d)),\rho_{i_1}\circ\rho_{i_2}(\Pi(\pmb{G}_d)),\cdots\}$ shrinks to a single point of $\mathbb{S}^1\cong \R/\Z$, which we denote by $\theta(i_1, i_2, \cdots)$. We define a \emph{$\pmb{G}_d$-ray at angle $\theta(i_1, i_2, \cdots)$} to be the concatenation of hyperbolic geodesics (in $\D$) connecting the consecutive points of the sequence $\{0,\rho_{i_1}(0),\rho_{i_1}\circ\rho_{i_2}(0),\cdots\}$. 
\end{definition}
\noindent (See Figure~\ref{deltoid_intro_fig} (right).)
We note that in general there may be more than one $\pmb{G}_d$-rays at a given angle.

\begin{definition}\label{dyn_ray_schwarz}
For $\sigma\in\cS_{\pmb{\cN}_d}$, the image of a $\pmb{G}_d$-ray at angle $\theta$ under the map $\psi_\sigma:\D \longrightarrow T^\infty(\sigma)$ (that conjugates $\pmb{\cN}_d$ to $\sigma$) is called a \emph{$\theta$-dynamical ray of $\sigma$}.
\end{definition}

We  denote the set of all (pre)periodic points of $\pmb{\cN}_d:\mathbb{S}^1\to\mathbb{S}^1$ by $\mathrm{Per}(\pmb{\cN}_d)$. As in the case for anti-polynomials with connected Julia set, the dynamical rays at angles in $\mathrm{Per}(\pmb{\cN}_d)$ land on $\Lambda(\sigma)$. This leads to the following analogue of rational laminations for the space~$\cS_{\pmb{\cN}_d}$.

\begin{definition}\label{def_preper_lami}
The \emph{preperiodic lamination} of $\sigma\in\cS_{\pmb{\cN}_d}$ is defined as the equivalence relation on $\mathrm{Per}(\pmb{\cN}_d)\subset\R/\Z$ such that $\theta, \theta'\in\mathrm{Per}(\pmb{\cN}_d)$ are related if and only if the dynamical rays of $\sigma$ at these angles land at the same point of $\Lambda(\sigma)$. We denote the preperiodic lamination of $\sigma$ by $\lambda(\sigma)$.
\end{definition}

\subsection{Coarse partition and topology}\label{polygonal_topology_subsec}

\begin{figure}[h!]
\captionsetup{width=0.96\linewidth}
\begin{tikzpicture}
\node[anchor=south west,inner sep=0] at (0,0) {\includegraphics[width=0.96\linewidth]{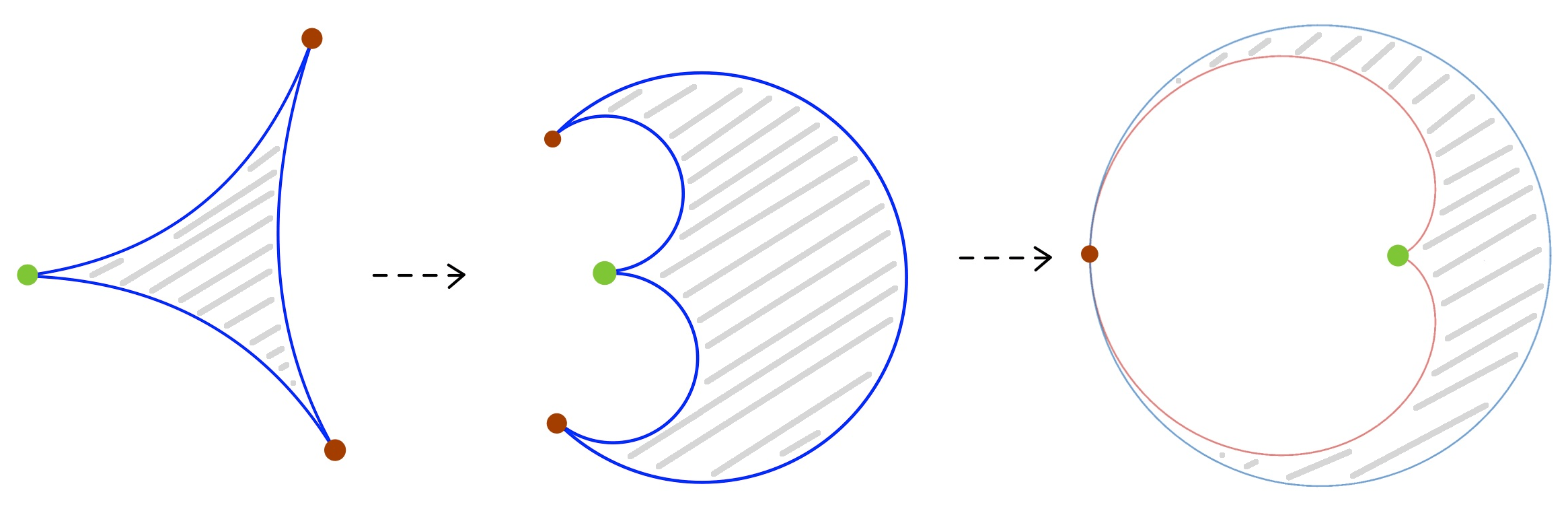}};
\end{tikzpicture}
\caption{Displayed is a schematic picture of a family of quadrature multi-domains defining maps in $\cS_{\pmb{\cN}_2}$. The leftmost quadrature domain (complement of the grey region) is the deltoid. It defines the mating of $\overline{z}^2$ and the Nielsen map $\pmb{\cN}_2$. The two brown cusps are the landing points of the $1/3$ and $2/3$ dynamical rays. The rightmost quadrature multi-domain belongs to the C\&C family, and defines the mating of the fat basilica anti-polynomial $\overline{z}^2-3/4$ and $\pmb{\cN}_2$. The $1/3$ and $2/3$ dynamical rays land at a common point (the touching point of the circle and the cardioid) for this map. The middle quadrature domain (complement of the grey region) defines the mating of some $\overline{z}^2+c$, $c\in(-3/4,0)$, with $\pmb{\cN}_2$.}
\label{deltoid_to_c_and_c_fig}
\end{figure}

Quadrature multi-domains defining the Schwarz reflections in $\cS_{\pmb{\cN}_2}$ come in essentially two different flavors: the deltoid-like quadrature domains, and the union of the cardioid and the exterior of a circumscribing circle. The first kind of maps are characterized by the fact that all of their fixed dynamical rays (at angles $0,1/3$ and $2/3$) land at distinct points, while the second type of maps are characterized by co-landing of their $1/3$ and $2/3$-dynamical rays (see Figure~\ref{deltoid_to_c_and_c_fig}). In higher degrees, there are many possibilities for the structures of the quadrature domains, and a book-keeping of these patterns is necessary to study convergence in the parameter space. These are captured by the following combinatorial data.

A degree $d$ polygonal Schwarz reflection $\sigma$ with connected limit set induces the action of the Nielsen map $\pmb{\cN}_d$ on the ideal boundary $I(T^\infty(\sigma))\cong \R/\Z$ of the escaping set, and hence has exactly $d+1$ fixed points $0,\cdots,\frac{d}{d+1}$ on the ideal boundary. We denote the set of these fixed points by $\mathrm{Fix}(\pmb{\cN}_d)$.
As in Definition~\ref{fixed_ray_lamination_def}, the \emph{fixed ray lamination} of $\sigma\in\cS_{\pmb{\cN}_d}$ is the equivalence relation $\mathfrak{L}_\sigma$ on $\mathrm{Fix}(\pmb{\cN}_d)$ such that $\theta_1\sim\theta_2$ if and only if the dynamical rays of $\sigma$ at angles $\theta_1,\theta_2$ land at the same point of $\Lambda(\sigma)$.

Thus, we have a coarse partition
$$
\cS_{\pmb{\cN}_d}= \bigsqcup_{\mathfrak{L}} \cS_{\pmb{\cN}_d, \mathfrak{L}},
$$
where $\cS_{\pmb{\cN}_d, \mathfrak{L}}$ consists of maps in $\cS_{\pmb{\cN}_d}$ with fixed ray lamination $\mathfrak{L}$, and the union is over all possible fixed ray laminations. Arguments similar to the ones employed in Subsection~\ref{conf_mating_identify_subsec} show the following properties of $(\sigma : \overline{\Omega} \longrightarrow \widehat{\C}) \in \cS_{\pmb{\cN}_d, \mathfrak{L}}$ (see \cite[\S 6.1]{LLM23} for details).
\begin{itemize}
\item The components of $\Omega$ are in one-to-one correspondence with the gaps of $\mathfrak{L}$; i.e.,
$$
\displaystyle\Omega=\bigsqcup_{\substack{\mathcal{G}\\ \textrm{gaps\ of}\ \mathfrak{L}}} \Omega_{\mathcal{G}}.
$$
\item Two quadrature domains $\Omega_{\mathcal{G}}, \Omega_{\mathcal{G}'}$ share at most one boundary point; moreover, they do so if and only if $\mathcal{G}$ and $\mathcal{G}'$ are adjacent.
\item The touching point of two quadrature domains $\Omega_{\mathcal{G}}, \Omega_{\mathcal{G}'}$ is a non-singular point for both $\partial\Omega_{\mathcal{G}}$ and $\partial\Omega_{\mathcal{G}'}$ (i.e., it is a double point of~$\partial\Omega$).
\item The degree of the uniformizing rational map of $\Omega_{\mathcal{G}}$ is equal to the number of arcs of $\mathbb{S}^1\setminus\mathrm{Fix}(\pmb{\cN}_d)$ on $\partial\mathcal{G}\cap\mathbb{S}^1$.
\item The number of cusps on $\partial\Omega_{\mathcal{G}}$ equals the number of points of $\mathrm{Fix}(\pmb{\cN}_d)$ in $\Int{(\partial\mathcal{G}\cap\mathbb{S}^1)}$.
\end{itemize}

We now define a topology on the space $\cS_{\pmb{\cN}_d}$ in terms of convergence of sequences. The definition given below is essentially an adaptation of Carath{\'e}odory convergence where we record the various Carath{\'e}odory limits for book-keeping purposes.
Since there are only finitely many fixed ray laminations, we can assume (after possibly passing to a subsequence) that any sequence in $\cS_{\pmb{\cN}_d}$ lies entirely in $\cS_{\pmb{\cN}_d, \mathfrak{L}}$, where $\mathfrak{L}$ is a particular fixed ray lamination. We enumerate the gaps of $\mathfrak{L}$ as $\mathcal{G}_1,\cdots,\mathcal{G}_k$.
For maps $\sigma:\overline{\Omega}\to\widehat{\C}$ in $\cS_{\pmb{\cN}_d, \mathfrak{L}}$, we denote the corresponding components of $\Omega$ as $\Omega_1,\cdots,\Omega_k$.

\begin{figure}[h!]
\captionsetup{width=0.96\linewidth}
\begin{tikzpicture}
\node[anchor=south west,inner sep=0,rotate=45] at (-0.6,-1.5) {\includegraphics[width=0.24\textwidth]{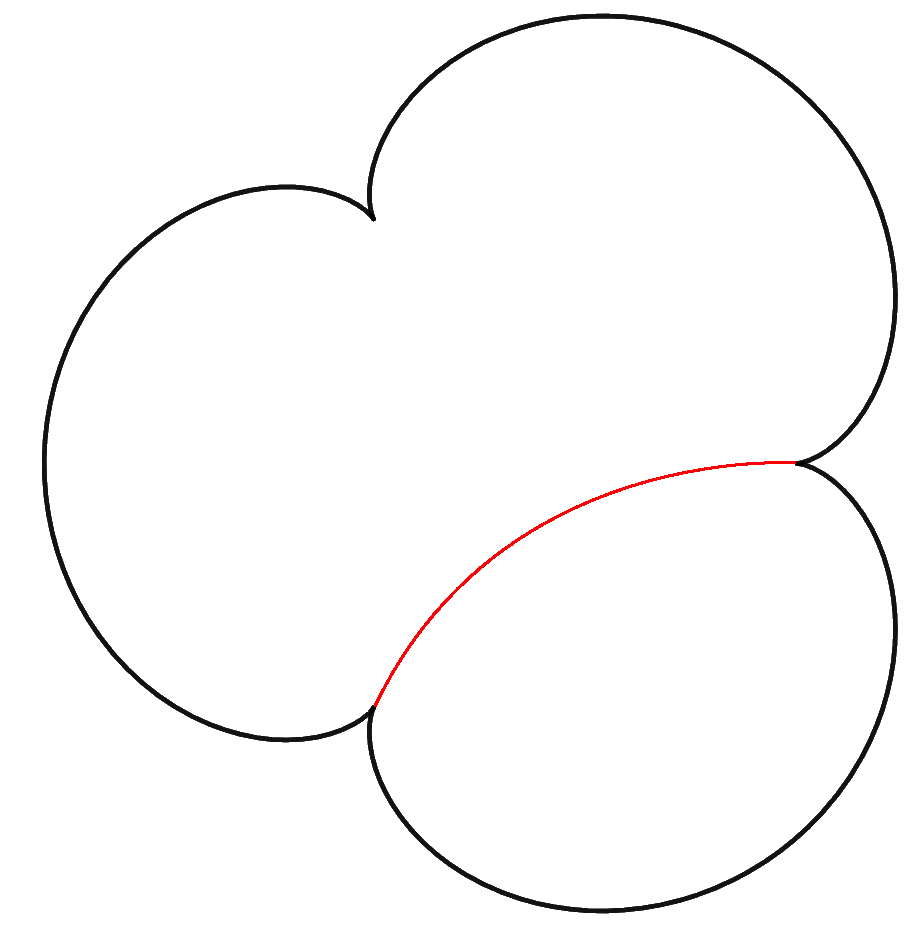}}; 
\node[anchor=south west,inner sep=0,rotate=-40] at (0.58,0.6) {\includegraphics[width=0.24\textwidth]{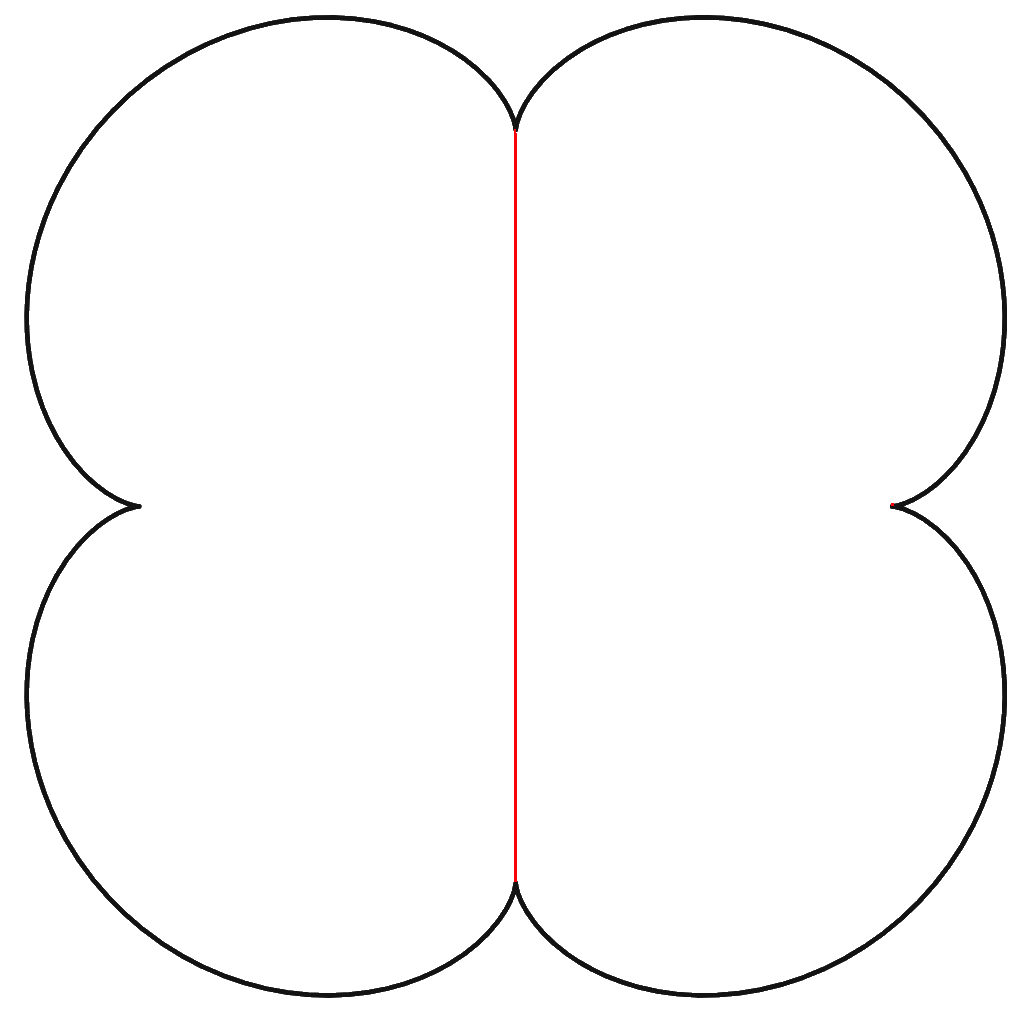}};
\node[anchor=south west,inner sep=0] at (-3.6,-5) {\includegraphics[width=0.4\textwidth]{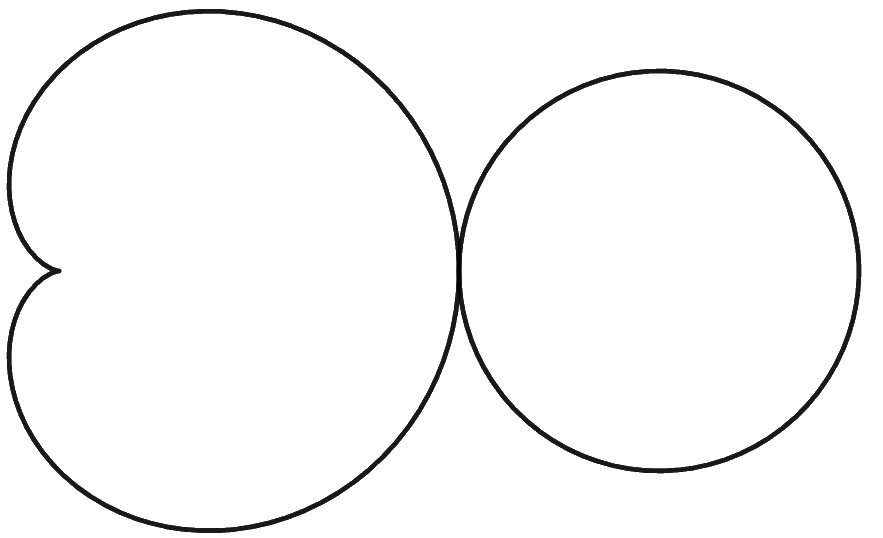}}; 
\node[anchor=south west,inner sep=0,rotate=120] at (4.6,-3.3) {\includegraphics[width=0.2\textwidth]{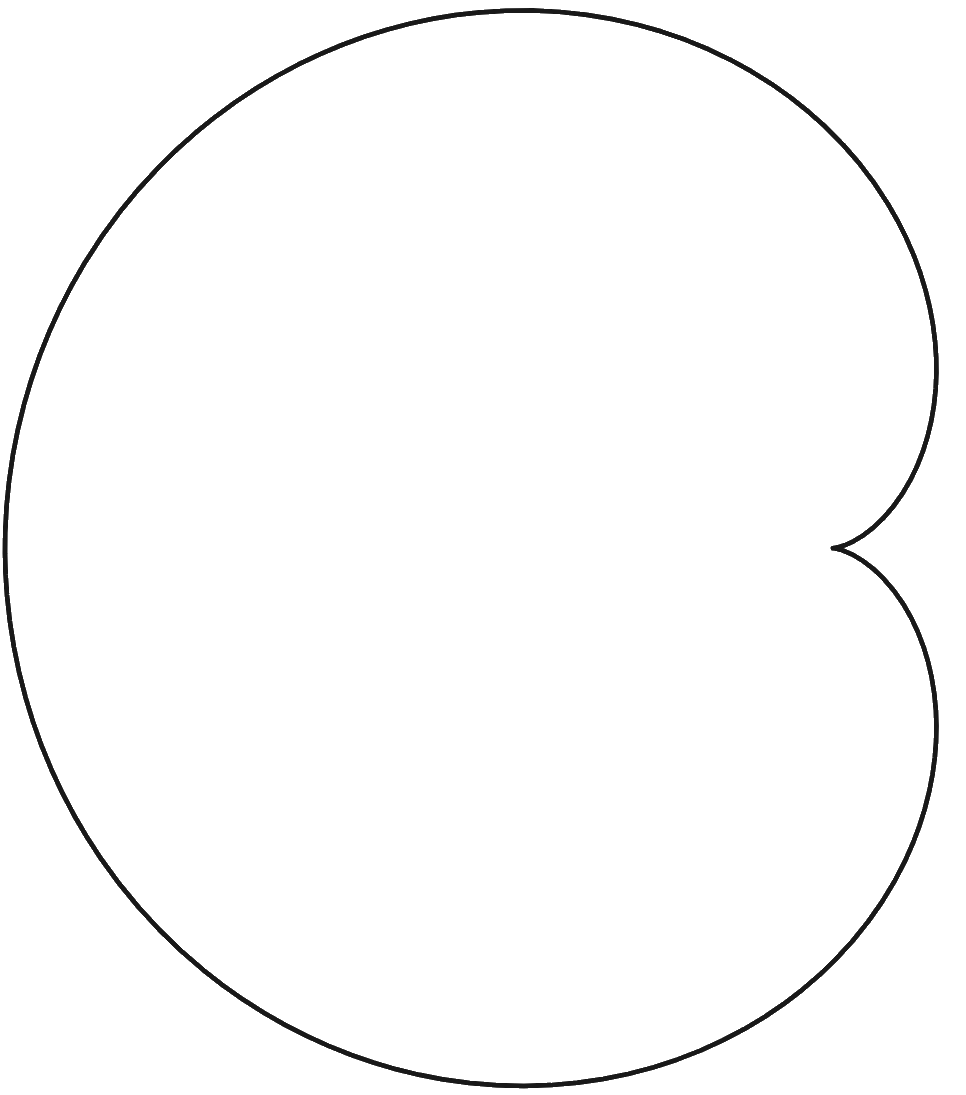}};
\node[anchor=south west,inner sep=0,rotate=-40] at (3.45,-3.2) {\includegraphics[width=0.2\textwidth]{postpinch_2.png}};

\node at (-0.8,0.6) {\begin{large}$\Omega_1^n$\end{large}};
\node at (3,0.1) {\begin{large}$\Omega_2^n$\end{large}};
\node at (-2.3,-3.5) {\begin{large}$\Omega_{1,1}$\end{large}};
\node at (0.2,-3.5) {\begin{large}$\Omega_{1,2}$\end{large}};
\node at (2.7,-3) {\begin{large}$\Omega_{2,1}$\end{large}};
\node at (5.25,-3) {\begin{large}$\Omega_{2,2}$\end{large}};
\end{tikzpicture}
\caption{The sequence of quadrature multi-domains $\{\Omega_1^n,\Omega_2^n\}$ depicted on the top converges to the quadrature multi-domain $\{ \Omega_{1,1},\Omega_{1,2}, \Omega_{2,1}, \Omega_{2,2}\}$ shown at the bottom. Each $\Omega_r^n$, $r\in\{1,2\}$, gets pinched into a pair of quadrature domains $\Omega_{r,1},\Omega_{r,2}$ in the limit. The red line segments in $\Omega_r^n$ indicate which cusps are pinched in the limit.}
\label{pinching_fig}
\end{figure}

\begin{definition}\label{snd_top_def}
We say that a sequence $\displaystyle\{\sigma_n:\overline{\Omega^n}=\bigcup_{r=1}^k\overline{\Omega_{r}^{n}}\to\widehat{\C}\}\subset\cS_{\pmb{\cN}_d, \mathfrak{L}}$ converges to $\displaystyle\sigma:\overline{\Omega}=\bigcup_{r=1}^{k}\bigcup_{j=1}^{l_r}\overline{\Omega_{r,j}}\to\widehat{\C}$ if 
\begin{enumerate}\upshape
\item $\{\Omega_{r,j}:\ j\in\{1,\cdots,l_r\} \}$ is the collection of all Carath{\'e}odory limits of the sequence of domains $\{\Omega_{r}^{n}\}_n$, $r\in\{1,\cdots,k\}$, and 

\item for $n$ large enough, the antiholomorphic maps $\sigma_n$ converge uniformly to $\sigma$ on compact subsets of $\Omega$.
\end{enumerate}
\end{definition}
\noindent (See Figures~\ref{pinching_fig},~\ref{deltoid_to_c_and_c_fig}.)

The pieces $\cS_{\pmb{\cN}_d, \mathfrak{L}}$ are not necessarily closed. For instance, matings of $\pmb{\cN}_2$ with appropriately chosen quadratic anti-polynomials in the principal hyperbolic component of the Tricorn (each of which arises from a single deltoid-like quadrature domain) can converge to the mating of $\pmb{\cN}_2$ with the fat Basilica anti-polynomial $\overline{z}^2-3/4$ (which lies in the C\&C family). This is illustrated in Figure~\ref{deltoid_to_c_and_c_fig}. In general, if a sequence in $\cS_{\pmb{\cN}_d, \mathfrak{L}_1}$ converges to some map in  $\cS_{\pmb{\cN}_d, \mathfrak{L}_2}$, then the latter fixed ray lamination must be stronger than the former; i.e., $\mathfrak{L}_1\subset\mathfrak{L}_2$ (see \cite[\S 6]{LLM23}).

The total space $\cS_{\pmb{\cN}_d}$, however, turns out to be compact. This is essentially a consequence of the fact that all maps in this space have the same external class.

\begin{proposition}\cite{LLM23}\label{polygonal_compact_prop}
$\cS_{\pmb{\cN}_d}$ is compact.
\end{proposition}

\subsection{Relation between geometrically finite maps}\label{polygonal_geom_fin_subsec}
The first result on the intimate relations between the spaces $\cS_{\pmb{\cN}_d}$ and $\mathscr{C}_d$ asserts that there is a natural bijection between their geometrically finite parameters. We remark that the $d=2$ case of this result was proved in \cite{LLMM2} using parameter tessellation of the escape locus of the C\&C family and rigidity of geometrically finite maps (see Subsections~\ref{c_and_c_straightening_mating_subsubsec} and~\ref{c_and_c_geom_fin_bijection_subsubsec}). It turns out that David surgery techniques and conformal removability results allow one to simplify the above arguments, and this strategy was adopted in \cite{LLM23} to handle the general case.

\begin{definition}\label{geom_fin_polygonal_def}
A Schwarz reflection map $\sigma\in\cS_{\pmb{\cN}_d}$ is said to be \emph{geometrically finite} if every critical point of $\sigma$ in the
limit set $\Lambda(\sigma)$ is preperiodic.
\end{definition}

By the classification of Fatou components of maps in $\cS_{\pmb{\cN}_d}$ (i.e., components of $\Int{K(\sigma)}$) and their relations with critical points, each periodic Fatou component of a geometrically finite map in $\cS_{\pmb{\cN}_d}$ is an attracting/parabolic basin, or a basin of attraction of some singular point on $\partial T(\sigma)$. Moreover, such a map has no Cremer cycle. 

Let us denote by $\cS_{\pmb{\cN}_d}^{gf}, \mathscr{C}_d^{gf}$ the collection of geometrically finite maps in $\cS_{\pmb{\cN}_d}, \mathscr{C}_d$ (respectively).

\begin{theorem}\cite{LLM23}\label{geom_finite_bij_polygonal_thm}
There is a bijection 
$$
\Phi: \mathscr{C}_{d}^{gf} \longrightarrow \cS_{\pmb{\cN}_d}^{gf},
$$ 
such that for each $f\in \mathscr{C}_{d}^{gf}$, the corresponding Schwarz reflection $\Phi(f)$ is a conformal mating of $f$ with $\pmb{\cN}_d$.
However, both $\Phi$ and $\Phi^{-1}$ are discontinuous.
\end{theorem}

\begin{figure}[h!]
\captionsetup{width=0.96\linewidth}
\includegraphics[width=0.45\linewidth]{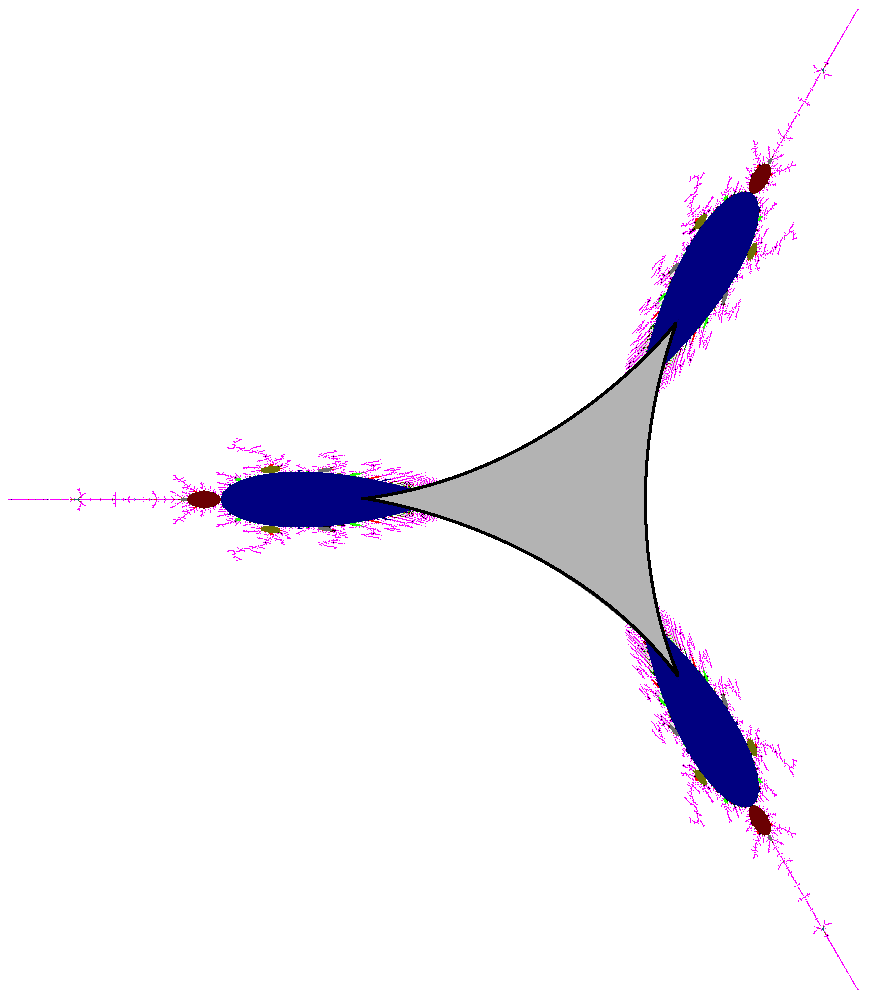}\hspace{8mm} \includegraphics[width=0.45\linewidth]{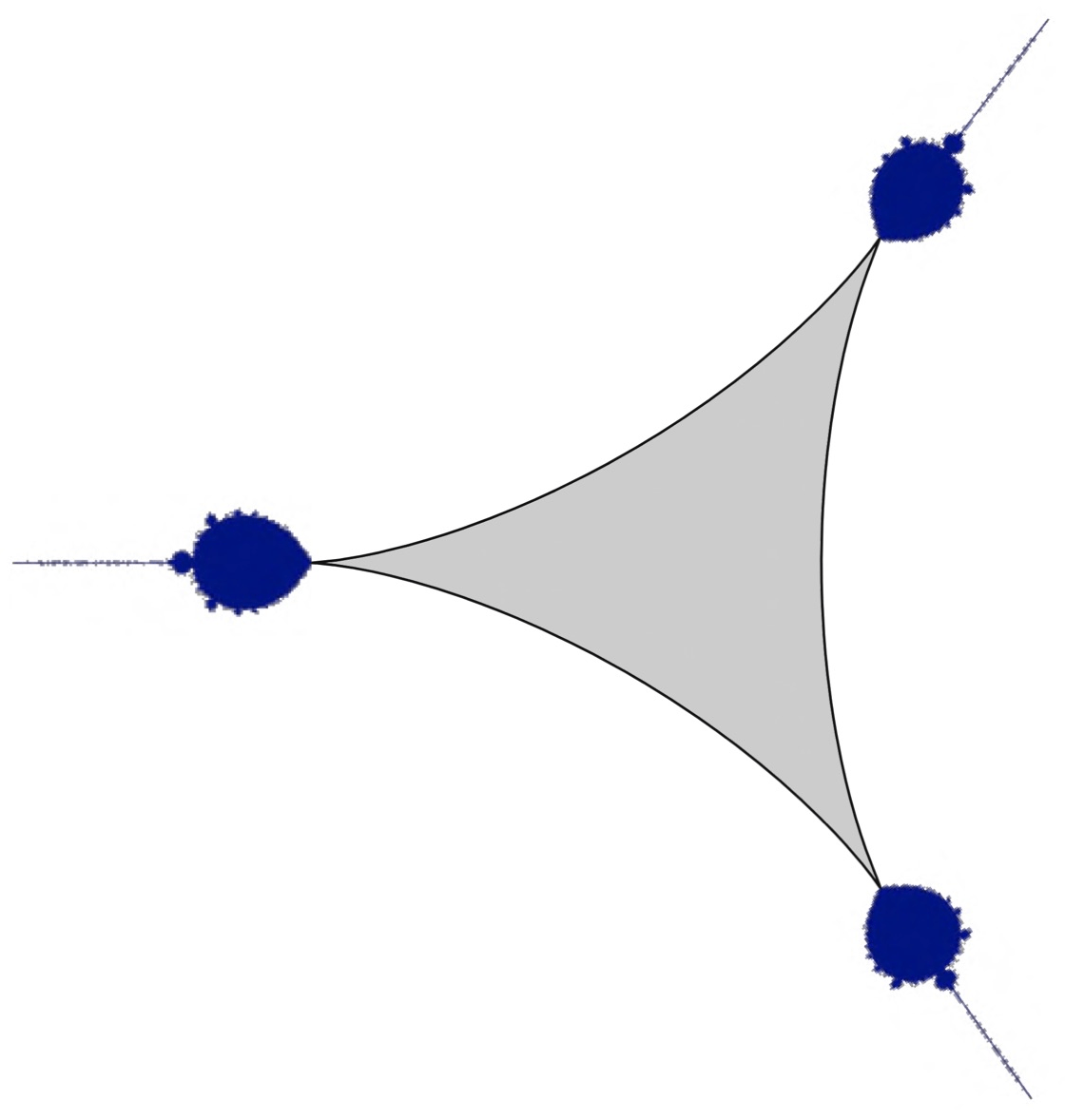}
\caption{The Tricorn $\mathscr{C}_2$ is displayed on the left, and a schematic picture of $\cS_{\pmb{\cN}_2}$ is depicted on the right. In both pictures, the central grey region is the period one hyperbolic component and the three symmetric blue bulbs are the period two hyperbolic components. The bifurcation structures of the period two hyperbolic components from the period one component are different in the two connectedness loci, and this causes discontinuity of $\Phi, \Phi^{-1}$.}
\label{discont_fig}
\end{figure}

The proof of the above result can essentially be split into the following parts. 
\begin{enumerate}\upshape
\item \textbf{Definition of the map $\Phi$}: The David surgery of Lemma~\ref{david_surgery_lemma} and its higher period variants allow one to turn geometrically finite anti-polynomials $f$ into geometrically finite Schwarz reflections with $\pmb{\cN}_d$ as their external class. Subsequently, conformal removability of limit sets of geometrically finite maps in $\cS_{\pmb{\cN}_d}$ imply uniqueness of the conformal mating of $f$ and $\pmb{\cN}_d$. This defines the map $\Phi$.
\item \textbf{Injectivity of $\Phi$}: This is a consequence of conformal removability of Julia sets of geometrically finite anti-polynomials with connected Julia set (see Section~\ref{conf_removable_subsec}). 
\item \textbf{Surjectivity of $\Phi$}: This is based on the combinatorial structure of preperiodic laminations of geometrically finite maps in $\cS_{\pmb{\cN}_d}$. Specifically, the proof uses 
\begin{enumerate}
\item the fact that the push-forward of the preperiodic lamination of a map in $\cS_{\pmb{\cN}_d}^{gf}$ under $\pmb{\mathcal{E}}_d$ satisfies the properties of rational laminations listed in Proposition~\ref{rat_lami_prop} (see Definition~\ref{push_lami_def} for the notion of push-forward of laminations), and
\item the existence of anti-polynomials with prescribed rational laminations.
\end{enumerate}
\end{enumerate}

It is worth pointing out that the source of discontinuity in the above theorem is quasiconformally deformable parabolic parameters (the same phenomenon causes discontinuity for other straightening maps in antiholomorphic dynamics too, cf. \cite{IM1}).

\subsection{Relation between periodically repelling maps}\label{polygonal_per_rep_subsec}
An anti-polynomial is called \emph{periodically repelling} if all of its periodic points in the plane are repelling. On the other hand, we say that a Schwarz reflection $\sigma\in\cS_{\pmb{\cN}_d}$ is \emph{periodically repelling} if no singular point of its droplet has an attracting direction in the non-escaping set $K(\sigma)$ and if all other periodic points of $\sigma$ are repelling.
The subspaces of $\cS_{\pmb{\cN}_d}, \mathscr{C}_d$ consisting of periodically repelling maps are denoted by $\cS_{\pmb{\cN}_d}^{r}, \mathscr{C}_{d}^{r}$ (respectively).

By Theorem~\ref{geom_finite_bij_polygonal_thm}, there exists a bijection between post-critically finite maps of $\mathscr{C}_{d}^{r}$ and $\cS_{\pmb{\cN}_d}^{r}$. Thanks to certain continuity properties for rational laminations of maps in $\mathscr{C}_{d}^{r}$ and preperiodic laminations of maps in $\cS_{\pmb{\cN}_d}^{r}$, one can extend the above bijection to a homeomorphism between periodically repelling combinatorial classes of $\mathscr{C}_{d}^{r}$ and $\cS_{\pmb{\cN}_d}^{r}$. {

\begin{theorem}\cite{LLM23}\label{comb_class_homeo_polygonal_thm}
There is a homeomorphism 
$$
\Phi: \faktor{\mathscr{C}_{d}^{r}}{\sim} \longrightarrow \faktor{\cS_{\pmb{\cN}_d}^{r}}{\sim}\ ,
$$
where $\sim$ identify maps with the same rational/preperiodic lamination.
\end{theorem}

Combining Theorem~\ref{comb_class_homeo_polygonal_thm} with combinatorial rigidity of at most finitely renormalizable maps in $\mathscr{C}_{d}^{r}$ and $\cS_{\pmb{\cN}_d}^{r}$ (i.e., such maps are uniquely determined by their rational/preperiodic laminations), we obtain the following corollary. We denote the collection of at most finitely renormalizable maps in $\mathscr{C}_{d}^{r}, \cS_{\pmb{\cN}_d}^{r}$ by $\mathscr{C}_{d}^{fr}, \cS_{\pmb{\cN}_d}^{fr}$, respectively.

\begin{corollary}\cite{LLM23}\label{polygonal_fin_renorm_cor}
There is a homeomorphism 
$$
\Phi: \mathscr{C}_{d}^{fr} \longrightarrow \cS_{\pmb{\cN}_d}^{fr}
$$
such that for each $f\in \mathscr{C}_{d}^{fr}$, the corresponding Schwarz reflection $\Phi(f)$ is a conformal mating of $f$ with $\pmb{\cN}_d$.
\end{corollary}

\section{Correspondences as matings: systematic theory via Schwarz dynamics}\label{general_mating_corr_sec}

In Subsection~\ref{chebyshev_gen_subsec}, we described the dynamics and parameter space of a family of $2$:$2$ antiholomorphic correspondences that arise as matings of quadratic parabolic anti-rational maps and the antiholomorphic analog $\mathbbm{G}_2$ of the modular group.
Such families, in the holomorphic setting \cite{BP,BuLo1,BuLo2,BuLo3} as well as in the antiholomorphic setting \cite{LLMM3}, were constructed by looking at explicit algebraic correspondences of bi-degree $2$:$2$. While this allows for a complete understanding of the dynamics and parameter planes of these correspondences, a shortcoming of this approach is that one is forced to rely heavily on the structure of low-dimensional parameter spaces (real two-dimensional) to realize matings as correspondences. 

In general, it seems hard to come up with explicit algebraic correspondences of bi-degree $d$:$d$ ($d\geq 2$) directly that exhibit similar mating phenomena. However, the intimate connection between antiholomorphic correspondences and Schwarz reflection maps (expounded in Subsections~\ref{chebyshev_subsec},~\ref{chebyshev_gen_subsec}) suggests that the task of constructing  antiholomorphic correspondences with prescribed hybrid dynamical behavior is strongly related to manufacturing Schwarz reflection maps with suitable dynamical properties.
In \cite{LMM3}, this approach was successfully adopted to generalize the family of $2$:$2$ antiholomorphic correspondences arising from cubic Chebyshev polynomials to arbitrary degree. We will now collect the main results of that paper, which settled a suitably modified version of a conjecture by Bullett and Freiberger in the antiholomorphic setting \cite[\S 3, p. 3926]{BuFr}. We will also touch briefly upon a slightly more general construction of antiholomorphic correspondences given in \cite{LLM23}.

We recall that the space $\pmb{\mathcal{B}}_d$ of parabolic anti-rational maps was introduced in Subsection~\ref{para_anti_rat_gen_subsubsec}. By definition, the \emph{anti-Hecke group} $\mathbbm{G}_d$ is a group of conformal and anti-conformal automorphisms of $\D$ generated by the rigid rotation by angle $\frac{2\pi}{d+1}$ around the origin and the reflection in the hyperbolic geodesic of $\D$ connecting $1$ to $e^{\frac{2\pi i}{d+1}}$. It is an index $d+1$ extension of the ideal $(d+1)$-gon reflection group $\pmb{G}_d$ (cf. Subsection~\ref{nielsen_first_return_external_map_subsubsec} and \cite[\S 3.1]{LMM3}). The name anti-Hecke is justified by the observation that replacing the reflection map with a M{\"o}bius inversion (fixing the geodesic connecting $1$ to $e^{\frac{2\pi i}{d+1}}$) in the above definition of $\mathbbm{G}_d$ yields the standard Hecke group (cf. \cite[\S 11.3]{Bea95}).

\begin{definition}\label{corr_mating_defn}
Let $R\in\pmb{\mathcal{B}}_d$ (respectively, let $p\in\mathscr{C}_d$) and let $G$ be a discrete subgroup of the group of conformal and anti-conformal automorphisms of the unit disk $\D$.
Let $\mathfrak{C}$ be an anti-holomorphic correspondence on a compact, simply connected (possibly noded) Riemann surface $\mathfrak{W}$.

We say that $\mathfrak{C}$ a {\em mating} of $R$ (respectively, $p$) and $G$ if there is a $\mathfrak{C}-$invariant partition $\mathfrak{W}=\mathcal{T}\sqcup\mathcal{K}$ such that the following hold.
\begin{enumerate}
	\item On $\cT$, the dynamics of $\mathfrak{C}$ is equivalent to the action of a group of (anti-)conformal automorphisms acting properly discontinuously. Further, $\cT/\mathfrak{C}$ is biholomorphic to~$\D/G$.
	
	\item $\cK$ can be written as the union of two copies $\cK_1, \cK_2$ of $\cK(R)$ (respectively, of $\cK(p)$), where $\cK(R)$ (respectively, $\cK(p)$) is the filled Julia set of $R$ (respectively, of $p$), such that $\cK_1$ and $\cK_2$ intersect at finitely many points. Furthermore, $\mathfrak{C}$ has a forward (respectively, backward) branch carrying $\cK_1$ (respectively, $\cK_2$) onto itself with degree $d$, and this branch is conformally (respectively, anti-conformally) conjugate to $R\vert_{\mathcal{K}(R)}$ or $p\vert_{\mathcal{K}(p)}$.\end{enumerate}
\end{definition}

We denote by $\mathfrak{U}_{d+1}$ the space of degree $d+1$ polynomials $f$ such that $f\vert_{\overline{\D}}$ is injective and $f$ has a unique (non-degenerate) critical point on $\mathbb{S}^1$.

\begin{theorem}\cite[Theorem~A]{LMM3}\label{general_mating_corr_existence_thm}
Let $R\in\pmb{\mathcal{B}}_d$. Then, there exists $f\in\mathfrak{U}_{d+1}$ such that the associated reversible antiholomorphic correspondence $\mathfrak{C}$ on $\widehat{\C}$ defined as
\begin{equation}
(z,w)\in\mathfrak{C}\iff \frac{f(w)-f(\eta(z))}{w-\eta(z)}=0
\label{corr_eqn}
\end{equation}
is the mating of the anti-Hecke group $\mathbbm{G}_d$ and $R$. 

\noindent Moreover, this mating operation yields a bijection between $\ \faktor{\pmb{\mathcal{B}}_d}{\mathrm{Aut}(\C)}$ and the connectedness locus of the moduli space of antiholomorphic correspondences generated by deck transformations of polynomials $f\in\mathfrak{U}_{d+1}$ and reflection in the unit disk.
\end{theorem}

\subsection{Sketch of proof of the realization theorem}\label{gen_mating_corr_existence_proof_subsec}

In order to motivate the proof of Theorem~\ref{general_mating_corr_existence_thm}, let us recall that in bi-degree $2$:$2$, the dynamical properties of the correspondences were obtained essentially from the parallel study of the dynamics of the associated Schwarz reflection maps. 

\subsubsection{Motivation from the cubic Chebyshev family}\label{learning_from_examples_subsubsec}

The following two relations between Schwarz reflections and correspondences described in Subsections~\ref{chebyshev_subsec},~\ref{chebyshev_gen_subsec} are of particular importance.
\smallskip

\noindent$\bullet$ The fact that the conformal class of the associated Schwarz reflections on their tiling sets is given by the map $\pmb{\cF}_2$ was instrumental in the $\mathbbm{G}_2$-structure of the correspondences on their lifted tiling sets (see Subsections~\ref{nielsen_first_return_external_map_subsubsec} and~\ref{chebyshev_center_corr_subsubsec}).
\smallskip

\noindent$\bullet$ The fact that the non-escaping set dynamics of the associated Schwarz reflections are hybrid conjugate to the actions of anti-rational maps in $\pmb{\mathcal{B}}_2$ on their filled Julia sets implied that suitable branches of the correspondences on their lifted tiling sets are hybrid conjugate to anti-rational maps.
\smallskip

Guided by this analogy, we will proceed to construct a space of Schwarz reflections exhibiting the above two features. 

\subsubsection{The degree $d$ anti-Farey map}\label{anti_farey_subsubsec}

To implement the above strategy, a degree $d$ generalization of the map $\pmb{\cF}_2$ was given in \cite[\S 3.1]{LMM3}. Succinctly, as the Nielsen map $\pmb{\cN}_d$ of the regular ideal $(d+1)$-gon reflection group $\pmb{G}_d$ commutes with $M_\omega(z)=\omega z$ (where $\omega:=e^{\frac{2\pi i}{d+1}}$), the quotient map $\overline{\D}\to\faktor{\overline{\D}}{\langle M_\omega\rangle}$ semi-conjugates $\pmb{\cN}_d:\overline{\D}\setminus \Int{\Pi(\pmb{G}_d)}\to\overline{\D}$ to a well-defined factor map 
$$
\pmb{\cF}_d:\faktor{\left(\overline{\D}\setminus \Int{\Pi(\pmb{G}_d)}\right)}{\langle M_\omega\rangle}\longrightarrow \faktor{\overline{\D}}{\langle M_\omega\rangle}.
$$
\begin{figure}[h!]
\captionsetup{width=0.96\linewidth}
\begin{tikzpicture}
\node[anchor=south west,inner sep=0] at (0,0) {\includegraphics[width=0.41\textwidth]{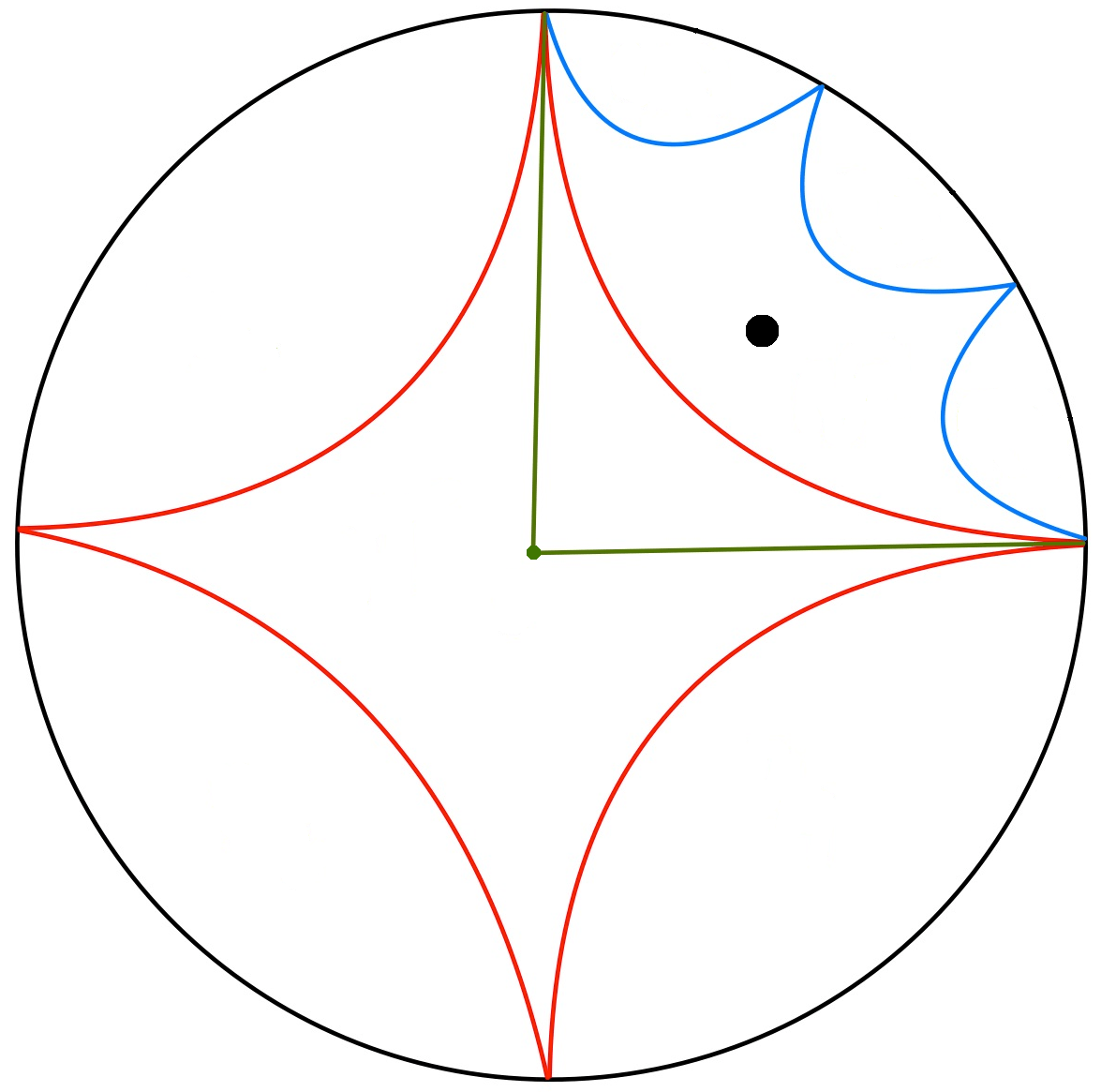}}; 
\node[anchor=south west,inner sep=0] at (6,0) {\includegraphics[width=0.4\textwidth]{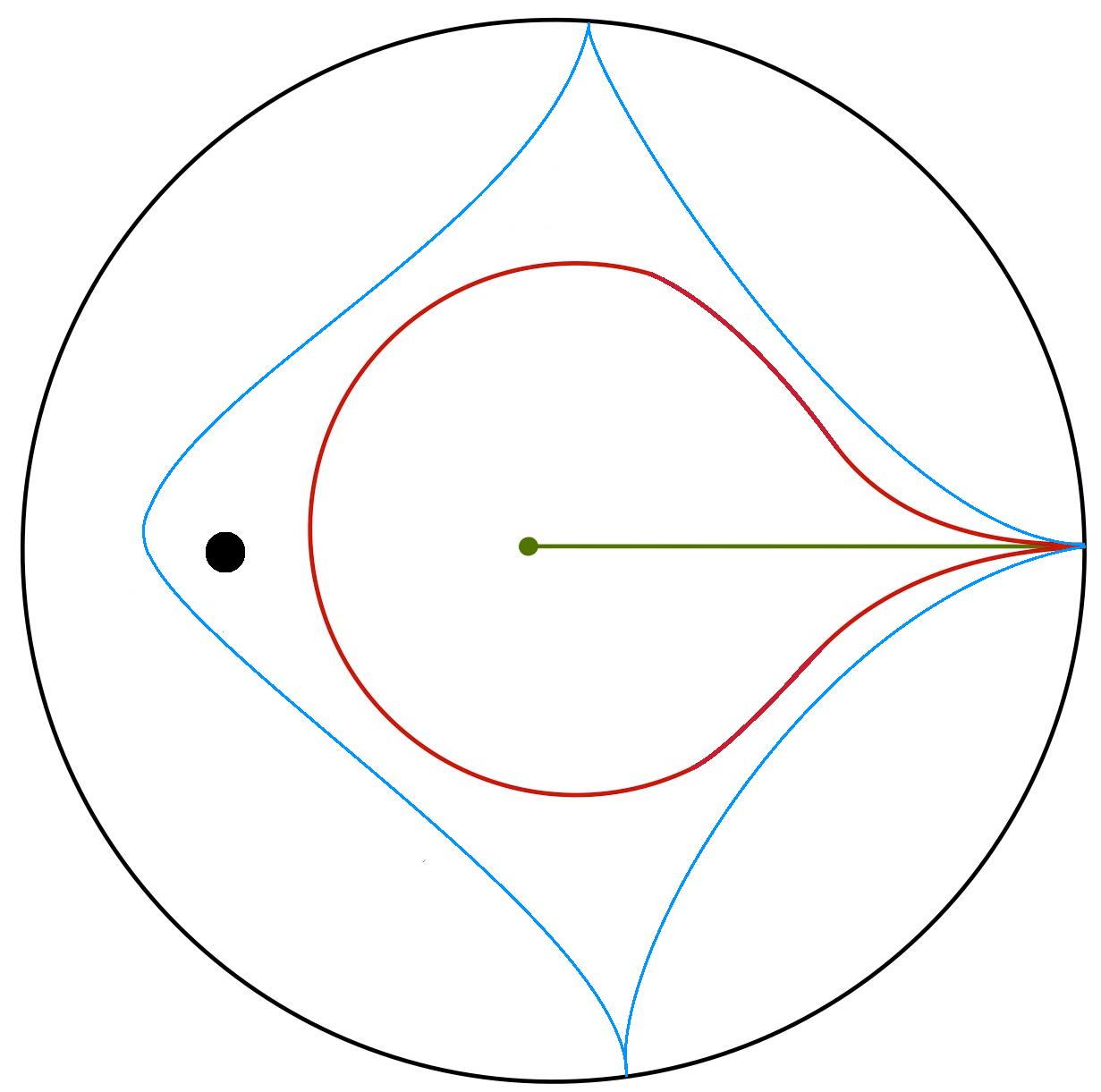}};
\end{tikzpicture}
\caption{Left: the fundamental domain $\Pi(\pmb{G}_3)\cap\D$ is the regular quadrilateral bounded by the red geodesics. A fundamental domain for the action of $\mathbbm{G}_3$ on $\D$ is given by the triangle having the two green lines and the red geodesic connecting $1$ to $i$ as its edges. Right: When the Riemannian orbifold $\faktor{\D}{\langle M_i\rangle}$ is uniformized by the disk, the map $\pmb{\cF}_3$ has the region bounded by the black unit circle and the red monogon as its domain of definition. The map $\pmb{\cF}_3$ fixes the red monogon pointwise, and acts as an orientation-reversing, degree three circle covering with a unique parabolic fixed point at $1$. The heavy black dot represents the triple critical point of $\pmb{\cF}_3$.}
\label{anti_farey_fig}
\end{figure}
\noindent This map $\pmb{\cF}_d$, which coincides with the map $\pmb{\cF}_2$ defined in Subsection~\ref{nielsen_first_return_external_map_subsubsec} for $d=2$, is called the \emph{degree $d$ anti-Farey map}. We note that $\pmb{\cF}_d$ has a unique critical point (of multiplicity $d$) with associated critical value at $0$ (shown in black in Figure~\ref{anti_farey_fig} (right)). A crucial feature of $\pmb{\cF}_d$ is that it acts as the identity map on the inner boundary of its domain of definition (shown in red in Figure~\ref{anti_farey_fig} (right)). Moreover, when the Riemannian orbifold $\faktor{\D}{\langle M_\omega\rangle}$ is uniformized by the unit disk, the restriction of $\pmb{\cF}_d$ on the outer boundary of its domain of definition becomes an expansive, orientation-reversing, degree $d$ circle covering with a unique parabolic fixed point.

\subsubsection{First step: constructing Schwarz reflections}\label{core_step_subsubsec}

Roughly speaking, the first step in the construction of an antiholomorphic correspondence that is a mating of $\mathbbm{G}_d$ and $R\in\pmb{\mathcal{B}}_d$ is to cook up a Schwarz reflection map $\sigma$ having $\pmb{\cF}_d$ as the conformal model of its tiling set dynamics and $R$ as the hybrid class of its non-escaping set dynamics. This is achieved by a careful \emph{quasiconformal} surgery procedure that glues in $\pmb{\cF}_d$ outside the filled Julia set of $R\in\pmb{\mathcal{B}}_d$ (see \cite[Theorem~5.1]{LMM3}). Specifically, such a surgery is facilitated by \cite[Lemma~4.10]{LMM3}, which states that the external dynamics $B_d$ of maps in $\pmb{\mathcal{B}}_d$ is quasiconformally compatible with the map $\pmb{\cF}_d$ that one needs to insert (see Subsection~\ref{para_anti_rat_gen_subsubsec} for the definition of $B_d$). The fact that the resulting hybrid conformal dynamical system $\sigma$ is indeed given by a Schwarz reflection map follows from the triviality of the action of $\pmb{\cF}_d$ on part of the boundary of its domain of definition.

Following \cite[Definition~3.4]{LMM3}, let us denote by $\mathcal{S}_{\pmb{\cF}_d}$ the space of Schwarz reflection maps $\sigma:\overline{\Omega}\to\widehat{\C}$ such that $\Omega$ is a Jordan quadrature domain and the tiling set dynamics of $\sigma$ is conformally conjugate to $\pmb{\cF}_d$. The above discussion can be summarized as follows.

\begin{proposition}\label{anti_farey_para_rat_mating_prop}
Let $R\in\pmb{\mathcal{B}}_d$. Then, there exist $\sigma:\overline{\Omega}\to\widehat{\C}$ in $\mathcal{S}_{\pmb{\cF}_d}$ and a global quasiconformal homeomorphism $\mathfrak{H}$ that is conformal on $\mathcal{K}(R)$,
such that $\mathfrak{H}$ conjugates $R\vert_{\mathcal{K}(R)}$ to $\sigma\vert_{K(\sigma)}$.
\end{proposition}

\subsubsection{Second step: constructing correspondences from Schwarz reflections}\label{upgrade_step_subsubsec}

It turns out that the Schwarz reflection $\sigma$ of Proposition~\ref{anti_farey_para_rat_mating_prop} arises from a Jordan quadrature domain $\Omega$ whose uniformizing map can be chosen to be a polynomial $f\in\mathfrak{U}_{d+1}$. This essentially follows from the existence of a multiplicity $d$ critical point of $\sigma$ in its tiling set (recall that $\pmb{\cF}_d$ has such a critical point). Moreover, this property characterizes the space $\cS_{\pmb{\cF}_d}$ of Schwarz reflections.

\begin{proposition}\cite[Proposition 3.3]{LMM3}\label{mating_equiv_cond_prop}
Let $f$ be a rational map of degree $d+1$ that is injective on $\overline{\D}$, $\Omega:=f(\D)$, and $\sigma$ the Schwarz reflection map associated with $\Omega$. Then the following are equivalent.
\begin{enumerate}\upshape
\item There exists a conformal conjugacy $\psi$ between 
$$
\quad \pmb{\cF}_d:\faktor{\left(\D\setminus \Int{\Pi(\pmb{G}_d)}\right)}{\langle M_\omega\rangle}\longrightarrow \faktor{\D}{\langle M_\omega\rangle}\quad \mathrm{and}\quad \sigma:T^\infty(\sigma)\setminus\Int{T^0(\sigma)}\longrightarrow T^\infty(\sigma).
$$
In particular, $T^\infty(\sigma)$ is simply connected.

\item After possibly conjugating $\sigma$ by a M{\"o}bius map and pre-composing $f$ with an element of $\mathrm{Aut}(\D)$, the uniformizing map $f$ can be chosen to be a polynomial with a unique critical point on $\mathbb{S}^1$. Moreover, $K(\sigma)$ is connected. 

\item $\Omega$ is a Jordan domain with a unique conformal cusp on its boundary. Moreover, $\sigma$ has a unique critical point in its tiling set $T^\infty(\sigma)$, and this critical point maps to $\Int{T^0(\sigma)}$ with local degree $d+1$.
\end{enumerate}
\end{proposition}

The polynomial $f\in\mathfrak{U}_{d+1}$ giving rise to the Schwarz reflection $\sigma$ (produced in Proposition~\ref{anti_farey_para_rat_mating_prop}) is precisely the one that appears in the statement of Theorem~\ref{general_mating_corr_existence_thm}. Indeed, the arguments of Subsection~\ref{chebyshev_center_corr_subsubsec}, combined with the dynamical properties of $\sigma$, apply mutatis mutandis to the present setting and imply that the antiholomorphic correspondence $\mathfrak{C}$ defined by Equation~\eqref{corr_eqn}
\smallskip

\noindent$\bullet$ is equivalent to the action of the anti-Hecke group $\mathbbm{G}_d$ on its lifted tiling set $f^{-1}(T^\infty(\sigma))$, and
\smallskip

\noindent$\bullet$ has a $d:1$ forward branch on `half' of its lifted non-escaping set (i.e., on $f^{-1}(K(\sigma))\cap\overline{\D}$) that is hybrid conjugate to $R$.
\smallskip

We refer the reader to \cite[Propositions~2.18,~3.12]{LMM3} for details.

\subsubsection{The bijection statement of Theorem~\ref{general_mating_corr_existence_thm}}\label{corr_para_anti_rat_bijection_sketch_subsubsec}

By Proposition~\ref{mating_equiv_cond_prop}, the space $\mathcal{S}_{\pmb{\cF}_d}$ is precisely the space of Schwarz reflections $\sigma:\overline{\Omega}\to\widehat{\C}$ with connected non-escaping set $K(\sigma)$ associated with Jordan quadrature domains $\Omega$ such that the Riemann uniformization of $\Omega$ can be chosen to be a polynomial in $\mathfrak{U}_{d+1}$.

In light of the defining equation~\eqref{corr_eqn} of the correspondences, it now follows that the connectedness locus of the space of antiholomorphic correspondences generated by deck transformations of polynomials in $\mathfrak{U}_{d+1}$ and the reflection map $\eta$ can be identified with the space $\faktor{\mathcal{S}_{\pmb{\cF}_d}}{\mathrm{PSL}_2(\C)}$. Here, we used the facts that
\smallskip

\noindent$\bullet$ the local branches $z\mapsto w$ of the correspondences defined by Equation~\eqref{corr_eqn} are given by $f^{-1}\circ f\circ\eta$ (i.e., composition of $\eta$ with local deck transformations of $f$), and
\smallskip

\noindent$\bullet$ M{\"o}bius conjugating a Schwarz reflection map $\sigma$ amounts to post-composing the Riemann uniformization of $\Omega$ with the same M{\"o}bius map, which leaves the associated correspondence unaltered.
\smallskip

According to \cite[Lemma~4.4]{LMM3}, M{\"o}bius conjugacy classes of maps in $\pmb{\mathcal{B}}_d$ and $\cS_{\pmb{\cF}_d}$ are completely determined by their hybrid classes (this is a standard consequence of the fact that all maps in these spaces have the same external class). Therefore, the quasiconformal surgery (or mating) operation explicated in Subsection~\ref{core_step_subsubsec} gives rise to a well-defined map from $\faktor{\pmb{\mathcal{B}}_d}{\mathrm{Aut}(\C)}$ to $\faktor{\mathcal{S}_{\pmb{\cF}_d}}{\mathrm{PSL}_2(\C)}$. For the bijection statement of Theorem~\ref{general_mating_corr_existence_thm}, one needs to argue that this map admits an inverse. This is the content of \cite[Theorem~4.12]{LMM3}, that reverses the construction of \cite[Theorem~5.1]{LMM3}. More precisely, one can use the quasiconformal compatibility of the anti-Farey map $\pmb{\cF}_d$ and the anti-Blaschke product $B_d$ (\cite[Lemma~4.10]{LMM3}) to start with a Schwarz reflection map $\sigma\in\faktor{\cS_{\pmb{\cF}_d}}{\mathrm{PSL}_2(\C)}$ and glue the anti-Blaschke product $B_d$ outside its non-escaping set. This produces a well-defined parabolic anti-rational map $R_\sigma\in\faktor{\pmb{\mathcal{B}}_d}{\mathrm{Aut}(\C)}$ that is hybrid conjugate to $\sigma$. Clearly, the association $\sigma\mapsto R_\sigma$ is the desired inverse map. Thus, we have a bijection
$$
\chi:\ \faktor{\mathcal{S}_{\pmb{\cF}_d}}{\mathrm{PSL}_2(\C)}\ \longrightarrow\  \faktor{\pmb{\mathcal{B}}_d}{\mathrm{Aut}(\C)},\quad [\sigma] \mapsto [R_\sigma].
$$ 
The map $\chi$ is called the \emph{straightening map}.

\subsubsection{A variation: mating $\mathbbm{G}_d$ with anti-polynomials}\label{most_poly_anti_hecke_mating_subsubsec}

The original conjecture of Bullett and Freiberger, translated to the antiholomorphic setting, asks whether any degree $d$ anti-polynomial with connected Julia set can be mated with the anti-Hecke group $\mathbbm{G}_d$ as a correspondence (cf. \cite[\S 3, p. 3926]{BuFr}). A modification of the proof of Theorem~\ref{general_mating_corr_existence_thm} combined with combinatorial continuity and rigidity arguments (similar to the ones alluded to in Subsection~\ref{polygonal_per_rep_subsec}) yield the following partial answer to this conjecture.

\begin{theorem}\cite[Theorem~C]{LMM3}\cite[\S 12]{LLM23}\label{semi_hyp_poly_mating_corr_thm}
Let $p\in\mathscr{C}_d$ be either
\begin{itemize}
\item geometrically finite; or 
\item periodically repelling and finitely renormalizable.
\end{itemize}
Then, there exists $f\in\mathfrak{U}_{d+1}$ such that the associated reversible correspondence $\mathfrak{C}$ on $\widehat{\C}$ given by Equation~\eqref{corr_eqn} is the mating of the anti-Hecke group $\mathbbm{G}_d$~and~$p$.
\end{theorem}

One of the main difficulties in carrying out the strategy of the proof of Theorem~\ref{general_mating_corr_existence_thm} in this setting is the hyperbolic-parabolic mismatch that we have encountered many times in this survey. Indeed, the external map $\overline{z}^d$ of a degree $d$ anti-polynomial (with connected Julia set) is not quasiconformally compatible with the map $\pmb{\cF}_d$. Thus, one is compelled to take the route of David surgery as described in Section~\ref{david_surgery_sec}. Specifically, Theorem~\ref{david_extension_general_thm} can be used to prove the existence of a circle homeomorphism that conjugates $\overline{z}^d$ to $\pmb{\cF}_d$ and admits a David extension to $\D$ \cite[Lemma~3.2]{LMM3}. This allows us to glue the map $\pmb{\cF}_d$ outside the filled Julia set of a degree $d$ anti-polynomial (with connected Julia set), \emph{provided} one has good control on the geometry of the basin of infinity of the anti-polynomial (cf. Subsection~\ref{integrable_thms_subsec}). This is indeed the case for a hyperbolic (or more generally, semi-hyperbolic) anti-polynomial $p\in\mathscr{C}_d$, which enables us to construct a Schwarz reflection map $\sigma$ as a conformal mating of $\pmb{\cF}_d$ and a hyperbolic anti-polynomial $p$ (see \cite[Propositions~3.7]{LMM3} for details of the construction). 

The next step of the proof is to take limits of conformal matings of $\pmb{\cF}_d$ with postcritically finite maps in $\mathscr{C}_d$ to manufacture conformal matings of $\pmb{\cF}_d$ with periodically repelling, finitely renormalizable anti-polynomials. As in the proofs of Theorem~\ref{comb_class_homeo_polygonal_thm} and Corollary~\ref{polygonal_fin_renorm_cor}, this procedure involves 
\begin{itemize}
\item combinatorial continuity results (specifically, continuity of rational laminations of maps in $\mathscr{C}_d$ and continuity of preperiodic laminations of maps in~$\cS_{\pmb{\cF}_d}$);
\item combinatorial rigidity of periodically repelling, finitely renormalizable maps in $\mathscr{C}_d, \cS_{\pmb{\cF}_d}$; and
\item local connectivity of Julia/limits sets of periodically repelling, finitely renormalizable maps in $\mathscr{C}_d, \cS_{\pmb{\cF}_d}$.
\end{itemize}

Finally, for a map $p\in\mathscr{C}_d$ satisfying one of the conditions of Theorem~\ref{semi_hyp_poly_mating_corr_thm}, the promotion of the conformal mating $\sigma$ of $\pmb{\cF}_d$ and $p$ to the desired correspondence follows the scheme of Subsection~\ref{upgrade_step_subsubsec}.

\subsection{Regularity of the mating surgery}\label{mating_regularity_subsec}

The bijection between the moduli space $\faktor{\pmb{\mathcal{B}}_d}{\mathrm{Aut}(\C)}$ of parabolic anti-rational maps and the moduli space $\faktor{\mathcal{S}_{\pmb{\cF}_d}}{\mathrm{PSL}_2(\C)}$ of Schwarz reflections (see Theorem~\ref{general_mating_corr_existence_thm}) is continuous at an abundant set of parameters. However, since the surgery procedure described in Subsection~\ref{core_step_subsubsec} involves the Riemann maps of the marked parabolic basins of the anti-rational maps, one needs to study continuity properties of Riemann maps to conclude continuity of the mating operation (cf. Remark~\ref{straightening_regularity_rem}). This is a non-trivial task since filled Julia sets do not move continuously in general. 

To circumvent this issue, a different mating surgery avoiding Riemann maps (but with equivalent outcome) was developed in \cite{LMM3}, one that restricts members of $\pmb{\mathcal{B}}_d^{\mathrm{simp}}$ (maps with a simple parabolic fixed point at $\infty$) to pinched anti-polynomial-like maps and replaces the external maps of such pinched anti-polynomial-like restrictions with the anti-Farey map $\pmb{\cF}_d$. According to \cite[Theorem~5.2]{LMM3}, this can be performed with continuous control over the dilatations of the hybrid conjugacies and the domains of definition of these conjugacies. 

The same can be done for the inverse surgery from $\faktor{\mathcal{S}_{\pmb{\cF}_d}^{\mathrm{simp}}}{\mathrm{PSL}_2(\C)}$ (where $\mathcal{S}_{\pmb{\cF}_d}^{\mathrm{simp}}$ consists of Schwarz reflections in $\mathcal{S}_{\pmb{\cF}_d}$ with a simple cusp on the associated quadrature domain boundary) to $\faktor{\pmb{\mathcal{B}}_d^{\mathrm{simp}}}{\mathrm{Aut}(\C)}$. Specifically, one can restrict $\sigma\in\mathcal{S}_{\pmb{\cF}_d}^{\mathrm{simp}}$ to a pinched anti-polynomial-like map and replace its external map with the anti-Blaschke product $B_d$, thus producing a parabolic anti-rational map $R_\sigma\in \faktor{\pmb{\mathcal{B}}_d^{\mathrm{simp}}}{\mathrm{Aut}(\C)}$ that is hybrid equivalent to $\sigma$.
This generalizes the straightening theorem for pinched anti-quadratic-like maps proved in \cite[Proposition~4.15, Theorem~5.4]{LLMM3} to arbitrary degree (see \cite[Theorem~4.8, Lemma~4.13]{LMM3}, cf. Subsections~\ref{chebyshev_center_hybrid_conj_subsubsec},~\ref{cubic_cheby_qc_straightening_subsubsec}). As in the degree two case, the general straightening theorem uses Warschawski's result on boundary behavior of conformal maps of topological strips in an essential way.

Utilizing these parameter dependencies of the mating surgeries, it was proved in \cite[\S 6]{LMM3} that the straightening map
$$
\chi:\ \faktor{\mathcal{S}_{\pmb{\cF}_d}^{\mathrm{simp}}}{\mathrm{PSL}_2(\C)}\ \longrightarrow\  \faktor{\pmb{\mathcal{B}}_d^{\mathrm{simp}}}{\mathrm{Aut}(\C)},\quad [\sigma] \mapsto [R_\sigma]
$$ 
(respectively, its inverse) is continuous at hyperbolic and quasiconformally rigid parameters of $\faktor{\mathcal{S}_{\pmb{\cF}_d}^{\mathrm{simp}}}{\mathrm{PSL}_2(\C)}$ (respectively, of $\faktor{\pmb{\mathcal{B}}_d^{\mathrm{simp}}}{\mathrm{Aut}(\C)}$). 

To conclude this section, we remark that as in the previous instances of straightening maps appearing in this survey, the  map $\chi$ is not necessarily continuous at quasiconformally deformable parabolic parameters (cf. Theorem~\ref{Straightening_discontinuity_Tricorn} and Subsection~\ref{polygonal_geom_fin_subsec}).

\subsection{Shabat polynomial slices in the space of correspondences}\label{shabat_slices_subsec}

Theorem~\ref{general_mating_corr_existence_thm}, which is a general existence theorem for correspondences that are matings of the anti-Hecke group $\mathbbm{G}_d$ and degree $d$ parabolic anti-rational maps, combined with the regularity of the mating surgery discussed in Subsection~\ref{mating_regularity_subsec}, paves the way for studying natural one-parameter slices of correspondences. Such slices generalize the one-parameter family of antiholomorphic correspondences arising from (injective restrictions of) the cubic Chebyshev polynomial (see Subsection~\ref{chebyshev_gen_subsec}). 
The following account is based on \cite{LMM4}.

\subsubsection{Shabat polynomials and their role}\label{shabat_subsubsec}

Generically, the complex dimension of a natural family of holomorphic/antiholomorphic maps equals the number of free/active critical orbits of the maps in the family. Guided by this philosophy, one aims at constructing one-parameter slices in $\ \faktor{\mathcal{S}_{\pmb{\cF}_d}}{\mathrm{PSL}_2(\C)}$ such that 
\begin{itemize}
\item the families are closed under quasiconformal deformation, and
\item the corresponding Schwarz reflections have a unique free critical value.
\end{itemize}

Let $\sigma:\overline{\Omega}\to\widehat{\C}$ be a Schwarz reflection in $\mathcal{S}_{\pmb{\cF}_d}$ and let $f\in\mathfrak{U}_{d+1}$ be a polynomial such that $f$ maps $\overline{\D}$ homeomorphically onto $\overline{\Omega}$ (see Subsection~\ref{upgrade_step_subsubsec}). Pre-composing $f$ with a rotation, we can assume that the unique critical point of $f$ on $\mathbb{S}^1$ is at $1$. 
Recall that the set of critical values of $\sigma$ is contained in the set of critical values of $f$, and the difference, if non-empty, is the cusp $f(1)$.
Note further that the critical value $f(\infty)=\infty$ of $\sigma$ lies inside the droplet, while the point $f(1)$ is a fixed point of $\sigma$ since it lies on the boundary of $\Omega$. Hence, none of these points is an active critical value. Thus, the condition that $\sigma$ has a unique free critical value is equivalent to the requirement that $f$ has exactly one critical value other than $f(1)$ and $\infty$; i.e., $f$ has three critical values in $\widehat{\C}$.  This leads us to the space of \emph{Shabat polynomials} \cite{LZ04,Sch94}, of which the cubic Chebyshev polynomial is the simplest non-trivial example.

\begin{definition}[Shabat and Belyi maps]\label{shabat_def}
1) A polynomial $g:\widehat{\C}\to\widehat{\C}$ (of degree at least three) is called a \emph{Shabat polynomial} if it has exactly two finite critical values. Two Shabat polynomials $g_1$ and $g_2$ are called \emph{equivalent} if there exist affine maps $A_1, A_2$ such that $g_2=A_1\circ g_1\circ A_2$.

2) A rational or anti-rational map $R:\widehat{\C}\to\widehat{\C}$ is said to be \emph{Belyi} if it has at most three critical values.
\end{definition}

There is an important combinatorial invariant associated with a Shabat polynomial $g$; namely, the \emph{dessin d'enfants} $\mathfrak{T}(g)$. Let us briefly recall the definition and basic properties of this invariant (see \cite[\S 2]{LZ04}, \cite[Appendix~A]{LMM4} for more background). 
For an arc $\gamma\subset\C$ connecting the two finite critical values $y_1$ and $y_2$ of a Shabat polynomial $g$, the preimage 
$$
\mathfrak{T}_\gamma(g):=g^{-1}(\gamma)
$$ 
is a plane tree with vertices at $g^{-1}(\{y_1, y_2\})$ (the set $\C\setminus\mathfrak{T}_\gamma(g)$ is connected since
 $g$ has no pole in $\C$ and this set is a topological annulus since $g:\C\setminus\mathfrak{T}_\gamma(g)\to\C\setminus\gamma$ is a covering map). One colors the pre-images of $y_1$ and $y_2$ black and white, respectively. Then $\mathfrak{T}_\gamma(g)$ has the structure of a bicolored plane tree. The tree $\mathfrak{T}_\gamma(g)$ has $\deg{(g)}$ many edges, and the valence of a vertex of $\mathfrak{T}_\gamma(g)$ is equal to the local degree of $g$ at that point. 
It is easily checked that the isotopy type of $\mathfrak{T}_\gamma(g)$ (relative to the vertices) is independent of the choice of the arc $\gamma$ connecting $y_1$ and $y_2$. In particular, the various $\mathfrak{T}_\gamma(g)$ are isomorphic as combinatorial bicolored plane trees. This combinatorial bicolored plane tree is called the dessin d'enfants of $g$, and is denoted by $\mathfrak{T}(g)$. Moreover, the isomorphism class of $\mathfrak{T}(g)$ (as a bicolored plane tree) remains unaltered if $g$ is replaced by a Shabat polynomial equivalent to $g$. The following classical result states that dessins d'enfants are complete invariants of Shabat polynomials.

\begin{theorem}\cite[Theorem~I.5]{Sch94}, \cite[Theorem~2.2.9]{LZ04}\label{shabat_classification_theorem}
The map $g\mapsto\mathfrak{T}(g)$ induces a bijection between the set of equivalence classes of Shabat polynomials and the set of isomorphism classes of bicolored plane trees (with at least one black and one white vertex of valence greater than one).
\end{theorem}

\subsubsection{One-parameter slices in $\mathcal{S}_{\pmb{\cF}_d}$}\label{one_para_slice_shabat_subsubsec}
We now return to the construction of one-parameter slices in $\ \faktor{\mathcal{S}_{\pmb{\cF}_d}}{\mathrm{PSL}_2(\C)}$. As mentioned before, a Schwarz reflection $\sigma:\overline{\Omega}\to\widehat{\C}$ in $\mathcal{S}_{\pmb{\cF}_d}$ arising from a Shabat polynomial $f\in\mathfrak{U}_{d+1}$ has at most three critical values. Hence, the pullback $\sigma^{-1}(\gamma')$ of an arc $\gamma'$ connecting $f(1)$ to the free critical value of $\sigma$ also carries the structure of a bicolored plane tree. We denote this combinatorial bicolored plane tree by $\mathfrak{T}(\sigma)$, and call it the dessin d'enfants of $\sigma$. Let us now assume that the free critical value of $\sigma$ lies in $\Omega$. It is not hard to see from the relation $\sigma\equiv f\circ\eta\circ(f\vert_{\overline{\D}})^{-1}$ that one can describe $\mathfrak{T}(\sigma)$ purely in terms of $\mathfrak{T}(f)$, and vice versa. 
\begin{figure}[h!]
\captionsetup{width=0.96\linewidth}
\begin{tikzpicture}
\node[anchor=south west,inner sep=0] at (7,0) {\includegraphics[width=0.4\textwidth]{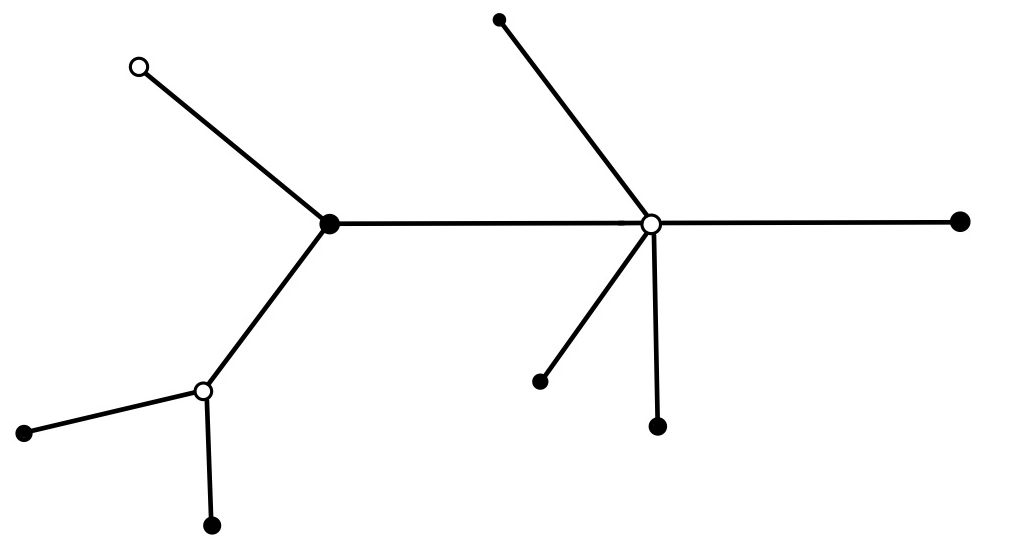}};
\node[anchor=south west,inner sep=0] at (-0.5,0) {\includegraphics[width=0.54\textwidth]{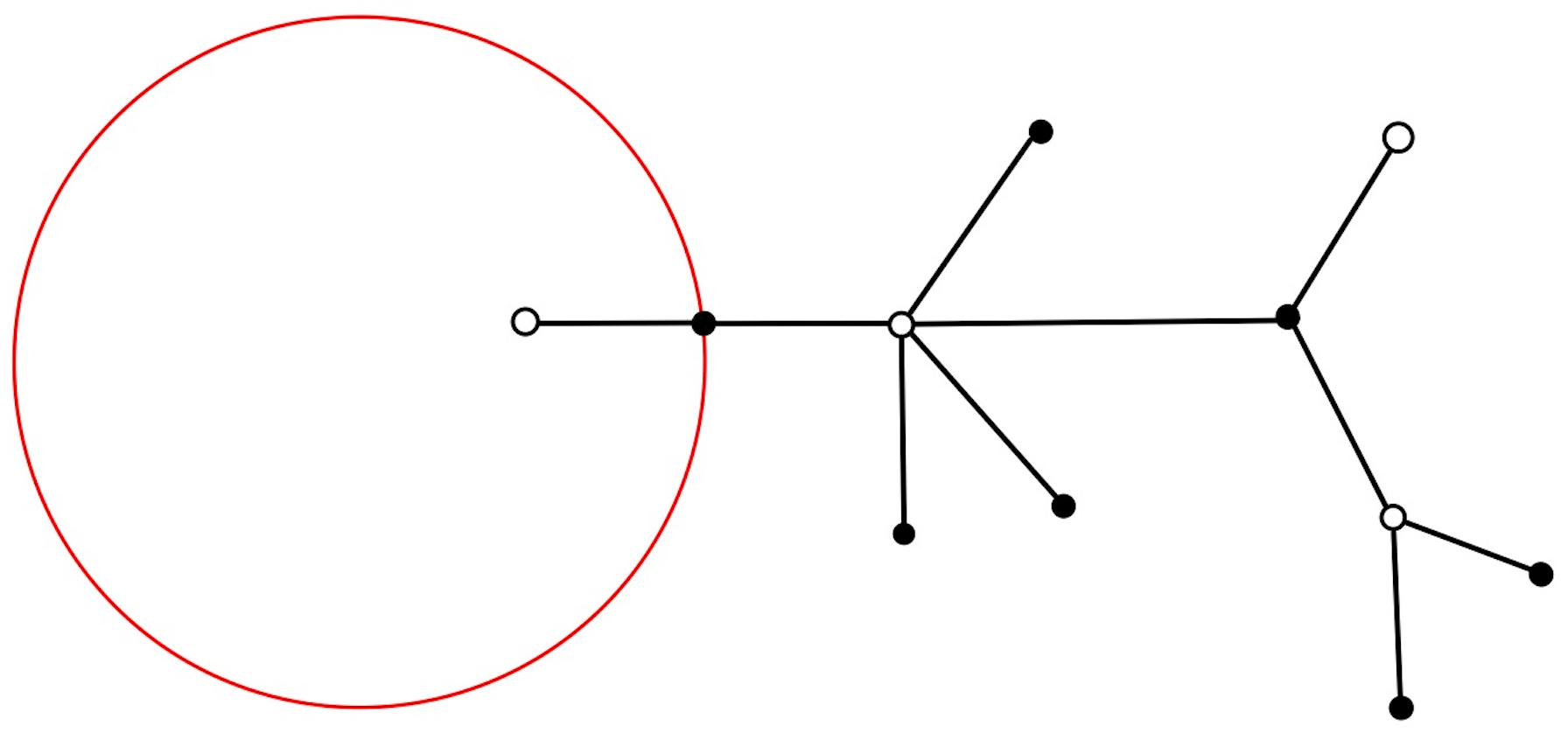}}; 
\node at (4.8,2.8) {$\mathfrak{T}(f)$};
\node at (10.8,2.5) {$\mathfrak{T}(\sigma)$};
\node at (0.4,1.6) {\begin{large}$\D$\end{large}};
\end{tikzpicture}
\caption{Relation between the dessin d'enfants $\mathfrak{T}(f)$ and $\mathfrak{T}(\sigma)$.}
\label{augmented_dessin_fig}
\end{figure}
Specifically, $\mathfrak{T}(\sigma)$ is obtained by pruning a distinguished peripheral edge from $\mathfrak{T}(f)$ and reversing the cyclic order of the edges around each vertex of the resulting tree (see \cite{LMM4}). We also note that if $\sigma$ and $\sigma_1$ (of the above form) are Hurwitz equivalent (in particular, if they are quasiconformally conjugate), then their dessins d'enfants $\mathfrak{T}(\sigma)$ and $\mathfrak{T}(\sigma_1)$ are isomorphic. The above discussion now implies that the Shabat polynomials $f$, $f_1$ associated with $\sigma$, $\sigma_1$ also have isomorphic dessins d'enfants. Thus, in light of Theorem~\ref{shabat_classification_theorem}, a one-parameter family of Schwarz reflections in $\ \faktor{\mathcal{S}_{\pmb{\cF}_d}}{\mathrm{PSL}_2(\C)}$ is closed under quasiconformal deformation precisely when the associated Shabat polynomials are equivalent in the sense of Definition~\ref{shabat_def} (see \cite{LMM4}).

After possibly conjugating the Schwarz reflections $\sigma$ by affine maps (which amounts to post-composing the associated Shabat polynomials with the same affine maps), we can require that all such $\sigma$ have the same marked critical values and cusps. Then, the corresponding Shabat polynomials $f$ only differ by pre-composition with affine maps.
Instead of fixing the domain of univalence $\D$ and varying the Shabat polynomials (that differ by pre-composition with affine maps and are injective on $\D$), it is more convenient to fix a Shabat polynomial $\pmb{f}$ with dessin d'enfants $\mathfrak{T}(\pmb{f})\cong\mathfrak{T}$, and restrict it to all possible disks of univalence such that the disks contain a marked critical critical point of $\pmb{f}$ on their boundary. This leads to the following space of Schwarz reflections.

Let us fix a degree $d+1$ Shabat polynomial $\pmb{f}$ with dessin d'enfants $\mathfrak{T}(\pmb{f})\cong\mathfrak{T}$ such that $\mathfrak{T}(\pmb{f})$ has a valence two `black' vertex $\pmb{v_b}$ with an adjacent valence one `white' vertex $\pmb{v_w}$.

\begin{definition}\label{shabat_schwarz_family_def}
We define the parameter space
$$
S_{\mathfrak{T}}:=\{ a\in\C: \pmb{v_w}\in\Delta_a:=B(a,\vert \pmb{v_b}-a\vert)\ \textrm{and}\ \pmb{f}\vert_{\overline{\Delta_a}}\ \textrm{is\ injective}\},
$$ 
and the associated space of Schwarz reflections
$$
\mathcal{S}_{\mathfrak{T}}:=\{\sigma_a\equiv\pmb{f}\circ\eta_a\circ (\pmb{f}\vert_{\overline{\Delta_a}})^{-1}:\Omega_a:=\pmb{f}(\Delta_a)\longrightarrow\widehat{\C}: a\in S_{\mathfrak{T}}\},
$$
where $\eta_a$ stands for reflection in the circle $\partial\Delta_a$.
\end{definition}

\subsubsection{Disks of univalence of Shabat polynomials}\label{shabat_disk_univ_subsubsec}

For the cubic Chebyshev polynomial $\pmb{f}(u)=u^3-3u$ discussed in Subsection~\ref{chebyshev_subsec}, an explicit description of $S_{\mathfrak{T}}$ was given in \cite[\S 3]{LLMM3} using exact numerical computation. 
In \cite{LMM4}, a detailed analysis of univalence properties of Shabat polynomials is carried out to provide a precise qualitative description of the parameter space $S_{\mathfrak{T}}$ for any Shabat polynomial $\pmb{f}$ with dessin d'enfants $\mathfrak{T}(\pmb{f})$ such that $\mathfrak{T}(\pmb{f})$ has a valence two vertex $\pmb{v_b}$ with an adjacent valence one vertex $\pmb{v_w}$.
We summarize the main results below.

\begin{theorem}[Topology of the parameter space $S_{\mathfrak{T}}$]\label{shabat_schwarz_para_space_thm_into}
\noindent\begin{enumerate}\upshape
\item $\Int{S_{\mathfrak{T}}}$ is a bounded Jordan domain, and $\overline{S_{\mathfrak{T}}}=\overline{\Int{S_{\mathfrak{T}}}}$. 

\item $\Int{S_{\mathcal{T}}}= \{a\in\C: \pmb{v_w}\in\Delta_a,\ f\vert_{\overline{\Delta_a}}\ \textrm{is\ injective,\ and}\ \pmb{f}(\pmb{v_b})\ \textrm{is\ a}\ (3,2)\ \textrm{cusp\ on}\ \partial\Omega_a \}.$

\item The boundary $\partial S_{\mathfrak{T}}$ has the structure of a topological quadrilateral such that
\begin{enumerate}\upshape
\item two \emph{horizontal} sides of $\partial S_{\mathfrak{T}}$ are characterized by the existence of a double point on the boundary of $\partial\Omega_a$,
\item one \emph{vertical} side of $\partial S_{\mathfrak{T}}$ is characterized by the condition that $\pmb{f}(\pmb{v_b})$ is a $(\nu,2)-$cusp on $\partial\Omega_a$, for $\nu\geq 5$, and
\item the other \emph{vertical} side of $\partial S_{\mathfrak{T}}$ is characterized by the condition that $\pmb{v_w}\in\partial\Delta_a$.
\end{enumerate}
\end{enumerate}
\end{theorem}

The proof of the above theorem essentially depends on the local dynamics of Schwarz reflection maps near conformal cusps and double points, and the relation between such singular points and the critical orbits of the corresponding Schwarz reflection maps. 

\subsubsection{Connectedness locus of $\mathcal{S}_{\mathfrak{T}}$}\label{shabat_slice_conn_locus_subsubsec}

As usual, the connectedness locus $\mathcal{C}(\mathcal{S}_{\mathfrak{T}})$ of $\mathcal{S}_{\mathfrak{T}}$ is defined as the collection of Schwarz reflections in the family with connected non-escaping set. Recall that for all maps $\sigma_a\in\mathcal{S}_{\mathfrak{T}}$, there is a passive critical value (of multiplicity $d$) in the tiling set that escapes to the fundamental tile in one iterate. A map $\sigma_a$ lies in the connectedness locus $\cC(\mathcal{S}_{\mathfrak{T}})$ if and only if the unique free critical value of $\sigma_a$ lies in the non-escaping set. The existence of a unique fully ramified critical point in the tiling set of maps in $\mathcal{C}(\mathcal{S}_{\mathfrak{T}})$ allows one to show that the tiling set dynamics of each map $\sigma_a\in\mathcal{C}(\mathcal{S}_{\mathfrak{T}})$ is conformally conjugate to the anti-Farey map $\pmb{\cF}_d$ (see Subsections~\ref{nielsen_first_return_external_map_subsubsec} and~\ref{cubic_cheby_qc_straightening_subsubsec} for discussions on the same result in the $d=2$ setting, cf. \cite{LMM4}). In other words, the one-parameter family $\mathcal{S}_{\mathfrak{T}}$ of Schwarz reflections meets the space $\mathcal{S}_{\pmb{\cF}_d}$ in its connectedness~locus. 

\begin{proposition}\label{shabat_slice_conn_locus_srd_prop}
$\mathcal{S}_{\mathfrak{T}}\cap \mathcal{S}_{\pmb{\cF}_d} = \mathcal{C}(\mathcal{S}_{\mathfrak{T}}).$
\end{proposition}

\subsubsection{Image of $\mathcal{C}(\mathcal{S}_{\mathfrak{T}})$ under the straightening map $\chi$}

We will now describe the image of the connectedness locus $\mathcal{C}(\mathcal{S}_{\mathfrak{T}})$ under the straightening map $\chi$ defined in Subsection~\ref{corr_para_anti_rat_bijection_sketch_subsubsec}. 

As each $\sigma_a\in\mathcal{C}(\mathcal{S}_{\mathfrak{T}})$ has at most two critical values in its non-escaping set (namely, the free critical value $\pmb{f}(\pmb{v_w})$ and possibly the conformal cusp $\pmb{f}(\pmb{v_b})$), it follows from the definition of $\chi$ that the straightened map $\chi(\sigma_a)$ has at most two critical values in its filled Julia set and exactly one (fully ramified) critical value in its completely invariant parabolic basin. In other words, each $R\in\chi(\mathcal{C}(\mathcal{S}_{\mathfrak{T}}))$ is a Belyi anti-rational map of $\widehat{\C}$. Moreover, if $\chi(\sigma_a)$ has three critical values, then the parabolic fixed point $\infty$ is one of them. The dessin d'enfants $\mathfrak{T}(R)$ of the Belyi map $R$ is defined as the combinatorial plane bicolored tree isomorphic to $R^{-1}(\gamma')$, where $\gamma'$ is an arc connecting the parabolic fixed point $\infty$ to the free critical value of $R$ (which lies in the filled Julia set of $R$).

Note further that the dessin d'enfants of any $\sigma_a$ has a distinguished valence one vertex at the cusp $\pmb{f}(\pmb{v_b})\in\partial\Omega_a$, which can be regarded as the root of the tree. We denote this abstract rooted bicolored plane tree by $(\mathfrak{T}^{\textrm{del}}, O)$ (the isomorphism class of this tree is independent of the parameter $a\in\mathcal{S}_{\mathfrak{T}}$). Here, the superscript `del' is chosen to indicate that the tree $\mathfrak{T}^{\textrm{del}}$ is obtained by deleting the edge $[\pmb{v_w},\pmb{v_b})$ from $\mathfrak{T}(\pmb{f})\cong\mathfrak{T}$ and reversing the cyclic order of the edges around each vertex of the resulting tree (see Subsection~\ref{one_para_slice_shabat_subsubsec}). 
Evidently, the dessin d'enfants of each $R\in\chi(\mathcal{C}(\mathcal{S}_{\mathfrak{T}}))$ is isomorphic to $\mathfrak{T}^{\textrm{del}}$. We also note that for $R\in\chi(\mathcal{C}(\mathcal{S}_{\mathfrak{T}}))$, the parabolic fixed point $\infty$ is a distinguished vertex of valence one on $\mathfrak{T}(R)\cong\mathfrak{T}^{\textrm{del}}$, and this vertex corresponds to the root vertex $O$ of the dessin d'enfants $\mathfrak{T}^{\textrm{del}}$ of the corresponding Schwarz reflection under the hybrid conjugacy. It follows that $\chi(\mathcal{C}(\mathcal{S}_{\mathfrak{T}}))$ is contained in the following space of parabolic anti-rational maps.

\begin{definition}\label{fd_slice_def}
We define
\begin{align*}
\mathfrak{F}_{\mathfrak{T}} & :=\bigg\{ R\in\pmb{\mathcal{B}}_d: R \textrm{ is Belyi, if } R \textrm{ has three critical values, then the parabolic }\\
&\qquad \textrm{ fixed  point } \infty\ \textrm{ is one of them, and }  \left(\mathfrak{T}(R),\infty\right)\cong\left(\mathfrak{T}^{\textrm{del}},O\right)\bigg\}/\mathrm{Aut}(\C),
\end{align*}
where the isomorphism is required to preserve the roots and the bicolored plane structures.
\end{definition}

The above discussion, combined with the fact that $\chi$ is a bijection between $\ \faktor{\mathcal{S}_{\pmb{\cF}_d}}{\mathrm{PSL}_2(\C)}$ and $\  \faktor{\pmb{\mathcal{B}}_d}{\mathrm{Aut}(\C)}$, implies that $\chi(\mathcal{C}(\mathcal{S}_{\mathfrak{T}}))=\ \faktor{\mathfrak{F}_{\mathfrak{T}}}{\mathrm{Aut}(\C)}$.

\subsubsection{Combinatorial model of $\mathcal{C}(\mathcal{S}_{\mathfrak{T}})$ and homeomorphism between models}\label{shabat_conn_locus_model_subsubsec}

The conformal position of the free critical value of $\sigma_a$ (in the tiling set) can be used to define a dynamically natural uniformization of the \emph{escape locus} $\mathcal{S}_{\mathfrak{T}}\setminus\mathcal{C}(\mathcal{S}_{\mathfrak{T}})$ (this is analogous to the dynamical uniformization of the escape locus of the C\&C family, see Subsection~\ref{c_and_c_geom_fin_bijection_subsubsec}). This uniformization gives a tessellation structure in the escape locus which allows us to construct parameter rays. The co-landing/co-accumulation patterns of these parameter rays is used in \cite{LMM4} to construct a model of the connectedness locus of $\mathcal{S}_{\mathfrak{T}}$ as a pinched disk.

Moreover, the continuity properties of the straightening map $\chi$ explicated in Subsection~\ref{mating_regularity_subsec} imply that $\chi$ induces a homeomorphism between the above pinched disk model of $\mathcal{C}(\mathcal{S}_{\mathfrak{T}})$ and a similar model for $\ \faktor{\mathfrak{F}_{\mathfrak{T}}}{\mathrm{Aut}(\C)}$. 
We remark that progress in combinatorial rigidity problems for the above parameter spaces would bring the pinched disk models closer to the actual connectedness loci.

\subsection{Correspondences as matings of $\pmb{G}_d$ and anti-polynomials}\label{deltoid_generalization_corr_subsec}

Every correspondence discussed in this survey thus far (in Subsections~\ref{deltoid_corr_subsubsec},~\ref{chebyshev_subsec},~\ref{chebyshev_gen_subsec}, and earlier in the current section) arises from the uniformizing rational map of a single quadrature domain. In \cite{LLM23}, a recipe for obtaining antiholomorphic correspondences from piecewise Schwarz reflections (e.g. Schwarz reflections in the Circle-and-Cardioid family discussed Subsection~\ref{c_and_c_general_subsec}, or more generally, polygonal Schwarz reflections studied in Section~\ref{mating_para_space_sec}) was given, and a combination theorem for the resulting correspondences was established.

\begin{theorem}\cite{LLM23}\label{polygonal_corr_thm}
Let $f$ be a degree $d$ anti-polynomial with connected Julia set which is either
\begin{itemize}
\item geometrically finite; or 
\item periodically repelling and finitely renormalizable.
\end{itemize}
Then there exists a reversible antiholomorphic correspondence $\mathfrak{C}$ on a compact, simply connected (possibly noded) Riemann surface $\mathfrak{W}$ which is a mating of $f$ and $\pmb{G}_d$.
\end{theorem}

The main new ingredient in the proof of Theorem~\ref{polygonal_corr_thm} is the construction of the nodal Riemann surface $\mathfrak{W}$ which is the phase space of the desired correspondence. We now outline how the surface $\mathfrak{W}$ naturally appears in the context of the theorem.

Let $f$ be as in the statement of Theorem~\ref{polygonal_corr_thm}. According to Theorem~\ref{geom_finite_bij_polygonal_thm} and Corollary~\ref{polygonal_fin_renorm_cor}, there exists a piecewise Schwarz reflection $(\sigma:\overline{\Omega}\to\widehat{\C})\in\mathcal{S}_{\pmb{\cN}_d}$ that is a conformal mating of $f$ and $\pmb{\cN}_d$. Recall that the quadrature multi-domain $\Omega$ may be disconnected with (Jordan) components $\Omega_1,\cdots,\Omega_k$. Let us denote the uniformizing rational maps of the quadrature domains by $\phi_j:\overline{\D}\to\overline{\Omega_j}$, $j\in\{1,\cdots,k\}$. For book-keeping purposes, we will denote the domain of $\phi_j$ (which is a copy of the Riemann sphere) by $\widehat{\C}_j$ and the unit disk in $\widehat{\C}_j$ by $\D_j$. Let us also set $\mathcal{I}:=\{1,\cdots,k\}$.

The limit set $\Lambda(\sigma)$ intersects $\partial\Omega$ in a finite set $S$, which consists of the singular points of $\partial\Omega$.
We define 
$$
\widetilde{S}_j\ := (\phi_j\vert_{\partial\D_j})^{-1}\left(S \cap \partial\Omega_j\right),\quad j\in\mathcal{I}.
$$

Let us now consider the disjoint union 
$$
\mathfrak{U}\ :=\ \bigsqcup_{i\in\mathcal{I}} \widehat{\C}_{i}\ \cong\ \widehat{\C}\times\mathcal{I}.
$$
Our goal is to capture the touching patterns of the various $\Omega_j$, and to construct the space $\mathfrak{W}$ by making the various spheres $\widehat{\C}_j$ touch exactly in the same fashion. To this end, we first define the following finite equivalence relation on $\mathfrak{U}$:
\begin{center}
For $z\in\widetilde{S}_i\subset\widehat{\C}_i$ and $w\in\widetilde{S}_j\subset\widehat{\C}_j$,
\smallskip

$(z,i)\sim (w,j)\iff \phi_i(z)=\phi_j(w)$.
\end{center}
Finally, we set  
$$
\mathfrak{W}\ :=\ \faktor{\mathfrak{U}}{\sim}.
$$
The space $\mathfrak{W}$ can naturally be regarded as a compact, simply connected, noded Riemann surface (triviality of $\pi_1(\mathfrak{W})$ is a consequence of simple connectedness of $\overline{\Omega}$).

Having constructed the nodal Riemann surface $\mathfrak{W}$, we now proceed to put all the rational maps $\phi_j$ together into a single branched cover $\mathfrak{W}\longrightarrow\widehat{\C}$. This is done simply by defining
$$
\pmb{R}:\ \mathfrak{W}\longrightarrow \widehat{\C},\quad (z,j)\mapsto R_{j}(z).
$$
(The definition of the equivalence relation $\sim$ readily shows that $\pmb{R}$ is a well-defined branched covering of degree $d+1$.)
To define the correspondence $\mathfrak{W}$, we also need an antiholomorphic involution (which acts as $\eta(w)=1/\overline{w}$ on each $\widehat{\C}_j$):
$$
\pmb{\eta}\ : \mathfrak{W}\longrightarrow \mathfrak{W},\quad (z,j)\mapsto (\eta(z),j).
$$

Generalizing the definition of the correspondence given by Equation~\eqref{corr_eqn}, we now define a correspondence $\mathfrak{C}$ on $\mathfrak{W}$ given by the equation: 
\begin{equation}
\{(\mathfrak{u}_1,\mathfrak{u}_2)\in\mathfrak{W}\times\mathfrak{W}: \frac{\pmb{R}(\mathfrak{u}_2)-\pmb{R}(\pmb{\eta}(\mathfrak{u}_1))}{\mathfrak{u}_2-\pmb{\eta}(\mathfrak{u}_1)}=0\}.
\label{corr_gen_eqn}
\end{equation}

Unsurprisingly, the dynamics of the correspondence $\mathfrak{C}$ on $\mathfrak{W}$ is intimately related to the dynamics of $\sigma$. Specifically, $\mathfrak{C}$ can be interpreted as lifts of $\sigma^{\pm 1}$ under the branched covering $\pmb{R}:\mathfrak{W}\to\widehat{\C}$. This allows one to translate the mating structure of $\sigma$ to a mating structure for $\mathfrak{C}$, and complete the proof of Theorem~\ref{polygonal_corr_thm}.

\section{Analytic applications}\label{anal_app_sec}

\subsection{Conformal removability}\label{conf_removable_subsec}

Conformal removability of various fractal sets, such as limit and Julia sets of Kleinian groups and rational maps, is an important question in geometric function theory. Using the fact that boundaries of John domains (more generally, H{\"o}lder domains) are conformally removable, Carleson, Jones, Smirnov, and Yoccoz deduced conformal removability of connected Julia sets of semi-hyperbolic (more generally, Collet--Eckmann) polynomials \cite{CJY,Jon95,JS00}. The situation is more subtle for Julia sets of parabolic polynomials due to the presence of cusps: indeed, the existence of cusps implies that the basin of infinity of a parabolic polynomial is not a John domain and hence the above results do not apply. In the same vein,  since the limit sets of necklace groups also have infinitely many cusps, it is natural to ask whether such limit sets are conformally removable.
It turns out that the David surgery techniques of Section~\ref{david_surgery_sec} can be sometimes used to address the above questions.

\begin{theorem}\cite[Theorems~9.1,9.2]{LMMN}\label{limit_julia_conf_removable_thm}
\begin{enumerate}\upshape
\item Let $P$ be a geometrically finite polynomial with connected Julia set $\mathcal{J}(P)$. Then $\mathcal{J}(P)$ is conformally removable.
\item  The limit set of a necklace reflection group is conformally removable.
\end{enumerate}
\end{theorem}

The proofs of both removability results are based on the fact that global David homeomorphisms carry $W^{1,1}$-removable compact sets to conformally removable ones (see Theorem~\ref{w11_removable_thm}).

For part (1) of Theorem~\ref{limit_julia_conf_removable_thm}, one first appeals to standard realization theorems in holomorphic dynamics to construct a postcritically finite polynomial $Q$ whose Julia dynamics is conjugate to that of $P$. Subsequently, one replaces suitable basins of attraction of $Q$ with parabolic basins using the David Surgery Lemma~\ref{david_surgery_lemma} (also see the David extension result of Example~\ref{example_2}) to recover $P$. This shows that $\mathcal{J}(P)$ is the image of the Julia set of the postcritically finite polynomial $Q$ under a global David homeomorphism. Since the basin of infinity of a postcritically finite polynomial is a John domain, it follows from Theorem~\ref{w11_removable_thm} that $\mathcal{J}(P)$ is conformally removable. 

The same strategy also yields part (2) of Theorem~\ref{limit_julia_conf_removable_thm} since the limit set of a necklace group is the image of the Julia set of a critically fixed anti-polynomial under a global David homeomorphism (see Subsection~\ref{david_regularity_subsubsec}). 

\subsection{Welding homeomorphisms}\label{welding_subsec}

A homeomorphism $h\colon \mathbb S^1 \to \mathbb S^1$ is called a \emph{welding homeomorphism} if there exists a Jordan curve $J$ and conformal homeomorphisms $H_1:\D\to U_1$ and $H_2:\D^*\to U_2$ (where $U_1, U_2$ are the interior and exterior of $J$, respectively) so that $h=\widetilde{H_2}^{-1}\circ \widetilde{H_1}$, where $\widetilde{H_1}$ and $\widetilde{H_2}$ are the homeomorphic extensions of $H_1$ and $H_2$ to the closures of $\D$ and $\D^*$, respectively. The Jordan curve $J$ is called a \emph{welding curve} corresponding to $h$. 
Note that if there exists a conformally removable welding curve $J$ corresponding to the welding homeomorphism $h$, then any other welding curve (corresponding to $h$) is a M\"obius image of $J$.

It is a straightforward consequence of the Ahlfors-Beurling Extension Theorem and the Measurable Riemann Mapping Theorem that every quasisymmetric homeomorphism of $\mathbb{S}^1$ is a welding homeomorphism, and the associated Jordan curve, which is a quasi-circle, is unique (up to M{\"o}bius transformations). This has applications to several important constructions in conformal dynamics; such as mating two Fuchsian groups to obtain a quasi-Fuchsian group, mating two Blaschke products to obtain a quasi-Blaschke rational map, etc.
\begin{figure}[h!]
\captionsetup{width=0.96\linewidth}
\begin{center}
\includegraphics[width=0.506\linewidth]{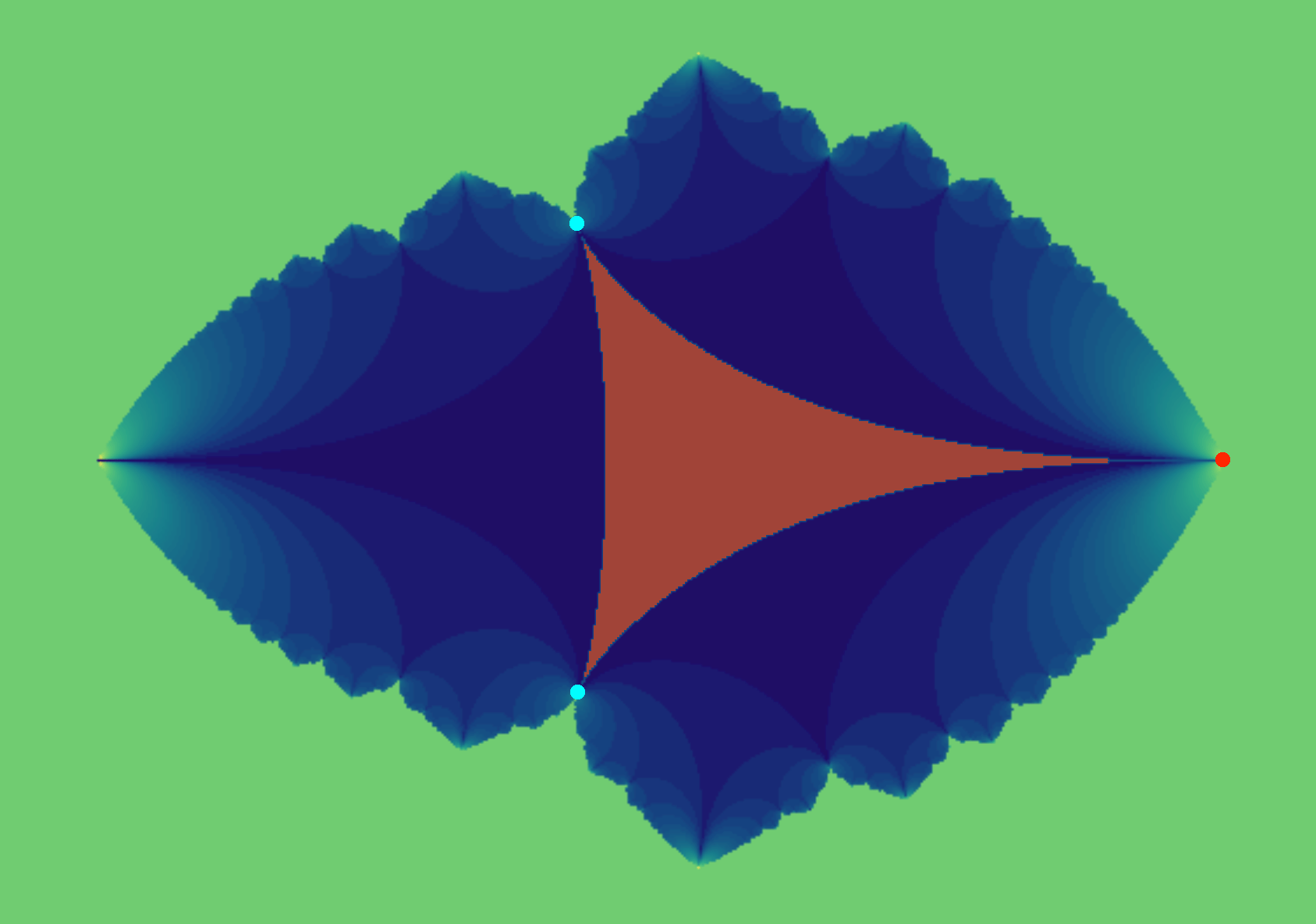}\quad \includegraphics[width=0.355\linewidth]{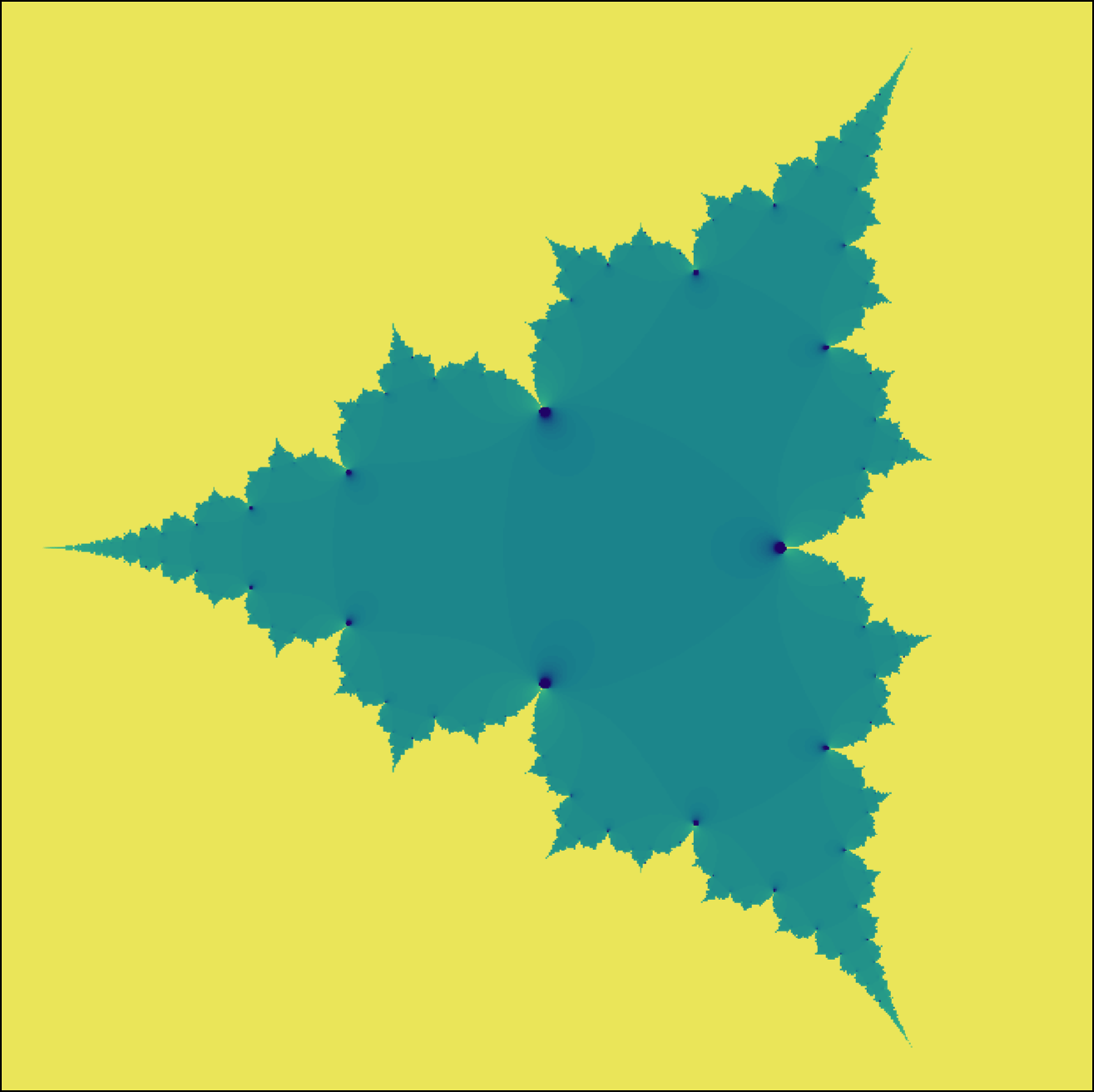}
\end{center}
\caption{Left: The dynamical plane of the Schwarz reflection map that is the unique conformal mating of $P(z)=\overline{z}^2+\frac14$ and $\pmb{\cN}_2$ is displayed. Its limit is a conformally removable Jordan curve with both cusps and sectors. Right: The pine tree shaped fractal is the Julia set of a cubic quasi-Blaschke product that has two parabolic basins. The Julia set, which has both inward and outward cusps, is a conformally removable Jordan curve.}
\label{pine_tree_fig}
\end{figure}

The David Extension Theorem~\ref{david_extension_general_thm} and the accompanying David surgery technique described in Subsection~\ref{david_surgery_hyp_para_subsec} lead to a general realization theorem for dynamically arising circle homeomorphisms as welding homeomorphisms.
 
\begin{theorem}\cite[Theorem~5.1]{LMMN}\label{welding_thm}
Let $f,g\colon \mathbb{S}^1\to \mathbb{S}^1$ be $C^1$, expansive, covering maps of the same degree and the same orientation, and $\mathcal P(f;\{a_0,\dots,a_r\})$, $\mathcal P(g;\{b_0,\dots,b_s\})$ be Markov partitions satisfying conditions \eqref{condition:uv} and \eqref{condition:holomorphic} of Section~\ref{david_surgery_sec} . Assume that each periodic point $a\in \{a_0,\dots,a_r\}$ of $f$ and each periodic point $b\in \{b_0,\dots,b_s\}$ of $g$ is either hyperbolic or symmetrically parabolic. Then any conjugating homeomorphism $h\colon \mathbb S^1\to \mathbb S^1$ between $f$ and $g$ is a welding homeomorphism and the corresponding welding curve is unique up to a M\"obius transformation.  
\end{theorem}

We point out that in the above theorem, the circle homeomorphism $h$ itself does not necessarily have a David extension to $\D$ (because it can move parabolic points to hyperbolic points). However, thanks to Corollary~\ref{power_map_cor}, the hypotheses of the theorem guarantee the existence of a pair of circle homeomorphisms conjugating $z^d$ or $\overline{z}^d$ (depending on the orientation) to $f, g$ (respectively) such that these conjugacies extend as David homeomorphisms of $\D$. One can replace the dynamics of the power map on $\D$ and $\D^*$ with the dynamics of $f$ and $g$ (respectively) using these David homeomorphisms, and then appeal to the David Integrability Theorem to conjugate this map to a holomorphic or antiholomorphic map (defined on a subset of $\widehat{\C}$) via a David homeomorphism $\Psi$. Since the resulting map is conformally conjugate to $f, g$ on the closures of the two complementary components of $\Psi(\mathbb{S}^1)$, it follows that $\Psi(\mathbb{S}^1)$ is the desired welding curve. It is worth mentioning that this construction produces first examples of dynamically defined conformally removable welding curves with infinitely many inward and outward cusps (see Figure~\ref{pine_tree_fig}).

\begin{remark}
We note that there is a classical quasiconformal surgery procedure 
(known as the Douady--Ghys surgery, see \cite{Ghy84}, \cite[\S 7.2]{BF14})
of mating a Blaschke product with a bounded type disk rotation to obtain a Siegel Julia set. 
It  was generalized by Petersen and Zakeri, by means of a David surgery, to almost all rotation numbers (see \cite{PZ04}, \cite[\S 9.2]{BF14}).
\end{remark}

\subsection{Extremal problems}\label{extremal_prob_subsec}

In classical complex analysis, the problem of coefficient maximization in families of univalent holomorphic maps played a central role. The celebrated De Brange's Theorem (earlier known as the Bieberbach conjecture) is the most prominent result in this area. However, for the class of \emph{external univalent} maps 
$$
\Sigma := \left\{ f(z)= z+\frac{a_1}{z} + \cdots +\frac{a_d}{z^d}+\cdots :\ f\vert_{\D^*} \textrm{ is univalent}\right\},
$$
the coefficient maximization problem is still unresolved (see \cite[\S 4.7]{Dur83}, \cite{HS05} for known results). It tuns out that the union of the rational subfamilies $\Sigma_d^*\subset\Sigma$ introduced in Subsection~\ref{sigma_d_subsec} is dense in $\Sigma$ (cf. \cite[Theorem~10]{Suf72}). 
Thus, it is natural to study the coefficient maximization problem for the spaces $\Sigma_d^*$. Recall that the set $H(\C\setminus\overline{\D})$ of all analytic functions on $\C\setminus\overline{\D}$ is a locally convex topological vector space, and $\Sigma_d^*$ is a compact subset of a finite-dimensional vector subspace of $H(\C\setminus\overline{\D})$. An \emph{extreme point} of $\Sigma_d^*$ is an element of $\Sigma_d^*$ which cannot be represented as a proper convex combination of two distinct elements of $\Sigma_d^*$.
By the Krein-Milman theorem, $\Sigma_d^*$ is contained in the convex hull of its extreme points. Hence, it is enough to investigate the coefficients of the extreme points of $\Sigma_d^*$ (see \cite{Bri70,Suf72} for an implementation of this strategy for the classical class $\mathcal{S}$ of normalized univalent holomorphic functions on $\D$).
By \cite[Theorem~2.5]{LM1}, an extreme point $f$ of $\Sigma_d^*$ is a Suffridge map (or a vertex of $\Sigma_d^*$) in the sense of Subsection~\ref{sigma_d_cell_structure_subsec}.
In other words, for an extreme point $f\in\Sigma_d^*$, the compact $\C\setminus f(\D^*)$ is a \emph{tree of triangles} and hence can be modeled by so-called \emph{bi-angled trees} \cite[\S 2.5, Table~1]{LMM1}. Such a tree is essentially the rooted adjacency/contact graph of the tree of triangles $\C\setminus f(\D^*)$ with the recording of whether the exit from a particular triangle is on the left or right.
The above discussion shows that the coefficient maximization problem in $\Sigma_d^*$ naturally leads one to the question of classifying Suffridge maps. Specifically, it motivates the following question: can any topological type of tree of triangles (or equivalently, bi-angled tree) be realized as a vertex of $\Sigma_d^*$?

The main difficulty in constructing Suffridge maps in $\Sigma_d^*$ is that in general, it is hard to check univalence of a rational map on a round disk. To circumvent this issue, one needs to look at the above problem through a different lens. Specifically, thanks to Proposition~\ref{simp_conn_quad_prop}, finding Suffridge maps in $\Sigma_d^*$ is equivalent to constructing Schwarz reflection maps with appropriate dynamical properties. One way of formalizing this statement is that the Suffridge maps in $\Sigma_d^*$ correspond to Schwarz reflection maps that are matings of $\overline{z}^d$ with maximally cusped necklace groups; i.e., necklace groups $G$ for which $\overline{\Pi(G)}\cap K(G)$ is a tree of triangles.
This perspective was adopted in \cite[\S 12]{LMMN} to prove the following classification theorem. 

\begin{theorem}\cite[Theorem~A]{LMM1}\cite[Theorem~12.7]{LMMN}
Let $d\geq2$. There is a canonical bijection between the following classes of objects.
\begin{itemize}
\item $\left\{ f \in \Sigma_d^* : f(\mathbb{T}) \textrm{ has } d+1 \textrm{ cusps and } d-2 \textrm{ double points} \right\} \big/ \hspace{1mm} \mathbb{Z}_{d+1}. $\\

\item $\{ \textrm{Bi-angled trees with } d-1 \textrm{ vertices up to isomorphism respecting the } \\  \textrm{ angular structure}\}.$
\end{itemize}
\end{theorem}
A different proof of the above theorem was originally given in \cite{LMM1} using quasiconformal deformation of Schwarz reflections and compactness of $\Sigma_d^*$. We refer the reader to Subsection~\ref{talbot_pinching_subsubsec} for an illustration of this approach.

Note that uniqueness of a member of $\Sigma_d^*$ realizing a given topological type of tree of triangles follows from conformal removability arguments.

\subsubsection{Zeroes of harmonic polynomials}\label{harmonic_poly_zero_subsubsec} 

The classification of vertices of $\Sigma_d^*$ is also related to construction of harmonic polynomials with maximal number of zeroes. It was observed by Crofoot and Sarason that a harmonic polynomial $p(z)-\overline{z}$ has the maximal number of zeroes (namely, $3d-2$ when $\deg{p}=d$) if $\overline{p}$ has $d-1$ distinct fixed critical points in $\C$. This can be seen by viewing the zeroes of $p(z)-\overline{z}$ as fixed points of the anti-polynomial $\overline{p(z)}$ and invoking the Lefschetz-Hopf Fixed Point Theorem combined with Fatou's count of the number of attracting fixed points of a polynomial. Crofoot and Sarason also conjectured the existence of such polynomials $p$, which was later proved by Geyer \cite{Gey08}. We term such polynomials $p$ as \emph{Crofoot-Sarason polynomials}. The correspondence between maps in $\Sigma_d^*$ and critically fixed anti-polynomials given in Theorem~\ref{sigma_d_bers_homeo_thm} improves the conclusion of \cite{Gey08}, and shows that Crofoot-Sarason polynomials of degree $d$ bijectively correspond to the vertices in $\Sigma_d^*$ (cf. \cite[Theorem~A]{LMM1}).

\subsection{Region of univalence for complex polynomials}\label{univalent_poly_subsec}

Regions of univalence of complex polynomials and rational maps is a well-studied problem in classical complex analysis (see \cite[\S 7.4]{She00} for general results, and \cite{Bra67,CR68,Suf72} for univalence loci of special families of polynomials). While much of the classical development on this topic is based on geometric function theory, it turns out that many questions on univalence loci of rational maps can be answered using the iteration theory of Schwarz reflection maps. For instance, if a rational map $f$ is known to be univalent on $\D$, it is natural to ask what sort of perturbations of $f$ would continue to be univalent on $\D$. This is particularly subtle when $f(\partial\D)$ contains a double point or a conformal cusp. However, it is often possible to quasiconformally deform the associated Schwarz reflection map $\sigma\equiv f\circ\eta\circ(f\vert_{\overline{\D}})^{-1}:f(\overline{\D})\to\widehat{\C}$ to nearby Schwarz reflections and then appeal to the characterization of simply connected quadrature domains (Proposition~\ref{simp_conn_quad_prop}) to construct rational maps close to $f$ that carry $\D$ injectively onto the deformed quadrature domains.
This strategy is successfully implemented in \cite{LMM2} to study the space $\Sigma_d^*$ of rational maps which are univalent on $\D^*$ and in \cite{LMM4} to study the region of univalence of Shabat polynomials.

\appendix

\section{Quasiconformal and David homeomorphisms}\label{qc_david_appendix}

\subsection{Basic definitions}\label{qc_david_def_subsec}

An orientation-preserving homeomorphism $H\colon U\to V$ between domains in the Riemann sphere $\widehat{\C}$ is called $K$-\emph{quasiconformal} for some constant $1\leq K<\infty$ if it lies in the Sobolev class $W^{1,2}_{\mathrm{loc}}(U)$ (i.e., the partial derivatives $\partial H/\partial z$ and $\partial H/\partial\overline{z}$ exist in the sense of distributions and belong to $L^2_{\mathrm{loc}}$) and satisfies $\vert\vert\mu_H\vert\vert_\infty\leq (K-1)/(K+1)$, where $\mu_H=  \frac{\partial H/ \partial \overline{z}}{\partial H/\partial z}$ is the Beltrami coefficient of $H$. Note that the constant $(K-1)/(K+1)$ is always contained in $[0,1)$.

An orientation-preserving homeomorphism $H\colon U\to V$ between domains in $\widehat{\C}$ is called a \emph{David homeomorphism} if it lies in the Sobolev class $W^{1,1}_{\mathrm{loc}}(U)$ and there exist constants $C,\alpha,\varepsilon_0>0$ with
\begin{align}\label{david_def}
\sigma(\{z\in U: |\mu_H(z)|\geq 1-\varepsilon\}) \leq Ce^{-\alpha/\varepsilon}, \quad \varepsilon\leq \varepsilon_0.
\end{align}
Here $\sigma$ is the spherical measure. By Condition~\eqref{david_def}, the Beltrami coefficient of a David homeomorphism takes values in $\D$ a.e.

\subsection{Integrability theorems and basic properties}\label{integrable_thms_subsec}

The Measurable Riemann Mapping Theorem \cite[Theorem 5.3.4, p.~170]{AIM09} states that if $\mu$ is a measurable function on $\widehat{\C}$ with $\|\mu\|_\infty<1$, then there exists a quasiconformal homeomorphism $H\colon \widehat{\C} \to \widehat{\C}$, unique up to postcomposition with M{\"o}bius maps, that solves the Beltrami equation
\begin{align*}
\frac{\partial H}{\partial \overline{z}}= \mu \frac{\partial H}{\partial z}.
\end{align*}

The following integrability result is a generalization of the Measurable Riemann Mapping Theorem, and makes David homeomorphisms useful in holomorphic dynamics. If $U$ is an open subset of $\widehat {\C}$ and $\mu \colon U \to \D$ is a measurable function satisfying Condition~\eqref{david_def} on $U$ for some constants $C,\alpha,\varepsilon_0>0$, then $\mu$ is called a \emph{David coefficient} on $U$. Such a function $\mu$ defines a measurable field of ellipses on $U$ (up to scaling), which we call a \emph{David conformal structure} on $U$.

\begin{theorem}[David Integrability Theorem]\cite{David}, \cite[Theorem~20.6.2, p.~578]{AIM09}\label{david_integrability_thm}
Let $\mu\colon \widehat{\C} \to \D$ be a David coefficient.  Then there exists a homeomorphism $H\colon \widehat{\C} \to \widehat{\C}$ of class $\in W^{1,1}(\widehat \C)$ that solves the Beltrami equation
\begin{align*}
\frac{\partial H}{\partial \overline{z}}= \mu \frac{\partial H}{\partial z}.
\end{align*}
Moreover, $H$ is unique up to postcomposition with M\"obius transformations.
\end{theorem}

The next theorem is a local version of the uniqueness part in Theorem~\ref{david_integrability_thm}, and plays an important role in applications of David homeomorphisms in dynamics.

\begin{theorem}\cite[Theorem~20.4.19, p.~565]{AIM09}\label{stoilow_thm}
Let $\Omega\subset \widehat{\C}$ be an open set and $f,g\colon \Omega\to \widehat{\C}$ be David embeddings with
$\mu_f=\mu_g$ almost everywhere. Then $f\circ g^{-1}$ is a conformal map on $g(\Omega)$.  
\end{theorem}

Let $\mu$ be a David conformal structure on an open subset $\Omega$ of $\widehat{\C}$. The local David homeomorphisms (defined on open subsets of $\Omega$) that straighten $\mu$ to the standard conformal structure (given by the circle field) define an atlas on $\Omega$, and by Theorem~\ref{stoilow_thm}, the corresponding change of coordinates are conformal maps. Thus, the David coefficient $\mu$ gives rise to a complex structure on $\Omega$.

While the composition of two quasiconformal homeomorphisms is always quasiconformal, the situation is more delicate for compositions of quasiconformal and David homeomorphisms. It turns out that post-composing a David homeomorphism $f:U\to V$ with a quasiconformal homeomorphism $g:V\to W$ always results in a David homeomorphism $g\circ f:U\to W$. However, the pre-composition of a David homeomorphism $f:U\to V$ with a quasiconformal homeomorphism $g:W\to U$ is not necessarily David since the David property of $f\circ g$ crucially depends on area distortion properties of $g$. The map $f\circ g: W\to V$ is indeed David if one has control on the map $g$ (for instance, if $g$ extends to a quasiconformal homeomorphism of an open neighborhood of $\overline{W}$ onto an open neighborhood of $\overline{U}$), or control over the geometry of the domains $U, W$ (for instance, if $U$ is a quasidisk and $W$ is a John domain). The proofs of these facts are given in \cite[Proposition~2.5]{LMMN}.

While the inverse of a $K$-quasiconformal map is also $K$-quasiconformal, the inverse of a David homeomorphism is not necessarily David (cf. \cite[p.~123]{Zak04}).

We direct the reader to \cite{Ahl06}, \cite[Chapters~3, 5, 20]{AIM09}, \cite[\S 2]{LMMN} for more background on the theory of quasiconformal and David homeomorphisms.

\subsection{Quasiconformal and David extensions of circle homeomorphisms}\label{qc_david_extension_subsec}

For an orientation-preserving homeomorphism $h\colon \mathbb S^1\to \mathbb S^1$, the \emph{distortion function} of $h$ is defined~as
\begin{align*}
\rho_h(z,t)= \max\left\{  \frac{|h(e^{2\pi i t}z)-h(z)| }{ |h(e^{-2\pi i t}z)-h(z)| } ,  \frac{|h(e^{-2\pi i t}z)-h(z)| }{ |h(e^{2\pi i t}z)-h(z)| }\right\},
\end{align*}
where $z\in \mathbb S^1$ and $0<t<1/2$. One further defines the \emph{scalewise distortion function} of $h$ to be  
\begin{align*}
\rho_h(t)= \max_{z\in \mathbb{S}^1}\rho_h(z,t),  
\end{align*}
where $0<t<1/2$. If $\rho_h(t)$ is bounded above, then $h$ is a quasisymmetric homeomorphism and the classical Ahlfors-Beurling Extension Theorem asserts that such an $h$ extends to a homeomorphism of $\overline{\D}$ that is quasiconformal on $\D$ \cite{BA56}. We will state a theorem, due to Chen-Chen-He and Zakeri, which says that if one has appropriate control on the growth of $\rho_h(t)$, then $h$ admits a David extension to the disk.

\begin{theorem}\cite[Theorem 3]{CCH96}, \cite[Theorem 3.1]{Zak08}\label{david_extension_criterion_thm}
Let $h\colon \mathbb S^1\to \mathbb S^1$ be an orientation-preserving homeomorphism and suppose that
\begin{align*}
\rho_h(t) = O(\log(1/t))\quad \textrm{as} \quad t\to 0.
\end{align*}
Then $h$ has an extension to a David homeomorphism $\widetilde h\colon \D\to \D$. 
\end{theorem}
 
A stronger David extension result for circle homeomorphisms was recently proved in~\cite{KN22}.

\subsection{David maps and removability}\label{removable_david_subsec}

A compact set $E\subset \widehat{\C}$ is said to be \emph{conformally removable} if every homeomorphism $f\colon \widehat{\C}\to \widehat{\C}$ that is conformal on $\widehat{\C} \setminus E$ is a M\"obius transformation. A compact set $E\subset \widehat{\C}$ is removable for $W^{1,1}$ functions if every continuous function $f\colon \widehat{\C}\to \R$ that lies in $W^{1,1}(\widehat{\C}\setminus E)$ in fact lies in $W^{1,1}(\widehat{\C})$. 

A domain $\Omega\subset \widehat{\C}$ is called a \emph{John domain} if for each base point $z_0\in \Omega$ there exists a constant $c>0$ such that for each point $z_1\in \Omega$ there exists an arc $\gamma$ joining $z_0$ to $z_1$ in $\Omega$ with the property that for each point $z$ on the path $\gamma$ we have 
\begin{align*}
\mathrm{dist}(z, \partial \Omega) \geq c\cdot \mathrm{length}(\gamma|_{[z,z_1]}),
\end{align*}
where $\gamma|_{[z,z_1]}$ denotes the subpath of $\gamma$ whose endpoints are $z$ and $z_1$ (here $\mathrm{dist}$ and $\mathrm{length}$ denote the spherical distance and spherical length). Roughly speaking, John domains are generalizations of quasidisks that allow for inward cusps but not outward cusps. By \cite[Theorem~4]{JS00}, boundaries of John domains are removable for $W^{1,1}$ functions.

\begin{theorem}\cite[Theorems~2.7,~2.12]{LMMN}\label{w11_removable_thm}
Suppose that $E\subset \widehat{\C}$ is a compact set that is removable for $W^{1,1}$ functions and $f\colon \widehat{\C}\to \widehat{\C}$ is a David homeomorphism. Then $f(E)$ is conformally removable. In particular, if $\Omega \subset {\widehat{\C}}$ is a John domain and $f\colon \widehat{\C}\to \widehat{\C}$ is a David homeomorphism, then $f(\partial \Omega)$ is conformally removable. 
\end{theorem}

\section*{List of notation} 
\begin{itemize}
\item $\D^*=\widehat{\C}\setminus\overline{\D}$.

\item $\eta(z)=1/\overline{z}$.

\item $B(a,r)=\{\vert z-a\vert<r\}$, where $a\in\C$ and $r>0$.

\noindent $\overline{B}(a,r)=\{\vert z-a\vert\leq r\}$, where $a\in\C$ and $r>0$.

\item $\textrm{Aut}(\widehat{\C})$ = Group of all M{\"o}bius automorphisms of $\widehat{\C}$,\\ 
$\textrm{Aut}^\pm(\widehat{\C})$ = Group of all M{\"o}bius and anti-M{\"o}bius automorphisms of $\widehat{\C}$.

\item $\textrm{Aut}(\D)$ = Group of all M{\"o}bius automorphisms of $\D$,\\ 
$\textrm{Aut}^\pm(\D)$ = Group of all M{\"o}bius and anti-M{\"o}bius automorphisms of $\D$.

\item $X^c =\widehat{\C}\setminus X$, for $X\subset\widehat{\C}$.

\item $m_{-d}:\mathbb{S}^1\to\mathbb{S}^1,\ \theta\mapsto -d\theta$.

\item $\mathscr{C}_d$ = Connectedness locus of monic, centered antiholomorphic polynomials of degree $d$.

\item $\pmb{G}_d$ = Regular ideal $(d+1)-$gon reflection group.

\noindent $\pmb{\cN}_d$ = Nielsen map of $\pmb{G}_{d}$ (Definition~\ref{regular_ideal_polygon_ref_group_def}).

\item $\cS_{\pmb{\cN}_d}$ = Space of normalized piecewise Schwarz reflection maps with $\pmb{\cN}_d$ as their external class (Section~\ref{mating_para_space_sec}).

\item $\mathbbm{G}_d$ = Anti-Hecke group isomorphic to $\Z/2\ast\Z/(d+1)$ (\S~\ref{nielsen_first_return_external_map_subsubsec}, \S~\ref{general_mating_corr_sec}).

\noindent $\pmb{\cF}_d$ = Anti-Farey map associated with $\mathbbm{G}_d$ (\S~\ref{nielsen_first_return_external_map_subsubsec},~\ref{anti_farey_subsubsec}).

\item $\cS_{\pmb{\cF}_d}$ = Space of Schwarz reflections having $\pmb{\cF}_d$ as their external class (\S~\ref{core_step_subsubsec}).

\item $\pmb{\mathcal{E}}_d$ = Minkowski circle homeomorphism conjugating $\pmb{\cN}_{d}$ to $\overline{z}^d$ (\S~\ref{question_mark_subsubsec}).

\item $\beta(\pmb{G}_d)$ = Bers slice of $\pmb{G}_d$ (\S~\ref{kissing_group_deform_space_subsubsec}).

\item $B_d(z) = \frac{(d+1)\overline z^d + (d-1)}{(d-1)\overline z^d + (d+1)}$ (unicritical parabolic antiholomorphic Blaschke product with critical {\'E}calle height $0$). 

\item $\pmb{\mathcal{B}}_d\ =\ $ Space of antiholomorphic rational maps having $B_d$ as their external class (\S~\ref{para_anti_rat_gen_subsubsec}).

\item $\mathcal{J}(R), \mathcal{F}(R)$ = Julia, Fatou set of a rational/anti-rational map $R$.

\noindent $\mathcal{K}(P), \mathcal{B}_\infty(P)$ = Filled Julia set, basin of infinity of a polynomial/anti-polynomial~$P$ (\S~\ref{anti_poly_dyn_general_subsubsec}).

\noindent $\mathcal{K}(R), \mathcal{B}(R)$ = Filled Julia set, marked parabolic basin of an anti-rational map $R\in\pmb{\mathcal{B}}_d$ (\S~\ref{para_anti_rat_gen_subsubsec}).

\item $G_{\mathcal{P}}$ = Kissing reflection group associated with a circle packing $\mathcal{P}$ (\S~\ref{kissing_group_subsubsec}).

\item $\Lambda(G), \Omega(G)$ = Limit set, domain of discontinuity of a reflection/Kleinian group.

\item $\Pi(G)$ = Canonical fundamental domain for the $G$-action on $\Omega(G)$ (\S~\ref{fund_dom_subsubsec}).

\item $\Omega_\infty(G)$ = Marked invariant component of the domain of discontinuity of a necklace group $G$.

\noindent $K(G)$ = The \emph{filled limit set} $\widehat{\C}\setminus \Omega_\infty(G)$ of a necklace group $G$.

\item $T(\sigma), T^0(\sigma), T^\infty(\sigma), K(\sigma)$ = Droplet, fundamental tile, tiling set, and non-escaping set of a piecewise Schwarz reflection map (\S~\ref{inv_partition_subsubsec}).

\item $\mathfrak{U}_{d+1}$ = Space of degree $d+1$ polynomials $f$ such that $f\vert_{\overline{\D}}$ is injective and $f$ has a unique (non-degenerate) critical point on $\mathbb{S}^1$.

\item $\mathscr{T}(R)$ = Tischler graph of a critically fixed anti-rational map $R$ (\S~\ref{group_map_bijection_subsubsec}).

\end{itemize}

\section*{List of figures}

\begin{itemize}
\item Figure~\ref{deltoid_intro_fig}: The ideal triangle reflection group and the deltoid Schwarz reflection map.

\item Figure~\ref{c_and_c_para_fig}: The Basilica limb of the Tricorn and the connectedness locus of the C\&C family.

\item Figure~\ref{limit_julia_intro_fig}: The Apollonian gasket limit set and the Apollonian Julia set.

\item Figure~\ref{itg_nielsen_fig}: The Nielsen map of the ideal triangle reflection group.

\item Figure~\ref{not_2_conn_fig}: Disconnected limit set of a kissing reflection group.

\item Figure~\ref{kissing_nielsen_fig}: Limit sets of various kissing reflection groups and their canonical fundamental domains.

\item Figure~\ref{question_mark_tree_fig}: The Farey and dyadic trees.

\item Figure~\ref{necklace_fig}: Limits sets of various necklace groups.

\item Figure~\ref{tricorn_fig}: The Tricorn and wiggling phenomena.

\item Figure~\ref{cardioid_disk_fig}: The Schwarz reflection map of the cardioid.

\item Figure~\ref{pre_deltoid_disk_fig}: Perturbation of the deltoid reflection map to anti-quadratic-like maps.

\item Figure~\ref{deltoid_disk_fig}: Mapping properties of the deltoid reflection map.

\item Figure~\ref{deltoid_corr_fig}: Dynamical planes of the deltoid reflection map and the associated antiholomorphic correspondence. 

\item Figure~\ref{basilica_schwarz_fig}: The rank zero and one tiles of a C\&C Schwarz reflection map.

\item Figure~\ref{basilica_anti_poly_fig}: The Basilica anti-polynomial Julia set and the limit set of the corresponding C\&C Schwarz reflection.

\item Figure~\ref{itg_nielsen_first_return_fig}: The degree two anti-Farey map.

\item Figure~\ref{anti_quad_like_fig}: A pinched anti-quadratic-like map.

\item Figure~\ref{chebyshev_center_fig}: Dynamical planes of a cubic Chebyshev Schwarz reflection map and the associated antiholomorphic correspondence.

\item Figure~\ref{talbot_schwarz_fig}: The dynamical planes of the Talbot Schwarz reflection, the Julia necklace group, and the Julia anti-polynomial.

\item Figure~\ref{apollo_cousins_fig}: The Apollonian Julia set and the limit set of the Deltoid-and-Circle Schwarz reflection.

\item Figure~\ref{affine_model_fig}: A quasiregular (partially affine) model of the Apollonian anti-rational map.

\item Figure~\ref{c_and_c_conn_cantor_limit_fig}: Limit sets of various Schwarz reflections in the C\&C family.

\item Figure~\ref{c_and_c_conn_locus_fig}: The connectedness locus of the C\&C family and the tessellation of its escape locus.

\item Figure~\ref{cheby_conn_locus_fig}: The connectedness locus of the cubic Chebyshev family of Schwarz reflections and the tessellation of its escape locus.

\item Figure~\ref{sigma_extremal_fig}: Droplets corresponding to two degree six Suffridge maps.

\item Figure~\ref{sigma_d_mating_fig}: Dynamical planes of two Schwarz reflection maps arising from $\Sigma_4^*$ and the corresponding necklace groups.

\item Figure~\ref{non_gasket_fig} and \ref{hamiltonian_polyhedral_fig}: Homeomorphic limit and Julia sets of various kissing reflection groups and critically fixed anti-rational maps.

\item Figure~\ref{enrichments_fig}: Two enrichments of the Tischler graph of $\overline{z}^3$.

\item Figure~\ref{ellipse_disk_fig}: The dynamical plane of the Schwarz reflection map in an ellipse and a pair of inscribed disks.

\item Figure~\ref{pinching_fig}: Pinching of quadrature multi-domains.

\item Figure~\ref{anti_farey_fig}: The degree three anti-Farey map.

\item Figure~\ref{augmented_dessin_fig}: Relation between the dessin d'enfants of Shabat polynomials and associated Schwarz reflections.

\item Figure~\ref{pine_tree_fig}: Conformally removable Jordan curves with infinitely many corners and inward/outward cusps.
\end{itemize}

\medskip

\noindent\textbf{Acknowledgments.} Part of this work was done during the authors' visits to MSRI (Simons Laufer Mathematical Sciences Institute), the Institute for Mathematical Sciences at Stony Brook, and the Urgench State University, Uzbekistan. The authors thank these institutes for their hospitality and support.
We would also like to thank Yusheng Luo and Dimitrios Ntalampekos for helpful conversations and useful comments.

\newcommand{\etalchar}[1]{$^{#1}$}

\end{document}